\documentclass[11pt,reqno]{amsart}
\usepackage[T1]{fontenc}
\usepackage[utf8]{inputenc}
\usepackage[a4paper,margin=25mm]{geometry}
\usepackage{amsmath,amssymb,amsthm,mathtools,array,longtable,booktabs}
\usepackage{microtype}
\usepackage{xurl}
\usepackage[section]{placeins}
\usepackage[colorlinks=true,linkcolor=blue,citecolor=blue,urlcolor=blue]{hyperref}
\hypersetup{pdftitle={Classification of G-g.o. metrics on compact homogeneous spaces G/H with G simple},pdfauthor={Huibin Chen, Zhiqi Chen, Zaili Yan and Fuhai Zhu},pdfsubject={Research draft}}
\newtheorem{maintheorem}{Theorem}

\newtheorem{maincorollary}[maintheorem]{Corollary}
\newtheorem{theorem}{Theorem}[section]
\newtheorem{lemma}[theorem]{Lemma}
\newtheorem{proposition}[theorem]{Proposition}
\newtheorem{corollary}[theorem]{Corollary}
\theoremstyle{remark}
\newtheorem{remark}[theorem]{Remark}
\newtheorem{example}[theorem]{Example}
\numberwithin{equation}{section}

\DeclareMathOperator{\ad}{ad}

\DeclareMathOperator{\SU}{SU}
\DeclareMathOperator{\PSU}{PSU}
\DeclareMathOperator{\SO}{SO}
\DeclareMathOperator{\Sp}{Sp}
\DeclareMathOperator{\Spin}{Spin}
\DeclareMathOperator{\U}{U}

\DeclareMathOperator{\SL}{SL}

\title[Classification of $G$-g.o. metrics]{Classification of $G$-g.o. metrics on compact homogeneous spaces $G/H$ with $G$ simple}
\author{Huibin Chen}
\author{Zhiqi Chen}
\author{Zaili Yan}
\author{Fuhai Zhu}
\date{}

\begin{document}
\begin{abstract}
We classify the $G$-g.o. metrics on compact homogeneous spaces $G/H$, where $G$ is a compact connected simple Lie group and $H$ is an arbitrary closed subgroup. We give the complete list of nonnormal metrics, including their parameter conditions and the possible mixing between equivalent isotropy summands. The proof uses rigidity criteria for arbitrary cometric operators, reductions through centralizers, and the stabilizers of the relevant representations. For every fixed pair $(G,H)$, the space spanned by the $G$-g.o. cometrics has dimension at most three, and we determine the equality cases. We also compare the classification with the compact weakly symmetric pairs and explain why the choice of the transitive group is essential.
\end{abstract}
\maketitle
\setcounter{tocdepth}{1}
\tableofcontents
\clearpage

\section{Introduction}
\label{sec:introduction}

Let $G/H$ be a compact homogeneous space, where $G$ is a compact connected simple Lie group and $H$ is a closed subgroup. A $G$-invariant Riemannian metric on $G/H$ is called a $G$-geodesic orbit metric, or a $G$-g.o. metric, if every geodesic is an orbit of a one-parameter subgroup of $G$~\cite{KoVa1991}. We keep the homogeneous presentation $G/H$ fixed. Here and throughout the paper, a connected Lie group is called simple if its Lie algebra is simple. We classify the $G$-g.o. metrics on $G/H$ without further restrictions on $H$.

Several classifications are known under additional assumptions on the isotropy. These include simple isotropy, abelian isotropy, maximal rank, regular embeddings, and certain homogeneous fibrations. For a general closed subgroup, equivalent isotropy modules may occur, and different irreducible summands of the defining representation may be coupled by the same simple ideal. A classification based only on diagonal metric operators does not cover these possibilities. Our purpose is to determine all $G$-g.o. metrics without such restrictions on $H$.

Set $Q=-B_{\mathfrak g}$ and $\mathfrak m=\mathfrak h^{\perp_Q}$, where $\mathfrak g$ and $\mathfrak h$ are the Lie algebras of $G$ and $H$. We represent a metric by its positive, $Q$-self-adjoint cometric operator $C$ through $g(u,v)=Q(C^{-1}u,v)$. A metric is called $G$-normal if $C=cI$ for some $c>0$. We call $C$ a $G$-g.o. cometric when its corresponding metric is $G$-g.o. Normality and natural reductivity are also understood with respect to $G$.

Our main result is the following. The group notation in the tables specifies standard models of Lie-algebra pairs. Finite central covers, effective central quotients, and closed isotropy extensions with the same identity component are included.

\begin{maintheorem}\label{thm11}
Let $G$ be a compact connected simple Lie group, and let $H\subsetneq G$ be closed. Every $G$-g.o. metric on $G/H$ is normal unless the Lie-algebra pair is one of the entries in Tables~\ref{tab:unitary-symplectic-list} and~\ref{tab:orthogonal-exceptional-list}, up to an automorphism of $\mathfrak g$. For each entry, all $G$-g.o. cometrics are precisely the operators described in \eqref{formula1}--\eqref{formula10}. The first entry has a nonnormal metric exactly when $w\ne0$. In every other entry, the displayed family contains nonnormal metrics.
\end{maintheorem}

\begin{table}[htbp]
\centering
\caption{The unitary and symplectic families.}
\label{tab:unitary-symplectic-list}
\small
\begin{tabular}{@{}p{.07\textwidth}p{.23\textwidth}p{.28\textwidth}p{.21\textwidth}p{.12\textwidth}@{}}
\toprule
& $\mathfrak g$ & $\mathfrak h$ & Range & Cometric \\
\midrule
1 & $\mathfrak{su}(d+q)$ & $\mathfrak{su}(d)\oplus\mathfrak r$ & $d>q\ge1$ & \eqref{formula1} \\
2 & $\mathfrak{su}(2n+1)$ & $\mathfrak{sp}(n)$ & $n\ge2$ & \eqref{formula2} \\
3 & $\mathfrak{su}(2n+1)$ & $\mathfrak{sp}(n)\oplus\mathbb R$ & $n\ge2$ & \eqref{formula3} \\
4 & $\mathfrak{sp}(n)$ & $\mathfrak{sp}(n-1)$ & $n\ge2$ & \eqref{formula4} \\
5 & $\mathfrak{sp}(n)$ & $\mathfrak{sp}(n-1)\oplus\mathbb R$ & $n\ge2$ & \eqref{formula4} \\
\bottomrule
\end{tabular}
\end{table}

In entry 1, the first summand acts on $\mathbb C^d$ and the isotropy satisfies $\SU(d)\subset H^0\subset S(\U(d)\times \U(q))$. The residual algebra $\mathfrak r$ is arbitrary within this inclusion. Put $z=i\operatorname{diag}(qI_d,-dI_q)$ and $w=P_{\mathfrak m}z$. When $w\ne0$, the family is
\begin{equation}\label{formula1}
C=bI+(a-b)P_{\mathbb Rw},\qquad a,b>0.
\end{equation}
When $w=0$, the family consists only of $C=bI$.

For entries 2 and 3, set $\mathfrak v=\mathfrak{su}(2n)\ominus\mathfrak{sp}(n)$, $\mathfrak p=(\mathbb C^{2n})_{\mathbb R}$ and $T=i\operatorname{diag}(I_{2n},-2n)$. In entry 2, the three coefficients are independent and
\begin{equation}\label{formula2}
C=aI_{\mathfrak v}\oplus bI_{\mathfrak p}\oplus tI_{\mathbb RT},\qquad a,b,t>0.
\end{equation}
The central circle in entry 3 is generated by $T$, and its cometrics are
\begin{equation}\label{formula3}
C=aI_{\mathfrak v}\oplus bI_{\mathfrak p},\qquad a,b>0.
\end{equation}
For entries 4 and 5, let $\mathfrak v$ and $\mathfrak p$ be the vertical and horizontal spaces of the fibration over $\Sp(n)/(\Sp(n-1)\Sp(1))$. The vertical dimension is three in entry 4 and two in entry 5. In both cases,
\begin{equation}\label{formula4}
C=aI_{\mathfrak v}\oplus bI_{\mathfrak p},\qquad a,b>0.
\end{equation}

\begin{table}[htbp]
\centering
\caption{The orthogonal and exceptional families.}
\label{tab:orthogonal-exceptional-list}
\small
\begin{tabular}{@{}p{.07\textwidth}p{.23\textwidth}p{.28\textwidth}p{.21\textwidth}p{.12\textwidth}@{}}
\toprule
& $\mathfrak g$ & $\mathfrak h$ & Range & Cometric \\
\midrule
6 & $\mathfrak{so}(2n+1)$ & $\mathfrak u(n)$ & $n\ge2$ & \eqref{formula5} \\
7 & $\mathfrak{so}(2n+1)$ & $\mathfrak{su}(n)$ & $n\ge2$ & \eqref{formula6} \\
8 & $\mathfrak{so}(2n)$ & $\mathfrak{su}(n)$ & $n\ge3$ odd & \eqref{formula7} \\
9 & $\mathfrak{so}(8+q)$ & $\mathfrak{spin}(7)\oplus\mathfrak r$ & $1\le q\le3$, $\mathfrak r\subset\mathfrak{so}(q)$ & \eqref{formula8} \\
10 & $\mathfrak{so}(9)$ & $\mathfrak g_2\oplus\mathfrak r$ & $\mathfrak r\subset\mathfrak{so}(2)$ & \eqref{formula8} \\
11 & $\mathfrak{so}(8)$ & $\mathfrak g_2$ & Standard embedding & \eqref{formula9} \\
12 & $\mathfrak e_6$ & $\mathfrak{so}(10)$ & Standard embedding & \eqref{formula10} \\
\bottomrule
\end{tabular}
\end{table}

For entries 6--8, write $\mathfrak v=(\Lambda^2\mathbb C^n)_{\mathbb R}$ and $J=iI_n$. In the odd-dimensional orthogonal group, put $\mathfrak p=(\mathbb C^n)_{\mathbb R}$. Entry 6 has the family
\begin{equation}\label{formula5}
C=aI_{\mathfrak v}\oplus bI_{\mathfrak p},\qquad a,b>0.
\end{equation}
In entry 7, the central tangent line has a prescribed coefficient,
\begin{equation}\label{formula6}
\begin{aligned}
C&=aI_{\mathfrak v}\oplus bI_{\mathfrak p}\oplus tI_{\mathbb RJ},\\
t&=\begin{cases}a,&n\text{ even},\\ ((n-1)a+b)/n,&n\text{ odd},\end{cases}
\qquad a,b>0.
\end{aligned}
\end{equation}
Entry 8 has the family
\begin{equation}\label{formula7}
C=aI_{\mathfrak v}\oplus bI_{\mathbb RJ},\qquad a,b>0.
\end{equation}
For entry 9, use the eight-dimensional spin embedding of $\Spin(7)$. For entry 10, use the seven-dimensional defining representation of $G_2$. With $(d,\mathfrak s)=(8,\mathfrak{spin}(7))$ or $(7,\mathfrak g_2)$, respectively, both families are
\begin{equation}\label{formula8}
C=aI_{\mathfrak{so}(d)\ominus\mathfrak s}\oplus bI_{(\mathfrak{so}(q)\ominus\mathfrak r)\oplus(\mathbb R^d\otimes\mathbb R^q)},\qquad a,b>0,
\end{equation}
where $q=2$ in entry 10. In entry 11, choose a $Q$-orthogonal equivariant identification $\mathfrak m=\mathbb R^7\otimes\mathbb R^2$. Every positive invariant cometric is allowed,
\begin{equation}\label{formula9}
C=I_{\mathbb R^7}\otimes P,\qquad P\in\operatorname{Sym}_2^+(\mathbb R).
\end{equation}
Finally, entry 12 has the decomposition $\mathfrak m=\mathbb R\oplus(\Delta^+_{10})_{\mathbb R}$, where the half-spin summand has real dimension $32$, and
\begin{equation}\label{formula10}
C=aI_{\mathbb R}\oplus bI_{(\Delta^+_{10})_{\mathbb R}},\qquad a,b>0.
\end{equation}

In the scalar families, normality means equality of all displayed coefficients. In entry 11, it means that $P$ is scalar. The table contains no nonnormal metric for $G_2,F_4,E_7$ or $E_8$. Low-rank identifications may identify entries as Lie-algebra pairs, even when their displayed global quotients are different.

\begin{maincorollary}\label{cor12}
For a fixed pair $(G,H)$ as in Theorem~\ref{thm11}, the vector space spanned by the $G$-g.o. cometrics has dimension at most three. Equality holds exactly when the Lie-algebra pair is equivalent to $(\mathfrak{su}(2n+1),\mathfrak{sp}(n))$, $n\ge2$, or $(\mathfrak{so}(8),\mathfrak g_2)$. This dimension includes the overall scale and is unchanged by finite isotropy extensions.
\end{maincorollary}

We describe the main points of the proof. The g.o. equation is linear in the cometric operator at each fixed tangent vector. However, an invariant cometric may mix equivalent isotropy summands, and its compression to a subspace need not be its restriction. We first prove rigidity statements which retain these mixed terms. Following T\'{o}th's Hamiltonian approach~\cite{To2010}, we combine constancy along paths with the connectedness of compact moment-map fibers. The energy of a $G$-g.o. cometric is then constant on every intersection of a coadjoint orbit with $\mathfrak h^\perp$. This gives a global proof of invariance under the full normalizer. We also show that two orthogonal invariant subspaces with free infinitesimal isotropy actions force normality, without requiring the cometric to preserve those subspaces.

For the classical groups, we apply these criteria to the irreducible summands of the defining representation. The remaining coupled terms are excluded by abelian generic stabilizers. In the orthogonal case, two odd-dimensional blocks require an additional argument. We use singular adjoint orbits and a six-dimensional Pfaffian calculation to prove that the cometric preserves the off-diagonal block. For the exceptional groups, we choose a simple ideal of $\mathfrak h$ and study its full centralizer. The resulting criteria allow arbitrary residual isotropy. Complete inventories of the embedded simple ideals then give the required exhaustion.

We recall the main developments related to the classification problem. Kowalski and Vanhecke introduced the term geodesic orbit space and classified these spaces in dimensions at most six~\cite{KoVa1991}. Naturally reductive spaces form a basic class of g.o. spaces. Kaplan's nilpotent examples and the seven-dimensional examples of Du\v{s}ek, Kowalski and Nik\v{c}evi\'{c} show that the g.o. condition is strictly more general than natural reductivity~\cite{Ka1983,DuKoNi2004}. Weakly symmetric Riemannian manifolds form another important class~\cite{BeKoVa1997}. Their compact classification is closely related to spherical pairs~\cite{Kr1979,Ng2000,Ya2004}. These results concern the isometry group of the metric, and the choice of a prescribed transitive group remains essential. We refer to~\cite{Ar2017,BeNi2020} for general accounts of the Riemannian theory.

The general structure of g.o. manifolds is closely related to the radical and nilradical of the isometry algebra. Gordon proved that a Riemannian g.o. nilmanifold is at most two-step nilpotent and that a g.o. manifold with nonpositive Ricci curvature is symmetric~\cite{Go1996}. Nikonorov obtained further restrictions from Killing fields of constant length, totally geodesic submanifolds and principal isotropy representations~\cite{Ni2013a,Ni2017}. Gordon and Nikonorov described the structure of g.o. manifolds diffeomorphic to Euclidean space~\cite{GoNi2018}. Within the nilpotent class, classifications are known for groups of Heisenberg type, modified H-type groups and groups associated with graphs~\cite{Ri1984,La1999,Ni2020}. Further constructions and the classification of nonsingular g.o. nilmanifolds show the range of examples which remain after the two-step reduction~\cite{Ni2024a,NiZi2026a}.

For compact spaces, many classification results use restrictions on the isotropy group or its representation. Alekseevsky and Arvanitoyeorgos classified $G$-g.o. metrics on generalized flag manifolds~\cite{AlAr2007}, and Alekseevsky and Nikonorov treated compact simply connected spaces of positive Euler characteristic~\cite{AlNi2009}. For a simple transitive group, their nonnormal metrics occur on two flag families. Further classifications concern spheres~\cite{Ni2013b}, spaces with two irreducible isotropy summands~\cite{ChNi2019}, and spaces with simple isotropy~\cite{ChNiNi2023}. Souris treated abelian isotropy, compact transitive groups of rank two, and semisimple regular isotropy~\cite{So2021,So2022,So2026}. In the last setting, the natural reductivity conclusion may involve an enlarged transitive group. It therefore does not imply normality with respect to the original simple group.

Homogeneous fibrations provide another approach to the compact classification. Tamaru classified the fibrations over irreducible symmetric spaces for which the two-scalar canonical variations are $G$-g.o.~\cite{Ta1999}. Classifications and metric restrictions have also been obtained for homogeneous bundles over Stiefel manifolds~\cite{ArSoSt2021a,ArSoSt2021b}, generalized Wallach spaces~\cite{ArWa2017}, and $M$-spaces and generalized $C$-spaces~\cite{ArWaZh2017,ArWaZh2018,ArQiWaZh2019}. Related constructions use strongly isotropy irreducible spaces and selected ideals in the isotropy of a flag manifold~\cite{ChChZh2021,ChChYaZh2023,ChChCh2025}. The intermediate-subgroup method in~\cite{ChChZh2026} gives criteria for passing to the base and the fiber and produces nonnaturally reductive g.o. metrics. These results motivate reductions through intermediate groups, with the metric-invariance assumptions checked at each step.

Equivalent isotropy summands require additional attention because an invariant metric may contain mixed terms. Souris derived algebraic restrictions on these terms and applied them to compact spaces including complex Stiefel manifolds~\cite{So2018}. Such terms also occur in the classifications of $G$-g.o. metrics on real flag manifolds of split classical type and split type $G_2$~\cite{GrGrNe2020,GrRoSa2025}. The compact transitive group in the latter case is locally $\SO(4)$. For Ledger--Obata spaces, Nikolayevsky and Nikonorov proved that the $G$-g.o. condition for the prescribed product group $G$ is equivalent to natural reductivity and described the corresponding metrics~\cite{NiNi2019}. These examples show why a general classification must retain the full invariant metric operator.

For left-invariant metrics on compact simple Lie groups, the g.o. condition is usually taken with respect to the connected isometry group, which can include both left and right translations. D'Atri and Ziller classified the naturally reductive metrics in this setting~\cite{DaZi1979}. Natural reductivity of g.o. metrics has been proved for families arising from flag manifolds and for metrics on the compact group $G_2$~\cite{ChChWo2018,So2025b}. Regular subgroup methods give further restrictions~\cite{So2025a}. The construction of Einstein metrics which are not g.o. also separates the Einstein condition from the g.o. condition~\cite{ChChDe2018,Ni2019,XuTa2019,LuTaXu2025}. These results address a different choice of acting group from the classification of $G$-g.o. metrics on $G/H$.

The present classification allows arbitrary closed isotropy in a compact group with simple Lie algebra. The additional cases include nonregular embeddings, repeated irreducible modules and simple ideals acting on several summands of the defining representation. We use the known classifications under their stated hypotheses and retain the attribution of the known metric families. Our task is to exclude all other embeddings while allowing every invariant cometric, including its mixed terms.

The choice of $G$ cannot be suppressed. A weakly symmetric Riemannian manifold may admit a compact simple transitive group $G$ for which its metric is not $G$-g.o. We compare the tables with the compact simple weakly symmetric pairs and give an explicit example on $\SU(6)/(\SU(3)\times \SU(3))$ in Section~\ref{sec:weak-symmetry-comparison}.

The paper is organized as follows. We recall the algebraic formulation and the external inputs in Section~\ref{sec:unified-preliminaries}, and prove the rigidity criteria in Section~\ref{sec:unified-rigidity}. Sections~\ref{sec:unified-a}, \ref{sec:unified-c} and~\ref{sec:unified-bd} treat the classical groups. The subsequent sections prove the exceptional classification by centralizer reductions and complete embedding inventories. We then prove Theorem~\ref{thm11} and Corollary~\ref{cor12}, and compare the result with weakly symmetric spaces. The appendices contain the high-rank embedding tables, the rank-one inventories and their order of dependence, and the branching calculations used in the proof.

\section{Preliminaries}
\label{sec:unified-preliminaries}

In this section, we recall the geodesic lemma and introduce the representation-theoretic facts used in the classification.

Let $G$ be a compact connected Lie group with simple Lie algebra $\mathfrak g$, and let $H\subsetneq G$ be a closed subgroup with Lie algebra $\mathfrak h$. Put $Q=-B_{\mathfrak g}$ and $\mathfrak m=\mathfrak h^{\perp_Q}$. We identify $\mathfrak m$ with the tangent space of $G/H$ at the origin. A $G$-invariant Riemannian metric is determined by a positive, $Q$-self-adjoint, $\operatorname{Ad}(H)$-equivariant operator $A$ through $g(u,v)=Q(Au,v)$. We use its inverse $C=A^{-1}$, which we call the cometric operator. The $G$-g.o. condition is expressed in terms of $C$ by the following geodesic lemma.

\begin{lemma}\label{lem21}
The metric determined by $C$ is $G$-g.o. if and only if
\begin{equation}\label{formula11}
 [CX,X]=[Z,X]\quad\text{for some }Z\in\mathfrak h,\qquad X\in\mathfrak m.
\end{equation}
The equality in \eqref{formula11} holds in $\mathfrak g$.
\end{lemma}
\begin{proof}
The geodesic lemma says that the velocity $u\in\mathfrak m$ is generated by $u+z\in\mathfrak g$ if and only if $Q(Au,[u+z,v])=0$ for every $v\in\mathfrak m$. Set $X=Au$, so $u=CX$. This condition is equivalent to the vanishing of the $\mathfrak m$-component of $[u+z,X]$. Its $\mathfrak h$-component also vanishes. Indeed, for $T\in\mathfrak h$, equivariance of $C$ and skew-adjointness of $\operatorname{ad}T$ give $Q([CX,X],T)=-Q(CX,[T,X])=0$, while $[z,X]\in\mathfrak m$. Thus $[CX+z,X]=0$. Replacing $z$ by $-Z$ proves the assertion.
\end{proof}

A metric is normal if it is induced by a bi-invariant metric on $G$. Since $\mathfrak g$ is simple, this is equivalent to $C=cI$ for some $c>0$. Natural reductivity will also refer to the fixed group $G$.

\begin{proposition}\label{prop22}
Suppose that a linear map $P\colon\mathfrak m\to\mathfrak h$ satisfies $[P(X),X]=[CX,X]$ for every $X\in\mathfrak m$. Then $C$ is scalar. In particular, a $G$-invariant metric is $G$-naturally reductive if and only if it is normal.
\end{proposition}
\begin{proof}
Averaging $P$ over $H$ preserves the identity and makes $P$ equivariant. In velocity coordinates $u=CX$, the map $u\mapsto-P(C^{-1}u)$ is a linear equivariant geodesic graph. Its graph is an $\operatorname{Ad}(H)$-invariant complement $\mathfrak v$ to $\mathfrak h$, and every vector of $\mathfrak v$ generates a geodesic. Polarizing the geodesic lemma on this complement gives natural reductivity. This is also the argument of \cite[Lemma~3]{ChNiNi2023}, which does not require simple isotropy.

For a naturally reductive complement, Kostant's extension theorem \cite[Theorem~4]{Ko1956} gives an invariant nondegenerate symmetric bilinear form on the ideal $\mathfrak v+[\mathfrak v,\mathfrak v]$, extending the metric on $\mathfrak v$ and making $\mathfrak v$ orthogonal to its intersection with $\mathfrak h$. This ideal is nonzero and therefore equals $\mathfrak g$. Every invariant symmetric bilinear form on compact simple $\mathfrak g$ is a multiple of $Q$. Its restriction to $\mathfrak v$ is positive, so the multiple is positive. Consequently $\mathfrak v=\mathfrak m$ and the metric is normal. Conversely, a normal metric is naturally reductive.
\end{proof}

For an orthogonal representation of a compact Lie algebra $\mathfrak l$ on $V$, write $\mathfrak l_v=\{Z\in\mathfrak l:Zv=0\}$. We say that $v$ is free if $\mathfrak l_v=0$. Thus freeness always means infinitesimal freeness. A free tuple is a tuple whose joint infinitesimal stabilizer is zero. When a free vector exists, free vectors form a nonempty Zariski-open subset, since injectivity of $Z\mapsto Zv$ is detected by maximal minors.

\begin{lemma}\label{lem23}
Let a real Lie algebra $\mathfrak l$ act on a real vector space $V$ of dimension $n$. Suppose that $\operatorname{rank}\rho(Z)\ge d$ for every nonzero $Z\in\mathfrak l$. If $qd>\dim\mathfrak l-1$, then $V$ has a free $q$-tuple. If $q\le n$, the tuple can be chosen linearly independent. For an orthogonal representation it can then be orthogonalized. In particular, $d\ge\dim\mathfrak l$ implies the existence of a free vector.
\end{lemma}
\begin{proof}
Parametrize nonzero elements of $\mathfrak l$, up to positive scaling, by its unit sphere. The incidence set of pairs consisting of $Z$ on this sphere and a $q$-tuple in $\ker\rho(Z)$ has dimension at most $\dim\mathfrak l-1+q(n-d)<qn$. To justify this bound, partition the sphere into finitely many semialgebraic strata on which the rank is constant and apply the dimension bound to each stratum. The projection of the incidence set cannot fill $V^q$. The free tuples form a nonempty Zariski-open set. When $q\le n$, intersect this set with the open dense set of independent tuples. Replacing a tuple by an orthonormal basis of its span preserves its pointwise stabilizer.
\end{proof}

\begin{lemma}\label{lem24}
Let $\mathfrak a$ and $\mathfrak b$ act faithfully and orthogonally on real spaces $U$ and $V$. Put $r=\min(\dim U,\dim V)$. If the representation on the larger space has a free tuple of at most $r$ vectors, then $U\otimes V$ has a free vector for $\mathfrak a\oplus\mathfrak b$.
\end{lemma}
\begin{proof}
Suppose that $\dim U\le\dim V$. View a tensor as a map $T\colon U\to V$. Choose $T$ injective with distinct positive singular values and with image containing a free tuple for $\mathfrak b$. Such a map exists by extending the span of the tuple to a subspace of dimension $\dim U$. If $BT=TA$, where $A$ and $B$ are the skew-adjoint representation operators, then $[A,T^*T]=0$. Since $T^*T$ has simple spectrum, $A=0$. Faithfulness kills the $\mathfrak a$ component. Then $B$ vanishes on $\operatorname{im}T$, so the $\mathfrak b$ component also vanishes. Transpose $T$ for the other dimension ordering.
\end{proof}

\begin{lemma}\label{lem25}
The following orthogonal representations have free tuples of the indicated lengths.
\begin{enumerate}
\item The adjoint representation of a compact semisimple Lie algebra has a free pair.
\item The conjugation representation of $\mathfrak{so}(d)$ on $\operatorname{Sym}^2_0\mathbb R^d$ has a free vector.
\item The conjugation representation of $\mathfrak{sp}(n)$ on the real space $W_n$ of traceless quaternionic Hermitian matrices has a free triple for $n\ge3$ and a free quadruple for $n=2$. Its real dimension is $2n^2-n-1$.
\end{enumerate}
\end{lemma}
\begin{proof}
For the first assertion, choose a regular element of a Cartan subalgebra and an element with nonzero components in every simple-root plane. A common centralizer belongs to the Cartan subalgebra and is annihilated by every simple root, so it is zero. For the second, a real symmetric matrix with distinct eigenvalues has zero skew-symmetric centralizer.

For the third, a real diagonal matrix with distinct entries reduces its infinitesimal centralizer to diagonal imaginary quaternions. A Hermitian matrix whose off-diagonal entries are all $1$ makes these diagonal entries equal. If $n\ge3$, a third Hermitian matrix with entries $i$ and $j$ in positions $(1,2)$ and $(1,3)$ forces their common value to commute with both $i$ and $j$, and hence to vanish. If $n=2$, use two additional matrices with entries $i$ and $j$ in position $(1,2)$. These matrices can all be chosen traceless. Finally, a quaternionic Hermitian matrix has $n$ real diagonal entries and $4\binom n2$ real off-diagonal coordinates, giving the stated dimension after removing the trace.
\end{proof}

\begin{lemma}\label{lem26}
Let $d\ge3$, let $V_d=\operatorname{Sym}^{d-1}\mathbb C^2$, and let $\chi$ be a nonzero character of $\U(1)$. The underlying real representation of $V_d\otimes\mathbb C_\chi$ has a free vector for $\mathfrak{su}(2)\oplus\mathbb R$.
\end{lemma}
\begin{proof}
Write $m=d-1$ and take the binary form $v=x^m+2y^m$. Choose parameters $t\in\mathbb R$ and $b\in\mathbb C$ so that an infinitesimal element of $\mathfrak{su}(2)$ acts by $Ax=itx-\overline b y$ and $Ay=bx-ity$. Write the scalar phase as $iqs$, where $q\ne0$. For $m\ge3$, the coefficients of the distinct mixed monomials in $(A+iqsI)v$ are $-m\overline b$ and $2mb$, so $b=0$. For $m=2$, the combined mixed coefficient is $2(2b-\overline b)$, which again forces $b=0$. The two remaining coefficients give $mt+qs=0$ and $-mt+qs=0$. Hence $t=s=0$.
\end{proof}

\section{Rigidity of \texorpdfstring{$G$-g.o.}{G-g.o.} cometrics}
\label{sec:unified-rigidity}

In this section, we prove the common rigidity statements used for the classical and exceptional groups. We first prove invariance under the full normalizer and then consider invariant complements and infinitesimal stabilizers.

The Hamiltonian observation underlying the next lemma appears in \cite[Lemma~3.9 and Theorem~3.10]{To2010}. We combine it with connectedness of the relevant moment-map fiber to obtain the global statement needed here.

\begin{lemma}\label{lem31}
Let $G$ be compact connected semisimple, let $H\subset G$ be closed, and let $Q$ be an invariant positive inner product on $\mathfrak g$. Put $\mathfrak m=\mathfrak h^{\perp_Q}$. If $C$ is a $G$-g.o. cometric, then $F(X)=\frac12Q(CX,X)$ is constant on every set $\operatorname{Ad}(G)X\cap\mathfrak m$.
\end{lemma}
\begin{proof}
Use $Q$ to identify $\mathfrak g$ with its dual. The compact connected coadjoint orbit $\mathcal O=\operatorname{Ad}(G)X$ carries its Kirillov--Kostant--Souriau symplectic form. The action of $H^0$ is Hamiltonian with moment map $\mu(Y)=Q(Y,\cdot)|_{\mathfrak h}$. Since $\mathcal O$ is compact, $\mu$ is proper. The connected-fiber theorem \cite[Theorem~1.1(ii)]{LeMeToWo1998} implies that $\mu^{-1}(0)=\mathcal O\cap\mathfrak m$ is connected. This theorem has no regular-value assumption.

We also need paths in this possibly singular fiber. The group $\operatorname{Ad}(G)=\operatorname{Aut}(\mathfrak g)^0$ is semialgebraic, since the automorphism group is real algebraic and its identity component is semialgebraic. Its orbit $\mathcal O$ and the intersection $\mathcal O\cap\mathfrak m$ are therefore semialgebraic. Every connected semialgebraic set is path connected by paths that admit a finite subdivision into $C^1$ arcs.

Along such an arc $Y(t)$, write $\dot Y=[T,Y]\in\mathfrak m$. By \eqref{formula11}, choose $Z\in\mathfrak h$ with $[CY,Y]=[Z,Y]$. Then
\begin{equation}\label{formula12}
 \frac{d}{dt}F(Y)=Q(CY,[T,Y])=Q([Y,CY],T)=-Q([Z,Y],T)=Q(Z,\dot Y)=0.
\end{equation}
By \eqref{formula12}, $F$ is constant on each arc. Continuity at the subdivision points proves constancy on the whole fiber.
\end{proof}

\begin{corollary}\label{cor32}
Under the assumptions of Lemma~\ref{lem31}, $C$ commutes with $\operatorname{Ad}(n)|_{\mathfrak m}$ for every $n\in N_G(H^0)$.
\end{corollary}
\begin{proof}
Every such $n$ preserves $\mathfrak h$ and $\mathfrak m$. Since $\operatorname{Ad}(n)X\in\operatorname{Ad}(G)X\cap\mathfrak m$, the lemma gives $F(\operatorname{Ad}(n)X)=F(X)$. Polarization proves the assertion.
\end{proof}

The conclusion of Corollary~\ref{cor32} is also stated in \cite[Corollary~4]{Ni2017}. The preceding proof includes the disconnected components of the normalizer.

\begin{corollary}\label{cor33}
For compact connected semisimple $G$ and closed $H$, the $G$-g.o. cometric cones of $G/H^0$ and $G/H$ coincide under their common tangent-space identification. Passing to a central covering group also leaves this cone unchanged.
\end{corollary}
\begin{proof}
Every element of $H$ normalizes $H^0$. Hence Corollary~\ref{cor32} makes every $G$-g.o. cometric on $G/H^0$ invariant under $H$. It descends and satisfies the same equation \eqref{formula11}. Conversely, a metric on $G/H$ lifts to $G/H^0$. Central coverings identify the Lie algebras and the equation, and their finite isotropy extensions have just been treated.
\end{proof}

\begin{corollary}\label{cor34}
Let $G=\SO(4k+2)$ and let $r\in\mathrm O(4k+2)$ normalize $\mathfrak h$. Then every $G$-g.o. cometric commutes with $\operatorname{Ad}(r)|_{\mathfrak m}$.
\end{corollary}
\begin{proof}
Only $\det r=-1$ requires an argument. Write $4k+2=2N$, where $N$ is odd. For every $X\in\mathfrak{so}(2N)$, the elements $X$ and $-\operatorname{Ad}(r)X$ are $\SO(2N)$-conjugate. To see this, put $X$ in a Cartan subalgebra with $N$ real rotation blocks. Conjugation by an orientation-reversing orthogonal matrix changes one coordinate sign, modulo the type-$D_N$ Weyl group. Subsequent negation changes the other $N-1$ signs. Since this number is even, the resulting change belongs to the Weyl group. As $r$ preserves $\mathfrak m$, Lemma~\ref{lem31} and the evenness of $F$ imply $F(\operatorname{Ad}(r)X)=F(-\operatorname{Ad}(r)X)=F(X)$. Polarization completes the proof.
\end{proof}

\begin{corollary}\label{cor35}
Suppose that $H^0\subset\Sp(N)$ preserves an orthogonal quaternionic decomposition $\mathbb H^N=V_1\oplus\cdots\oplus V_r$. Put $\mathfrak k=\bigoplus_i\mathfrak{sp}(V_i)$ and $\mathfrak f=\mathfrak k\ominus\mathfrak h$. Let $\mathfrak p_{ij}$ be the off-diagonal skew-Hermitian block joining $V_i$ and $V_j$. Every $\Sp(N)$-g.o. cometric preserves $\mathfrak f$ and each $\mathfrak p_{ij}$.
\end{corollary}
\begin{proof}
The element equal to $-I$ on $V_i$ and $I$ on its orthogonal complement belongs to $\Sp(N)$ and centralizes $H^0$. All these signs act trivially on $\mathfrak f$, while their character on $\mathfrak p_{ij}$ is $\varepsilon_i\varepsilon_j$. These characters are nontrivial and distinct for distinct unordered pairs. Apply Corollary~\ref{cor32} to their common eigenspaces.
\end{proof}

We next use free subspaces without assuming that the cometric preserves them.

\begin{lemma}\label{lem36}
Let $\mathfrak g$ be compact simple and suppose that $\mathfrak m=\bigoplus_{i=1}^sV_i$ is an orthogonal $\mathfrak h$-invariant decomposition. If every complement $W_i=\bigoplus_{j\ne i}V_j$ has a free vector for $\mathfrak h$, then every $G$-g.o. metric on $G/H$ is normal.
\end{lemma}
\begin{proof}
Put $B(X)=[CX,X]$. By \eqref{formula11}, $B(X)\in U$ whenever $X\in U$ and $U$ is $\mathfrak h$-invariant. Each free locus in $W_i$ is nonempty and Zariski open. Hence there is a nonempty open set $\Omega\subset\mathfrak m$ on which all complement projections are free. The compensator $Z(X)$ in \eqref{formula11} is unique on $\Omega$.

Fix $i$ and write $X=x+y$, where $x\in W_i$ is free and $y\in V_i$. Projecting \eqref{formula11} to $W_i$, we obtain
\begin{equation}\label{formula13}
 [Z(x+y),x]=B(x)+P_{W_i}\bigl([Cx,y]+[Cy,x]\bigr).
\end{equation}
The map $Z\mapsto[Z,x]$ is injective. Thus \eqref{formula13} makes the unique compensator affine in $y$. Choose a product of open balls contained in $\Omega$. On this product $Z$ is affine in each variable separately, so it is a multiaffine polynomial $P$. The polynomial identity $[P(X),X]=B(X)$ holds on an open set and therefore holds everywhere. Its homogeneous degree-two component is $[P_1(X),X]=[CX,X]$, where $P_1$ is linear. Proposition~\ref{prop22} proves normality. In particular, no regularity of a geodesic graph at singular vectors is used.
\end{proof}

\begin{proposition}\label{prop37}
Let $L\subset G$ be closed, with $\mathfrak g$ compact simple. Suppose that $\mathfrak l^\perp$ contains two orthogonal $L^0$-invariant subspaces, each with a free vector for $\mathfrak l$. Then every $G$-g.o. metric on $G/H$ is normal for every closed subgroup $H\subset L$.
\end{proposition}
\begin{proof}
The same vectors are free for $\mathfrak h$. Use one subspace as the first summand of $\mathfrak m$ and its orthogonal complement as the second. The second summand contains the other free subspace. Lemma~\ref{lem36} applies. Corollary~\ref{cor33} includes disconnected isotropy.
\end{proof}

\begin{lemma}\label{lem38}
Let $\mathfrak h\subset\mathfrak k\subset\mathfrak g$ and put $\mathfrak p=\mathfrak k^{\perp_Q}$. If a $G$-g.o. cometric $C$ preserves $\mathfrak p$, then $C|_{\mathfrak p}$ commutes with $\operatorname{ad}(\mathfrak k)|_{\mathfrak p}$ and satisfies the g.o. equation for $(\mathfrak g,\mathfrak k)$. In particular, it defines a $G$-g.o. metric on $G/K$ when $K$ is the connected closed subgroup with Lie algebra $\mathfrak k$.
\end{lemma}
\begin{proof}
For $x\in\mathfrak p$, equation \eqref{formula11} gives $[C x,x]=[Z,x]\in\mathfrak p$, since $Z\in\mathfrak h\subset\mathfrak k$. For every $T\in\mathfrak k$, pairing with $T$ gives $Q(Cx,[T,x])=0$. Since $C|_{\mathfrak p}$ is self-adjoint and $\operatorname{ad}T|_{\mathfrak p}$ is skew-adjoint, polarization gives their commutation. The same $Z\in\mathfrak h\subset\mathfrak k$ supplies the descended g.o. equation.
\end{proof}

The metric version of Lemma~\ref{lem38} is contained in \cite[Lemmas~4.1--4.2]{ChChZh2026}. The proof records explicitly the equivariance needed for descent.

\begin{lemma}\label{lem39}
Let $\mathfrak h\subset\mathfrak k\subsetneq\mathfrak g$, with $\mathfrak g$ compact simple, and put $\mathfrak p=\mathfrak k^\perp$. Suppose that $C$ preserves $\mathfrak p$, that $C|_{\mathfrak p}=cI$, and that $\mathfrak p$ has a free vector for $\mathfrak k$. Then $C=cI$ on $\mathfrak m$. The scalarity assumption on $C|_{\mathfrak p}$ holds, in particular, if $\mathfrak p$ is an irreducible $\mathfrak k$-module.
\end{lemma}
\begin{proof}
Put $\mathfrak w=\mathfrak k\ominus\mathfrak h$ and $E=C|_{\mathfrak w}$. At $y+x$, with $y\in\mathfrak w$ and $x\in\mathfrak p$, the $\mathfrak p$-component of \eqref{formula11} gives $[(E-cI)y-Z,x]=0$ for some $Z\in\mathfrak h$. Choose a free $x$. Then $(E-cI)y=Z$, and the two sides lie in the orthogonal spaces $\mathfrak w$ and $\mathfrak h$. Thus both vanish. The last assertion follows from Lemma~\ref{lem38} and self-adjointness.
\end{proof}

\begin{proposition}\label{prop310}
Put $\mathfrak n=\mathfrak n_{\mathfrak g}(\mathfrak h)=\mathfrak h\oplus\mathfrak l$, where $\mathfrak l=\mathfrak c_{\mathfrak m}(\mathfrak h)$, and put $\mathfrak p=\mathfrak n^\perp$. Suppose that the descended $G$-g.o. metric on $G/N_G(H^0)$ is normal. Then
\begin{equation}\label{formula14}
 C=E|_{\mathfrak l}\oplus bI_{\mathfrak p},\qquad b>0,
\end{equation}
where $E$ is positive, self-adjoint and commutes with $\operatorname{ad}(\mathfrak l)$. Conversely, an isotropy-invariant operator of the form \eqref{formula14} is $G$-g.o. if and only if
\begin{equation}\label{formula15}
 [DY,X]\in[\mathfrak h,X]\quad(Y\in\mathfrak l,\ X\in\mathfrak p),\qquad D=E-bI.
\end{equation}
Its image $\mathfrak i=\operatorname{im}D$ is an ideal of $\mathfrak n$. The connected group generated by $H^0$ and $\exp\mathfrak i$ has the same orbits on $\mathfrak p$ as $H^0$.
\end{proposition}
\begin{proof}
The subspace $\mathfrak l=\mathfrak m^{\mathfrak h}$ is preserved by every isotropy-equivariant self-adjoint operator. Its orthogonal complement in $\mathfrak m$ is $\mathfrak p$. Lemma~\ref{lem38} proves descent, and Corollary~\ref{cor32} gives invariance under the full normalizer. The hypothesis therefore gives \eqref{formula14}. The same corollary implies that $E$ commutes with $\operatorname{ad}(\mathfrak l)$. Hence $[EY,Y]=0$, since $\operatorname{ad}Y(EY)=E\operatorname{ad}Y(Y)=0$. Also $[\mathfrak h,\mathfrak l]=0$. At $Y+X$, equation \eqref{formula11} becomes exactly \eqref{formula15}.

The image of the equivariant map $D$ is an ideal of $\mathfrak l$, and it commutes with $\mathfrak h$. Thus it is an ideal of $\mathfrak n$. By \eqref{formula15}, every vector field generated by $\mathfrak i$ is tangent to every $H^0$-orbit. These orbits are compact, so the flows preserve them. This proves the orbit assertion.
\end{proof}

We now use a semisimple ideal of the isotropy algebra to obtain the invariant complement required for descent.

\begin{lemma}\label{lem311}
Let $0\ne\mathfrak s\triangleleft\mathfrak h\subsetneq\mathfrak g$ be semisimple, and put $\mathfrak k=\mathfrak s\oplus\mathfrak c_{\mathfrak g}(\mathfrak s)$. Then $\mathfrak k$ is the Lie algebra of the proper closed connected subgroup $K=N_G(S)^0$, where $S$ has Lie algebra $\mathfrak s$. Moreover, $\mathfrak h\subset\mathfrak k$ and every isotropy-equivariant self-adjoint operator preserves $\mathfrak m=\mathfrak w\oplus\mathfrak p$, where $\mathfrak w=\mathfrak k\ominus\mathfrak h$ and $\mathfrak p=\mathfrak k^\perp$.
\end{lemma}
\begin{proof}
Every derivation of a semisimple Lie algebra is inner, so $\mathfrak n_{\mathfrak g}(\mathfrak s)=\mathfrak s\oplus\mathfrak c_{\mathfrak g}(\mathfrak s)$. The ideal assumption gives $\mathfrak h\subset\mathfrak k$. A connected compact semisimple subgroup is closed. Thus $K$ is closed, and it is proper because otherwise $\mathfrak s$ would be a nonzero ideal of the simple algebra $\mathfrak g$. Finally, $\mathfrak w=\mathfrak m^{\mathfrak s}$. Indeed, $\mathfrak w\subset\mathfrak c_{\mathfrak g}(\mathfrak s)$, whereas $\mathfrak p$ has no $\mathfrak s$-fixed vectors. Equivariance and self-adjointness prove the asserted preservation.
\end{proof}

\begin{lemma}\label{lem312}
Use the notation of Lemma~\ref{lem311}. Suppose every $G$-g.o. metric on $G/K$ is normal and, for every $T\in\mathfrak w$, there is $x\in\mathfrak p$ such that $\mathfrak k_x\subset T^\perp$. Then every $G$-g.o. metric on $G/H$ is normal.
\end{lemma}
\begin{proof}
By Lemmas~\ref{lem311} and \ref{lem38}, $C=E|_{\mathfrak w}\oplus cI_{\mathfrak p}$. For $Y\in\mathfrak w$ put $T=(E-cI)Y$. The $\mathfrak p$-component of the g.o. equation at $Y+x$ gives
\begin{equation}\label{formula16}
 T-Z\in\mathfrak k_x\quad\text{for some }Z\in\mathfrak h.
\end{equation}
Choose $x$ as in the hypothesis. Since $T\perp\mathfrak h$, taking the inner product of \eqref{formula16} with $T$ gives $Q(T,T)=0$. Thus $E=cI$.
\end{proof}

\begin{theorem}\label{thm313}
Let $0\ne\mathfrak s\triangleleft\mathfrak h\subsetneq\mathfrak g$ be semisimple and put $\mathfrak l=\mathfrak s\oplus\mathfrak c_{\mathfrak g}(\mathfrak s)$, $L=N_G(S)^0$ and $\mathfrak p=\mathfrak l^\perp$. Suppose every $G$-g.o. metric on $G/L$ is normal. If some $x_0\in\mathfrak p$ has an abelian stabilizer $\mathfrak l_{x_0}\subset[\mathfrak l,\mathfrak l]$, then every $G$-g.o. metric on $G/H$ is normal.
\end{theorem}
\begin{proof}
For $T\in\mathfrak l\ominus\mathfrak h$, decompose $T=T_z+T_s$ according to $\mathfrak l=\mathfrak z(\mathfrak l)\oplus[\mathfrak l,\mathfrak l]$. The abelian subalgebra $\mathfrak l_{x_0}$ lies in a Cartan subalgebra of $[\mathfrak l,\mathfrak l]$. Every vector of a compact semisimple Lie algebra is orthogonal to some conjugate of a given Cartan subalgebra. This follows from Kostant's linear convexity theorem \cite{Ko1973}, since the average of a Weyl orbit is zero. Conjugate $x_0$ by an element of the connected derived subgroup of $L$ so that its stabilizer is orthogonal to $T_s$. It is also orthogonal to $T_z$, because it lies in $[\mathfrak l,\mathfrak l]$. Lemma~\ref{lem312} now applies.
\end{proof}

\begin{remark}\label{rem314}
For reductive $\mathfrak l$, abelianness of a principal stabilizer is not sufficient for Theorem~\ref{thm313}. Its projection to $\mathfrak z(\mathfrak l)$ must also vanish. If $\mathfrak l$ is semisimple, this extra condition is automatic.
\end{remark}

\begin{corollary}\label{cor315}
Under the first assumption of Theorem~\ref{thm313}, if $\mathfrak p$ has a free vector for $\mathfrak l$, then every $G$-g.o. metric on $G/H$ is normal. No semisimplicity assumption on $\mathfrak l$ is required.
\end{corollary}
\begin{proof}
The zero stabilizer satisfies the second assumption of the theorem. Alternatively, apply Lemma~\ref{lem39} after Lemma~\ref{lem311}.
\end{proof}

\begin{lemma}\label{lem316}
Let $\mathfrak l$ be compact reductive and act orthogonally on $V$. If $[\mathfrak l,\mathfrak l]$ has an abelian principal stabilizer on $V$, then $\mathfrak l$ also has an abelian principal stabilizer.
\end{lemma}
\begin{proof}
Choose $x$ in the common principal strata. Since $\mathfrak l/[\mathfrak l,\mathfrak l]$ is abelian, $[\mathfrak l_x,\mathfrak l_x]\subset[\mathfrak l,\mathfrak l]_x$. The right side is abelian by hypothesis. The left side is the derived algebra of a compact Lie algebra and is therefore semisimple. It must be zero.
\end{proof}

\begin{lemma}\label{lem317}
Let $L\subset G$ be closed and connected, where $\mathfrak g$ is compact simple. Put $\mathfrak l=\mathfrak d\oplus\mathfrak z$, where $\mathfrak d=[\mathfrak l,\mathfrak l]$ and $\mathfrak z=\mathfrak z(\mathfrak l)$, and put $\mathfrak p=\mathfrak l^\perp$. Let $D$ be the connected subgroup with Lie algebra $\mathfrak d$. If every $G$-g.o. metric on $G/D$ is normal, then the projection of a principal $L$-stabilizer on $\mathfrak p$ to $\mathfrak z$ is zero.
\end{lemma}
\begin{proof}
At a principal point $x_0$, set $\mathfrak e=P_{\mathfrak z}\mathfrak l_{x_0}$. This subspace is independent of the principal point, because principal stabilizers are conjugate and $L$ acts trivially on its center. We claim that $\mathfrak e\subset P_{\mathfrak z}\mathfrak l_x$ for every $x\in\mathfrak p$. Choose principal points $x_j\to x$ and elements $k_j\in L$ with $\mathfrak l_{x_j}=\operatorname{Ad}(k_j)\mathfrak l_{x_0}$. Pass to a convergent subsequence. For every $U\in\mathfrak l_{x_0}$, the limit of $\operatorname{Ad}(k_j)U$ stabilizes $x$ and has the same central projection as $U$. This proves the claim.

If $\mathfrak e\ne0$, choose a nonzero positive semidefinite self-adjoint map $A\colon\mathfrak z\to\mathfrak z$ with image in $\mathfrak e$. On $\mathfrak d^\perp=\mathfrak z\oplus\mathfrak p$, define $C=(cI+A)|_{\mathfrak z}\oplus cI_{\mathfrak p}$ with $c>0$. This is $D$-equivariant. For $y\in\mathfrak z$ and $x\in\mathfrak p$, the claim provides $Z\in\mathfrak d$ with $Ay-Z\in\mathfrak l_x$. Since $\mathfrak z$ is central in $\mathfrak l$, we obtain $[C(y+x),y+x]=[Ay,x]=[Z,y+x]$. Thus $C$ is a nonscalar $G$-g.o. cometric on $G/D$, contrary to the hypothesis.
\end{proof}

The next criterion will be used when a principal stabilizer is nonabelian.

\begin{lemma}\label{lem318}
Let $K$ be compact connected semisimple and let $S\subset K$ be connected and closed. Suppose that every fundamental representation of the simply connected cover of $K$ has a nonzero vector fixed by the connected lift of $S$. Then every adjoint orbit in $\mathfrak k$ meets $\mathfrak s^\perp$ for every invariant positive inner product.
\end{lemma}
\begin{proof}
Replace $K$ by its simply connected cover. If $V_\lambda$ and $V_\mu$ have nonzero $S$-fixed vectors $w_\lambda,w_\mu$, then their Cartan component $V_{\lambda+\mu}$ also has one. Indeed, for highest-weight vectors $v_\lambda,v_\mu$, the two matrix coefficients $\langle w_\lambda,kv_\lambda\rangle$ and $\langle w_\mu,kv_\mu\rangle$ are nonzero real analytic functions on connected $K$. Their product is not identically zero. Hence the projection of $w_\lambda\otimes w_\mu$ onto the submodule generated by $v_\lambda\otimes v_\mu$ is nonzero and $S$-fixed. The hypothesis therefore gives an $S$-fixed vector in every dominant integral highest-weight module.

Fix such a weight $\lambda$ and $0\ne w\in V_\lambda^S$. On the compact highest-weight projective orbit $K\cdot[v_\lambda]$, maximize $|\langle w,v\rangle|^2/\|v\|^2$. The maximum is positive. This orbit is also the closed complex-group orbit of the highest-weight line. Thus, for $A\in\mathfrak s$, the curve $[\exp(itA)v]$ stays in it. At a maximizing point, the numerator along this curve is constant because $Aw=0$. Differentiating the denominator gives $\langle iAv,v\rangle=0$ for all $A\in\mathfrak s$.

The functional $A\mapsto\langle iAv,v\rangle/\|v\|^2$ traces the coadjoint orbit of $\lambda$, up to the overall sign in the infinitesimal-weight convention. At the highest-weight line this follows from the vanishing of root-operator expectations and the defining Cartan weight, and the general statement follows by equivariance. Hence this coadjoint orbit meets the annihilator of $\mathfrak s$. Positive rational multiples of dominant integral weights are dense in a Weyl chamber. Approximation and compactness of $K$ extend the conclusion to every orbit. Identifying the dual with $\mathfrak k$ by the chosen invariant form finishes the proof.
\end{proof}

\section{The unitary groups}
\label{sec:unified-a}\label{sec:classical-a-complete}
In this section, we classify the $G$-g.o. metrics for the compact special unitary groups. We first exclude irreducible defining representations and then reduce every nonnormal metric to a unique irreducible block of dimension greater than half the defining dimension.

We use skew-Hermitian matrices. The form $-\operatorname{tr}(AB)$ is a positive multiple of $Q$ and may be used in all orthogonality calculations below.

\begin{lemma}
\label{lem41}
Let $V=V_1\otimes\cdots\otimes V_r$, where $r\geq2$ and $d_i=\dim_{\mathbb C}V_i\geq2$, and let $\mathfrak k=\bigoplus_{i=1}^r\mathfrak{su}(V_i) \subset\mathfrak{su}(V)$ be the tensor product embedding.  Every $\SU(V)$-g.o. cometric on $\SU(V)/K$ is normal, where $K$ is the connected image of $\prod_i\SU(V_i)$.
\end{lemma}

\begin{proof}
Let $E_i=\operatorname{Herm}_0(V_i)$ be the real vector space of traceless Hermitian endomorphisms of $V_i$.  For every nonempty subset $S\subset\{1,\ldots,r\}$, define $M_S=i\left(\bigotimes_{i\in S}E_i\right) \otimes\left(\bigotimes_{i\notin S}\mathbb RI_{V_i}\right).$ The factors are always placed in their original order.  These spaces give an orthogonal real decomposition $\mathfrak{su}(V)=\bigoplus_{S\ne\varnothing}M_S, \qquad \mathfrak k=\bigoplus_{|S|=1}M_S, \qquad \mathfrak p=\bigoplus_{|S|\geq2}M_S.$ Each $M_S$ is an irreducible real $\mathfrak k$-module. It is an outer tensor product of adjoint modules, all of real type.  Modules associated with different supports $S$ are inequivalent.  Thus every positive self-adjoint $\mathfrak k$-equivariant cometric on $\mathfrak p$ has the form $C|_{M_S}=c_S I$.

We use the following immediate consequence of the g.o. equation.  If $U,W\subset\mathfrak p$ are $\mathfrak k$-invariant scalar subspaces with unequal eigenvalues, then $[U,W]\subset U+W$.  Indeed, applying the equation to $u+w$ gives $(c_U-c_W)[u,w]=[z,u]+[z,w]\in U+W$.

If $r=2$, there is only one summand, and the assertion follows. Suppose $r\geq3$.  For distinct $i,j,k$, the bracket $[M_{\{i,j\}},M_{\{i,k\}}]$ has a nonzero component in $M_{\{i,j,k\}}$. Choose two noncommuting traceless Hermitian matrices in the $i$th factor and nonzero matrices in factors $j,k$.  Therefore $c_{\{i,j\}}=c_{\{i,k\}}$.  It follows that all two-element supports have the same coefficient $c$.

Let $|S|\geq3$, and choose distinct $i,j\in S$.  Choose a nonzero $A_i\in E_i$, noncommuting $B_j,D_j\in E_j$, and nonzero $A_\ell\in E_\ell$ for $\ell\in S\setminus\{i,j\}$.  Use these matrices to form $x\in M_{\{i,j\}}$ and $y\in M_S$, with the same $A_i$ in both.  Up to a nonzero scalar, their bracket is $A_i^2\otimes[B_j,D_j]\otimes \bigotimes_{\ell\in S\setminus\{i,j\}}A_\ell.$ Its projection onto $M_{S\setminus\{i\}}$ is nonzero, since the scalar component of $A_i^2$ is $d_i^{-1}\operatorname{tr}(A_i^2)I_{V_i}\ne0$. The target support differs from both $S$ and $\{i,j\}$, so the g.o. equation forces $c_S=c_{\{i,j\}}=c$.  Consequently $C=cI$.
\end{proof}

\begin{theorem}
\label{thm42}
Let $H$ be a proper closed subgroup of $\SU(N)$ and assume that its Lie algebra $\mathfrak h$ acts irreducibly on $\mathbb C^N$ as a complex representation.  Then every $\SU(N)$-g.o. Riemannian metric on $\SU(N)/H$ is normal. Equivalently, every positive cometric satisfying the g.o. equation is a scalar endomorphism.
\end{theorem}

\begin{proof}
By Schur's lemma, an element of the center of $\mathfrak h$ acts on $V=\mathbb C^N$ as an imaginary scalar.  Since it is traceless, it is zero.  Thus $\mathfrak h$ is semisimple.  Write $\mathfrak h=\mathfrak h_1\oplus\cdots\oplus\mathfrak h_r, \qquad V=V_1\otimes\cdots\otimes V_r,$ where the $\mathfrak h_i$ are simple and each $V_i$ is a nontrivial complex irreducible $\mathfrak h_i$-module.  In particular, $d_i=\dim_{\mathbb C}V_i\geq2$ and $\mathfrak h_i\subset\mathfrak{su}(V_i)$.

First suppose $r=1$. By \cite[Theorems~1--2 and Table~1]{ChNiNi2023}, the only nonnormal cases with ambient algebra $\mathfrak{su}(N)$ have standard isotropy $\mathfrak{su}(n)\subset\mathfrak{su}(n+p)$, where $n\ge2$ and $1\le p<n$, or $\mathfrak{sp}(n)\subset\mathfrak{su}(2n+1)$, where $n\ge2$. Both defining representations are reducible. Thus the metric is normal. We apply the cited theorem to the simply connected space $\SU(N)/H^0$ and use Corollary~\ref{cor33} for $H$.

Now suppose $r\geq2$.  Introduce the connected tensor subgroup $K$ with Lie algebra $\mathfrak k=\bigoplus_i\mathfrak{su}(V_i), \qquad \mathfrak w_i=\mathfrak{su}(V_i)\ominus\mathfrak h_i, \qquad \mathfrak p=\mathfrak k^{\perp_Q}.$ Then $\mathfrak m=\left(\bigoplus_i\mathfrak w_i\right)\oplus\mathfrak p, \qquad \mathfrak p=\bigoplus_{|S|\geq2}M_S.$ Here $M_S$ is as in Lemma~\ref{lem41}.  Schur's lemma also gives $E_i^{\mathfrak h_i}=0$.  Thus every irreducible real component of $E_i$ is nontrivial for $\mathfrak h_i$.  Consequently, every irreducible component of $\mathfrak w_i$ has support exactly $\{i\}$, while every component of $M_S$ has support exactly $S$. Components with different supports cannot be equivalent.  Since $C$ is $\mathfrak h$-equivariant, it preserves every $\mathfrak w_i$ and also $\mathfrak p$.  This argument permits arbitrary multiplicities and mixing between equivalent components within an individual support.

By Lemma~\ref{lem38}, $C|_{\mathfrak p}$ is a $\mathfrak k$-equivariant $G$-g.o. cometric for $(\mathfrak g,\mathfrak k)$. Lemma~\ref{lem41} now gives $C|_{\mathfrak p}=cI$ for some $c>0$.

It remains to determine $C|_{\mathfrak w_i}$.  Let $0\ne y\in\mathfrak w_i$ be a $C$-eigenvector of eigenvalue $a$. The centralizer of $\mathfrak h_i$ in $\mathfrak{su}(V_i)$ is zero, so there is $u\in\mathfrak h_i$ with $[y,u]\ne0$.  Choose $j\ne i$ and $0\ne v\in\mathfrak h_j$, and put $x=i\,u\otimes v$, with identity factors elsewhere.  Thus $x\in M_{\{i,j\}}\subset\mathfrak p$. The factor $i$ ensures that $x$ is skew-Hermitian.  Applying the g.o. equation to $y+x$ yields $(a-c)[y,x]=[z,y]+[z,x]\qquad(z\in\mathfrak h).$ The nonzero vector $[y,x]$ belongs to the real subspace $i\,\mathfrak w_i\otimes\mathfrak h_j$ of $M_{\{i,j\}}$. On the right, $[z,y]\in\mathfrak w_i$, whereas $[z,x]\in i\,\mathfrak h_i\otimes\mathfrak h_j$. These spaces are orthogonal to $i\,\mathfrak w_i\otimes\mathfrak h_j$. Projection therefore gives $a=c$.  This applies to every eigenvector of each $C|_{\mathfrak w_i}$, proving $C=cI$ on $\mathfrak m$.
\end{proof}

The irreducibility hypothesis concerns $H^0$. We next allow arbitrary invariant summands, including equivalent modules and coupled block projections.

\begin{theorem}\label{thm43}
Let $H\subset \SU(N)$ be a closed subgroup, and suppose that its Lie algebra preserves an orthogonal decomposition $\mathbb C^N=V_1\oplus\cdots\oplus V_r, \qquad r\geq2,\qquad d_i=\dim_{\mathbb C}V_i.$ The summands need not be irreducible or pairwise inequivalent, and the projections of $\mathfrak h$ onto the different blocks may be coupled.  If $\max_i d_i\leq \frac N2,$ then every $\SU(N)$-g.o. metric on $\SU(N)/H$ is normal.  The conclusion also holds when $H$ is disconnected.
\end{theorem}

\begin{proof}
Work with skew-Hermitian matrices and $Q(A,B)=-\operatorname{tr}(AB)$. Put $K=S\bigl(\mathrm \U(V_1)\times\cdots\times\mathrm \U(V_r)\bigr), \qquad \mathfrak f=\mathfrak k\ominus\mathfrak h, \qquad \mathfrak p=\mathfrak k^\perp.$ The space $\mathfrak p$ consists of the off-diagonal blocks, and $\mathfrak m=\mathfrak f\oplus\mathfrak p$.  Let $C$ be a $G$-g.o. cometric.  The block-scalar torus of $K$ centralizes $H^0$. By Corollary~\ref{cor32}, $C$ commutes with this torus. Its zero-weight space is $\mathfrak f$, whereas every off-diagonal block has a nonzero weight. Hence $C$ preserves $\mathfrak f\oplus\mathfrak p$. By Lemma~\ref{lem38}, $C|_{\mathfrak p}$ is $K$-equivariant and is a $\SU(N)$-g.o. cometric on $\SU(N)/K$.

For $i<j$, let $\mathfrak p_{ij}$ denote the real off-diagonal block corresponding to $\operatorname{Hom}_{\mathbb C}(V_j,V_i)$. These are mutually inequivalent irreducible real $K$-modules. Their block-scalar torus weights distinguish the unordered pairs, and the unitary block factors give irreducibility.  Thus $C|_{\mathfrak p_{ij}}=b_{ij}I$.  If $r=2$, write $b=b_{12}$. If $r\geq3$, choose distinct $i,j,k$ and $x\in\mathfrak p_{ij}$, $y\in\mathfrak p_{jk}$ with $0\ne[x,y]\in\mathfrak p_{ik}$.  Such matrices exist by choosing one matrix unit in each block.  The $\mathfrak p_{ik}$ component of the g.o. equation for $x+y$ is $(b_{ij}-b_{jk})[x,y]=0,$ because $[\mathfrak k,x+y]\subset \mathfrak p_{ij}\oplus\mathfrak p_{jk}$.  Applying this to all triangles gives a common scalar.  In either case, $C|_{\mathfrak p}=bI.$

Set $D=(C-bI)|_{\mathfrak f}$.  For any $Y\in\mathfrak f$ and $X\in\mathfrak p$, the off-diagonal component of the g.o. equation for $Y+X$ gives a $Z\in\mathfrak h$ such that
\begin{equation}\label{formula17}
 [DY-Z,X]=0.
\end{equation}
It remains to show that this forces $D=0$.

First suppose $r\geq3$.  Relabel so that $d_r=b_0=\max_i d_i$. Then $d_i\leq b_0$ for $i<r$ and $\sum_{i<r}d_i\geq b_0$. For each $i<r$, choose an orthonormal family $v_{i1},\ldots,v_{id_i}\in\mathbb C^{b_0}$ such that the union of these families spans $\mathbb C^{b_0}$ and
\begin{equation}\label{formula18}
 \langle v_{ia},v_{jc}\rangle\ne0
 \qquad(i\ne j).
\end{equation}
These simultaneous choices exist.  The product of the complex Stiefel manifolds of the families is connected and real analytic. Each specified cross inner product is not identically zero. Full spanning is also possible. Choose coordinate orthonormal families whose union covers the $b_0$ coordinate vectors, using $\sum_{i<r}d_i\geq b_0$.  Hence the full-spanning condition and the finitely many nonvanishing conditions in \eqref{formula18} are open dense conditions in this product, and have nonempty intersection.

Choose positive, pairwise distinct numbers $s_{i1},\ldots,s_{id_i}$ for each $i$, and let $X_i$ be the $d_i\times b_0$ matrix whose $a$-th row is $s_{ia}v_{ia}^*$.  Define $X\in\mathfrak p$ by $X_{ir}=X_i$, $X_{ri}=-X_i^*$, with all other blocks zero. We claim its infinitesimal stabilizer in $\mathfrak k$ is zero. Indeed, if $A=\operatorname{diag}(A_1,\ldots,A_r)\in\mathfrak k$ commutes with $X$, then $A_iX_i=X_iA_r \qquad(i<r).$ Since $A_i$ and $A_r$ are skew-Hermitian, $A_i$ commutes with $X_iX_i^*=\operatorname{diag}(s_{ia}^2)$.  Distinctness gives $A_i=\operatorname{diag}(\mathrm i\theta_{ia})$, and then $A_rv_{ia}=\mathrm i\theta_{ia}v_{ia}$. Skew-Hermitian matrices have orthogonal eigenspaces for different eigenvalues, so \eqref{formula18} forces all these eigenvalues to agree.  There are at least two leaf blocks, so this also connects any two vectors in the same family.  Since the vectors span the central block, $A_r=\mathrm i\theta I$, and $A_i=\mathrm i\theta I$ for every $i$.  The total trace condition gives $N\mathrm i\theta=0$, proving the claim. Use this fixed $X$ in \eqref{formula17}.  It follows that $DY=Z$. Orthogonality of $\mathfrak f$ and $\mathfrak h$ gives $DY=0$ for every $Y$.

It remains to consider $r=2$.  The dimension hypothesis forces $d_1=d_2=d$.  Fix $Y\in\mathfrak f$ and write $T=DY=\operatorname{diag}(T_1,T_2)$.  Choose unitary matrices $U,V$ such that the diagonals of $U^*T_1U$ and $V^*T_2V$ are constant.  To obtain such a basis, diagonalize the traceless part of the skew-Hermitian matrix and then conjugate by the unitary Fourier matrix. All diagonal entries of the traceless part become their mean, which is zero. Let $S=\operatorname{diag}(s_1,\ldots,s_d)$ have distinct positive entries, and choose the off-diagonal block of $X\in\mathfrak p$ to be $USV^*$.  Its stabilizer in $\mathfrak k$ is
\begin{equation*}
\mathfrak k_X=
 \left\{
 \operatorname{diag}
 \bigl(U\operatorname{diag}(\mathrm i\theta_a)U^*,
       V\operatorname{diag}(\mathrm i\theta_a)V^*\bigr):
       \sum_a\theta_a=0
 \right\}.
\end{equation*}
The constant-diagonal choices imply $\mathfrak k_X\subset T^\perp$. By \eqref{formula17}, $T-Z\in\mathfrak k_X$ for some $Z\in\mathfrak h$.  Since $T\perp\mathfrak h$, pairing with $T$ gives $Q(T,T)=0$.  Thus again $D=0$, and $C=bI$ on all of $\mathfrak m$.

\end{proof}

\begin{corollary}\label{cor44}
If a closed subgroup $H\subset\SU(N)$ admits a nonnormal $\SU(N)$-g.o. metric, then the defining $\mathfrak h$-representation is reducible and has a unique irreducible summand of complex dimension greater than $N/2$.
\end{corollary}

\begin{proof}
The irreducible case is excluded by Theorem~\ref{thm42}.  Apply Theorem~\ref{thm43} to an orthogonal decomposition into irreducible $\mathfrak h$-modules.  At least one summand must have dimension greater than $N/2$, and two such summands cannot occur.
\end{proof}

The following compression argument does not require the cometric to preserve an individual internal complement.

\begin{proposition}\label{prop45}
Let $H\subset\SU(N)$ be closed, and decompose the defining $\mathfrak h$-representation into orthogonal irreducible complex summands $\mathbb C^N=V_1\oplus\cdots\oplus V_r.$ The summands may be equivalent, and the block projections of $\mathfrak h$ need not be independent.  Let $\mathfrak h_i\subset\mathfrak u(V_i)$ be the image of $\mathfrak h$ on $V_i$, and put $\mathfrak s_i=[\mathfrak h_i,\mathfrak h_i], \qquad \mathfrak w_i=\mathfrak{su}(V_i)\ominus\mathfrak s_i.$ Regard $\mathfrak w_i$ as a single diagonal block in $\mathfrak{su}(N)$.  Then $\mathfrak w_i\subset\mathfrak m$. For every $\SU(N)$-g.o. cometric $C$, its orthogonal compressions satisfy
\begin{align}
 P_{\mathfrak w_i}C|_{\mathfrak w_i}&=a_iI
       &&\text{for some }a_i>0\text{ if }\mathfrak w_i\ne0,
       \label{formula19}\\
 P_{\mathfrak w_i}C|_{\mathfrak w_j}&=0
       &&(i\ne j).
       \label{formula20}
\end{align}
The statement does not assert that $C$ preserves each $\mathfrak w_i$.
\end{proposition}

\begin{proof}
Use skew-Hermitian matrices and the trace inner product.  Since $V_i$ is irreducible, the center of $\mathfrak h_i$ acts by scalar matrices.  Thus $\mathfrak h_i\subset\mathfrak s_i+\mathrm i\mathbb RI_{V_i},$ and $\mathfrak w_i$ is orthogonal to the image of $\mathfrak h$ in its block.  This proves $\mathfrak w_i\subset\mathfrak m$.  It is an $\mathfrak h$-invariant subspace.

Let $A_i=P_{\mathfrak w_i}C|_{\mathfrak w_i}$.  This compression is positive, self-adjoint, and $\mathfrak h_i$-equivariant. For $Y\in\mathfrak w_i$, write the $i$-th diagonal block of $CY$ as $(CY)_{ii}=A_iY+L_i(Y)+\mathrm i\alpha(Y)I, \qquad L_i(Y)\in\mathfrak s_i.$ Other blocks of $CY$ do not contribute to the $i$-th diagonal block of $[CY,Y]$.  Projecting the g.o. equation for $Y$ onto this block gives $[A_iY,Y]=[Z_i^{\mathrm{ss}}-L_i(Y),Y]$ for some $Z\in\mathfrak h$, where $Z_i^{\mathrm{ss}}\in\mathfrak s_i$ is the semisimple part of its $i$-th block.  Hence $A_i$ is a $\SU(V_i)$-g.o. cometric on $\SU(V_i)/S_i$, where $S_i$ is the connected compact subgroup with Lie algebra $\mathfrak s_i$.

If $\dim_{\mathbb C}V_i\geq2$, the defining representation of $\mathfrak s_i$ is irreducible. Adjoining the scalar center does not change invariant subspaces.  The irreducible defining exclusion theorem, Theorem~\ref{thm42}, therefore gives $A_i=a_iI$.  If $\dim_{\mathbb C}V_i=1$, then $\mathfrak w_i=0$, so no assertion is needed.  This proves \eqref{formula19}.

Fix $i\ne j$, $Y\in\mathfrak w_i$, and $W\in\mathfrak w_j$, and put $T=P_{\mathfrak w_i}CW$.  Apply the g.o. equation to $Y+W$ and project first onto its $i$-th diagonal block and then onto $\mathfrak s_i$.  All brackets between $\mathfrak s_i$ and $\mathfrak w_i$ belong to $\mathfrak w_i$, scalar parts commute, and $[P_{\mathfrak w_i}CY,Y]=[a_iY,Y]=0$.  The other diagonal block $W$ does not enter this projection.  We obtain $P_{\mathfrak s_i}[T,Y]=0 \qquad\text{for every }Y\in\mathfrak w_i.$ Pairing with $A\in\mathfrak s_i$ gives $Q(T,[Y,A])=0$.  On the other hand, $\operatorname{span}[\mathfrak s_i,\mathfrak w_i] =\mathfrak w_i.$ Indeed, its orthogonal complement in $\mathfrak w_i$ is the fixed subspace $\mathfrak w_i^{\mathfrak s_i}$, which is zero by complex irreducibility and Schur's lemma. A traceless skew-Hermitian matrix commuting with $\mathfrak s_i$ is zero. It follows that $T=0$, proving \eqref{formula20}.
\end{proof}

Although the compressions between different $\mathfrak w_i$ vanish, mixing between an internal complement and the remaining diagonal space has not yet been excluded. We retain these terms in the next argument.

\begin{lemma}
\label{lem46}
Let $\mathfrak t\subset\mathfrak{su}(d)$ be abelian, and let $1\le q<d$. There exists a complex $q$-plane $E\subset\mathbb C^d$ whose infinitesimal stabilizer in $\mathfrak t$ is zero.
\end{lemma}
\begin{proof}
Simultaneously diagonalize $\mathfrak t$. Choose a rank-$q$ Hermitian orthogonal projection $P$ all of whose off-diagonal entries are nonzero. Such projections form a nonempty open dense subset of the complex Grassmannian. Each condition $P_{ij}=0$ is a proper real-analytic condition, since one can choose a $q$-plane containing $(e_i+e_j)/\sqrt2$ and $q-1$ coordinate vectors different from $e_i,e_j$. There are only finitely many pairs, and the Grassmannian is connected. If a diagonal skew-Hermitian matrix $A=\operatorname{diag}(ia_j)$ preserves $E=\operatorname{im}P$, it also preserves $E^\perp$, and therefore $[A,P]=0$. Thus $(a_i-a_j)P_{ij}=0$ for every $i\ne j$, so all $a_i$ are equal. The trace-zero condition gives $A=0$. This proof does not require $\mathfrak t$ to exponentiate to a closed torus.
\end{proof}

\begin{lemma}
\label{lem47}
Let $\mathfrak s\subsetneq\mathfrak{su}(d)$ be semisimple and complex-irreducible on $\mathbb C^d$. Unless $d=2n$ and $\mathfrak s=\mathfrak{sp}(n)$ in its standard representation, the stabilizer $\mathfrak s_Y=\{A\in\mathfrak s:[A,Y]=0\}$ is abelian for a dense set of $Y\in\mathfrak w=\mathfrak{su}(d)\ominus\mathfrak s$.
\end{lemma}
\begin{proof}
Apply \cite[Theorems~1.6 and~3.1, Table~1]{BeKo2021} to $\mathfrak s_{\mathbb C}\subset\mathfrak{sl}_d(\mathbb C)$. The nonabelian-generic-stabilizer alternatives in type A are standard block algebras containing a dominant $\mathfrak{sl}_p$ factor, and the standard $\mathfrak{sp}_{d/2}\subset\mathfrak{sl}_d$. The block alternatives act reducibly on $\mathbb C^d$. Thus only the stated symplectic exception can remain. The passage to the compact real form follows from their Lemma~2.12(1). The compact representation has reductive stabilizers, and abelian generic stabilizers are equivalent for the real representation and its complexification.
\end{proof}

\begin{theorem}
\label{thm48}
Let $d>q\ge1$, and let $H\subset S(\U(d)\times \U(q))\subset \SU(d+q)$ be closed. Assume that the projection of $\mathfrak h$ on the first block acts complex-irreducibly. Let $\mathfrak s$ be its semisimple part, and assume $\mathfrak s\ne\mathfrak{su}(d),\qquad \mathfrak s\ne\mathfrak{sp}(d/2) \quad\text{in the standard representation}.$ Then every $\SU(d+q)$-g.o. metric on $\SU(d+q)/H$ is normal. The projection of $\mathfrak h$ to the two blocks may be coupled, and the action on the smaller block is arbitrary.
\end{theorem}
\begin{proof}
Use skew-Hermitian matrices and $Q(A,B)=-\operatorname{tr}(AB)$. Put
\begin{equation*}
\mathfrak k=\mathfrak s(\mathfrak u(d)\oplus\mathfrak u(q)),
 \quad\mathfrak f=\mathfrak k\ominus\mathfrak h,
 \quad\mathfrak p=\mathfrak k^\perp,
 \quad\mathfrak w=\mathfrak{su}(d)\ominus\mathfrak s.
\end{equation*}
The space $\mathfrak w$, supported only on the first block, is contained in $\mathfrak f$ and is $\mathfrak h$-invariant. Let $C$ be a positive self-adjoint $G$-g.o. cometric. The block-scalar circle centralizes $\mathfrak h$, so normalizer invariance makes $C$ preserve $\mathfrak f$ and $\mathfrak p$. Descent to the complex Grassmannian $\SU(d+q)/K$ gives $C=E\oplus cI_{\mathfrak p},\qquad D=E-cI_{\mathfrak f}.$ In particular, $D$ is self-adjoint on $\mathfrak f$. For $y\in\mathfrak f$ and $P(X)=\left(\begin{smallmatrix}0&X\\-X^*&0\end{smallmatrix}\right) \in\mathfrak p$, the two components of the g.o. equation are
\begin{equation}\label{formula21}
 [Z-Dy,y]=0,\qquad
 (Dy-Z)_{11}X=X(Dy-Z)_{22}
 \quad\text{for some }Z\in\mathfrak h.
\end{equation}
The same $Z$ occurs in both equations.

We use the already proved compression fact $P_{\mathfrak w}C|_{\mathfrak w}=aI \quad(a>0).$ For clarity, it holds without $C$ preserving $\mathfrak w$. Project the g.o. equation at $Y\in\mathfrak w$ onto its first block, absorb the first-block $\mathfrak s$ component of $CY$ into the compensator, and obtain a $\SU(d)$-g.o. cometric on $\SU(d)/S$. The irreducible-defining-representation exclusion theorem makes that compression scalar. Consequently, for every $Y\in\mathfrak w$, $(DY)_{11}=(a-c)Y+L(Y)+i\alpha(Y)I_d, \qquad L(Y)\in\mathfrak s.$

First we prove $D\mathfrak w=0$, including all possible mixed components of $DY$. Choose nonzero $Y\in\mathfrak w$ with $\mathfrak s_Y$ abelian, as provided by Lemma~\ref{lem47}. Write $Z_{11}=A+i\eta I_d$, with $A\in\mathfrak s$. The first equation in \eqref{formula21}, for $y=Y$, gives $A'=A-L(Y)\in\mathfrak s_Y.$ The algebra $\mathfrak t_Y=\mathfrak s_Y+\mathbb RY\subset\mathfrak{su}(d)$ is abelian. By Lemma~\ref{lem46}, choose a $q$-plane $E_Y$ with zero infinitesimal $\mathfrak t_Y$-stabilizer, and a full-column-rank $X$ whose image is $E_Y$. The second equation in \eqref{formula21} implies that $(a-c)Y-A'+i(\alpha(Y)-\eta)I_d$ preserves $E_Y$. Its traceless part lies in $\mathfrak t_Y$, and hence vanishes. Since $Y\perp\mathfrak s$ and $Y\ne0$, we obtain $a=c$ and $A'=0$. Now $(DY-Z)_{11}=i\delta I_d$. Because $X$ has full column rank, the horizontal equation gives $(DY-Z)_{22}=i\delta I_q$. Total trace zero gives $(d+q)\delta=0$, hence $DY=Z$. But $DY\in\mathfrak f$ and $Z\in\mathfrak h$, so $DY=0$. The set of such $Y$ is dense and $D$ is linear. Therefore
\begin{equation}\label{formula22}
 D\mathfrak w=0.
\end{equation}
Self-adjointness also gives $D(\mathfrak f)\subset\mathfrak f\cap\mathfrak w^\perp$. Thus all previously possible mixing with $\mathfrak w$ has been eliminated.

It remains to kill the full diagonal remainder. Fix $r\in\mathfrak f\cap\mathfrak w^\perp$. Both $r_{11}$ and $(Dr)_{11}$ lie in $\mathfrak s+i\mathbb RI_d$. Choose a fixed $Y\in\mathfrak w$ with abelian $\mathfrak s_Y$, and apply \eqref{formula21} to $y=r+Y$. By \eqref{formula22}, $Dy=Dr$. Write $(Z-Dr)_{11}=A+i\delta I_d, \qquad A\in\mathfrak s.$ Projecting the first-block vertical equation $[(Z-Dr)_{11},r_{11}+Y]=0$ onto $\mathfrak w$ gives $[A,Y]=0$. The bracket with $r_{11}\in\mathfrak s+i\mathbb RI_d$ lies in $\mathfrak s$, whereas $[A,Y]\in\mathfrak w$. Thus $A\in\mathfrak s_Y$. Choose a full-column-rank $X$ whose image has zero $\mathfrak s_Y$-stabilizer. The horizontal equation forces $A=0$, then forces $(Z-Dr)_{22}=i\delta I_q$. Total trace zero gives $\delta=0$ and $Dr=Z$. Orthogonality gives $Dr=0$. Since $r$ was arbitrary, $D=0$ and $C=cI$.
\end{proof}

It remains to determine the metrics when the dominant block contains a full special unitary or standard symplectic factor. We treat the special unitary case first.

Let $N=d+q$ with $d>q\ge1$, and let $H$ be any closed subgroup satisfying $\SU(d)\subset H\subset K=S\bigl(\U(d)\times \U(q)\bigr)\subset\SU(N),$ where $\SU(d)$ acts on the first summand of $\mathbb C^N=\mathbb C^d\oplus\mathbb C^q$. Write $\mathfrak g=\mathfrak k\oplus\mathfrak p,\qquad \mathfrak k=\mathfrak h\oplus\mathfrak v,\qquad \mathfrak m=\mathfrak v\oplus\mathfrak p,$ with orthogonality taken for $Q=-B_{\mathfrak{su}(N)}$. Identify $\mathfrak p$ with $M_{d,q}(\mathbb C)$ by
\begin{equation*}
P(X)=\begin{pmatrix}0&X\\-X^*&0\end{pmatrix}.
\end{equation*}
Set $z=i\operatorname{diag}(qI_d,-dI_q),\qquad h_0=P_{\mathfrak h}z,\qquad w=P_{\mathfrak m}z=z-h_0\in\mathfrak v.$

\begin{theorem}\label{thm49}
If $w=0$, every $\SU(N)$-g.o. metric on $\SU(N)/H$ is normal.  If $w\ne0$, its complete family of $G$-g.o. cometrics is
\begin{equation}\label{formula23}
 C=bI+(a-b)P_{\mathbb Rw},\qquad a,b>0.
\end{equation}
Here the metric represented by $C$ is $g(u,v)=Q(C^{-1}u,v)$.

\end{theorem}

\begin{proof}
We use the g.o. equation $[Cu,u]=[Z,u],\qquad Z\in\mathfrak h,$ for every $u\in\mathfrak m$.

Since $z$ centralizes $\mathfrak h$, we have $h_0\in\mathfrak z(\mathfrak h)$ and $w\in\mathfrak m^{\mathfrak h}$. By Corollary~\ref{cor32}, $C$ commutes with $\operatorname{ad}(z)$. This operator vanishes on $\mathfrak v$ and acts as the nonzero scalar complex structure $iN$ on $\mathfrak p$. Hence $C$ preserves $\mathfrak v\oplus\mathfrak p$. By Lemma~\ref{lem38} and the irreducibility of the complex Grassmannian, $C=E\oplus bI_{\mathfrak p}$ for some $b>0$.

Introduce the linear subspace $\mathcal L=\{T\in\mathfrak k: [T,P(X)]\in[\mathfrak h,P(X)] \text{ for every }X\in M_{d,q}(\mathbb C)\}.$ It is a linear subspace because, for each fixed $X$, the condition is membership in the inverse image of a vector subspace under a linear map.  It contains $\mathfrak h$. It also contains $z$.  Indeed, for $F=\operatorname{im}X$ and $k=\dim F\le q<d$, the matrix $A_X=iN P_F-iN\frac{k}{d-k}(I_d-P_F)\in\mathfrak{su}(d)$ satisfies $A_XX=iNX$ and hence $[A_X,P(X)]=[z,P(X)]$. It follows that $w=z-h_0\in\mathcal L$.

For $Y\in\mathfrak v$, put $T=(E-bI)Y$. The off-diagonal part of the g.o. equation at $Y+P(X)$ gives $[T,P(X)]=[Z,P(X)]$ for some $Z\in\mathfrak h$. Therefore
\begin{equation}\label{formula24}
 \operatorname{im}(E-bI)\subset\mathfrak v\cap\mathcal L.
\end{equation}
We claim
\begin{equation}\label{formula25}
 \mathfrak v\cap\mathcal L=\mathbb Rw.
\end{equation}
The right side is interpreted as zero when $w=0$.

Let $T\in\mathfrak v\cap\mathcal L$.  Orthogonality to $\mathfrak{su}(d)\subset\mathfrak h$ makes its first diagonal block a scalar $itI_d$.  Since $\operatorname{tr}T=0$, $Q(T,z)$ is a positive constant times $dNt$. Moreover $Q(T,z)=Q(T,w)$ because $T\perp\mathfrak h$. If $w\ne0$, subtract the orthogonal projection of $T$ onto $\mathbb Rw$. If $w=0$, leave $T$ unchanged.  The resulting $T_0\in\mathfrak v\cap\mathcal L$ has the form $T_0=\operatorname{diag}(0,B),\qquad B\in\mathfrak{su}(q).$

There is an orthonormal basis of $\mathbb C^q$ in which all diagonal entries of $B$ vanish.  For example, diagonalize the traceless Hermitian matrix $-iB$ and then apply the discrete Fourier unitary.  Choose a positive definite Hermitian matrix $S$ with distinct eigenvalues and with this eigenbasis, and choose a full-column-rank $X$ satisfying $X^*X=S$. For $q=1$, the conclusion $B=0$ is immediate and this construction is unnecessary.

Since $T_0\in\mathcal L$, some $Z\in\mathfrak h$ satisfies $[T_0-Z,P(X)]=0$.  Write its two skew-Hermitian diagonal blocks as $K_1,K_2$.  The intertwining identity $K_1X=XK_2$ and its adjoint imply $[K_2,S]=0.$ Thus $K_2$ is diagonal in the selected eigenbasis of $S$. As $B$ has zero diagonal there, $Q(T_0,T_0-Z)=Q(\operatorname{diag}(0,B), \operatorname{diag}(K_1,K_2))=0.$ But $T_0\perp\mathfrak h$, so the same expression is $Q(T_0,T_0)$.  Hence $T_0=0$, proving \eqref{formula25}. By \eqref{formula24} and self-adjointness of $E$, the only possibilities are exactly the operators in \eqref{formula23}, or $C=bI$ when $w=0$.

Suppose $C$ has the stated form and write $u=Y+P(X)$ with $Y\in\mathfrak v$.  There is a real number $\lambda$ such that $(E-bI)Y=\lambda w$. Let $F=\operatorname{im}X$ and $k=\dim F<d$, and define $A=i\lambda N P_F-i\lambda N\frac{k}{d-k}(I_d-P_F) \in\mathfrak{su}(d)\subset\mathfrak h, \qquad Z=A-\lambda h_0\in\mathfrak h.$ The first block of $Y$ is scalar, so $[A,Y]=0$. Also $[A,P(X)]=\lambda[z,P(X)]$ and $[z,Y]=0$. Using $w=z-h_0$, we obtain $[Z,Y]=\lambda[w,Y],\qquad [Z,P(X)]=\lambda[w,P(X)].$ Their sum is precisely $[Z,u]=[Cu,u]$. The construction works also when $X=0$ or has deficient rank. By Corollary~\ref{cor33}, these cometrics also descend to every closed isotropy extension.
\end{proof}

We now consider the standard symplectic projection. The calculation below also excludes arbitrary residual isotropy when the complementary dimension is at least two.

\begin{theorem}\label{thm410}
Let $n\ge2$, $1\le q<2n$, and let $H\subset S(\U(2n)\times \U(q))$ have Lie algebra $\mathfrak h=\mathfrak{sp}(n)\oplus\mathfrak l, \qquad \mathfrak l\subset\{\operatorname{diag}(i\eta I_{2n},B):2ni\eta+\operatorname{tr}B=0\}.$ If $q\ge2$, every $\SU(2n+q)$-g.o. metric is normal. If $q=1$, put $T=i\operatorname{diag}(I_{2n},-2n)$. For $\mathfrak l=0$, all $G$-g.o. cometrics have three independent positive coefficients on $\mathfrak{su}(2n)\ominus\mathfrak{sp}(n)$, $(\mathbb C^{2n})_{\mathbb R}$ and $\mathbb RT$. For $\mathfrak l=\mathbb RT$, the last summand is absent and the other two positive coefficients are independent.
\end{theorem}
\begin{proof}
Put $\mathfrak k=\mathfrak s(\mathfrak u(2n)\oplus\mathfrak u(q))$, $\mathfrak w=\mathfrak{su}(2n)\ominus\mathfrak{sp}(n)$ and $\mathfrak p=\mathfrak k^\perp$. Write $\mathfrak k\ominus\mathfrak h=\mathfrak w\oplus\mathfrak v_0$. The first block of every element of $\mathfrak v_0$ is scalar. The block-scalar circle centralizes $\mathfrak h$, so $C$ preserves $\mathfrak k\ominus\mathfrak h$ and $\mathfrak p$. The factor $\Sp(n)$ acts trivially on $\mathfrak v_0$ and has no fixed vectors in $\mathfrak w$. Thus $C$ also preserves $\mathfrak w$ and $\mathfrak v_0$. By descent and the irreducibility of $\SU(2n)/\Sp(n)$ and the complex Grassmannian, we have $C=aI_{\mathfrak w}\oplus E_{\mathfrak v_0}\oplus bI_{\mathfrak p},\qquad a,b>0.$

Suppose $q\ge2$. Choose $Y=i\operatorname{diag}(y_1I_2,\ldots,y_nI_2)\in\mathfrak w$, where the $y_j$ are distinct and have zero sum. A compensator at $Y+P(X)$ has first block $A+i\eta I_{2n}$, where $A\in\mathfrak{sp}(n)$ commutes with $Y$. Hence $A=\operatorname{diag}(A_1,\ldots,A_n)$ with $A_j\in\mathfrak{su}(2)$. Write $X_j$ for the two-row blocks of $X$. The horizontal equation gives
\begin{equation}\label{formula26}
 D_jX_j=X_jB,\qquad D_j=A_j-i\bigl((a-b)y_j-\eta\bigr)I_2,\qquad B\in\mathfrak u(q).
\end{equation}
Choose full-row-rank $X_j$ with distinct singular values such that all entries of $X_jX_k^*$ are nonzero in their Gram eigenbases for $j\ne k$, and the joint row span is $\mathbb C^q$. Such choices exist. We may choose independent orthonormal two-row frames, then multiply their rows by distinct positive constants. Every forbidden cross inner product defines a proper real-analytic subset of the product of Stiefel manifolds. The full-span condition is nonempty because $q<2n$, and it is open dense. Thus all these conditions hold simultaneously.

Taking adjoints in \eqref{formula26}, we obtain $[D_j,X_jX_j^*]=0$ and $D_jX_jX_k^*=X_jX_k^*D_k$. The simple Gram spectra make the $D_j$ diagonal in their respective Gram eigenbases. The nonzero entries of every cross Gram matrix force all their eigenvalues to coincide. Comparing traces, we find that $(a-b)y_j-\eta$ is independent of $j$. Since the $y_j$ are distinct, $a=b$.

Fix $R_0\in\mathfrak v_0$ and set $(E-bI)R_0=T_0=\operatorname{diag}(itI_{2n},T_q)$. Apply the g.o. equation at $Y+R_0+P(X)$ with the same $Y$ and $X$. The first-block vertical equation again gives $A\in\mathfrak{sp}(n)_Y$. If $Z=\operatorname{diag}(A+i\eta I_{2n},Z_q)$ is the compensator and $\delta=t-\eta$, the horizontal equation becomes $(i\delta I_2-A_j)X_j=X_j(T_q-Z_q)\qquad(1\le j\le n).$ The preceding Gram calculation makes all $i\delta I_2-A_j$ equal scalar matrices. Since $\operatorname{tr}A_j=0$, we have $A_j=0$. The joint row span then gives $T_q-Z_q=i\delta I_q$. Taking the total trace of $T_0-Z$, we obtain $(2n+q)i\delta=0$. Therefore $T_0=Z$. Since $T_0\perp\mathfrak h$, we have $T_0=0$. Thus $E=bI$, proving normality for $q\ge2$.

Suppose $q=1$. Then $\mathfrak l$ is either zero or $\mathbb RT$. The internal module and the horizontal module are inequivalent real irreducible $\Sp(n)$-modules. The central line, when present, is trivial. Hence every invariant positive self-adjoint operator has the scalar form stated in the theorem.

We prove sufficiency directly. Let the three coefficients be $a,b,c$ and write a tangent vector as $Y+P(v)+tT$. Diagonalize $Y$ by $\Sp(n)$, so $Y=i\operatorname{diag}(y_1I_2,\ldots,y_nI_2)$ with $\sum_jy_j=0$, and write $v=(v_1,\ldots,v_n)^t$. Put $\delta_j=(a-b)y_j+(c-b)t(2n+1)$. For $v_j\ne0$, set $A_j=i\delta_j(2P_{\mathbb Cv_j}-I_2)\in\mathfrak{su}(2)$, and put $A_j=0$ for $v_j=0$. Then $A=\operatorname{diag}(A_1,\ldots,A_n)\in\mathfrak{sp}(n)$ commutes with $Y$ and satisfies $Av=(a-b)Yv+i(c-b)t(2n+1)v$. Consequently $[A,Y+P(v)+tT]=[C(Y+P(v)+tT),Y+P(v)+tT]$. The construction also applies when eigenvalues repeat or some $v_j$ vanish. If the circle belongs to the isotropy, omit the term $tT$. By Corollary~\ref{cor33}, these formulas apply to every closed isotropy extension.
\end{proof}

\begin{theorem}
\label{thm411}
Let $H$ be a proper closed subgroup of $\SU(N)$, $N\ge2$. Every $\SU(N)$-g.o. metric is normal except for the following families, up to conjugacy and finite isotropy extensions.
\begin{enumerate}
\item There is a splitting $\mathbb C^N=\mathbb C^d\oplus\mathbb C^q$,
$d>q\ge1$, with $\SU(d)\subset H\subset S(\U(d)\times \U(q)).$ Let $Q=-B_{\mathfrak{su}(N)}$, $\mathfrak m=\mathfrak h^{\perp_Q}$, $z=i\operatorname{diag}(qI_d,-dI_q),\qquad w=P_{\mathfrak m}z.$ If $w=0$, only the normal metric occurs. Otherwise all $G$-g.o. cometrics, and no others, are $C=bI+(a-b)P_{\mathbb Rw},\qquad a,b>0.$ The action on the $q$-dimensional complement is unrestricted.
\item $N=2n+1$, $n\ge2$, and the identity-component isotropy
algebra is the standard $\mathfrak{sp}(n)$, or $\mathfrak{sp}(n)\oplus\mathbb RT$, where $T=i\operatorname{diag}(I_{2n},-2n)$. For $\mathfrak h=\mathfrak{sp}(n)$ all positive parameters on the three pairwise inequivalent summands $\mathfrak{su}(2n)\ominus\mathfrak{sp}(n),\qquad (\mathbb C^{2n})_{\mathbb R},\qquad\mathbb RT$ are independent. For the circle extension the last summand is absent and both remaining positive parameters are independent.
\end{enumerate}
In each family, normality means that all displayed coefficients coincide. The description is relative to the fixed acting group.
\end{theorem}
\begin{proof}
The irreducible-defining and balanced-decomposition theorems imply that every nonnormal example has a unique irreducible defining summand $V$ of dimension $d>N/2$. Put $q=N-d$. The whole subgroup $H$ preserves $V$, including when $H$ is disconnected. Uniqueness prevents a component of $H$ from permuting it with another irreducible $H^0$-summand. Thus $H\subset S(\U(d)\times \U(q))$.

Let $\mathfrak s$ be the semisimple part of the first-block projection. Its representation on $V$ is irreducible. Theorem~\ref{thm48} leaves only $\mathfrak s=\mathfrak{su}(d)$ or standard $\mathfrak s=\mathfrak{sp}(d/2)$. In either case, this simple projection is the image of a simple ideal of $\mathfrak h$. That ideal cannot act nontrivially on the smaller complement. The least dimension of a nontrivial complex representation is $d$ for both $\mathfrak{su}(d)$ and standard $\mathfrak{sp}(d/2)$, whereas $q<d$. It therefore acts independently on the dominant block. All other simple ideals and central directions have scalar first block. Theorem~\ref{thm49} gives family~(1), and Theorem~\ref{thm410} gives family~(2).

The explicit compensators in those classifications prove sufficiency at every vector, including singular matrix ranks. In family~(1), $w$ is fixed by $H$, since $H$ preserves the two blocks, fixes $z$, and preserves $\mathfrak h$ and its orthogonal complement. In family~(2), the normalizer preserves each listed summand. Hence all parameters stated above descend to the allowed disconnected isotropy extensions.
\end{proof}

\section{The symplectic groups}
\label{sec:unified-c}\label{sec:classical-c-progress}
In this section, we classify the $G$-g.o. metrics for the compact symplectic groups. We first treat quaternionically irreducible isotropy and then use a largest irreducible block to reduce the proof to a standard quaternionic projective fibration.

\begin{theorem}\label{thm51}
Every closed subgroup $H\subset \U(n)\subset \Sp(n)$ has only normal $\Sp(n)$-g.o. metrics.
\end{theorem}
\begin{proof}
Write $\mathfrak k=\mathfrak u(n)$, $\mathfrak f=\mathfrak k\ominus\mathfrak h$ and $\mathfrak p=\mathfrak k^\perp$. The central scalar circle of $\U(n)$ centralizes $H^0$, regardless of whether that circle is contained in $H$. The normalizer lemma makes a $G$-g.o. cometric commute with this circle. Its action on $\mathfrak f$ is trivial and its action on $\mathfrak p\cong(\operatorname{Sym}^2\mathbb C^n)_{\mathbb R}$ has nonzero scalar weight. Thus $C$ preserves $\mathfrak f\oplus\mathfrak p$. Descent makes $C_{\mathfrak p}$ a $\Sp(n)$-g.o. cometric on $\Sp(n)/\U(n)$. This symmetric space is isotropy irreducible, so $C_{\mathfrak p}=cI$.

The $\U(n)$ action on $\operatorname{Sym}^2\mathbb C^n$ has a locally free vector. Indeed, take a positive diagonal symmetric matrix $D$ with distinct diagonal entries. For $A\in\mathfrak u(n)$, the equation $AD+DA^t=0$ forces $A=0$. The diagonal entries vanish, and the off-diagonal equations, combined with skew-Hermitian symmetry, force both their real and imaginary parts to vanish. For $y\in\mathfrak f$ and this $x=D\in\mathfrak p$, the horizontal g.o. equation is $[(C-cI)y-Z,x]=0,\qquad Z\in\mathfrak h.$ Local freeness gives $(C-cI)y=Z$, and orthogonality to $\mathfrak h$ gives both sides zero. Hence $C=cI$.
\end{proof}

For the standard tensor inclusion $K=\Sp(p)\SO(q)\subset \Sp(pq)$, let $A_p=\mathfrak{sp}(p),\quad W_p=\operatorname{Herm}_0(p,\mathbb H),\quad B_q=\mathfrak{so}(q),\quad S_q=\operatorname{Sym}_0^2\mathbb R^q.$ Exterior and symmetric squaring of the complex defining tensor representation give the real orthogonal decomposition
\begin{equation}\label{formula27}
 \mathfrak{sp}(pq)\ominus(\mathfrak{sp}(p)+\mathfrak{so}(q))
   =(A_p\otimes S_q)\oplus(W_p\otimes B_q).
\end{equation}
For $p=1$, the second summand is zero.

We use Lemmas~\ref{lem24} and~\ref{lem25}.

\begin{theorem}\label{thm52}
Let $p\ge2$, $q\ge3$, and $(p,q)\ne(2,3)$. Every closed subgroup $H\subset \Sp(p)\SO(q)\subset \Sp(pq)$ has only normal $\Sp(pq)$-g.o. metrics.
\end{theorem}
\begin{proof}
Both modules in \eqref{formula27} are locally free for $K$. For $A_p\otimes S_q$, the required base sizes are two and one, while both dimensions are at least three. For $W_p\otimes B_q$, if $p\ge3$, the required sizes are three and two and $\dim B_q\ge3$, $\dim W_p\ge14$. If $p=2$ and $q\ge4$, the required sizes are four and two, while $\dim B_q\ge6$ and $\dim W_2=5$. Proposition~\ref{prop37} applies to these orthogonal submodules of $\mathfrak k^\perp$. The argument also covers $q=4$. The adjoint representation of $\mathfrak{so}(4)$ is reducible but faithful and still has a two-vector base.
\end{proof}

\begin{theorem}\label{thm53}
Let $H$ be closed and suppose $H^0$ acts irreducibly on the quaternionic defining module $\mathbb H^n$ of $\Sp(n)$. Then every $\Sp(n)$-g.o. metric on $\Sp(n)/H$ is normal.
\end{theorem}
\begin{proof}
First suppose the complex defining module $\mathbb C^{2n}$ is reducible. The standard Schur decomposition for compact quaternionic representations says that a quaternionically irreducible but complex-reducible module has an invariant complex Lagrangian. Thus, after conjugation, $H^0\subset \U(n)$. Theorem~\ref{thm51} applies at Lie algebra level.

Now suppose the complex defining module is irreducible. Schur's lemma and preservation of its symplectic form imply that $\mathfrak h$ has zero center. If $\mathfrak h$ is simple, apply \cite[Theorems~1--2 and Table~1]{ChNiNi2023}. The only nonnormal case with ambient $\Sp(n)$ has the standard $\Sp(n-1)$, whose defining representation is reducible, so it does not occur here. If $\mathfrak h=\bigoplus_j\mathfrak h_j$ has at least two simple ideals, its complex defining representation factors as $\bigotimes_j V_j$. Each $V_j$ is self-dual, and an odd number have quaternionic type, since the product is symplectic. Choose one quaternionic factor of dimension $2p$. The product of the others has real type and dimension $q\ge3$. Thus $H^0\subset \Sp(p)\SO(q)$ in its standard tensor realization. Theorem~\ref{thm52} treats all cases except $p=1$ and $(p,q)=(2,3)$.

If $p=1$, the corresponding simple ideal is the full $\mathfrak{sp}(1)$. The complement $\mathfrak f=\mathfrak k\ominus\mathfrak h$ is fixed by this ideal, whereas $\mathfrak p=\mathfrak k^\perp=A_1\otimes S_q$ has no fixed vector. Hence $C$ preserves $\mathfrak f\oplus\mathfrak p$. The $K$-module $\mathfrak p$ is irreducible, so descent makes $C_{\mathfrak p}=cI$. It is locally free by Lemmas~\ref{lem24} and~\ref{lem25}. The horizontal g.o. equation then gives $C=cI$.

In the remaining case $(p,q)=(2,3)$, the real three-dimensional factor is the full $\SO(3)$. Its fixed vectors in $\mathfrak h^\perp$ are exactly $\mathfrak f$, while the two horizontal modules have types $(10,5)$ and $(5,3)$ and hence no $\SO(3)$-fixed vectors. Again $C$ preserves the split. The two-isotropy-summand classification in \cite[Theorem~2]{ChNi2019} implies that every $\Sp(6)$-g.o. cometric on $\Sp(6)/(\Sp(2)\SO(3))$ is scalar. The only nonnormal symplectic entry in that classification is $\Sp(r+1)/(\Sp(r)\U(1))$. The $(10,5)$ module is locally free. Scalar descent therefore finishes this case as well.
\end{proof}

Let $\mathbb H^n=\bigoplus_{i=1}^r V_i$, $r\ge2$, be an orthogonal decomposition into quaternionically irreducible $H^0$-modules. Put $d_i=\dim_{\mathbb H}V_i$, $\mathfrak k=\bigoplus_i\mathfrak{sp}(V_i)$ and $\mathfrak f=\mathfrak k\ominus\mathfrak h$. Every block sign, equal to $-I$ on one $V_i$ and $I$ elsewhere, centralizes $H^0$. The full normalizer lemma therefore separates $\mathfrak f$ from the off-diagonal modules $\mathfrak p_{ij}=\operatorname{Hom}_{\mathbb H}(V_j,V_i)$, and separates these modules from one another. Descent to the maximal-rank space $\Sp(n)/\prod_iSp(d_i)$, followed by \cite[Theorem~6.2]{AlNi2009}, gives
\begin{equation}\label{formula28}
 C=E\oplus cI_{\mathfrak p},\qquad
 \mathfrak p=\bigoplus_{i<j}\mathfrak p_{ij},\quad c>0.
\end{equation}
The quotient is never the exceptional presentation $\Sp(n)/(\Sp(n-1)\U(1))$. Set $D=(E-cI)|_{\mathfrak f}$ and extend it by zero on $\mathfrak p$.

Let $\mathfrak h_i$ be the projection of $\mathfrak h$ to $\mathfrak{sp}(V_i)$ and put $\mathfrak w_i=\mathfrak{sp}(V_i)\ominus\mathfrak h_i$. We regard $\mathfrak w_i$ as a subspace of $\mathfrak f$.

\begin{lemma}
\label{lem54}
The positive self-adjoint compression $A_i=P_{\mathfrak w_i}C|_{\mathfrak w_i}$ is a $\Sp(d_i)$-g.o. cometric on $\Sp(d_i)/H_i$. If every $\Sp(d_i)$-g.o. cometric on this quotient is scalar, write $A_i=a_iI$. For $i\ne j$, under this scalar-compression assumption, $P_{\mathfrak w_i}C(\mathfrak w_j) \subset\mathfrak w_i^{\mathfrak h_i}.$ Consequently $P_{\mathfrak w_j}CY=0$ whenever $Y\in\mathfrak w_i$ is orthogonal to the fixed subspace.
\end{lemma}
\begin{proof}
For $Y\in\mathfrak w_i$, write $(CY)_i=A_iY+B_iY$ with $B_iY\in\mathfrak h_i$. The $i$th block of the g.o. equation gives $[A_iY,Y]=[Z_i-B_iY,Y]$ with $Z_i-B_iY\in\mathfrak h_i$. Equivariance and positivity of the compression are immediate, so this is precisely the asserted g.o. condition. For $V\in\mathfrak w_j$, apply the g.o. equation to $Y+V$ and project its $i$th block to $\mathfrak h_i$. All terms except $[P_{\mathfrak w_i}CV,Y]_{\mathfrak h_i}$ vanish. Pairing with arbitrary elements of $\mathfrak h_i$ shows that $P_{\mathfrak w_i}CV$ is orthogonal to $[\mathfrak h_i,\mathfrak w_i]$, hence is fixed by $\mathfrak h_i$. The last claim follows by self-adjointness.
\end{proof}

We also use an elementary torus fact. If a compact connected abelian group acts effectively at Lie algebra level on a connected manifold, its generic infinitesimal stabilizer is zero. Its principal stabilizers are identical, since the group is abelian, and therefore fix an open dense set and hence the entire manifold. For an abelian subalgebra whose analytic subgroup is not closed, apply this fact to its compact torus closure. Every nontrivial conjugacy orbit of a quaternionic Hermitian matrix is an effective $\Sp(d)$-space at Lie algebra level for $d\ge2$, since $\mathfrak{sp}(d)$ is simple. In particular, for any abelian $\mathfrak t\subset\mathfrak{sp}(d)$ and $1\le e\le d$, one can choose a positive semidefinite rank-$e$ Hermitian matrix $M$ with
\begin{equation}\label{formula29}
 \{T\in\mathfrak t:[T,M]=0\}=0.
\end{equation}
When $e=d$, choose distinct positive eigenvalues to ensure that the conjugacy orbit is nontrivial.

\begin{proposition}
\label{prop55}
Suppose $d_i=\max_jd_j\ge2$ and $\mathfrak h_i\subsetneq\mathfrak{sp}(d_i)$. Then $C=cI$.
\end{proposition}
\begin{proof}
Theorem~\ref{thm53} and Lemma~\ref{lem54} give $A_j=a_jI$ on every nonzero $\mathfrak w_j$. First suppose the complex representation of $\mathfrak h_i$ on $\mathbb C^{2d_i}$ is irreducible. Then $\mathfrak h_i$ is semisimple, and $\mathfrak w_i^{\mathfrak h_i}=0$. By \cite[Theorems~1.6 and~3.1, Table~1]{BeKo2021}, applied to the complexification, we see that its generic stabilizer on $\mathfrak w_i$ is abelian. The only symplectic-series exceptions have a reducible defining representation with a full standard symplectic block, or are two full symplectic blocks. The compact real form has the same generic abelianness by \cite[Lemma~2.12(1)]{BeKo2021}. Choose a generic nonzero $Y\in\mathfrak w_i$ with abelian $\mathfrak a=(\mathfrak h_i)_Y$ and set $\mathfrak t=\mathfrak a+\mathbb RY$. By Lemma~\ref{lem54}, the other block components of $CY$ lie in $\mathfrak h_j$. Choose a star-shaped off-diagonal vector $X$, with each block $X_{ij}:V_j\to V_i$ injective, and one Gram matrix $X_{ij}X_{ij}^*$ satisfying \eqref{formula29} for $\mathfrak t$. This is possible because $d_i\ge d_j$.

Apply the g.o. equation to $Y+X$. Its $i$th diagonal component says $A:=Z_i-B_iY\in\mathfrak a$. Its $(i,j)$ component is $((a_i-c)Y-A)X_{ij}=X_{ij}((CY)_j-Z_j).$ Both coefficient operators are skew-Hermitian, so the operator on the left commutes with $X_{ij}X_{ij}^*$. It belongs to $\mathfrak t$ and must vanish. Orthogonality of $Y$ to $\mathfrak h_i$ then gives $a_i=c$ and $A=0$. Injectivity of all $X_{ij}$ gives $(CY)_j=Z_j$ for all $j\ne i$. Thus $DY=Z\in\mathfrak f\cap\mathfrak h=0$. Density of the generic set proves $D\mathfrak w_i=0$.

For arbitrary $r\in\mathfrak f$, self-adjointness implies $(Dr)_i\in\mathfrak h_i$. Add $Y\in\mathfrak w_i$ so that $P_{\mathfrak w_i}r_i+Y$ has abelian $\mathfrak h_i$-stabilizer. In the g.o. equation at $r+Y+X$, the $i$th diagonal equation puts $(Dr)_i-Z_i$ in an abelian subalgebra of $\mathfrak h_i$. It centralizes $r_i+Y$, hence also its $\mathfrak w_i$-projection. Choose one star Gram matrix freely for that abelian algebra. The same argument gives $(Dr)_i=Z_i$, followed by $(Dr)_j=Z_j$ in every other block. Thus $Dr=Z=0$.

Now suppose the complex defining representation of $\mathfrak h_i$ is reducible. Quaternionic irreducibility puts $\mathfrak h_i\subset\mathfrak u(d_i)$ after conjugation. Let $P=(\operatorname{Sym}^2\mathbb C^{d_i})_{\mathbb R}$ inside the $i$th block. The scalar circle $J_i$ on this block centralizes the entire $\mathfrak h$. Its adjoint weights on $P$, the incident cross blocks, and all remaining diagonal blocks are respectively $2,1,0$. Thus $C$ preserves $P$, and its scalar compression gives $C|_P=a_iI$. Choose $Y\in P$ with zero $\mathfrak u(d_i)$-stabilizer. The diagonal g.o. equation at $Y+X$ gives $Z_i=0$. Choose one star Gram matrix which does not commute with $Y$. This is possible for $d_i\ge2$ by \eqref{formula29} applied to $\mathbb RY$. The mixed equation forces $a_i=c$, so $DP=0$. For arbitrary $r$, self-adjointness now gives $(Dr)_i\in\mathfrak u(d_i)$. Add $Y\in P$ so that $P r_i+Y$ has zero full $\mathfrak u(d_i)$-stabilizer. The diagonal equation at $r+Y+X$, projected to $P$, gives $(Dr)_i=Z_i$. Injectivity of the star blocks again gives $Dr=Z=0$.
\end{proof}

\begin{proposition}
\label{prop56}
Suppose a simple ideal of $\mathfrak h$ projects onto $\mathfrak{sp}(d)$ on at least two quaternionically irreducible blocks of maximal dimension $d$. Then $C=cI$, including $d=1$.
\end{proposition}
\begin{proof}
The least dimension of a nontrivial quaternionic representation of $\mathfrak{sp}(d)$ is $d$, with equality only for the standard module. Thus this ideal acts diagonally on a collection of equally sized maximal blocks, and trivially on all remaining blocks. After identifying the standard copies, its commutant contains the real orthogonal rotations of the multiplicity space. For two such blocks, $Y=\operatorname{diag}(A,-A,0,\ldots,0),\qquad A\in\mathfrak{sp}(d),$ lies in $\mathfrak f$. A real rotation through $\pi/4$ in the two-copy multiplicity space conjugates $Y$ into the $(1,2)$ off-diagonal block. Normalizer invariance and \eqref{formula28} give $CY=cY$. By self-adjointness, for arbitrary $r\in\mathfrak f$ the components of $Dr$ on all these repeated blocks are equal. Orthogonality to the diagonal $\mathfrak{sp}(d)$ ideal then makes all of them zero.

Choose such a difference $Y$ so that $r_1+A$ is regular in $\mathfrak{sp}(d)$. In the g.o. equation at $r+Y+X$, the first diagonal component puts $Z_1$ in a Cartan subalgebra $\mathfrak t$. On the two repeated blocks $Z_1=Z_2$. Their cross equation is therefore $[Z_1,X_{12}]=0$. The conjugation action of $\mathfrak t$ on $\operatorname{Mat}_d(\mathbb H)$ is faithful at Lie algebra level. Its generic stabilizer is zero. Choose $X_{12}$ invertible with this property, so $Z_1=0$. Choose all remaining star blocks $X_{1j}$ injective. Their equations give $(Dr)_j=Z_j$ for every $j$. Thus $Dr=Z\in\mathfrak f\cap\mathfrak h=0$. This commutator argument also works for $d=1$. A Gram-matrix argument alone would not, since quaternionic Hermitian one-by-one matrices are real scalars.
\end{proof}

\begin{theorem}
\label{thm57}
Let $G$ be compact connected with Lie algebra $\mathfrak{sp}(n)$ and let $H$ be any closed subgroup. Every $G$-g.o. metric on $G/H$ is normal, unless $n\ge2$ and, up to conjugacy, $\mathfrak h=\mathfrak{sp}(n-1) \quad\text{or}\quad \mathfrak h=\mathfrak{sp}(n-1)\oplus\mathfrak u(1)$ in their standard block embeddings. In each exceptional case all $G$-g.o. cometrics are exactly $C=aI_{\mathfrak v}\oplus bI_{\mathfrak p}$, $a,b>0$, where $\mathfrak v=\mathfrak{sp}(1)$ or $\mathfrak{sp}(1)\ominus\mathfrak u(1)$ and $\mathfrak p\cong\mathbb H^{n-1}$. The nonnormal ones have $a\ne b$.
\end{theorem}
\begin{proof}
Induct on $n$. For $n=1$, proper connected isotropy is trivial or a circle, and only the normal metrics occur. Quaternionically irreducible isotropy is covered by Theorem~\ref{thm53}. Otherwise take the irreducible quaternionic block decomposition used above and a block of largest dimension $d$. If $d\ge2$ and a maximal block has proper projected isotropy, Proposition~\ref{prop55} gives normality. If $d=1$ and all projected algebras are proper, conjugate them into a product of circles. Then $H^0\subset \U(n)$ and Theorem~\ref{thm51} applies. In all other cases some maximal block has full projected $\mathfrak{sp}(d)$. This projection arises from one simple ideal of $\mathfrak h$. All other ideals and the center act trivially on that block. If that ideal acts on another block, Proposition~\ref{prop56} applies.

We are left with an independent full standard block $\mathfrak h=\mathfrak{sp}(d)\oplus\mathfrak r, \qquad \mathfrak r\subset\mathfrak{sp}(q),\qquad d+q=n.$ For $K=\Sp(d)\Sp(q)$, the vertical space $\mathfrak f=\mathfrak{sp}(q)\ominus\mathfrak r$ is precisely the fixed subspace of the first ideal. The horizontal module has no fixed vectors, so $C$ splits as $E\oplus cI$. Descent gives the scalar horizontal part because $\Sp(n)/K$ is isotropy irreducible. Restriction of the g.o. equation to the vertical space makes $E$ a $\Sp(q)$-g.o. cometric on $\Sp(q)/R$.

If $E$ is nonnormal, induction gives $\mathfrak r=\mathfrak{sp}(q-1)$ or $\mathfrak{sp}(q-1)\oplus\mathfrak u(1)$, with $q\ge2$. In the first case the original isotropy is a standard product of two nonzero symplectic blocks. The quaternionic Stiefel-bundle theorem \cite[Theorem~1.1]{ArSoSt2021b} makes the original metric normal. In the second case the original isotropy has maximal rank, and \cite[Theorem~6.2]{AlNi2009} likewise makes it normal. Two nonzero symplectic blocks cannot be the single $\Sp(n-1)$ factor of its exceptional presentation.

It remains that $E=aI$. If $a=c$, the metric is normal. If $a\ne c$, it is a nontrivial canonical variation over the irreducible symmetric base $\Sp(n)/(\Sp(d)\Sp(q))$. By \cite[Proposition~2.1 and the following remark]{Ta1999}, one nontrivial g.o. canonical variation implies that all positive parameter choices are $G$-g.o. \cite[Theorem~4.1, Table~1, rows (1.8)--(1.9)]{Ta1999} has, in the symplectic series, only the triples $(\mathfrak g,\mathfrak k,\mathfrak h) =(\mathfrak{sp}(n),\mathfrak{sp}(n-1)+\mathfrak{sp}(1), \mathfrak{sp}(n-1))$ and the same triple with an additional $\mathfrak u(1)$ in $\mathfrak h$. Since the full first block $\mathfrak{sp}(d)$ lies in $\mathfrak h$ and $\mathfrak r\subsetneq\mathfrak{sp}(q)$, these are exactly $q=1$, $\mathfrak r=0$ or $\mathfrak u(1)$. The case $\mathfrak r=\mathfrak{sp}(q)$ has no vertical space and is normal. This proves necessity.

For the sphere presentation, normalizer invariance makes the three-dimensional vertical block scalar and isotropy irreducibility makes the horizontal block scalar. For the projective presentation the two blocks are inequivalent real irreducibles. Their arbitrary positive coefficients satisfy the g.o. equation directly. For $x\in\mathbb H^{n-1}\setminus\{0\}$ and $q\in\operatorname{Im}\mathbb H$, the skew-Hermitian matrix $xqx^*/\|x\|^2\in\mathfrak{sp}(n-1)$ sends $x$ to $xq$. Given a vertical vector $y$, choose $q=-(a-b)y$. This compensator solves the mixed equation, and the vertical equation vanishes because the vertical coefficient is scalar. If $x=0$, use the zero compensator. The same formula applies to the two-dimensional complement of $\mathfrak u(1)$. Finite components of $H$ normalize its identity component and preserve these two scalar blocks, so the same classification descends to the given global quotient.
\end{proof}

\section{The orthogonal groups}
\label{sec:unified-bd}\label{sec:classical-bd-complete}
In this section, we classify the $G$-g.o. metrics for compact simple orthogonal groups. The proof separates two defining blocks, reduces a nonnormal metric to one dominant irreducible block, and determines the metrics for the remaining unitary, spin and $G_2$ projections.

\begin{theorem}
\label{thm61}
Each of the following subgroups has the property that every closed subgroup $H$ contained in it admits only normal $G$-g.o. metrics.
\begin{enumerate}
\item $\SO(a)\otimes \SO(b)\subset \SO(ab)\subset \SO(N)$,
where $a,b\ge3$ and $N\ge ab$. \item $\Sp(m)\otimes \Sp(n)\subset \SO(4mn)\subset \SO(N)$, where $m,n\ge2$ and $N\ge4mn$.
\end{enumerate}
The last inclusions fix the extra coordinates pointwise. Finite kernels of the tensor homomorphisms do not affect the statement.
\end{theorem}
\begin{proof}
In the real tensor case, the orthogonal complement of the product algebra in $\mathfrak{so}(ab)$ is $E=(\Lambda^2\mathbb R^a)\otimes \operatorname{Sym}^2_0\mathbb R^b, \qquad F=\operatorname{Sym}^2_0\mathbb R^a\otimes (\Lambda^2\mathbb R^b).$ Each is locally free by Lemmas~\ref{lem24} and~\ref{lem25}. The symmetric representation needs one vector and the adjoint representation needs two, whereas every smaller dimension here is at least three. This also covers $a=4$ or $b=4$. Simplicity of each orthogonal factor is unnecessary.

In the quaternionic tensor case, the tensor product of the two standard quaternionic representations has a real form of dimension $4mn$. Put $W_k=\operatorname{Herm}_0(k,\mathbb H)$. Its skew endomorphisms decompose as the two product ideals and $E=\mathfrak{sp}(m)\otimes W_n, \qquad F=W_m\otimes\mathfrak{sp}(n).$ These are real representations, not the underlying real spaces of the displayed complexifications. The adjoint factors have dimension at least $10$, the $W$ factors have dimension at least $5$, and Lemma~\ref{lem25} gives free-tuple sizes two and at most four, respectively. Lemma~\ref{lem24} applies to both. Proposition~\ref{prop37} proves both assertions. The same two subspaces remain in $\mathfrak l^\perp$ after adding fixed coordinates to the ambient representation.
\end{proof}

\begin{theorem}
\label{thm62}
Let $0\ne\mathfrak s\subsetneq\mathfrak{so}(d)$ be semisimple, and let $L$ be the image of $S\times \SO(2)$ on $\mathbb R^d\otimes\mathbb R^2$. For $N\ge2d$, every closed subgroup $H\subset L\subset \SO(N)$ admits only normal $\SO(N)$-g.o. metrics.
\end{theorem}
\begin{proof}
Put $R=\mathfrak{so}(d)\ominus\mathfrak s\ne0$ and $P=\operatorname{Sym}^2_0\mathbb R^2$. The circle acts on $P$ with nonzero weight two. The decomposition of $\Lambda^2(\mathbb R^d\otimes\mathbb R^2)$ provides the orthogonal $\mathfrak l$-subspaces $F=\mathfrak s\otimes P,\qquad E=\big(\operatorname{Sym}^2_0\mathbb R^d\big) \oplus(R\otimes P) \quad\subset\mathfrak l^\perp.$ The first summand of $E$ is tensored with $\Lambda^2\mathbb R^2$, so the circle acts trivially there. A simple-spectrum symmetric matrix kills $\mathfrak s$ in the stabilizer of a vector of $E$. Any nonzero component in $R\otimes P$ then kills the circle. For $F$, take a two-vector adjoint free tuple, make it orthogonal, and assign distinct positive singular values. The rectangular-matrix argument first kills the circle and then $\mathfrak s$. Apply Proposition~\ref{prop37}.
\end{proof}

\begin{theorem}
\label{thm63}
Let $n\ge2$, and use the standard action of $\Sp(n)\Sp(1)$ on $\mathbb H^n\simeq\mathbb R^{4n}$. For any $\mathfrak r\subset\mathfrak{sp}(n)$, every $\SO(4n)$-g.o. metric with isotropy $\mathfrak h=\mathfrak r\oplus\mathfrak{sp}(1)$ is normal.
\end{theorem}
\begin{proof}
The full ideal-centralizer envelope of the right scalar $\mathfrak{sp}(1)$ is $\mathfrak l=\mathfrak{sp}(n)\oplus\mathfrak{sp}(1)$. Its orthogonal complement is the irreducible real module $W_n\otimes\mathbb R^3$. The vertical space $\mathfrak{sp}(n)\ominus\mathfrak r$ is fixed by the scalar $\Sp(1)$, whereas this horizontal module has no fixed vectors. Hence $C$ preserves the two spaces. By Lemma~\ref{lem38}, the horizontal cometric is scalar. For $n\ge3$, a rectangular tensor whose three image vectors are a free tuple from Lemma~\ref{lem25}, with distinct singular values, has zero $\mathfrak l$-stabilizer. The horizontal g.o. equation then kills the vertical cometric defect. For $n=2$, this module is $\mathbb R^5\otimes\mathbb R^3$ for $\mathfrak{so}(5)\oplus\mathfrak{so}(3)$. A full-rank matrix with distinct singular values has the abelian stabilizer $\mathfrak{so}(2)$ in the first factor. More explicitly, a vertical defect lies in $\mathfrak{sp}(2)\ominus\mathfrak r\simeq \mathfrak{so}(5)\ominus\mathfrak r$. For any such skew matrix $T$, choose a two-plane on which the two-form $\langle T\cdot,\cdot\rangle$ vanishes, and choose the rectangular tensor with image its orthogonal three-plane and distinct singular values. Its stabilizer is orthogonal to $T$. Pairing the horizontal g.o. equation with $T$ gives $\|T\|^2=0$.
\end{proof}

\begin{corollary}
\label{cor64}
Let $\mathfrak s\subset\mathfrak{so}(d)$ be semisimple, and suppose its defining real representation is irreducible of real type, that is, $\operatorname{End}_{\mathfrak s}(\mathbb R^d)=\mathbb R$. Then every $\SO(d)$-g.o. metric is normal, whenever the ambient algebra is simple.
\end{corollary}
\begin{proof}
For simple $\mathfrak s$, apply \cite[Theorems~1--2 and Table~1]{ChNiNi2023}. None of its nonnormal orthogonal cases has a real-type irreducible defining isotropy representation. For several simple ideals, complexify the defining module. It is an irreducible tensor product of self-conjugate modules, each of real or quaternionic type, with an even number of quaternionic factors. If the tensor factors can be partitioned into two nonempty collections of real type, the isotropy lies in the first family of Theorem~\ref{thm61}. Otherwise there are exactly two quaternionic factors. If both have complex dimension at least four, the second family of that theorem applies. The only remaining nontrivial product consists of a quaternionic two-dimensional factor and one other quaternionic factor. It lies in $\Sp(1)\Sp(n)$ and contains the scalar $\Sp(1)$ ideal, so Theorem~\ref{thm63} applies. If $n=1$, the ambient algebra is $\mathfrak{so}(4)$, outside the hypothesis.
\end{proof}

\begin{theorem}
\label{thm65}
Let $p\ge2$, $q\ge1$, $p+q\ge5$, and $\mathfrak h=\mathfrak{so}(p)\oplus\mathfrak r \subset\mathfrak{so}(p)\oplus\mathfrak{so}(q) \subset\mathfrak{so}(p+q), \qquad\mathfrak r\subset\mathfrak{so}(q).$ If $p\ge\lfloor q/2\rfloor$, every $\SO(p+q)$-g.o. metric is normal. No irreducibility or semisimplicity assumption is imposed on $\mathfrak r$.
\end{theorem}
\begin{proof}
Put $\mathfrak k=\mathfrak{so}(p)\oplus\mathfrak{so}(q)$, $\mathfrak w=\mathfrak{so}(q)\ominus\mathfrak r$, and $\mathfrak p=\mathfrak k^\perp\simeq \operatorname{Hom}(\mathbb R^p,\mathbb R^q)$. The $\mathfrak{so}(p)$-fixed space of $\mathfrak m$ is precisely $\mathfrak w$, so $C$ preserves $\mathfrak w\oplus\mathfrak p$. The descent lemma and irreducibility of the real Grassmannian isotropy representation give $C|_{\mathfrak p}=cI$. This includes $p=2$. The excluded $(p,q)=(2,2)$ is exactly the nonsimple ambient exception.

Fix $Y\in\mathfrak w$ and put $T=(C-cI)Y\in\mathfrak w$. The horizontal component of the g.o. equation at $Y+X$ says $T-Z\in\mathfrak k_X \quad\text{for some }Z\in\mathfrak h, \qquad X\in\mathfrak p.$ If $p\ge q$, choose a full-rank matrix $X$ with distinct positive singular values. Its stabilizer in $\mathfrak k$ is $\mathfrak{so}(p-q)$ in the first block, orthogonal to $T$. Thus $Q(T,T-Z)=\|T\|^2=0$.

Suppose $p<q$. For any skew-adjoint $T$ on $\mathbb R^q$, there is a subspace $W$ of dimension at least $\lceil q/2\rceil$ such that $\langle Tu,v\rangle=0$ for all $u,v\in W$. To construct it, take one vector from each nonzero rotation plane of $T$ and the whole kernel. By the hypothesis, choose such a $W$ of dimension $q-p$. Choose $X$ with image $W^\perp$ and distinct positive singular values. Its stabilizer is exactly $\mathfrak{so}(W)$ in the second block. The defining property of $W$ means $\mathfrak{so}(W)\perp T$ for the trace form. Again $Q(T,T-Z)=\|T\|^2=0$. Thus $C=cI$.
\end{proof}

\begin{lemma}\label{lem66}
Let $d,q\ge3$ be odd, and let $\mathfrak h\subset\mathfrak{so}(d)\oplus\mathfrak{so}(q) \subset\mathfrak{so}(d+q)$. Every $G$-g.o. cometric preserves the off-diagonal block $\mathfrak p\simeq\operatorname{Mat}_{d\times q}(\mathbb R)$.
\end{lemma}
\begin{proof}
Put $\mathfrak f=(\mathfrak{so}(d)\oplus \mathfrak{so}(q))\ominus\mathfrak h$. Fix $Y=\operatorname{diag}(A,B)\in\mathfrak f$, and write
\begin{equation*}
P(X)=\begin{pmatrix}0&X\\-X^t&0\end{pmatrix},\qquad
 \ell_Y(X)=Q(CY,P(X)).
\end{equation*}
The block reflection $R=\operatorname{diag}(I_d,-I_q)$ sends $Y+P(X)$ to $Y-P(X)$. If $Y+P(X)$ is singular, these two matrices lie in the same $\SO(d+q)$ adjoint orbit. A reflection in a kernel direction centralizes the first matrix and corrects the determinant of $R$. Orbit-energy constancy therefore gives $0=F(Y+P(X))-F(Y-P(X))=2\ell_Y(X) \quad\text{whenever }Y+P(X)\text{ is singular}.$ If $A$ or $B$ has nullity at least three, every rank-one $X$ leaves a nonzero kernel vector. Thus $\ell_Y$ vanishes on all rank-one matrices and is zero.

Otherwise both nullities are one. Let $u,v$ be unit kernel vectors of $A,B$. For rank-one $X$, expansion of the Pfaffian shows that $Y+P(X)$ is singular exactly when $u^tXv=0$. The coefficient is the product of the nonzero Pfaffians on $u^\perp$ and $v^\perp$. Equivalently, this follows directly by solving the two block equations for a kernel vector. The rank-one matrices with $u^tXv=0$ span the hyperplane defined by that equation. Consequently $\ell_Y(X)=\lambda u^tXv$.

Choose a nonzero rotation two-plane of each of $A$ and $B$, with rotation parameters $a,b\ne0$. On these two-planes together with the kernel lines, take $X=\operatorname{diag}(t_1,t_2,t_3)$, zero elsewhere. The Pfaffian of the resulting six-dimensional block, up to an immaterial fixed sign, is $t_3(ab-t_1t_2).$ Take $t_1t_2=ab$ and $t_3\ne0$. The whole skew matrix is singular, but $u^tXv=t_3\ne0$. Hence $\lambda=0$. Thus $Q(CY,\mathfrak p)=0$ for every $Y\in\mathfrak f$, and self-adjointness proves preservation.
\end{proof}

\begin{corollary}
\label{cor67}
If $H$ preserves an orthogonal splitting $\mathbb R^N=V\oplus W$ with $\dim V,\dim W\ge2$, every $G$-g.o. cometric preserves its off-diagonal block. It also does so when one dimension is one and the other is even.
\end{corollary}
\begin{proof}
If one block is even-dimensional, its sign change belongs to $\SO(N)$ and centralizes $H^0$. Full normalizer invariance separates the two eigenspaces of its adjoint action. If both dimensions are odd and at least three, apply Lemma~\ref{lem66}.
\end{proof}

\begin{lemma}
\label{lem68}
Suppose $\mathfrak h\subset\mathfrak{so}(p)\oplus \mathfrak{so}(q)$, where $2\le p\le q\le2p+1$ and $p+q\ge5$. Every $\SO(p+q)$-g.o. metric is normal.
\end{lemma}
\begin{proof}
Two-block separation and descent to the real Grassmannian give $C=E\oplus cI_{\mathfrak p}$. For any vertical $Y$, put $T=(E-cI)Y\in\mathfrak h^\perp\cap (\mathfrak{so}(p)\oplus\mathfrak{so}(q))$. The horizontal g.o. equation gives $T-Z\in\mathfrak k_X$, where $Z\in\mathfrak h$ and $\mathfrak k=\mathfrak{so}(p)\oplus\mathfrak{so}(q)$. Here identify a cross block with a map $X:\mathbb R^p\to\mathbb R^q$ (the transpose of the upper right matrix block). Choose a $(q-p)$-dimensional subspace $W\subset\mathbb R^q$ isotropic for the two-form of the second block of $T$. This is possible because $q-p\le\lceil q/2\rceil$. Take a full-rank rectangular $X$ with image $W^\perp$ and distinct positive singular values. Then $\mathfrak k_X=\mathfrak{so}(W)$ is orthogonal to $T$. Pairing $T-Z$ with $T$ gives $\|T\|^2=0$.
\end{proof}

\begin{lemma}
\label{lem69}
If a $G$-g.o. metric is nonnormal and the defining representation of $H^0$ is reducible, it has a unique real irreducible block of dimension $d$, with complementary dimension $q$, such that $d>2q+1.$ The possible exception to two-block separation is $q=1$ and $d$ odd.
\end{lemma}
\begin{proof}
Let the dimensions of real irreducible blocks be $d_i$, and put $L=\lceil(N-1)/3\rceil$. If no $d_i>N-L$, one can group some blocks into a dimension $p\in[L,N-L]$. Either one block already lies in that interval, or successively add blocks of size at most $L-1$ until reaching $L$. The latter sum is at most $2L-2\le N-L$. Replace $p$ by $\min(p,N-p)$. For $N\ge5$, it satisfies $2\le p\le q=N-p\le2p+1$, contradicting Lemma~\ref{lem68}. Thus some $d>N-L=\lfloor(2N+1)/3\rfloor$, which is equivalent to $d>2(N-d)+1$. Uniqueness is immediate. The last assertion follows from two-block separation.
\end{proof}

We next determine which internal projections can survive on the dominant block. The unitary case uses the irreducible defining representation theorem from the first classical section.

\begin{lemma}
\label{lem610}
Let $n\ge3$, let $\mathfrak s\subsetneq\mathfrak{su}(n)$ be semisimple and complex-irreducible on $\mathbb C^n$, and let $\mathfrak l=\mathfrak s$ or $\mathfrak l=\mathfrak s\oplus\mathbb RJ$, where $J=iI_n$. Every $\SO(2n)$-g.o. cometric with isotropy $\mathfrak l$ is scalar. The generic stabilizer of $\mathfrak l$ on $\mathfrak{so}(2n)\ominus\mathfrak l$ is abelian.
\end{lemma}
\begin{proof}
Normalizer invariance under $J$ separates $\mathfrak w=\mathfrak{su}(n)\ominus\mathfrak s, \qquad\mathfrak v=(\Lambda^2\mathbb C^n)_{\mathbb R}, \qquad\mathbb RJ\quad\text{if }J\notin\mathfrak l.$ There is no mixing of $\mathfrak w$ with the last line, since complex irreducibility makes $\mathfrak w^{\mathfrak s}=0$. Descent to $\SO(2n)/\U(n)$ makes $C|_{\mathfrak v}=bI$. Compression to $\SU(n)/S$, and the irreducible defining type-A theorem, give $C|_{\mathfrak w}=aI$.

First assume $\mathfrak s$ is not standard $\mathfrak{sp}(n/2)$. The type-A generic-stabilizer theorem gives a generic $R\in\mathfrak w$ with $\mathfrak s_R$ abelian. At $R+Y$, $Y\in\mathfrak v$, the g.o. equation implies $[Z,R]=0,\qquad [(a-b)R-Z,Y]=0.$ The relevant algebra is contained in the abelian algebra $\mathfrak s_R+\mathbb RR+\mathbb RJ$. The representation of $\mathfrak u(n)$ on $\Lambda^2\mathbb C^n$ is faithful for $n\ge3$. A faithful orthogonal representation of an abelian compact algebra has vectors with zero infinitesimal stabilizer (choose nonzero components in all its weight spaces). Choose such a $Y$. Then $(a-b)R=Z$, and orthogonality gives $a=b$. If the central line belongs to the tangent space, write $C|_{\mathbb RJ}=tI$ and apply the g.o. equation at $R+J+Y$. The same argument with $\mathfrak s_R+\mathbb RJ$ gives $t=b$.

For $\mathfrak s=\mathfrak{sp}(k)$, $n=2k\ge4$, let $\Omega\in\Lambda^2\mathbb C^{2k}$ be the standard symplectic form. It is fixed by $\mathfrak{sp}(k)$, while $J$ acts on it by a nonzero scalar complex rotation. Choose $R=i\operatorname{diag}(r_1I_2,\ldots,r_kI_2)\in\mathfrak w$ with $\sum r_j=0$ and not all $r_j=0$. The $\mathfrak v$ component of the g.o. equation at $R+\Omega$ would give $(a-b)[R,\Omega]\in\mathbb R i\Omega$. But $[R,\Omega]$ has distinct block coefficients and is not in that line. Thus $a=b$. If $J$ is tangent, its full quaternionic scalar centralizer $\Sp(1)$ rotates $J,\operatorname{Re}\Omega, \operatorname{Im}\Omega$. Normalizer invariance gives the same coefficient $b$ on the remaining central line.

For the generic-stabilizer assertion, apply \cite[Theorems~1.6 and~3.1, Table~1]{BeKo2021}. Among the nonabelian-generic-stabilizer orthogonal embeddings, the only one commuting with a complex structure and acting complex-irreducibly is the full $\mathfrak{su}(n)\subset\mathfrak{so}(2n)$. The dominant real standard blocks and the $G_2$ and spin representations in that table cannot commute with a complex structure on their real-type summand of multiplicity one. The full $\mathfrak{su}(n)$ case is excluded here. The passage from the complex table to the compact real representation is \cite[Lemma~2.12(1)]{BeKo2021}. Adding $\mathbb RJ$ does not change abelianness. The $\mathfrak s$-action on that central line is trivial, and the derived algebra of a compact $\mathfrak l$-stabilizer lies in the corresponding abelian $\mathfrak s$-stabilizer.
\end{proof}

\begin{lemma}\label{lem611}
Let $\mathfrak l\subset\mathfrak{so}(d)$ act real-irreducibly, with $d\ge5$. Except for the following projections, the $\SO(d)$-g.o. cometrics on $\SO(d)/L$ are scalar and its generic internal stabilizers are abelian. $\mathfrak l=\mathfrak{so}(d),\qquad (d,\mathfrak l)=(7,\mathfrak g_2),(8,\mathfrak{spin}(7)), \qquad d=2n,\quad\mathfrak l=\mathfrak{su}(n)\text{ or } \mathfrak u(n).$
\end{lemma}
\begin{proof}
For real type, the projection is semisimple and Corollary~\ref{cor64} gives scalarity. The generic-stabilizer alternatives are precisely \cite[Propositions~3.4--3.5]{BeKo2021}. The only proper irreducible orthogonal exceptions are $G_2$ in dimension seven and the spin $\Spin(7)$ in dimension eight. For complex type, the projection is $\mathfrak s\oplus(0\text{ or }\mathbb RJ)$ inside $\mathfrak u(n)$, with $\mathfrak s$ complex-irreducible. Use Lemma~\ref{lem610} unless $\mathfrak s=\mathfrak{su}(n)$. For quaternionic type, the projection is semisimple, lies in $\mathfrak{sp}(n)\subsetneq\mathfrak{su}(2n)$, and is complex-irreducible on $\mathbb C^{2n}$. The same lemma applies for $n\ge2$. The dimension-four cases are treated separately below.
\end{proof}

\begin{lemma}
\label{lem612}
Let $d>q\ge1$, let $\mathfrak h\subset\mathfrak{so}(d)\oplus\mathfrak{so}(q)$, and let $\mathfrak l$ be its first projection. Suppose $C$ preserves the off-diagonal block, every $\SO(d)$-g.o. cometric on $\SO(d)/L$ is scalar, and $\mathfrak l_Y$ is abelian for generic $Y\in\mathfrak w=\mathfrak{so}(d)\ominus\mathfrak l\ne0$. Then $C$ is scalar. No independence of the two projections of $\mathfrak h$ is assumed.
\end{lemma}
\begin{proof}
By Lemma~\ref{lem38}, write $C=E\oplus cI_{\mathfrak p}$ and put $D=E-cI$ on $\mathfrak f=(\mathfrak{so}(d)\oplus\mathfrak{so}(q))\ominus\mathfrak h$. For $y\in\mathfrak f$ and $X\in\operatorname{Mat}_{d\times q}(\mathbb R)$, the same compensator $Z\in\mathfrak h$ satisfies
\begin{equation}\label{formula30}
 [Z-Dy,y]=0,\qquad (Dy-Z)_{11}X=X(Dy-Z)_{22}.
\end{equation}
The compression $P_{\mathfrak w}C|_{\mathfrak w}$ is a $\SO(d)$-g.o. cometric on $\SO(d)/L$. Indeed, project the equation at $Y\in\mathfrak w$ onto the first block and absorb the $\mathfrak l$-component of $CY$ into the compensator, exactly as in Proposition~\ref{prop45}. Thus this compression is $aI$, and $(DY)_{11}=(a-c)Y+L(Y)$ for some $L(Y)\in\mathfrak l$.

Choose nonzero $Y\in\mathfrak w$ with abelian stabilizer $\mathfrak l_Y$. The vertical equation in \eqref{formula30} gives $A=Z_{11}-L(Y)\in\mathfrak l_Y$. The algebra $\mathfrak l_Y+\mathbb RY$ is abelian. Its action on the real Grassmannian of $q$-planes in $\mathbb R^d$ is infinitesimally effective. A compact torus with an effective action on a connected manifold has zero generic infinitesimal stabilizer. Its principal stabilizers are equal, so a principal stabilizer fixes a dense open set and hence the whole manifold. For an abelian subalgebra that does not integrate to a closed torus, apply this argument to its torus closure. Choose a full-column-rank $X$ whose image has zero stabilizer in $\mathfrak l_Y+\mathbb RY$. The horizontal equation implies that $(a-c)Y-A$ preserves this image. Hence $(a-c)Y=A$, and by orthogonality we have $a=c$ and $A=0$. Since $X$ has full column rank, we obtain $DY=Z$. As $DY\perp\mathfrak h$, we have $DY=0$. By density, $D\mathfrak w=0$.

Since $D$ is self-adjoint, $D\mathfrak f\subset\mathfrak w^\perp$. Fix $r\in\mathfrak f\cap\mathfrak w^\perp$ and apply \eqref{formula30} to $y=r+Y$. Both $r_{11}$ and $(Dr)_{11}$ lie in $\mathfrak l$. Projecting the first-block vertical equation onto $\mathfrak w$, we obtain $(Z-Dr)_{11}\in\mathfrak l_Y$. Choose the image of $X$ to have zero stabilizer in this abelian algebra. The horizontal equation gives $(Z-Dr)_{11}=0$, then $(Z-Dr)_{22}=0$. Therefore $Dr=Z=0$. Thus $D=0$ and $C=cI$.
\end{proof}

\begin{lemma}
\label{lem613}
Suppose $q\ge2$, $q<n$, and $\mathfrak h=\mathfrak{su}(n)\oplus\mathfrak r, \qquad \mathfrak r\subset\mathbb RJ\oplus\mathfrak{so}(q), \qquad \mathfrak g=\mathfrak{so}(2n+q).$ Every $G$-g.o. metric is normal. The algebra $\mathfrak r$ may contain an arbitrary coupling with the circle.
\end{lemma}
\begin{proof}
Set $\mathfrak k=\mathfrak u(n)\oplus\mathfrak{so}(q)$. Two-block separation and normalizer invariance under $J$ make $C$ preserve $\mathfrak p=\mathfrak k^\perp = (\Lambda^2\mathbb C^n)_{\mathbb R} \oplus(\mathbb C^n\otimes\mathbb R^q)_{\mathbb R}.$ By descent and \cite[Theorem~6.2]{AlNi2009}, we have $C|_{\mathfrak p}=cI$. Indeed the only nonnormal orthogonal member of that classification is $\SO(2n+1)/\U(n)$, whereas $q\ge2$ here.

The $K$-module $\mathfrak p$ has a free vector. Choose a complex skew matrix $Y$ with distinct nonzero Takagi block parameters, and its generic maximum rank. Its stabilizer in $\U(n)$ is $\Sp(1)^{\lfloor n/2\rfloor}$, with an additional $\U(1)$ on the kernel line when $n$ is odd. Choose a real cross matrix $X$ whose $q$ columns are orthogonal with distinct lengths, and whose first column has a nonzero component in every Takagi two-dimensional complex block and in the kernel line if present. For a stabilizer pair $(A,B)$, the equation $AX=XB$ forces $B$ to commute with the simple-spectrum real Gram matrix $X^tX$, and hence $B=0$. Then $A$ kills the first column. The action of $\Sp(1)$ on a nonzero quaternionic scalar, and of $\U(1)$ on a nonzero complex scalar, has zero infinitesimal stabilizer. Thus $A=0$. The horizontal g.o. equation for any vertical defect now gives $T=Z\in\mathfrak h\cap\mathfrak h^\perp=0$.
\end{proof}

\begin{lemma}
\label{lem614}
For the spin embedding $\mathfrak{spin}(7)\subset\mathfrak{so}(8)$ and any $\mathfrak r\subset\mathfrak{so}(q)$, $1\le q\le3$, all $\SO(8+q)$-g.o. cometrics on $\SO(8+q)/(\Spin(7)\cdot R)$ are exactly $C=aI_{\mathfrak{so}(8)\ominus\mathfrak{spin}(7)} \oplus bI_{(\mathfrak{so}(q)\ominus\mathfrak r) \oplus(\mathbb R^8\otimes\mathbb R^q)}, \qquad a,b>0.$ For $G_2\subset \SO(7)$ and $\mathfrak r\subset\mathfrak{so}(2)$, the corresponding formula holds on $\SO(9)/(G_2\cdot R)$ with $\mathbb R^8$ replaced by $\mathbb R^7$.
\end{lemma}
\begin{proof}
Two-block separation and descent make the cross coefficient $b$ scalar. The internal seven-dimensional module is irreducible and nontrivial for the large simple factor, whereas the residual vertical space is trivial for it. Thus $C=aI\oplus E\oplus bI$. For $Y$ in the seven-dimensional internal module and $R_0\in\mathfrak{so}(q)\ominus\mathfrak r$, put $T=(E-bI)R_0$. The g.o. equations give $[A,Y]=0,\qquad (A-(a-b)Y)X=X(B-T),\quad A\in\mathfrak s, B\in\mathfrak r.$ Choose $X$ with full column rank and distinct real singular values. Since both end operators are skew-adjoint, $B-T$ commutes with $X^tX$, so $B-T=0$. Orthogonality gives $T=0$. Thus $E=bI$. Conversely, we construct a compensator in the large simple factor. For the spin embedding, use the Clifford decomposition $\mathfrak{so}(8)=\mathfrak{spin}(7)\oplus\operatorname{span}_{\mathbb R}\{J_1,\ldots,J_7\}$, where $J_iJ_j+J_jJ_i=-2\delta_{ij}I$. For a nonzero internal vector $Y$, we have $Y^2=-r^2I$ for some $r>0$. The operator $J=Y/r$ is a complex structure, and $\mathfrak{spin}(7)_Y=\mathfrak{su}(4)$ on $(\mathbb R^8,J)=\mathbb C^4$. Let $U$ be the complex span of the $q$ real columns of $X$, with $k=\dim_{\mathbb C}U\le q\le3$. The skew-Hermitian operator equal to $(a-b)Y$ on $U$ and $-(a-b)kY/(4-k)$ on $U^\perp$ has complex trace zero. It therefore belongs to $\mathfrak{su}(4)$, commutes with $Y$, and agrees with $(a-b)Y$ on every column of $X$.

For $G_2\subset \SO(7)$, identify the internal module with the cross-product operators $Y_y:x\mapsto y\times x$. For $y\ne0$, the stabilizer $(G_2)_y=\SU(3)$ fixes $y$ and acts complex-linearly on $y^\perp\cong\mathbb C^3$, where $Y_y/\|y\|$ is the complex structure. The projections of the two columns of $X$ span a complex subspace of dimension $k\le2$. The same trace correction, with $3-k$ in the denominator, gives $A\in\mathfrak{su}(3)$ satisfying $[A,Y_y]=0$ and $AX=(a-b)Y_yX$. The component along $y$ is killed by both operators.

In either case, use this $A$ and put $B=0$. If the internal vector is zero, use $A=0$. These choices solve both the horizontal and vertical equations for every $R_0\in\mathfrak{so}(q)\ominus\mathfrak r$. This proves sufficiency.
\end{proof}

\begin{theorem}
\label{thm615}
Let $\mathfrak g=\mathfrak{so}(N)$ be compact simple, $N\ge5$, and let $H$ be a proper closed isotropy subgroup. Up to automorphisms of the acting Lie algebra, the following list gives every possibility for a nonnormal $G$-g.o. metric. In all other cases $C=cI$, $c>0$. The formulas below also include the normal metrics.
\begin{enumerate}
\item $\mathfrak g=\mathfrak{so}(2n+1)$,
$\mathfrak h=\mathfrak u(n)$, $n\ge2$. On $\mathfrak m=\mathfrak v\oplus\mathfrak p$, where $\mathfrak v=(\Lambda^2\mathbb C^n)_{\mathbb R}$ and $\mathfrak p=(\mathbb C^n)_{\mathbb R}$, all cometrics are $C=aI_{\mathfrak v}\oplus bI_{\mathfrak p},\qquad a,b>0.$
\item $\mathfrak g=\mathfrak{so}(2n+1)$,
$\mathfrak h=\mathfrak{su}(n)$, $n\ge2$. Add the central tangent line $\mathbb RJ$ to the preceding decomposition. All $G$-g.o. cometrics are
      \begin{equation*}
C=aI_{\mathfrak v}\oplus bI_{\mathfrak p}
          \oplus tI_{\mathbb RJ},\qquad a,b>0,
       \qquad
       t=\begin{cases}
          a,& n\text{ even},\\[2pt]
          ((n-1)a+b)/n,& n\text{ odd}.
         \end{cases}
\end{equation*}
\item $\mathfrak g=\mathfrak{so}(2n)$,
$\mathfrak h=\mathfrak{su}(n)$, $n\ge3$ odd. All $G$-g.o. cometrics are $C=aI_{(\Lambda^2\mathbb C^n)_{\mathbb R}} \oplus bI_{\mathbb RJ},\qquad a,b>0.$
\item $\mathfrak g=\mathfrak{so}(8+q)$,
$\mathfrak h=\mathfrak{spin}(7)\oplus\mathfrak r$, $1\le q\le3$, $\mathfrak r\subset\mathfrak{so}(q)$. The embedding on $\mathbb R^8$ is the spin embedding, and the cometrics are exactly those in Lemma~\ref{lem614}.
\item $\mathfrak g=\mathfrak{so}(9)$,
$\mathfrak h=\mathfrak g_2\oplus\mathfrak r$, $\mathfrak r\subset\mathfrak{so}(2)$. The cometrics are exactly those in Lemma~\ref{lem614}.
\item $\mathfrak g=\mathfrak{so}(8)$,
$\mathfrak h=\mathfrak g_2$. Under a $Q$-orthogonal identification $\mathfrak m=\mathbb R^7\otimes\mathbb R^2$, every invariant positive cometric is $G$-g.o., and $C=I_{\mathbb R^7}\otimes P, \qquad P\in\operatorname{Sym}^+_2(\mathbb R).$
\end{enumerate}
In the first five entries, normality means that the two free parameters are equal. In the last, it means that $P$ is scalar.
\end{theorem}

\begin{proof}
First suppose the defining real representation is irreducible. The real-type and proper-unitary internal results, including the quaternionic case, leave only $\SU(n)\subset \SO(2n)$ and $\U(n)\subset \SO(2n)$. The latter quotient is irreducible symmetric. \cite[Theorems~1--2 and Table~1]{ChNiNi2023} gives precisely entry (3) for the former. It is normal when $n$ is even.

Suppose the defining representation is reducible. Lemma~\ref{lem69} gives a unique block of dimension $d>2q+1$. First consider the one situation where block separation was not proved, namely $q=1$ and $d$ odd. The isotropy algebra is semisimple and acts real-irreducibly on $\mathbb R^d$. If it is simple, \cite[Theorems~1--2 and Table~1]{ChNiNi2023} leaves only $\SO(8)/G_2$ among nonnormal metrics with this defining representation. If it has several simple ideals, the odd dimension forces all its tensor factors to be of real type and odd dimension. It lies in an orthogonal tensor subgroup $\SO(a)\otimes \SO(b)\subset \SO(ab)\subset \SO(ab+1)$, $a,b\ge3$. The two free tensor submodules inside $\mathfrak{so}(ab)\ominus\mathfrak h$ remain orthogonal submodules of $\mathfrak{so}(ab+1)\ominus\mathfrak h$. Thus the hereditary criterion, as included in Theorem~\ref{thm61}, gives normality. Thus this situation gives exactly entry (6).

In all other cases the off-diagonal block is preserved. Let $\mathfrak l$ be the isotropy projection on the large block. Compression supplies a $\SO(d)$-g.o. cometric on $\SO(d)/L$. The internal alternatives and Lemma~\ref{lem612} exclude every projection except full $\mathfrak{so}(d)$, the full unitary semisimple algebra, $G_2$ in dimension seven, and the spin $\Spin(7)$ in dimension eight.

For a full $\mathfrak{so}(d)$ projection, the inequality $q<d$ forces its simple ideal to be supported entirely on the first block. A homomorphism into $\mathfrak{so}(q)$ cannot be injective by dimension. The dimension-four case has $q=1$ and both $\mathfrak{su}(2)$ ideals are also independent. Theorem~\ref{thm65}, using the full large block as its first factor, gives normality.

For the unitary semisimple projection, write $d=2n$. Dominance gives $q<n$. The simple $\mathfrak{su}(n)$ ideal cannot map nontrivially to $\mathfrak{so}(q)$, again by dimension. Hence $\mathfrak h=\mathfrak{su}(n)\oplus\mathfrak r$ with $\mathfrak r\subset\mathbb RJ\oplus\mathfrak{so}(q)$. For $q\ge2$, Lemma~\ref{lem613} gives normality. For $q=1$, the residual algebra is zero or $\mathbb RJ$, giving entries (1) and (2). The dimension-four irreducible projections are $\mathfrak{so}(4)$, $\mathfrak{su}(2)$ and $\mathfrak u(2)$. The last two are the $n=2$ cases just listed, and dominance forces $q=1$.

For $G_2$, dominance gives $q\le2$. For spin $\Spin(7)$, it gives $q\le3$. The large simple ideals cannot couple to these small orthogonal algebras. The isotropy is therefore an independent product with an arbitrary residual $\mathfrak r\subset\mathfrak{so}(q)$. The $G_2$ case $q=1$ was already treated, and Lemma~\ref{lem614} gives entries (4) and (5). This exhausts every real defining decomposition. We now verify the coefficient conditions and sufficiency.

The scalar formulas in entries (4) and (5) were proved directly in Lemma~\ref{lem614}. Entry (1) follows from \cite[Theorem~6.2]{AlNi2009}. Entries (2), (3), and (6) are the explicit cometric forms of \cite[Theorem~2 and Table~1]{ChNiNi2023}. For odd $n$ in entry (2), its metric-eigenvalue equation $n\alpha_J^{-1}=(n-1)\alpha_{\mathfrak v}^{-1} +\alpha_{\mathfrak p}^{-1}$ is exactly the displayed linear cometric relation. The low-rank cases $\SO(5)/\SU(2)$ and $\SO(6)/\SU(3)$ are respectively $\Sp(2)/\Sp(1)$ and $\SU(4)/\SU(3)$ at the Lie-algebra level. Their formulas agree with the same list. In every listed family, the normal metrics are exactly the scalar cometrics.
\end{proof}

By Corollary~\ref{cor33}, all the displayed cometrics descend to closed isotropy groups with the same identity component. Finite central covers and quotients do not change the Lie-algebra g.o. equation.

\section{Regular isotropy and regular simple ideals in exceptional groups}\label{sec:unified-exceptional-regular}
In this section, we classify $G$-g.o. metrics with regular isotropy in an exceptional compact simple group $G$. We then treat arbitrary isotropy containing a regular simple ideal. The two root-$A_1$ cases require the stabilizer orthogonality criterion.

\begin{proposition}
\label{prop71}
Let $G$ be compact connected with exceptional simple Lie algebra. The following assertions hold for a closed proper subgroup $H$.
\begin{enumerate}
\item If $\mathfrak h$ is abelian, every $G$-g.o. metric is normal.
\item If $\operatorname{rank}\mathfrak h=\operatorname{rank}\mathfrak g$,
every $G$-g.o. metric is normal.
\item If $\mathfrak h$ is simple, a nonnormal $G$-g.o. metric exists
only for the standard pair $(\mathfrak g,\mathfrak h)=(\mathfrak e_6,\mathfrak{so}(10))$. For that pair,
$\mathfrak m=\mathbb R\oplus(\Delta^+_{10})_{\mathbb R},\qquad C=a\operatorname{Id}_{\mathbb R} \oplus b\operatorname{Id}_{(\Delta^+_{10})_{\mathbb R}}, \qquad a,b>0,$
and every such metric is $G$-g.o.  Here the half-spin module has complex dimension $16$ and real dimension $32$.
\item The same exceptional list holds when the real isotropy
representation has precisely two irreducible summands.
\end{enumerate}
\end{proposition}

\begin{proof}
The first assertion is the abelian-isotropy theorem of Souris~\cite{So2021}.  The second is the positive-Euler-characteristic classification of Alekseevsky--Nikonorov~\cite{AlNi2009}, Theorem~6.2. Its only nonnormal simple-ambient families have ambient algebras of type $B$ or $C$.  The third is the simple-isotropy classification of Chen--Nikolayevsky--Nikonorov~\cite{ChNiNi2023}, restricted to exceptional ambient algebras. No extra irreducibility assumption on the Riemannian manifold is needed here. For simple $\mathfrak g$ and nonzero $\mathfrak h$, no metric eigenspace in $\mathfrak h^\perp$ can contain a nonzero ideal of $\mathfrak g$. Thus the product-removal step in their Lemma~4 is absent, and the reduction at the beginning of their Section~5 gives the list in Theorem~2(A) directly. The fourth follows from the two-isotropy-summand classification of Chen--Nikonorov~\cite{ChNi2019}, Theorem~2, including the additional subgroup cases discussed in Chen--Nikolayevsky--Nikonorov, footnote~1. In the fourth assertion, two irreducible isotropy summands must not be replaced by the weaker assumption that $C$ has two eigenvalues.
\end{proof}

\begin{theorem}
\label{thm72}
Let $G$ be a compact connected Lie group with Lie algebra $\mathfrak g_2$. For every closed subgroup $H\subset G$, every $G$-g.o. Riemannian metric on $G/H$ is normal.
\end{theorem}

\begin{proof}
The cases of abelian and simple $\mathfrak h$ follow from Proposition~\ref{prop71}. Otherwise, compact reductivity gives
$\mathfrak h=\mathfrak z(\mathfrak h) \oplus\mathfrak h_1\oplus\cdots\oplus\mathfrak h_s,$
where $\mathfrak h_i$ are simple.  Since $\mathfrak h$ is nonabelian and nonsimple, either $s\geq2$, or $s=1$ and $\dim\mathfrak z(\mathfrak h)\geq1$.  Hence $\operatorname{rank}\mathfrak h\geq2$. The opposite bound follows from $\operatorname{rank}\mathfrak g_2=2$. Thus $\mathfrak h$ has maximal rank, and the second assertion of Proposition~\ref{prop71} applies.
\end{proof}

\begin{theorem}\label{thm73}
Let $G$ be compact connected simple of type $G_2,F_4,E_6,E_7$, or $E_8$, and let $H$ be a closed connected proper regular subgroup.  Choose a maximal torus $T$ normalizing $H$, and put
$K=HT,\qquad \mathfrak k=\mathfrak h+\mathfrak t, \qquad \mathfrak z=\mathfrak z(\mathfrak k), \qquad \mathfrak a=\mathfrak h\cap\mathfrak z.$
Then $K$ is a closed connected subgroup of maximal rank, and
\begin{equation*}
 \mathfrak h=[\mathfrak k,\mathfrak k]\oplus\mathfrak a,
 \qquad \mathfrak v=\mathfrak z\ominus\mathfrak a,
 \qquad \mathfrak p=\mathfrak k^\perp,
 \qquad \mathfrak m=\mathfrak v\oplus\mathfrak p.
\end{equation*}
All orthogonality is for $Q=-B_{\mathfrak g}$. Let $\mathfrak s$ be a principal stabilizer Lie algebra for the compact isotropy representation of $K$ on $\mathfrak p$, and define
$\mathfrak e=P_{\mathfrak z}\mathfrak s, \qquad L=P_{\mathfrak v}\mathfrak e.$
Then the complete set of $G$-g.o. cometrics is
\begin{equation}\label{formula31}
 C=cI_{L^\perp}\oplus S_L,
 \qquad c>0,\qquad S_L\text{ arbitrary positive self-adjoint on }L.
\end{equation}
In particular, the entire classification for regular $H$ reduces to the central projection of one principal stabilizer for each maximal-rank pair $(G,K)$.  Simple ideals of the normalizer outside $\mathfrak h$ do not provide additional metric parameters.
\end{theorem}

\begin{proof}
Since $T$ normalizes $H$, the product $HT$ is a compact subgroup. The root-space description of a regular subalgebra shows that $[\mathfrak k,\mathfrak k]\subset\mathfrak h$, and hence gives the stated decomposition.  Equivalently, $\mathfrak v=\mathfrak t\ominus(\mathfrak t\cap\mathfrak h)$. Every element of $\mathfrak v$ centralizes $\mathfrak h$.

Normalizer invariance of a $G$-g.o. cometric makes $C$ commute with $\operatorname{ad}\mathfrak t$.  The zero-weight subspace of $\mathfrak m$ is $\mathfrak v$, while $\mathfrak p$ is the sum of its nonzero root weights.  Thus $C$ preserves both subspaces. The descent lemma makes $C|_{\mathfrak p}$ a $G$-g.o. cometric on $G/K$.  This space has positive Euler characteristic. The Alekseevsky--Nikonorov classification of compact $G$-g.o. spaces of positive Euler characteristic has no nonnormal case with exceptional simple ambient group.  Consequently
$C|_{\mathfrak p}=cI,\qquad C|_{\mathfrak v}=S, \qquad D=S-cI.$
This step uses the positive-Euler-characteristic classification, not merely the classification of generalized flag manifolds. A maximal-rank subgroup need not be a torus centralizer.

For $u=y+x$, $y\in\mathfrak v$, $x\in\mathfrak p$, the g.o. equation is exactly
$[Dy,x]=[Z,x],\qquad Z\in\mathfrak h,$
because $\mathfrak v$ is abelian and centralizes $\mathfrak h$. For $T_0\in\mathfrak v$ and fixed $x$, this equation is solvable if and only if
$T_0\in\mathfrak a+P_{\mathfrak z}(\mathfrak k_x).$
Indeed, $T_0-Z\in\mathfrak k_x$ is the required stabilizer element. Conversely its semisimple part can be absorbed in $[\mathfrak k,\mathfrak k]\subset\mathfrak h$.

For principal $x$, the algebras $\mathfrak k_x$ are conjugate under $K$.  Since $K$ acts trivially on its center, $P_{\mathfrak z}(\mathfrak k_x)=\mathfrak e$ is independent of the principal point.  Moreover
\begin{equation}\label{formula32}
 \mathfrak e\subset P_{\mathfrak z}(\mathfrak k_x)
 \quad\text{for every }x\in\mathfrak p.
\end{equation}
To prove this, choose principal $x_j\to x$ and write $\mathfrak k_{x_j}=\operatorname{Ad}(k_j)\mathfrak s$. After taking a subsequence, compactness gives $k_j\to k$. For every $A\in\mathfrak s$, passing to the limit in $[\operatorname{Ad}(k_j)A,x_j]=0$ gives $\operatorname{Ad}(k)\mathfrak s\subset\mathfrak k_x$. Taking central projections proves \eqref{formula32}.

It follows that the g.o. condition is equivalent to
$\operatorname{im}D\subset \mathfrak v\cap(\mathfrak a+\mathfrak e) =P_{\mathfrak v}\mathfrak e=L.$
Self-adjointness of $D$ then makes it vanish on $L^\perp$. Conversely, every operator in \eqref{formula31} has its image condition, and \eqref{formula32} supplies the g.o. compensator for every vector, including nonprincipal vectors. This proves both necessity and sufficiency.
\end{proof}

\begin{proposition}
\label{prop74}
For a maximal-rank $K\subset G$ with center $\mathfrak z$, let $L=C_G(Z(K)^0)^0$.  Then $K\subset L$ and $\mathfrak z(\mathfrak l)=\mathfrak z$. If there is a vector $x\in\mathfrak l^\perp$ such that
$P_{\mathfrak z}(\mathfrak l_x)=0,$
then the subspace $\mathfrak e$ in Theorem~\ref{thm73} is zero. Thus a zero central-projection entry for a Levi subgroup disposes of all maximal-rank subgroups of that Levi subgroup with the same center.
\end{proposition}

\begin{proof}
In root notation, the roots of $\mathfrak k$ span $\mathfrak z^\perp\subset\mathfrak t^*$. The roots of $\mathfrak l$ are all roots vanishing on $\mathfrak z$. Their span is therefore the same space. This proves equality of the centers.  Since $\mathfrak k_x\subset\mathfrak l_x$ for the specified $x$, its central projection is zero.  By \eqref{formula32}, the principal central projection $\mathfrak e$ is then zero as well.
\end{proof}

\begin{corollary}\label{cor75}
Suppose $[\mathfrak k,\mathfrak k]$ is simple.  Then $\mathfrak e=0$, except for the standard maximal-rank pair
$(\mathfrak g,\mathfrak k) =(\mathfrak e_6,\mathfrak{so}(10)\oplus\mathbb R).$
In this exceptional pair $\mathfrak e=\mathfrak z$. Hence all regular intermediate central quotients in these cases are normal, apart from $E_6/\Spin(10)$, whose two positive parameters on the real spin module and the residual central line are arbitrary.
\end{corollary}

\begin{proof}
Apply Theorem~\ref{thm73} with isotropy algebra $[\mathfrak k,\mathfrak k]$ and no isotropy torus.  A nonzero $\mathfrak e$ would supply a nonnormal $G$-g.o. metric with simple isotropy.  The Chen--Nikolayevsky--Nikonorov classification leaves exactly $E_6/\Spin(10)$ when the ambient group is exceptional. For that space, all positive spin/central parameters are known to be $G$-g.o., so the one-dimensional central projection is the entire center.  Adding its central circle to the isotropy gives the irreducible symmetric space $E_6/(\Spin(10)\U(1))$.
\end{proof}

\begin{proposition}
\label{prop76}
Let $V_1,\ldots,V_s$ be complex $K$-submodules of $\mathfrak p$, viewed as real modules, such that the center acts on $V_j$ by the scalar $i\alpha_j$, with $\alpha_j\in\mathfrak z^*$. Suppose the semisimple group $[K,K]$ has a nonzero homogeneous complex polynomial invariant $f_j$ of positive degree on $V_j$. Then
$\mathfrak e\subset\bigcap_j\ker\alpha_j.$
In particular, if these weights span $\mathfrak z^*$, every regular quotient arising from this $K$ is normal.
\end{proposition}

\begin{proof}
Choose a principal $x$ whose component $v_j\in V_j$ satisfies $f_j(v_j)\ne0$ for every $j$.  These are finitely many nonempty open dense conditions.  If $A+Z\in\mathfrak k_x$, with $A\in[\mathfrak k,\mathfrak k]$ and $Z\in\mathfrak z$, then
$Av_j+i\alpha_j(Z)v_j=0.$
Apply the differential of $f_j$ at $v_j$. Semisimple invariance makes the first term zero, while homogeneity makes the second term $i\alpha_j(Z)\deg(f_j)f_j(v_j)$. Hence $\alpha_j(Z)=0$ for every $j$, proving the assertion.
\end{proof}

\begin{proposition}
\label{prop77}
Let $\mathfrak l=\mathfrak{so}(10)\oplus\mathbb Rz$ be the standard Levi subalgebra of $\mathfrak e_6$.  Suppose $\mathfrak k\subsetneq\mathfrak l$ is regular of full rank in $\mathfrak e_6$ and $\mathfrak z(\mathfrak k)=\mathbb Rz$. Then the principal central projection for $K$ acting on $\mathfrak k^\perp\subset\mathfrak e_6$ is zero. Consequently all regular isotropy algebras $[\mathfrak k,\mathfrak k]\oplus\mathfrak a$, $\mathfrak a\subset\mathbb Rz$, have only normal $G$-g.o. metrics.
\end{proposition}

\begin{proof}
The semisimple part of $\mathfrak k$ is a proper full-rank closed root subsystem of $D_5$.  Such a subsystem is $D_3+D_2$, up to the Weyl group. This particular root-system assertion can be seen in the model $D_5=\{\pm e_i\pm e_j:i\ne j\}$. Form the graph on the five coordinate indices, joining two indices whenever a root of the subsystem uses them. On a connected component with $r$ coordinates, closure under root reflections and root addition gives either $D_r$, of rank $r$, or a signed copy of $A_{r-1}$, of rank $r-1$. Full rank excludes every $A$ component and every unused coordinate.  Properness requires at least two $D$ components, each using at least two coordinates.  The only partition of five is $3+2$. Thus the semisimple group is locally $\SU(4)\times\SU(2)\times\SU(2)$.

The complement of $\mathfrak l$ in $\mathfrak e_6$ is the underlying real module of a complex half-spin module $\Delta^+_{10}$, and $z$ acts on this entire complex module by one nonzero scalar $i\alpha(z)$.  Its restriction is
$\Delta^+_{10}|_{\Spin(6)\times\Spin(4)} =(\mathbb C^4\otimes\mathbb C^2\otimes\mathbb C) \oplus((\mathbb C^4)^*\otimes\mathbb C\otimes\mathbb C^2).$
Indeed, the half-spin tensor formula for an orthogonal direct sum is $\Delta^+_{10}=(\Delta^+_6\otimes\Delta^+_4) \oplus(\Delta^-_6\otimes\Delta^-_4)$. The half-spin modules in dimensions six and four are respectively $\mathbb C^4,(\mathbb C^4)^*$ and the two defining $\SU(2)$ modules.

Write a vector as $(A,B)$ with $A\in\mathbb C^4\otimes\mathbb C^2$ and $B\in(\mathbb C^4)^*\otimes\mathbb C^2$, represented by two $4\times2$ column matrices.  The polynomial
$f(A,B)=\det(B^tA)$
is a nonzero invariant for $\SL(4,\mathbb C)\times\SL(2,\mathbb C) \times\SL(2,\mathbb C)$. The $\SL(4)$ transformations cancel in the contraction, and the two $\SL(2)$ factors act on its two sides with determinant one.  Nonvanishing is witnessed by taking both pairs of columns to be the first two coordinate vectors and covectors. The polynomial has degree four, and hence nonzero weight $4\alpha$ under the Levi center, which acts with the same character on both summands of $\Delta^+_{10}$. Proposition~\ref{prop76} now annihilates the one-dimensional central projection.
\end{proof}

\begin{proposition}
\label{prop78}
For nested connected maximal-rank subgroups $K\subset L\subset G$, write $\mathfrak e_K,\mathfrak e_L$ for the principal central projections of their isotropy representations.  Then
$P_{\mathfrak z(\mathfrak l)}\mathfrak e_K \subset\mathfrak e_L.$
\end{proposition}

\begin{proof}
Both subgroups contain the selected maximal torus, so their root descriptions give $\mathfrak z(\mathfrak l)\subset\mathfrak z(\mathfrak k)$. Choose $x\in\mathfrak l^\perp$ principal for $L$. For $T\in\mathfrak e_K$, \eqref{formula32} gives $A\in\mathfrak k_x$ with $P_{\mathfrak z(\mathfrak k)}A=T$. Since $\mathfrak k_x\subset\mathfrak l_x$ and the central subspaces are nested, $P_{\mathfrak z(\mathfrak l)}T =P_{\mathfrak z(\mathfrak l)}A\in\mathfrak e_L$.
\end{proof}

\begin{proposition}
\label{prop79}
Suppose $K\subset M\subsetneq G$, where $M$ is connected semisimple, and the compact isotropy representation of $M$ on $\mathfrak u^\perp$, with $\mathfrak u=\operatorname{Lie}(M)$, has an abelian principal stabilizer Lie algebra. Then $\mathfrak e_K=0$.
\end{proposition}

\begin{proof}
Let $T\in\mathfrak e_K\subset\mathfrak k\subset\mathfrak u$. For a compact semisimple Lie algebra, every element is orthogonal to a Cartan subalgebra.  Indeed, the linear convexity theorem of Kostant~\cite{Ko1973} identifies projection of an adjoint orbit onto a fixed Cartan with the convex hull of a Weyl orbit. The Weyl average is zero. An abelian subalgebra of a compact Lie algebra lies in a Cartan. Conjugating a principal vector in $\mathfrak u^\perp$ by $M$, we may therefore obtain $x$ with $\mathfrak u_x\subset T^\perp$. By \eqref{formula32}, some $A\in\mathfrak k_x$ has central projection $T$. Thus $A=T+A_0$ with $A_0\in[\mathfrak k,\mathfrak k]$. But $A\in\mathfrak u_x\subset T^\perp$, whereas $Q(T,A)=Q(T,T)$ because $T$ is central in $\mathfrak k$. Hence $T=0$.
\end{proof}

\begin{proposition}
\label{prop710}
For the standard maximal Levi subalgebra $\mathfrak l=\mathfrak{so}(10)\oplus\mathfrak{su}(2)\oplus\mathbb R$ of $\mathfrak e_7$, the principal central projection is zero.
\end{proposition}

\begin{proof}
The extended root diagram places this Levi in the semisimple subalgebra $\mathfrak m=\mathfrak{so}(12)\oplus\mathfrak{su}(2)$, with the standard inclusion $\mathfrak{so}(10)\oplus\mathfrak{so}(2)\subset\mathfrak{so}(12)$. The submodule $\mathfrak m\ominus\mathfrak l$ of $\mathfrak l^\perp\subset\mathfrak e_7$ is the underlying real module of $\mathbb C^{10}$.  Its semisimple action is the complexified standard $\Spin(10)$ representation. The other $\SU(2)$ factor acts trivially, and the central circle has nonzero scalar weight.  The polynomial $f(v)=v^tv$ is a nonzero quadratic invariant of the semisimple part.  Proposition~\ref{prop76} therefore kills the entire one-dimensional center.
\end{proof}

\begin{theorem}
\label{thm711}
Let $G$ be compact connected simple of exceptional type, and let $H$ be a closed subgroup whose Lie algebra is proper and regular. All $G$-g.o. metrics on $G/H$ are normal, except when the Lie algebra pair is the standard $(\mathfrak e_6,\mathfrak{so}(10))$. In that case the tangent space is the real half-spin module plus a central line, and every positive scalar on each of these two summands gives a $G$-g.o. cometric.  These are all $G$-g.o. metrics.
\end{theorem}

\begin{proof}
We use the complete root-deletion calculation in Table~\ref{tab:exceptional-levi-deletions}. For a deleted simple root $\alpha_i$, write $L_i$ for the maximal Levi subgroup and $M_i$ for the full-rank root subgroup obtained by replacing $\alpha_i$ with the negative highest root when the highest-root coefficient is greater than one. The methods A, B and C in the table refer to Proposition~\ref{prop79}, Corollary~\ref{cor75} and Proposition~\ref{prop710}, respectively. N denotes the nonzero $D_5$ entry in $E_6$.

The abelian-stabilizer assertion in all rows marked A follows from Benoist--Kobayashi~\cite[Theorems~1.6, 4.1 and Table~2]{BeKo2021}.  Among proper semisimple full-rank subgroups of exceptional groups, the only nonabelian generic-stabilizer pairs in that table are
$(G_2,A_2)$, $(F_4,B_4)$, $(F_4,D_4)$, $(E_7,D_6+A_1)$ and $(E_8,E_7+A_1)$.
None occurs in a row marked A.  Their result for complexified isotropy modules applies to compact real forms.  The complex generic-stabilizer stratum contains a nonempty Zariski-open set. The compact real form is Zariski dense and its principal stratum is open dense, so these two strata meet.  At a real point the stabilizer complexifies to the complex stabilizer. Hence the compact principal stabilizer is abelian as well. Proposition~\ref{prop79} completes each row A.  Methods B and C complete the other zero entries.

Let now $L_J$ be any proper Levi subgroup, with nonempty set $J$ of deleted simple roots.  Its center is the span of the fundamental coweights $\omega_i^\vee$, $i\in J$, and $L_J\subset L_i$ for every $i\in J$. If $G\ne E_6$, every maximal-Levi entry has zero central projection.  Proposition~\ref{prop78} therefore makes $\mathfrak e_{L_J}$ orthogonal to each $\omega_i^\vee$, $i\in J$.  Since it lies in their span, it is zero.

For $G=E_6$, if $J$ contains any of $2,3,4,5$, then $L_J\subset L_i\subset M_i$ for a row A, so Proposition~\ref{prop79} gives zero central projection directly.  If $J=\{1,6\}$, the derived Levi algebra is simple of type $D_4$, and Corollary~\ref{cor75} again gives zero.  The only nonzero Levi entries are consequently the two conjugate $D_5+\mathbb R$ entries.

Finally take the maximal-rank subgroup $K=H^0T$ from Theorem~\ref{thm73}. If its center is zero there are no possible deformations. Otherwise pass to $L=C_G(Z(K)^0)^0$. Proposition~\ref{prop74} excludes $K$ unless $L=D_5+\mathbb R$ in $E_6$. Inside that exceptional Levi, Proposition~\ref{prop77} excludes every proper $K$ with the same center. For $K=L$, the isotropy algebra is either $D_5$ or $D_5+\mathbb R$.  The latter is maximal rank and normal, while the former has exactly the known two-parameter $G$-g.o. family. This proves the connected classification.

Corollary~\ref{cor33} gives the same $G$-g.o. cone for $H$ and $H^0$.
\end{proof}

\begin{remark}
The abelian-stabilizer classification is \cite[Theorems~1.6, 4.1 and Table~2]{BeKo2021}, with numbering as in the authors' preprint. The orthogonal-Cartan assertion follows from Kostant's linear convexity theorem~\cite{Ko1973}. The argument uses a semisimple overgroup, its conjugate Cartan subalgebras and the central-projection characterization. Abelianness of the principal isotropy alone does not suffice for an arbitrary homogeneous pair.
\end{remark}

\begin{remark}
Suppose that $G/K$ has only normal $G$-g.o. metrics. For $H\subsetneq K$, descent first requires the metric operator to preserve $\mathfrak k^\perp$. This is not automatic for an arbitrary nonregular $H$.  Even when it does preserve that space, normality of the descended metric fixes only the horizontal restriction.  The fiber operator on $\mathfrak k\ominus\mathfrak h$ remains to be determined. The example $\Spin(10)\subset\Spin(10)\U(1)\subset E_6$ demonstrates this failure directly. The base is irreducible symmetric and normal, while the total space admits the nonnormal central deformation.
\end{remark}

\begin{lemma}\label{lem714}
Every adjoint orbit of compact $E_7$ meets the orthogonal complement of its standard regular $D_4$ subalgebra.
\end{lemma}
\begin{proof}
By Proposition~\ref{propC1}, all seven fundamental representations contain nonzero fixed vectors for this subgroup. Lemma~\ref{lem318} proves the assertion.
\end{proof}

\begin{theorem}
\label{thm715}
Let $G$ be compact simple of type $E_7$ or $E_8$. If $\mathfrak h$ has a simple ideal which is a root $A_1$ of $\mathfrak g$, then every $G$-g.o. metric on $G/H$ is normal. The complementary isotropy in its centralizer is arbitrary.
\end{theorem}

\begin{proof}
The two ideal-centralizer envelopes are, respectively,
$\mathfrak k=\mathfrak{so}(12)\oplus\mathfrak{su}(2), \qquad \mathfrak k=\mathfrak e_7\oplus\mathfrak{su}(2),$
where the specified root ideal is the last factor. Both base quotients are irreducible compact symmetric spaces, so their invariant cometrics are scalar. Moreover $\mathfrak h=\mathfrak{su}(2)\oplus\mathfrak r$ with $\mathfrak r$ contained in the first factor, and every vector $T\in\mathfrak w=\mathfrak k\ominus\mathfrak h$ lies in that first factor.

We identify a stabilizer in each horizontal representation by an explicit root calculation.  If $\beta_1,\ldots,\beta_4$ are strongly orthogonal noncompact roots, choose unit vectors $X_j$ in their compact real root planes inside $\mathfrak p$. Their root $\mathfrak{su}(2)$ algebras commute. Simultaneous rotations in these four algebras conjugate a generic $x=\sum t_jX_j$ to $\sum t_jH_{\beta_j}$. Choose the real coefficients outside the finitely many hyperplanes that arise from roots not orthogonal to every $\beta_j$. Then the dimension of $\mathfrak g_x$ is
$\operatorname{rank}\mathfrak g+ \#\{\gamma\in\Phi(\mathfrak g): (\gamma,\beta_j)=0\text{ for all }j\}.$

For $E_7$ use the roots in $\mathbb R^8$ orthogonal to $e_7+e_8$.  The compact roots are $\pm e_i\pm e_j$ ($1\le i<j\le6$) and $\pm(e_7-e_8)$. The noncompact roots are the half-roots $\tfrac12(s_1,\ldots,s_6,s_7,-s_7)$ with an even number of minus signs.  Take
\begin{equation*}
 \begin{array}{c|rrrrrrrr}
 2\beta_1&-1& 1& 1& 1& 1& 1&1&-1\\
 2\beta_2& 1&-1&-1&-1& 1& 1&1&-1\\
 2\beta_3& 1&-1& 1& 1&-1&-1&1&-1\\
 2\beta_4&-1& 1&-1&-1&-1&-1&1&-1.
 \end{array}
\end{equation*}
These roots are pairwise orthogonal and hence strongly orthogonal. The roots orthogonal to all four are exactly
$\pm(e_1+e_2),\qquad \pm(e_3-e_4),\qquad\pm(e_5-e_6).$
For compact roots this follows by comparing the six columns of the displayed matrix.  No half-root can be orthogonal to all four, because $\sum_j\beta_j=2(e_7-e_8)$. The corresponding root algebra $\mathfrak q=3A_1$ lies in $\mathfrak{so}(12)$ and centralizes every $X_j$. Since $\dim\mathfrak g_x=7+6=13$, while $\mathfrak q\subset\mathfrak k_x$ has dimension $9$ and $\operatorname{span}\{X_1,\ldots,X_4\}\subset\mathfrak p_x$ has dimension $4$, we obtain $\mathfrak k_x=\mathfrak q$.

This $3A_1$ preserves three mutually orthogonal real four-planes $A,B,C$ in $\mathbb R^{12}$, hence lies in $\mathfrak{so}(A)\oplus\mathfrak{so}(B)\oplus\mathfrak{so}(C)$. Choose orthonormal bases $a_i,b_i,c_i$, $1\le i\le4$. The six rotations supported on the disjoint two-planes
$(a_1,b_1),\ (a_2,b_2),\ (b_3,c_1),\ (b_4,c_2),\ (c_3,a_3),\ (c_4,a_4)$
span a Cartan subalgebra of $\mathfrak{so}(12)$ orthogonal to $\mathfrak q$. Every $T\in\mathfrak{so}(12)$ is conjugate into this Cartan. Equivalently some conjugate of $\mathfrak k_x$ is orthogonal to $T$. Lemma~\ref{lem312} proves the $E_7$ case.

For $E_8$ use the standard roots $\pm e_i\pm e_j$ together with the half-roots having an even number of minus signs.  The external $A_1$ is generated by $\alpha=e_7+e_8$. The $E_7$ roots are orthogonal to $\alpha$. Choose the four noncompact roots
$\beta_1=e_5+e_7,\quad\beta_2=e_5-e_7,\quad \beta_3=e_6+e_8,\quad\beta_4=e_6-e_8.$
They are strongly orthogonal, and their common orthogonal roots are exactly the $D_4$ roots $\pm e_i\pm e_j$ with $i,j\le4$. Thus $\dim\mathfrak g_x=8+24=32$. The corresponding compact root algebra $\mathfrak q=D_4$ has dimension $28$ and centralizes the four-dimensional span of the $X_j$. Hence again $\mathfrak k_x=\mathfrak q$. This is precisely the standard $D_4\subset D_6\subset E_7$ used in Lemma~\ref{lem714}. For every $T\in\mathfrak e_7$ that lemma gives a conjugate of $\mathfrak q$ orthogonal to $T$. Lemma~\ref{lem312} proves the $E_8$ case.
\end{proof}

\begin{corollary}
\label{cor716}
Let $G$ be compact simple of exceptional type and suppose $\mathfrak h$ contains a nonzero simple ideal $\mathfrak s$ which is regular in $\mathfrak g$.  The rest of $\mathfrak h$ need not be regular.  Then every $G$-g.o. metric is normal, except when $\mathfrak g=\mathfrak e_6$ and $\mathfrak h=\mathfrak{so}(10)$ in the standard embedding. That case has the known two-parameter $G$-g.o. family.
\end{corollary}

\begin{proof}
Choose a maximal torus normalizing $\mathfrak s$. The ideal-centralizer algebra $\mathfrak l=\mathfrak s\oplus\mathfrak c_{\mathfrak g}(\mathfrak s)$ contains this Cartan subalgebra, so it is regular and of full rank. The positive-Euler-characteristic g.o. classification gives normality of every metric on $G/L$. Theorem~\ref{thm711} and its principal-central-projection calculation give zero central projection for the $L$-isotropy module except when $(\mathfrak g,\mathfrak l)=(\mathfrak e_6,D_5+\mathbb R)$.

Let $\mathfrak l_s=[\mathfrak l,\mathfrak l]$. As an $\mathfrak l_s$-module, $\mathfrak g/\mathfrak l_s$ is the direct sum of $\mathfrak p=\mathfrak l^\perp$ and the trivial module $\mathfrak z(\mathfrak l)$, so generic stabilizers on these two modules agree.  Benoist--Kobayashi's exceptional list of nonabelian generic stabilizers therefore applies. Outside that list, Lemma~\ref{lem316} and Theorem~\ref{thm313} prove normality.

We check all entries of that list which can be regular. For $G_2$, the entry $A_2$ has full rank. For $F_4$, the entries $B_4$ and $D_4$ have full rank. Since $\mathfrak s$ is a simple ideal of $\mathfrak l_s$, in these cases it is the whole displayed simple algebra. Thus $H$ itself has full rank, and its $G$-g.o. metrics are normal.

For $E_6$, the entries $F_4$ and $B_4$ cannot be regular in the simply laced root system $E_6$. The remaining entry is $D_5$, whose centralizer is one-dimensional. Here $\mathfrak s=D_5$ and $\mathfrak h$ is either $D_5$ or $D_5+\mathbb R$, giving exactly the stated exception.

For $E_7$, the potentially regular entries are $D_6+A_1$, $D_6$, and $E_6$. The root centralizer of $D_6$ is $A_1$, so its full ideal-centralizer algebra is $D_6+A_1$. If $\mathfrak s=D_6$, then $\mathfrak h=D_6+\mathfrak r$ with $\mathfrak r\subset\mathfrak{su}(2)$. The latter algebra is zero, a toral line, or all of $\mathfrak{su}(2)$. The first case has simple isotropy and the others have full rank.  All are normal by the known classifications. If $\mathfrak s=E_6$, its centralizer is a line, so the same simple-isotropy/full-rank argument applies. The only alternative left by the abelian-stabilizer argument is the root $A_1$ ideal in $D_6+A_1$. It is closed, with arbitrary residual isotropy in $D_6$, by Theorem~\ref{thm715}.

For $E_8$, the entries are $E_7+A_1$ and $E_7$. The root centralizer of $E_7$ is $A_1$. If $\mathfrak s=E_7$, the simple-isotropy/full-rank argument again applies to $E_7+\mathfrak r$, $\mathfrak r\subset\mathfrak{su}(2)$. The remaining root $A_1$ is likewise closed, with arbitrary residual isotropy in $E_7$, by Theorem~\ref{thm715}. This exhausts the nonabelian-stabilizer list.
\end{proof}

\section{Subgroups of exceptional product subgroups}\label{sec:unified-exceptional-envelopes}
In this section, we prove normality for every isotropy subgroup in the product subgroups used below. These results apply to arbitrary subgroups of the indicated product subgroups and provide the cases used in centralizer propagation.

\begin{corollary}
\label{cor81}
Suppose $H\subset E_6$ has a simple ideal which is one of the factors in the standard subgroups
$A_5A_1,\qquad A_2A_2A_2,\qquad G_2A_2.$
Then every $E_6$-g.o. metric on $E_6/H$ is normal, even when the residual isotropy inside the complementary factor or factors is nonregular or diagonally coupled.
\end{corollary}

\begin{proof}
For each specified factor, its full centralizer is the product of the other displayed factors.  The corresponding isotropy modules are, respectively, $\mathbb C^2\otimes\Lambda^3\mathbb C^6$, the underlying real module of $\mathbb C^3\otimes\mathbb C^3\otimes\mathbb C^3$, and $\mathbb R^7\otimes\mathfrak{su}(3)$. Every factor acts nontrivially and these modules have no fixed vectors for any one factor. The corresponding ideal-centralizer algebras are semisimple. Each displayed isotropy module is real irreducible, so the base quotient has only scalar invariant cometrics. None of these pairs is in the exceptional nonabelian-generic-stabilizer list of Benoist--Kobayashi~\cite{BeKo2021}. Theorem~\ref{thm313} applies.
\end{proof}

\begin{remark}
Normality on $G/S$ does not generally imply normality on $G/L$ for $S\triangleleft L$, even when $G$ and $S$ are simple. For example, take $G=\SU(5)$, $S=\SU(2)$ on the smaller standard block, and $L=\SU(3)\times\SU(2)$. The simple-isotropy classification makes every metric on $G/S$ normal, whereas $G/L$ admits the nonnormal central deformation associated with the dominant $\SU(3)$ block. This $L$ is not the full ideal-centralizer algebra. That algebra also contains a scalar circle.  The example rules out a general ascending inference, and the first hypothesis of Theorem~\ref{thm313} is retained explicitly.
\end{remark}

\paragraph{\texorpdfstring{Two nonregular exclusion results in $F_4$}{Two nonregular exclusion results in F4}}
\label{sec:f4-g2-so3-verified}

Let $G$ be the compact group of type $F_4$, and let $K\subset G$ be its subgroup with Lie algebra $\mathfrak k=\mathfrak g_2\oplus\mathfrak{so}(3)$. Its isotropy representation is
\begin{equation}\label{formula33}
 \mathfrak g=\mathfrak g_2\oplus\mathfrak{so}(3)
      \oplus(V\otimes W),
 \qquad \dim V=7,\quad\dim W=5,
\end{equation}
where $V$ is the standard real $G_2$-module and $W=\operatorname{Sym}_0^2(\mathbb R^3)$ is the irreducible five-dimensional $\SO(3)$-module.  The tensor module is irreducible.  This branching is given in \cite[p.~21]{BeElMa2009}. All orthogonal complements below use an invariant inner product $Q$ on $\mathfrak g$.

\begin{proposition}\label{prop83}
For the inclusion $\mathfrak h=\mathfrak g_2\oplus\mathfrak{so}(2) \subset\mathfrak k$, every $F_4$-g.o. metric on $G/H$ is normal.
\end{proposition}

\begin{proof}
Put $\mathfrak f=\mathfrak{so}(3)\ominus\mathfrak{so}(2)$ and $\mathfrak p=V\otimes W$.  The $\mathfrak g_2$ action on $\mathfrak f$ is trivial, whereas $\mathfrak p$ has no $\mathfrak g_2$-fixed vector.  Thus an invariant cometric preserves $\mathfrak f$ and $\mathfrak p$. In particular, there is no mixing assumption here.  For a $G$-g.o. cometric, descent to $G/K$ and irreducibility of $\mathfrak p$ imply $C|_{\mathfrak p}=cI$.  The two-dimensional circle module $\mathfrak f$ is irreducible, so $C|_{\mathfrak f}=aI$.

Choose $0\ne Y\in\mathfrak f$ and $X=v\otimes w\in\mathfrak p$. Write the g.o. compensator as $Z=A+zT$, with $A\in\mathfrak g_2$ and $\mathfrak{so}(2)=\mathbb RT$. The $\mathfrak k$ component of the g.o. equation for $Y+X$ is $z[T,Y]=0$, so $z=0$.  Its $\mathfrak p$ component is
$(a-c)v\otimes(Y\cdot w)=(A\cdot v)\otimes w.$
Take $v\ne0$ and choose $w$ with $Y\cdot w\ne0$. The action is skew-adjoint, so $Y\cdot w\perp w$. Projecting onto $V\otimes w^\perp$ gives $a=c$.
\end{proof}

\begin{proposition}\label{prop84}
Let $L\subset K$ be the subgroup with Lie algebra $\mathfrak l=\mathfrak{so}(4)\oplus\mathfrak{so}(3)$, where $\mathfrak{so}(4)\subset\mathfrak g_2$ is the standard maximal-rank subalgebra.  For every closed subgroup $H\subset L$, all $F_4$-g.o. metrics on $G/H$ are normal.  In particular, the conclusion allows arbitrary toral directions, diagonal couplings between simple factors, and disconnected subgroups.
\end{proposition}

\begin{proof}
Write $\mathfrak{so}(4)=\mathfrak a\oplus\mathfrak b$, with $\mathfrak a\cong\mathfrak b\cong\mathfrak{su}(2)$, and let $\mathfrak c=\mathfrak{so}(3)$ be the external factor.  The standard branching in $G_2$ is
$V=V_3\oplus V_4, \qquad \mathfrak g_2=\mathfrak a\oplus\mathfrak b\oplus U_8.$
Here $V_3$ is the adjoint module of $\mathfrak a$, trivial on $\mathfrak b$, and $V_4$ is the standard real four-dimensional $\mathfrak{so}(4)$-module.  The complexification of $U_8$ is, up to interchanging the factors, $\mathbb C^2\otimes\operatorname{Sym}^3(\mathbb C^2)$. In particular, every nonzero element of either simple factor acts invertibly on $U_8$. The weights of its two-dimensional or four-dimensional factor are nonzero.  These statements also follow directly from the short and long orthogonal root $A_1$-subsystems of $G_2$.

Define orthogonal $L$-invariant subspaces
$E=V_4\otimes W, \qquad F=(V_3\otimes W)\oplus U_8.$
They exhaust $\mathfrak l^\perp$.  We show that both contain a vector with zero infinitesimal $\mathfrak l$-stabilizer.

We use the following elementary tensor observation.  Identify $U\otimes W$ with maps $T:U\to W$ using invariant Euclidean inner products.  A stabilizer satisfies $BT=TA$, where $A$ and $B$ are skew-adjoint infinitesimal actions.  If $T$ is injective and has distinct positive singular values, then $A$ commutes with $T^*T$ and is therefore zero.  Consequently $B$ kills $\operatorname{im}T$.  For the five-dimensional spin-two representation of $\mathfrak c$, every nonzero $B$ has a one-dimensional kernel. Its complex weights are $-2,-1,0,1,2$, up to a common nonzero scale.

Apply the observation first with $U=V_4$.  A full-rank $T:V_4\to W$ with distinct singular values forces the $\mathfrak{so}(4)$ action to vanish and forces the $\mathfrak c$ action to kill a four-dimensional space. Both are zero, since the former action is faithful and the latter has kernel dimension one when nonzero.  Hence $E$ contains a free vector for $\mathfrak l$.

Apply the same observation to a full-rank map $V_3\to W$. Its stabilizer in $\mathfrak a\oplus\mathfrak c$ is zero, so its stabilizer in $\mathfrak l$ is exactly $\mathfrak b$. Add any nonzero vector of $U_8$.  Invertibility of every nonzero $\mathfrak b$ action on $U_8$ makes the resulting vector in $F$ free for $\mathfrak l$.

For any $H\subset L$, these same vectors are free for $\mathfrak h$, and
$\mathfrak m=E\oplus \bigl(F\oplus(\mathfrak l\ominus\mathfrak h)\bigr)$
is an orthogonal $H$-invariant splitting in which both summands contain a vector with zero infinitesimal stabilizer.  The two-summand case of the deleting-any-summand lemma implies normality.  That lemma does not require $C$ to preserve the displayed summands, so possible isotypic mixing is fully included.
\end{proof}

\begin{remark}
These propositions close the specified nonregular families in $F_4$. They are used in the exhaustive argument of Theorem~\ref{thm91}.
\end{remark}

\paragraph{\texorpdfstring{All isotropy subgroups contained in the nonregular $G_2\SO(3)$ subgroup}{All isotropy subgroups contained in the nonregular G2\SO(3) subgroup}}
\label{sec:f4-g2-so3-all}

\begin{theorem}\label{thm86}
Let $K\subset F_4$ be the subgroup with Lie algebra $\mathfrak g_2\oplus\mathfrak{so}(3)$ and isotropy module $\mathfrak p=V\otimes W$ as in \eqref{formula33}.  For every closed subgroup $H\subset K$, all $F_4$-g.o. metrics on $F_4/H$ are normal.
\end{theorem}

\begin{proof}
By Corollary~\ref{cor33}, it suffices to treat connected $H$. Let $\pi:\mathfrak h\to\mathfrak g_2$ be the first projection. We use the following exhaustive subgroup fact. Every proper connected subgroup of $G_2$, except $\SU(3)$ and the principal $\SO(3)$, is contained in a conjugate of $\SO(4)$.  See \cite[Proposition~15, Table~5 and the following paragraph]{Ko2009}. Abelian subgroups are also covered, since a maximal torus of $G_2$ lies in $\SO(4)$.

If $\pi(\mathfrak h)$ is contained in a conjugate of $\mathfrak{so}(4)$, Proposition~\ref{prop84} applies.  This includes arbitrary couplings with the external $\mathfrak{so}(3)$.  It remains to treat the three projections $\mathfrak g_2$, $\mathfrak{su}(3)$ and a principal $\mathfrak a\cong\mathfrak{so}(3)$ in $\mathfrak g_2$.

If the projection is $\mathfrak g_2$ or $\mathfrak{su}(3)$, the elementary Lie-algebra form of Goursat's lemma gives
$\mathfrak h=\pi(\mathfrak h)\oplus\mathfrak l, \qquad \mathfrak l\subset\mathfrak{so}(3).$
Indeed, neither of these simple algebras has a nonzero quotient isomorphic to a quotient of a subalgebra of $\mathfrak{so}(3)$.  Up to conjugation, $\mathfrak l$ is $0$, $\mathfrak{so}(2)$ or $\mathfrak{so}(3)$. If the projection is $\mathfrak a$, the same three product possibilities occur, with the additional possibility of a graph of an isomorphism $\mathfrak a\to\mathfrak{so}(3)$. For the graph case, the tensor module $V\otimes W$ restricts to the diagonal $\mathfrak{so}(3)$ as the tensor product of spins $3$ and $2$, hence as the sum of spins $1,2,3,4,5$. Its orthogonal spin-two and spin-three submodules each contain a vector with zero infinitesimal stabilizer. A nonzero rotation has kernel dimension one, and the union of its kernels has dimension at most three, less than either $5$ or $7$.  The deletion lemma therefore gives normality in the graph case, even with arbitrary mixing of equivalent modules elsewhere in $\mathfrak m$. The following product cases exhaust the remaining work.

\medskip \noindent\emph{Projection $\mathfrak g_2$.} For $\mathfrak h=\mathfrak g_2\oplus\mathfrak l$, put $\mathfrak f=\mathfrak{so}(3)\ominus\mathfrak l$. The $\mathfrak g_2$ types distinguish $\mathfrak f$ from $\mathfrak p$, so $C$ preserves both.  Descent to $F_4/K$ gives $C|_{\mathfrak p}=cI$.  For $\mathfrak l=0$, normalizer invariance under the external $\SO(3)$ gives $C|_{\mathfrak f}=aI$.  For $\mathfrak l=\mathfrak{so}(2)$ the same statement follows from the real irreducibility of its two-dimensional complement. For $\mathfrak l=\mathfrak{so}(3)$ there is no vertical space.

In either nonzero vertical case, choose $0\ne Y\in\mathfrak f$. For $X=v\otimes w$, write the compensator as $A+B$, with $A\in\mathfrak g_2$, $B\in\mathfrak l$. The vertical equation gives $[B,Y]=0$, which forces $B=0$ when $\mathfrak l=\mathfrak{so}(2)$. It is already zero when $\mathfrak l=0$.  The horizontal equation is
$(a-c)v\otimes(Y\cdot w)=(A\cdot v)\otimes w.$
Take $v\ne0$ and $Y\cdot w\ne0$, and project onto $V\otimes w^\perp$.  This gives $a=c$.

\medskip \noindent\emph{Projection $\mathfrak{su}(3)$.} Write
$V=\mathbb Re_0\oplus V_6, \qquad \mathfrak g_2=\mathfrak{su}(3)\oplus U_6,$
where $V_6$ and $U_6$ are the underlying real modules of the standard complex $\SU(3)$-module, and $e_0$ is a fixed unit vector. For $\mathfrak l=\mathfrak{so}(3)$, the external factor separates $U_6$ from $\mathfrak p$, so descent gives $C|_{\mathfrak p}=cI$ and $C|_{U_6}=aI$. For $\mathfrak l=0$, normalizer invariance under the full external $\SO(3)$ gives the same conclusion and $C|_{\mathfrak{so}(3)}=bI$.  In this case the external modules $\mathbb R$ on $U_6$, the adjoint module on $\mathfrak{so}(3)$ and the spin-two module on $\mathfrak p=(\mathbb R\oplus V_6)\otimes W$ have distinct types, so this argument also excludes all mixing.

For either of these two choices of $\mathfrak l$, take $0\ne Y\in U_6$ and $X=e_0\otimes w$ with $w\ne0$. Since the stabilizer algebra of $e_0$ in $\mathfrak g_2$ is $\mathfrak{su}(3)$, one has $Y\cdot e_0\ne0$. For a compensator $A+B\in\mathfrak{su}(3)\oplus\mathfrak l$, the horizontal g.o. equation reads
$(a-c)(Y\cdot e_0)\otimes w=e_0\otimes(B\cdot w).$
Its two sides lie in orthogonal first-factor spaces, so $a=c$. When $\mathfrak l=0$, testing instead a nonzero $Y\in\mathfrak{so}(3)$ and $X=v\otimes w$ gives
$(b-c)v\otimes(Y\cdot w)=(A\cdot v)\otimes w,$
which implies $b=c$ by projecting onto $V\otimes w^\perp$.

The remaining case $\mathfrak l=\mathfrak{so}(2)$ must allow mixing, and we handle it by the deletion lemma. Decompose the spin-two module under this circle as
$W=W_0\oplus W_1\oplus W_2, \qquad \dim W_0=1,\quad\dim W_1=\dim W_2=2,$
where the two planes have nonzero weights $1$ and $2$. Put $\mathfrak q=\mathfrak{so}(3)\ominus\mathfrak{so}(2)$, a nonzero-weight two-plane.  The orthogonal invariant spaces
$E=U_6\oplus(V_6\otimes W_0)\oplus\mathfrak q, \qquad F=(V_6\otimes W_1)\oplus(\mathbb Re_0\otimes W_2)$
are disjoint and both contain an infinitesimally free vector. For $E$, a nonzero $\mathfrak q$ component kills the circle stabilizer, and two complex-linearly independent vectors in the two standard $\mathbb C^3$ modules then kill $\mathfrak{su}(3)$.  For $F$, a nonzero component in $\mathbb Re_0\otimes W_2$ first kills the circle. Choose the two $V_6$ coefficients of a vector in $V_6\otimes W_1$ to be complex-linearly independent to kill $\mathfrak{su}(3)$.  A traceless skew-Hermitian $3\times3$ matrix killing two complex-linearly independent vectors is zero, which justifies both assertions. Thus $E$ and its orthogonal complement in $\mathfrak m$ are both infinitesimally free, the latter because it contains $F$.  The deletion lemma gives normality without any assumption that $C$ preserves these spaces.

\medskip \noindent\emph{Projection a principal $\mathfrak a$.} Under the principal algebra $\mathfrak a\cong\mathfrak{so}(3)$,
$V=V_7,\qquad \mathfrak g_2=\mathfrak a\oplus U_{11},$
where $V_7$ and $U_{11}$ are the irreducible real modules of spins $3$ and $5$, respectively.  For $\mathfrak h=\mathfrak a\oplus\mathfrak l$, one has
$\mathfrak m=U_{11}\oplus (\mathfrak{so}(3)\ominus\mathfrak l) \oplus(V_7\otimes W).$
The principal-factor types have spins $5$, $0$ and $3$. Thus the three displayed spaces cannot mix. Descent gives $C|_{V_7\otimes W}=cI$, and irreducibility gives $C|_{U_{11}}=aI$. On the nonzero external complement $\mathfrak q$, write $C|_{\mathfrak q}=bI$. For $\mathfrak l=0$ this follows from normalizer invariance under the external $\SO(3)$, and for a circle from irreducibility.

There is a vector $Y\in U_{11}$ with zero infinitesimal $\mathfrak a$-stabilizer.  Indeed, every nonzero infinitesimal rotation on this spin-five module has a one-dimensional kernel.  The family of its kernels is parametrized by rotation axes, a two-dimensional space. Their union has dimension at most three and cannot fill the eleven-dimensional module. Fix such a $Y$.  For $X=v\otimes w$, the vertical g.o. equation forces the $\mathfrak a$ component of the compensator to vanish.  Hence its horizontal equation becomes
$(a-c)(Y\cdot v)\otimes w=v\otimes(B\cdot w), \qquad B\in\mathfrak l.$
The standard seven-dimensional representation of $\mathfrak g_2$ is faithful, so choose $v$ with $Y\cdot v\ne0$.  Since $Y\cdot v\perp v$, projection onto $v^\perp\otimes W$ implies $a=c$.

If $\mathfrak q\ne0$, take $0\ne Y\in\mathfrak q$. The vertical g.o. equation for $Y+v\otimes w$ forces the $\mathfrak l$ component $B$ to vanish, as in the first case. The horizontal equation is then
$(b-c)v\otimes(Y\cdot w)=(A\cdot v)\otimes w, \qquad A\in\mathfrak a.$
Choosing $Y\cdot w\ne0$ and projecting onto $V_7\otimes w^\perp$ gives $b=c$. This completes all product cases and hence the proof.
\end{proof}

\begin{remark}
The theorem exhausts subgroups of this particular nonregular maximal subgroup of $F_4$, including non-product subgroups. The other maximal-subgroup branches are incorporated in Theorem~\ref{thm91}, which completes the fixed-$F_4$ classification.
\end{remark}

\paragraph{Further subgroups of $F_4$.}\label{sec:f4-exhaustion}
We use the convention $V_d=\operatorname{Sym}^{d-1}\mathbb C^2$. When $d$ is odd, the same symbol may denote its invariant real form. A tensor of an even number of quaternionic factors denotes its real form. The notation $(\cdot)_{\mathbb R}$ denotes the underlying real representation when no real form is specified. In particular $(V_2\otimes V_2\otimes V_2)_{\mathbb R}$ has real dimension $16$.

\begin{lemma}
\label{lem88}
The following representations are locally free at some vector.
\begin{enumerate}
\item $V_2\otimes V_{2n}$ of $\mathfrak{su}(2)^2$ for every $n\geq2$.
\item $V_3\otimes V_5$ of $\mathfrak{so}(3)^2$.
\item $(V_2\otimes V_2\otimes V_2)_{\mathbb R}$ of
$\mathfrak{su}(2)^3$.
\end{enumerate}
\end{lemma}

\begin{proof}
We use the standard weights and conjugate an element of the compact Lie algebra into a fixed Cartan subalgebra.  It suffices to show that the union of the kernels of all nonzero infinitesimal operators has dimension smaller than the representation. For (1), the weights are $\pm a+(2j-1)b$. A nonzero operator with a kernel has one of finitely many projective ratios $a:b$, both entries nonzero. The real kernel has dimension $2$, and its conjugates form a family of dimension at most $4$. The union therefore has dimension at most $6<4n$. For (2), when both Cartan entries are nonzero a zero-weight space has dimension at most $3$, with conjugacy parameters of dimension at most $4$.  Generic projective ratios have only the constant zero weight, and add at most one parameter.  If one entry vanishes, the kernel has dimension at most $5$, with at most two conjugacy parameters. All resulting bounds are smaller than $15$.

For (3), the complex weights are $\pm a\pm b\pm c$. With all entries nonzero, a kernel has complex dimension $2$ and the projective parameters satisfying a zero-weight equation have dimension at most one.  Including the six conjugacy parameters, the corresponding union in the underlying real module has dimension at most $11$. With precisely two nonzero entries, a kernel occurs only when their absolute values agree. It has real dimension $8$ and at most four conjugacy parameters.  With one nonzero entry there is no kernel. Thus the union has dimension at most $12<16$.
\end{proof}

\begin{lemma}
\label{lem89}
All subgroups of each of the following connected subgroups of $F_4$ have only normal $F_4$-g.o. metrics.
\begin{enumerate}
\item $K=\Spin(3)^3\subset\Spin(9)$,
where $\mathbb R^9=\mathbb R^3\oplus\mathbb R^3\oplus\mathbb R^3$.
\item $K=A_1^{B_2}\cdot A_1\cdot A_1$, where $A_1^{B_2}$ is
principal in $\Spin(5)$ and $K\subset\Spin(5)\cdot\Spin(4) \subset\Spin(9)$.
\item the three-factor envelope of the $A_1$ with label
$A_1+\widetilde A_1$ in the nilpotent-orbit notation.
\end{enumerate}
\end{lemma}

\begin{proof}
For (1), the quotient representation splits as
$W_1=(V_3\otimes V_3\otimes1)\oplus(V_3\otimes1\otimes V_3) \oplus(1\otimes V_3\otimes V_3),\qquad W_2=(V_2\otimes V_2\otimes V_2)_{\mathbb R}.$
The tensor $\mathbb R^3\otimes\mathbb R^3$ is locally free for $\SO(3)^2$. A matrix with distinct nonzero singular values has zero infinitesimal stabilizer.  Thus two of the pair modules already make $W_1$ locally free for all three factors.  The preceding lemma treats $W_2$.

For (2), use
$W_1=(V_4\otimes V_2\otimes1)\oplus(V_4\otimes1\otimes V_2), \qquad W_2=(V_5\otimes V_2\otimes V_2)\oplus(V_7\otimes1\otimes1).$
The first summand of $W_1$ kills the first two infinitesimal factors. A nonzero vector in the second then kills the third. A generic vector in $V_7$ has zero $\mathfrak{su}(2)$ stabilizer. After the first factor is killed, $V_5\otimes V_2\otimes V_2$ is five copies of the real standard $\SO(4)$ module.  Three independent quaternions, for example $1,i,j$, kill its remaining two factors.

For (3), the quotient representation is
$W_1=V_5\otimes V_2\otimes V_2,\qquad W_2=(V_5\otimes V_3\otimes1)\oplus(1\otimes V_4\otimes V_2).$
The first module is the restriction of the standard $\SO(5)\times \SO(4)$ representation on real $5\times4$ matrices. A full-rank matrix with distinct nonzero singular values has zero Lie stabilizer even for this larger group. In $W_2$, the first summand kills the first two factors by Lemma~\ref{lem88}, and a nonzero vector in the second kills the third.  The two-free-module criterion and its hereditary form now apply in all three cases.
\end{proof}

\begin{lemma}
\label{lem810}
Let $S$ be principal $\SO(3)$ in the long-root $\SU(3)$ factor of the regular subgroup $\SU(3)\cdot \SU(3)\subset F_4$. Let $R$ be any connected subgroup of the other $\SU(3)$ factor. Then every $F_4$-g.o. metric with isotropy algebra $\mathfrak s\oplus\mathfrak r$ is normal.
\end{lemma}

\begin{proof}
For $K=\SU(3)\cdot S$, the quotient representation is
$\mathfrak f_4/\mathfrak k=V_5\oplus (V_3\otimes\operatorname{Sym}^2\mathbb C^3)_{\mathbb R},$
where $S$ acts on $V_3,V_5$ and the short-root $\SU(3)$ acts on $\mathbb C^3$.  This orientation follows from formula (4.6) in Shan~\cite{Sh2024}. The cross module for the short- and long-root factors is $\operatorname{Sym}^2\mathbb C^3_{\rm short}\otimes \mathbb C^3_{\rm long}$ and its dual. For $R=\SU(3)$ the quotient has precisely two real irreducible summands, of dimensions $5$ and $36$.  the two-isotropy-summand classification of Chen--Nikonorov~\cite{ChNi2019} therefore gives normality.

Every proper connected subgroup of $\SU(3)$ is contained in a conjugate of its standard $\U(2)$ or is principal $\SO(3)$. This follows as well by considering whether its defining three-dimensional complex representation is reducible.

For $R=\U(2)$, decompose $\mathbb C^3=\mathbb C^2\oplus\det^{-1}$. The two modules
$(V_3\otimes\operatorname{Sym}^2\mathbb C^2)_{\mathbb R},\qquad (V_3\otimes\mathbb C^2\otimes\det^{-1})_{\mathbb R}$
are both locally free for $S\times \U(2)$. For the first, the semisimple factors act as $\SO(3)\times \SO(3)$ on complex $3\times3$ matrices and the center acts by a nonzero scalar phase.  A real diagonal matrix with distinct positive entries has zero Lie stabilizer. Its imaginary part kills the center, and its singular values kill both rotation factors. For the second, use a real $3\times2$ matrix with distinct positive singular values.  The imaginary part kills the imaginary component of the $\mathfrak u(2)$ operator. The remaining $\mathfrak{so}(3)\oplus\mathfrak{so}(2)$ stabilizer is zero. The determinant twist is an automorphism of $\mathfrak u(2)$, so it does not affect this conclusion. These two free modules restrict to every subgroup of $\U(2)$. For principal $R=\SO(3)$, the tensor module contains two real copies of $V_3\otimes V_5$, both locally free.  This finishes all cases.
\end{proof}

\begin{lemma}
\label{lem811}
Let $\mathfrak h=\mathfrak{sp}(2)\oplus \operatorname{diag}\mathfrak{sp}(1)$ in $\mathfrak{sp}(2)\oplus\mathfrak{sp}(1)\oplus\mathfrak{sp}(1) \subset\mathfrak{so}(9)\subset\mathfrak f_4$. Then every $F_4$-g.o. metric on $F_4/H$ is normal.
\end{lemma}

\begin{proof}
There are two invariant free modules
$W_1=\mathbb R^5\oplus(\mathbb R^5\otimes\mathbb R^3),\qquad W_2=\mathbb R^3\oplus\mathbb H^2\oplus\mathbb H^2.$
For $W_1$, take a matrix with columns $e_1,2e_2,3e_3$ and the additional vector $e_4$. Its stabilizer in $\SO(5)\times \SO(3)$ has zero Lie algebra.  For $W_2$, take quaternionic columns $q_1=e_1$ and $q_2=e_1i+e_2$. If $(A,b)\in\mathfrak{sp}(2)\oplus\mathfrak{sp}(1)$ fixes both, then $Aq_j=q_jb$ and skew-adjointness forces $[b,i]=0$. An additional adjoint vector $j\in\mathfrak{sp}(1)$ forces $[b,j]=0$, hence $b=0$.  Since the two columns form a quaternionic basis, $A=0$. Apply the two-free-module criterion.
\end{proof}

\paragraph{Hereditary normality in the exceptional tensor chains}
\label{sec:exceptional-tensor-chains}

We work with the compact forms of the following Tits decompositions.
\begin{equation}\label{formula34}
\begin{array}{c|c|c|c}
\mathfrak g&\mathfrak k=\mathfrak g_2\oplus\mathfrak l
 &W&\mathfrak p=\mathfrak k^{\perp_Q}\\\hline
\mathfrak e_6&\mathfrak g_2\oplus\mathfrak{su}(3)
 &\operatorname{Herm}_0(3,\mathbb C)&\mathbb R^7\otimes W\\
\mathfrak e_7&\mathfrak g_2\oplus\mathfrak{sp}(3)
 &\operatorname{Herm}_0(3,\mathbb H)&\mathbb R^7\otimes W\\
\mathfrak e_8&\mathfrak g_2\oplus\mathfrak f_4
 &\operatorname{Herm}_0(3,\mathbb O)&\mathbb R^7\otimes W.
\end{array}
\end{equation}
Here $\dim W=8,14,26$, respectively, and all the displayed tensor modules are real irreducible.  The $14$-dimensional module of $\mathfrak{sp}(3)$ has highest weight $\omega_2$, not $\omega_3$. These decompositions and their realizations by trace-free Jordan algebras are given in \cite[pp.~21--22]{BeElMa2009}. We use the reductive decompositions and the indicated module structures. Finite central quotients of the product notation do not affect any of the arguments below.

\paragraph{Three vectors and a descent lemma}

\begin{lemma}\label{lem812}
In each of the three representations $(\mathfrak l,W)$ in \eqref{formula34}, there are three linearly independent vectors whose joint infinitesimal stabilizer is zero.  They may be chosen orthonormal.
\end{lemma}

\begin{proof}
For $\mathfrak{su}(3)$ take a real diagonal Hermitian matrix $D$ with distinct diagonal entries and trace zero, and a Hermitian matrix $S$ whose three off-diagonal entries are $1$. An element of $\mathfrak{su}(3)$ commuting with $D$ is diagonal. Commuting also with $S$ makes all its diagonal entries equal, so it is zero.  Add a third independent vector.

For $\mathfrak{sp}(3)$ use the same $D,S$ in the quaternionic Hermitian model.  An infinitesimal stabilizer of $D$ is $\operatorname{diag}(q_1,q_2,q_3)$ with $q_i\in\operatorname{Im}\mathbb H$. Fixing $S$ gives $q_1=q_2=q_3=q$.  A third Hermitian matrix with off-diagonal entries $i,j,0$ forces $[q,i]=[q,j]=0$, hence $q=0$.

For $\mathfrak f_4$ use the Albert algebra and write its three off-diagonal Peirce coordinates as an ordered triple in $\mathbb O^3$. Take
$D=\operatorname{diag}(1,0,-1),\qquad S=(1,1,1),\qquad T=(i,j,\ell),$
where $i,j$ generate a quaternion subalgebra and $\ell$ is a unit orthogonal to that subalgebra. Fixing $D$ fixes its three spectral idempotents, so its stabilizer is the frame stabilizer $\Spin(8)$, acting on the three Peirce spaces by triality.  Fixing the three entries of $S$ makes the related triple a common octonion automorphism.  At Lie algebra level this follows directly from the triality identity
$d_1(xy)=d_2(x)y+x\,d_3(y): \quad d_1(1)=d_2(1)=d_3(1)=0 \ \Longrightarrow\ d_1=d_2=d_3\in\operatorname{Der}(\mathbb O).$
Finally, a derivation fixing $i,j,\ell$ fixes the algebra they generate, namely $\mathbb O$, and therefore vanishes. The standard frame and triality description can also be read directly from the Peirce products in the Albert algebra. Orthogonalization changes neither the span nor its pointwise stabilizer, proving the last assertion.
\end{proof}

We record a useful elementary matrix observation.  Let $A\in\mathfrak{so}(U)$, $B\in\mathfrak{so}(W)$, and let $T:U\longrightarrow W$ be injective with distinct positive singular values.  If $BT=TA$, then $A$ commutes with $T^*T$. Hence $A=0$ and $B$ vanishes on $\operatorname{im}T$.  Consequently, if $\operatorname{im}T$ contains a pointwise-free frame for the acting algebra on $W$, the tensor represented by $T$ has zero stabilizer for the product algebra.

Lemma~\ref{lem812} and the matrix observation show that $\mathfrak k$ has a free vector in each tensor module $\mathfrak p=\mathbb R^7\otimes W$. Take an injective $T:\mathbb R^7\to W$ with distinct singular values whose image contains the three-vector base.

\begin{proposition}\label{prop813}
In each row of \eqref{formula34}, all $G$-g.o. metrics are normal for isotropy algebras of any of the following forms.
\begin{equation*}
 \mathfrak g_2\oplus\mathfrak r\quad(\mathfrak r\subset\mathfrak l),
 \qquad
 \mathfrak r\oplus\mathfrak l\quad(\mathfrak r\subset\mathfrak g_2),
 \qquad
 \mathfrak a_{1,\mathrm{pr}}\oplus\mathfrak r
       \quad(\mathfrak r\subset\mathfrak l),
\end{equation*}
where $\mathfrak a_{1,\mathrm{pr}}$ is principal in $\mathfrak g_2$.
\end{proposition}

\begin{proof}
In the first case $\mathfrak k/\mathfrak h$ is trivial under the full $\mathfrak g_2$ factor, whereas $\mathfrak p$ is a sum of standard seven-dimensional modules.  Thus equivariance prevents mixing between these two spaces.  In the second case use the full $\mathfrak l$ factor instead.  In the third case the principal $\mathfrak{su}(2)$ acts on $\mathfrak k/\mathfrak h$ only through the $11$-dimensional irreducible module and trivial modules, whereas its action on $\mathfrak p$ is a sum of $7$-dimensional irreducible modules.  These assertions follow from $\mathfrak g_2|_{\mathfrak a_{1,\mathrm{pr}}}=V_3\oplus V_{11}$ and $\mathbb R^7|_{\mathfrak a_{1,\mathrm{pr}}}=V_7$. In every case $C$ preserves $\mathfrak p$, and Lemma~\ref{lem39} applies.
\end{proof}

\paragraph{\texorpdfstring{A hereditary $\SO(4)$ branch in all three groups}{A hereditary \SO(4) branch in all three groups}}

\begin{theorem}\label{thm814}
In every row of \eqref{formula34}, every closed subgroup $H\subset \SO(4)\cdot L$ has only normal $G$-g.o. metrics. The embedding $\SO(4)\subset G_2$ is the standard maximal subgroup.
\end{theorem}

\begin{proof}
Write $\mathfrak{so}(4)=\mathfrak a\oplus\mathfrak b$, where $\mathfrak a$ is the short-root and $\mathfrak b$ the long-root $\mathfrak{su}(2)$ factor.  The standard restrictions are
$\mathbb R^7=V_3\oplus V_4,\qquad \mathfrak g_2=\mathfrak a\oplus\mathfrak b\oplus U_8.$
Here $V_3$ is the adjoint of $\mathfrak a$, trivial on $\mathfrak b$, and $V_4$ is the standard real $\SO(4)$ module. Moreover, $(U_8)_{\mathbb C}=\operatorname{Sym}^3\mathbb C^2_{\mathfrak a} \otimes\mathbb C^2_{\mathfrak b}$. Every nonzero element of $\mathfrak b$ therefore acts invertibly on $U_8$.

In the quotient by $\SO(4)\cdot L$ take the orthogonal invariant modules
$E=V_4\otimes W,\qquad F=(V_3\otimes W)\oplus U_8.$
For $E$, use an injective map $\mathbb R^4\to W$ with distinct singular values and image containing the three-vector base. It kills both the $\mathfrak{so}(4)$ and $\mathfrak l$ components of its stabilizer. For the tensor summand of $F$, the analogous map $\mathbb R^3\to W$ kills $\mathfrak a$ and $\mathfrak l$. Add any nonzero vector of $U_8$ to kill $\mathfrak b$. Thus both $E$ and $F$ are free for the whole envelope. They remain free for every subalgebra $\mathfrak h$ of that envelope and lie in $\mathfrak h^\perp$. The two-free-module criterion proves the claim.
\end{proof}

\paragraph{Two further hereditary envelopes}

\begin{lemma}\label{lem815}
On the standard real $G_2$ module $\mathbb R^7$, every nonzero infinitesimal operator has kernel dimension at most $3$. The three octonions $i,j,\ell$ used above have joint $\mathfrak g_2$-stabilizer zero. On $\mathbb R^3\otimes\mathbb R^3$, every nonzero infinitesimal operator of $\mathfrak{so}(3)^2$ has kernel dimension at most $3$.
\end{lemma}

\begin{proof}
The weights of the seven-dimensional module are $0,\pm\varepsilon_1,\pm\varepsilon_2,\pm\varepsilon_3$ with $\varepsilon_1+\varepsilon_2+\varepsilon_3=0$. On a nonzero Cartan element at most one of the three $\varepsilon_i$ vanishes.  The assertion follows by compact conjugacy.  The three-vector assertion follows because $i,j,\ell$ generate $\mathbb O$.  For the last assertion the nine weights are the pairwise sums of $0,\pm a$ and $0,\pm b$. For $(a,b)\ne(0,0)$ at most three are zero.
\end{proof}

\begin{theorem}\label{thm816}
The following two envelopes have hereditary normality.
\begin{enumerate}
\item $G_2\cdot \SO(3)\cdot \SO(3)\subset E_7$, with the last
two factors inside $\Sp(3)$ via $\mathbb C^6=V_3\otimes V_2$.
\item $G_2\cdot G_2\cdot \SO(3)\subset E_8$, with the last
two factors the standard $G_2\times \SO(3)\subset F_4$.
\end{enumerate}
That is, every closed subgroup of either envelope has only normal $G$-g.o. metrics for the displayed exceptional acting group.
\end{theorem}

\begin{proof}
For (1), restriction in the defining symplectic representation gives
$\Lambda^2_0(V_3\otimes V_2)=(V_3\otimes V_3)\oplus(V_5\otimes1), \qquad \mathfrak{sp}(3)/(\mathfrak{so}(3)\oplus\mathfrak{so}(3)) =V_5\otimes V_3.$
Consequently the quotient by the three-factor envelope is
$E=(V_7\otimes V_3\otimes V_3),\qquad F=(V_7\otimes V_5\otimes1)\oplus(1\otimes V_5\otimes V_3).$
For $E$, view a tensor as an injective map $\mathbb R^7\to\mathbb R^3\otimes\mathbb R^3$ with distinct singular values.  The matrix observation kills the $G_2$ action. The remaining $\mathfrak{so}(3)^2$ operator kills a seven-dimensional space and is zero by Lemma~\ref{lem815}. In the first summand of $F$, use an injective $\mathbb R^5\to\mathbb R^7$ with distinct singular values. It first kills the middle $\mathfrak{so}(3)$ and then kills $\mathfrak g_2$, by the kernel bound.  In the second summand choose two linearly independent coefficients in the last $\mathbb R^3$ to kill its remaining $\mathfrak{so}(3)$.

For (2), the Jordan-algebra restriction is
$\mathbb R^{26}|_{G_2\times \SO(3)} =(V_7\otimes V_3)\oplus(1\otimes V_5), \qquad \mathfrak f_4/(\mathfrak g_2\oplus\mathfrak{so}(3))=V_7\otimes V_5.$
Thus take
$E=V_7\otimes V_7\otimes V_3,\qquad F=(V_7\otimes1\otimes V_5)\oplus(1\otimes V_7\otimes V_5).$
The representation $V_7\otimes V_3$ has a free vector for $\mathfrak g_2\oplus\mathfrak{so}(3)$. As a map $\mathbb R^3\to\operatorname{Im}\mathbb O$, take the three columns $i,2j,3\ell$.  The distinct singular values kill $\mathfrak{so}(3)$, and the three-vector base kills $\mathfrak g_2$. Choose now an injective map $\mathbb R^7\to V_7\otimes V_3$ with distinct singular values and with this free vector in its image.  The matrix observation proves freeness of $E$. Each tensor in $F$ admits a free vector for its two nontrivial factors, by the $\mathbb R^5\to\mathbb R^7$ argument above. Taking such vectors in both summands makes $F$ free for all three factors.  The hereditary two-free-module criterion finishes both claims.
\end{proof}

\begin{proposition}
\label{prop817}
For every row of \eqref{formula34}, let $\mathfrak h\subset\mathfrak g_2\oplus\mathfrak l$ project surjectively onto the first factor.  Then all $G$-g.o. metrics are normal.
\end{proposition}

\begin{proof}
The ideal $\mathfrak h\cap\mathfrak g_2$ is either the whole $\mathfrak g_2$ or zero.  The first case is a product covered by Proposition~\ref{prop813}. In the second case projection onto $\mathfrak l$ identifies $\mathfrak h$ with a compact reductive subalgebra $\mathfrak r\subset\mathfrak l$ having a quotient $\mathfrak g_2$.  Thus $\mathfrak r$ contains a simple ideal isomorphic to $\mathfrak g_2$. This is impossible for $\mathfrak l=\mathfrak{su}(3)$ by dimension, and for $\mathfrak l=\mathfrak{sp}(3)$ because the defining complex six-dimensional representation cannot contain a nontrivial $\mathfrak g_2$ module. For $\mathfrak l=\mathfrak f_4$, the $G_2$ embedding is unique up to conjugacy and has centralizer $\mathfrak{so}(3)$.
~\cite[Sections~4.5.1.7 and~4.5.3]{Sh2024} Consequently $\mathfrak h$ lies in the second envelope of Theorem~\ref{thm816}, including any diagonal identification of its two $\mathfrak g_2$ factors.
\end{proof}

\paragraph{\texorpdfstring{Completion of the $E_6$ tensor subgroup}{Completion of the E6 tensor subgroup}}

For this subsection write $A=\SU(3)\subset G_2$ for the stabilizer of a unit imaginary octonion, and $B=\SU(3)$ for the other factor of the $E_6$ tensor subgroup.  Then
$\mathbb R^7|_A=\mathbb R\oplus U, \qquad \mathfrak g_2/\mathfrak a=U, \qquad U=(\mathbb C^3)_{\mathbb R}.$

\begin{lemma}\label{lem818}
Every subgroup of $A\cdot \U(2)\subset A\cdot B\subset E_6$ has only normal $E_6$-g.o. metrics.
\end{lemma}

\begin{proof}
For the standard $\U(2)\subset B$ write
$\mathfrak b=\mathfrak{su}(2)\oplus\mathfrak z\oplus P, \qquad \dim\mathfrak z=1,\quad \dim_{\mathbb R}P=4.$
The action on $P$ is the standard $\U(2)$ action with a determinant twist. At Lie algebra level it is equivalent to the standard action composed with an automorphism of $\mathfrak u(2)$. Inside the quotient by $A\cdot \U(2)$ choose
$E=U\oplus(U\otimes\mathfrak z)\oplus P \oplus(\mathbb R\otimes\mathfrak{su}(2)),\qquad F=U\otimes P.$
The first $U$ and $P$ in $E$ are the vertical complements in $\mathfrak g_2$ and $\mathfrak b$. The other terms come from $(\mathbb R\oplus U)\otimes\mathfrak b$. These are disjoint orthogonal invariant modules.

Choose two complex linearly independent vectors in the two copies of $U$ in $E$. They kill $\mathfrak a$. A nonzero vector of $P$ has a one-dimensional stabilizer in $\mathfrak u(2)$ whose projection onto $\mathfrak{su}(2)$ is nonzero.  Choose an adjoint vector not on that stabilizer's axis. It kills the remaining one-dimensional algebra.  Thus $E$ is free. For $F$, use an injective map $P\simeq\mathbb R^4\to U\simeq\mathbb R^6$ with distinct singular values. It kills $\mathfrak u(2)$, and the remaining $\mathfrak{su}(3)$ fixes its four-dimensional real image pointwise. That image contains two complex linearly independent vectors, whose stabilizer in $\mathfrak{su}(3)$ is zero. So $F$ is free as well.
\end{proof}

\begin{lemma}\label{lem819}
Every subgroup of $A\cdot \SO(3)\subset A\cdot B\subset E_6$, where $\SO(3)$ is principal in $B$, has only normal $E_6$-g.o. metrics.
\end{lemma}

\begin{proof}
The adjoint of $B$ restricts as $V_3\oplus V_5$. The two orthogonal modules $U\otimes V_3$ and $U\otimes V_5$ are each free for $\mathfrak{su}(3)\oplus\mathfrak{so}(3)$. Represent a tensor as an injective map $\mathbb R^d\to\mathbb C^3$ with distinct singular values, for $d=3$ or $5$, and choose the real image to contain two complex linearly independent vectors. The matrix observation kills $\mathfrak{so}(3)$, after which pointwise fixation of those two vectors kills $\mathfrak{su}(3)$.
\end{proof}

\begin{lemma}
\label{lem820}
The following isotropy algebras also have only normal $E_6$-g.o. metrics.
\begin{enumerate}
\item a graph $\mathfrak{su}(3)\subset\mathfrak a\oplus\mathfrak b$.
\item a graph between the principal $\mathfrak{su}(2)$ in
$\mathfrak g_2$ and a standard $\mathfrak{su}(2)\subset\mathfrak b$, with or without its commuting central circle.
\item a graph between the principal $\mathfrak{su}(2)$ in
$\mathfrak g_2$ and the principal $\mathfrak{so}(3)\subset\mathfrak b$.
\end{enumerate}
\end{lemma}

\begin{proof}
For (1), the complex tensor rule
$\mathbb C^3\otimes\mathfrak{sl}_3(\mathbb C) =\mathbb C^3\oplus(\operatorname{Sym}^2\mathbb C^3)^* \oplus V_{(2,1)}$
gives a real $12$-dimensional module $E=((\operatorname{Sym}^2\mathbb C^3)^*)_{\mathbb R}$ in $U\otimes\mathfrak b$. The complex symmetric matrix with distinct positive diagonal entries has zero $\mathfrak{su}(3)$ stabilizer under $S\mapsto USU^t$. Writing $A=X+iY$, its stabilizer equation $AS+SA^t=0$ makes both $X$ and $Y$ zero. Thus $E$ is free.  A second free module is the sum of the vertical $U=\mathfrak g_2/\mathfrak a$ and the copy of $(\mathbb C^3)_{\mathbb R}$ in the displayed tensor decomposition. Two complex independent vectors give zero stabilizer there. The argument is unchanged for an outer-automorphism graph, because the adjoint representation twisted by an automorphism is equivalent to itself.

For (2), under the diagonal $\mathfrak{su}(2)$ the tensor $V_7\otimes\mathbb C^2$ splits as $V_8\oplus V_6$. Here $V_d=\operatorname{Sym}^{d-1}\mathbb C^2$. The corresponding two underlying-real modules in $\mathbb R^7\otimes P$ both carry a nonzero scalar character of the central circle, if that circle is present. Each is infinitesimally free for $\mathfrak{su}(2)\oplus\mathbb R$. After conjugating an element into a Cartan, its complex eigenvalues are $i(am+t)$, where the $d$ weights $m$ are distinct. For $a\ne0$ its kernel has real dimension at most $2$. There are only finitely many projective ratios $a:t$ giving a kernel, and the conjugacy axes contribute at most two dimensions. The union of all nonzero-operator kernels therefore has dimension at most $4<2d$.  For $a=0$, $t\ne0$, there is no kernel. Without the circle the same conclusion follows by restriction. The two-free-module criterion applies to $V_8$ and $V_6$.

For (3), the orthogonal modules $\mathfrak g_2/\mathfrak a_{1,\mathrm{pr}}=V_{11}$ and $\mathfrak b/\mathfrak{so}(3)=V_5$ are already free at generic vectors for the diagonal $\mathfrak{su}(2)$. Indeed, for any odd-dimensional $V_d$, $d\ge5$, the kernel of a nonzero Cartan element is one-dimensional, and its conjugates have at most two parameters.  The union of these kernels has dimension at most $3<d$.
\end{proof}

\begin{theorem}
\label{thm821}
Every closed subgroup $H\subset G_2\cdot \SU(3)\subset E_6$ has only normal $E_6$-g.o. metrics.
\end{theorem}

\begin{proof}
It suffices to treat the Lie algebra of $H$. Every proper connected subgroup of $G_2$, except the standard $\SU(3)$ and the principal $\SO(3)$, lies in a conjugate $\SO(4)$.
~\cite[Proposition~15 and Table~5]{Ko2009} Consider the projection of $\mathfrak h$ onto $\mathfrak g_2$. A projection contained in $\mathfrak{so}(4)$ is handled by Theorem~\ref{thm814}. A full projection is handled by Proposition~\ref{prop817}.

If the first projection is $\mathfrak a=\mathfrak{su}(3)$, let $\mathfrak r\subset\mathfrak b$ be the second projection. Every proper connected subgroup of $\SU(3)$ is contained in its standard $\U(2)$ or is its principal $\SO(3)$. The first alternative is covered by Lemma~\ref{lem818}, and the second by Lemma~\ref{lem819}. If $\mathfrak r=\mathfrak b$, Goursat's lemma for two simple algebras gives either their direct product or a graph. The product is covered by Proposition~\ref{prop813}, and the graph by Lemma~\ref{lem820}(1).

The remaining first projection is principal $\mathfrak{su}(2)$. If $\mathfrak h$ contains this factor, it is a product and Proposition~\ref{prop813} applies. Otherwise Goursat's lemma says that its second projection has a simple quotient $\mathfrak{su}(2)$. Inside $\mathfrak{su}(3)$ it is therefore a standard $\mathfrak{su}(2)$, the corresponding $\mathfrak u(2)$, or principal $\mathfrak{so}(3)$. These are exactly Lemma~\ref{lem820}(2)--(3). This exhausts all Lie algebras, including central tori and coupled projections.  Corollary~\ref{cor33} includes disconnected $H$.
\end{proof}

\paragraph{\texorpdfstring{Completion of the $E_7$ tensor subgroup}{Completion of the E7 tensor subgroup}}

\begin{lemma}\label{lem822}
Every subgroup of either $G_2\cdot \Sp(2)\cdot \Sp(1)$ or $G_2\cdot \U(3)$ inside $G_2\cdot \Sp(3)\subset E_7$ has only normal $E_7$-g.o. metrics.
\end{lemma}

\begin{proof}
For $R=\Sp(2)\cdot \Sp(1)$, let $P$ be its standard real eight-dimensional module on $\mathbb H^2$. Restriction of the defining $\Sp(3)$ module gives
$W_{14}|_R=V_5\oplus P\oplus\mathbb R, \qquad \mathfrak{sp}(3)/\mathfrak r=P.$
Here $\Sp(2)$ acts on $V_5$ through $\SO(5)$ and $\Sp(1)$ acts trivially there.  Take
$E=V_7\otimes P,\qquad F=(V_7\otimes V_5)\oplus P.$
A full-rank map $\mathbb R^7\to\mathbb R^8$ with distinct singular values has zero stabilizer even in $\mathfrak{so}(7)\oplus \mathfrak{so}(8)$. The remaining skew operator on a one-dimensional orthogonal complement is zero.  Thus $E$ is free. For $F$, a full-rank map $\mathbb R^5\to\mathbb R^7$ first kills $\mathfrak{sp}(2)$ and then $\mathfrak g_2$, using Lemma~\ref{lem815}.  A nonzero vector of $P$ kills the remaining $\mathfrak{sp}(1)$, which acts by quaternionic scalar multiplication.

For $R=\U(3)$, write
$W_{14}|_R=\mathfrak{su}(3)\oplus U_6, \qquad \mathfrak{sp}(3)/\mathfrak u(3)=P_{12},$
where $U_6$ is the underlying real module of $\Lambda^2\mathbb C^3$ and $P_{12}$ that of $\operatorname{Sym}^2\mathbb C^3$.  Both have nonzero central circle weight. The action of $\mathfrak u(3)$ on $U_6$ is faithful. Take
$E=V_7\otimes U_6,\qquad F=(V_7\otimes\mathfrak{su}(3))\oplus P_{12}.$
For $E$, a full-rank map $\mathbb R^6\to\mathbb R^7$ with distinct singular values kills both factors. For the tensor part of $F$, a full-rank map $\mathbb R^7\to\mathbb R^8$ kills $\mathfrak g_2$ and the semisimple $\mathfrak{su}(3)$ factor.  Any nonzero vector of $P_{12}$ kills the remaining central circle. The two-free-module criterion proves both hereditary claims.
\end{proof}

\begin{lemma}\label{lem823}
Every proper connected subgroup of $\Sp(3)$ is contained in a conjugate of one of
$\Sp(2)\cdot \Sp(1),\qquad \U(3),\qquad \SO(3)\cdot \Sp(1),\qquad A_{1,\mathrm{pr}}.$
The third embedding is $\mathbb C^6=V_3\otimes V_2$, and the last is $\mathbb C^6=V_6$.
\end{lemma}

\begin{proof}
Consider the defining complex symplectic module $V=\mathbb C^6$. If its restriction is reducible, the usual decomposition into real, complex, and quaternionic isotypic components has one of two consequences.  Either there is a proper nonzero invariant nondegenerate symplectic summand, of complex dimension $2$ or $4$, or there is an invariant Lagrangian of complex dimension $3$. Indeed, a quaternionic irreducible summand can be chosen nondegenerate. Real-type multiplicity spaces are symplectic and complex-type summands occur paired with their duals. Unless such a nondegenerate block already fills $V$, it yields the first alternative.  A block filling $V$ with reducible complex restriction has two isotropic halves, yielding the second. The first alternative lies in $\Sp(2)\cdot \Sp(1)$. Compactness in the second gives an invariant Hermitian orthogonal Lagrangian splitting and containment in $\U(3)$.

If $V$ is irreducible, the center of its compact acting Lie algebra acts by scalars and is zero because the representation preserves a nondegenerate alternating form.  The acting algebra is semisimple. For a simple acting algebra, the six-dimensional irreducible modules of rank at most three leave only the standard $C_3$ module and $V_6$ of $A_1$ as symplectic possibilities. The six-dimensional modules of $A_2$ are not self-dual, and $\Lambda^2\mathbb C^4$ of $A_3$ is orthogonal. $B_2,B_3,G_2$ supply no six-dimensional irreducible module. The $C_3$ possibility is the whole $\Sp(3)$ and is excluded. For more than one simple factor, irreducibility forces a tensor product of dimensions $2$ and $3$.  Symplectic type forces the standard $A_1$ module times the real three-dimensional $A_1$ module, giving $\SO(3)\cdot \Sp(1)$.
\end{proof}

\begin{theorem}
\label{thm824}
Every closed subgroup $H\subset G_2\cdot \Sp(3)\subset E_7$ has only normal $E_7$-g.o. metrics.
\end{theorem}

\begin{proof}
Consider the projection $\mathfrak r$ of $\mathfrak h$ onto $\mathfrak{sp}(3)$.  If it is the full algebra, Goursat's lemma and $\dim\mathfrak g_2<\dim\mathfrak{sp}(3)$ imply that $\mathfrak h$ contains the full $\mathfrak{sp}(3)$ factor, so Proposition~\ref{prop813} applies. The first three alternatives of Lemma~\ref{lem823} are covered by Lemma~\ref{lem822} and Theorem~\ref{thm816}(1). Any proper connected subgroup of the remaining principal $A_1$ is a torus, and is also contained in a conjugate maximal torus of $\Sp(2)\cdot \Sp(1)$.  It remains to treat the exact principal $\mathfrak{su}(2)$ projection.

Its action on $W_{14}$ is $V_9\oplus V_5$, since $\Lambda^2 V_6=V_9\oplus V_5\oplus V_1$. If the first projection lies in $\mathfrak{so}(4)$ or is the full $\mathfrak g_2$, the earlier theorems apply. If it is the standard $\mathfrak{su}(3)\subset\mathfrak g_2$, the two projections have no common simple quotient and $\mathfrak h$ is their direct product.  Its two modules
$U\otimes V_9,\qquad U\otimes V_5$
are free.  For the first, use a full-rank map $\mathbb R^6\to \mathbb R^9$ with distinct singular values. It kills $\mathfrak{su}(3)$, and a nonzero $\mathfrak{su}(2)$ operator on $V_9$ has only a one-dimensional kernel. For the second use $\mathbb R^5\to U$ with image containing two complex independent vectors, as before. If the first projection is principal $\mathfrak{su}(2)$, the two projections are either a product, already covered, or a graph.  In the graph case the two orthogonal free modules $V_{11}\subset\mathfrak g_2$ and $V_7\subset\mathfrak{sp}(3)$ lie outside the graph algebra.
$\mathfrak g_2|_{A_{1,\mathrm{pr}}}=V_3\oplus V_{11},\qquad \mathfrak{sp}(3)|_{A_{1,\mathrm{pr}}}=V_3\oplus V_7\oplus V_{11}.$
This exhausts the possible first projections and proves the claim.
\end{proof}

\paragraph{\texorpdfstring{Completion of the $E_8$ tensor subgroup}{Completion of the E8 tensor subgroup}}

\begin{lemma}\label{lem825}
In $E_8$, every subgroup of each of
$G_2\cdot\Spin(9),\qquad G_2\cdot \Sp(3)\cdot \Sp(1),\qquad G_2\cdot \SU(3)\cdot \SU(3)$
inside $G_2\cdot F_4$ has only normal $E_8$-g.o. metrics.
\end{lemma}

\begin{proof}
For $R=\Spin(9)$ the standard branchings are
$W_{26}|_R=\mathbb R\oplus\mathbb R^9\oplus\Delta_{16}, \qquad \mathfrak f_4/\mathfrak{so}(9)=\Delta_{16}.$
Take
$E=V_7\otimes\Delta_{16},\qquad F=(V_7\otimes\mathbb R^9)\oplus\Delta_{16}.$
The spin weights are $\tfrac12(\pm a_1\pm a_2\pm a_3\pm a_4)$. If one $a_i\ne0$, pairing the weights differing only in its sign shows that at most eight weights vanish.  Thus every nonzero spin operator has rank at least eight. Lemma~\ref{lem23}, with $q=7$ and $\dim\mathfrak{so}(9)=36$, supplies a free seven-frame in $\Delta_{16}$.  A map $\mathbb R^7\to\Delta_{16}$ onto that frame, with distinct singular values, proves that $E$ is free. For the first summand of $F$, a full-rank map $\mathbb R^7\to\mathbb R^9$ with distinct singular values kills $\mathfrak g_2$ and leaves only the rotations of a two-dimensional orthogonal complement in $\mathfrak{so}(9)$. Every nonzero such plane rotation acts invertibly on the spin module.  Add any nonzero spinor in the second summand to kill it.

For $R=\Sp(3)\cdot \Sp(1)$ the branchings are
$W_{26}|_R=W_{14}\oplus P_{12},\qquad \mathfrak f_4/\mathfrak r=P_{28},$
where $W_{14}=\operatorname{Herm}_0(3,\mathbb H)$ is trivial under $\Sp(1)$, $P_{12}=\mathbb H^3$ has the standard left-right $\Sp(3)\cdot \Sp(1)$ action, and $(P_{28})_{\mathbb C}=\Lambda^3_0\mathbb C^6\otimes\mathbb C^2$. The five real vectors
$e_1,e_2,e_3,e_1i,e_1j\in\mathbb H^3$
have zero joint stabilizer. Indeed, $Ae_s=e_sb$ for $s=1,2,3$ makes $A=bI_3$, and the other two equations force $[b,i]=[b,j]=0$, hence $b=0$. Complete them to an orthonormal seven-frame. Therefore
$E=V_7\otimes P_{12},\qquad F=(V_7\otimes W_{14})\oplus P_{28}$
are free. For $E$ use that seven-frame. For the tensor part of $F$ choose a rank-seven map with distinct singular values and image containing the three-vector base for $\Sp(3)$ in $W_{14}$. After $\mathfrak g_2$ and $\mathfrak{sp}(3)$ vanish, a nonzero vector of $P_{28}$ kills $\mathfrak{sp}(1)$, whose standard two-dimensional complex factor has no zero weight.

For $R=\SU(3)_a\cdot \SU(3)_b$, choose $a$ to be the short-root factor of $F_4$.  Then
$W_{26}|_R=(\mathfrak{su}(3)_a,1)\oplus P_{18},\qquad \mathfrak f_4/\mathfrak r=P_{36},$
where $P_{18}=(\mathbb C^3_a\otimes(\mathbb C^3_b)^*)_{\mathbb R}$ and $P_{36}=(\operatorname{Sym}^2\mathbb C^3_a\otimes \mathbb C^3_b)_{\mathbb R}$. On $P_{18}\simeq\operatorname{Mat}_3(\mathbb C)$ the three matrices $I$, a real diagonal matrix with distinct entries, and a matrix with all off-diagonal entries nonzero have zero joint stabilizer for $\mathfrak{su}(3)_a\oplus\mathfrak{su}(3)_b$. Indeed, their successive stabilizers lie in the diagonal $\mathfrak{su}(3)$, its diagonal Cartan, and zero. Choose
$E=V_7\otimes P_{18},\qquad F=(V_7\otimes\mathfrak{su}(3)_a)\oplus P_{36}.$
The three-vector base makes $E$ free by the matrix observation. In $F$ a full-rank map $\mathbb R^7\to\mathbb R^8$ with distinct singular values first kills $\mathfrak g_2$ and $\mathfrak{su}(3)_a$.  Two complex independent coefficients of $\mathbb C^3_b$ in $P_{36}$ then kill $\mathfrak{su}(3)_b$. This proves all three hereditary claims. The displayed $F_4$ branchings are also recorded in Shan, \S4.3.1--4.3.3, loc.~cit.
\end{proof}

\begin{theorem}
\label{thm826}
Every closed subgroup $H\subset G_2\cdot F_4\subset E_8$ has only normal $E_8$-g.o. metrics.
\end{theorem}

\begin{proof}
If the second projection is the full $\mathfrak f_4$, dimension and Goursat's lemma imply that $\mathfrak h$ contains that full factor, so Proposition~\ref{prop813} applies. Otherwise its projection lies in one of the maximal connected proper subgroups
$\Spin(9),\quad \Sp(3)\cdot \Sp(1),\quad \SU(3)\cdot \SU(3),\quad G_2\times \SO(3),\quad A_{1,\mathrm{pr}}$
of $F_4$. See Shan, \S4.3, loc.~cit. The first four cases are covered by Lemma~\ref{lem825} and Theorem~\ref{thm816}(2). Proper subgroups of the final principal $A_1$ are toral and also lie in a conjugate $\Spin(9)$, which has maximal rank.  Thus only the exact principal $A_1$ projection remains.

Under this principal $A_1$,
$W_{26}=V_{17}\oplus V_9,\qquad \mathfrak f_4=V_3\oplus V_{11}\oplus V_{15}\oplus V_{23}.$
The first equality is Shan's partition $[17,9]$. The second follows from the exponents $1,5,7,11$. As in the $E_7$ argument, a first projection in $\SO(4)$ or equal to $G_2$ is already covered.  For first projection the standard $\SU(3)$, the projections have no common simple quotient and the isotropy is a direct product.  The two modules $U\otimes V_{17}$ and $U\otimes V_9$ are free. Maps $\mathbb R^6\to\mathbb R^{17}$ and $\mathbb R^6\to\mathbb R^9$ with distinct singular values kill $\mathfrak{su}(3)$, and a nonzero operator of the remaining $\mathfrak{su}(2)$ has only a one-dimensional kernel in either target. For first projection principal $A_1$, there is a product or a graph.  Products are covered by Proposition~\ref{prop813}. For a graph use the two orthogonal free modules $V_{11}\subset\mathfrak g_2$ and $V_{15}\subset\mathfrak f_4$. This exhausts the connected isotropy groups. Corollary~\ref{cor33} includes disconnected isotropy.
\end{proof}

\paragraph{\texorpdfstring{Hereditary normality inside $F_4\cdot \SO(3)\subset E_7$}{Hereditary normality inside F4  \SO(3) in  E7}}
\label{sec:e7-f4-a1-hereditary}

Consider the standard maximal tensor subgroup $K=F_4\cdot \SO(3)\subset E_7$.  Its isotropy representation is
$\mathfrak e_7=\mathfrak f_4\oplus\mathfrak{so}(3) \oplus(W_{26}\otimes\mathbb R^3), \qquad W_{26}=\operatorname{Herm}_0(3,\mathbb O).$
This is the other nonsymmetric irreducible tensor decomposition of $\mathfrak e_7$. See \cite[pp.~21--22]{BeElMa2009}. Throughout the section a free vector means a vector with zero infinitesimal stabilizer, and normality refers to the fixed acting group $E_7$.

We use the matrix observation from the tensor-chain section. If $T:U\to V$ has full column rank and distinct positive singular values, then a stabilizer equation $BT=TA$, with both operators skew-adjoint, forces $A=0$ and forces $B$ to vanish on $\operatorname{im}T$.

\begin{lemma}\label{lem827}
The real spin representation $\Delta_{16}$ of $\Spin(9)$ has three vectors with zero joint infinitesimal stabilizer. Consequently $\Delta_{16}\otimes\mathbb R^3$ has a free vector for $\mathfrak{so}(9)\oplus\mathfrak{so}(3)$.
\end{lemma}

\begin{proof}
The stabilizer of a unit spinor is $\Spin(7)$, and restriction of the spin representation to this stabilizer is
$\Delta_{16}=\mathbb R\oplus\mathbb R^7\oplus\Delta_8.$
Choose a second spinor with nonzero components in $\mathbb R^7$ and $\Delta_8$.  The stabilizer of a unit vector in the spin module $\Delta_8$ is $G_2$. Choose the seven-dimensional component nonzero, so that the common stabilizer is $\SU(3)$. Under this $\SU(3)$ one has
$\Delta_{16}=\mathbb R^4\oplus (\mathbb C^3)_{\mathbb R}\oplus (\mathbb C^3)_{\mathbb R}.$
Choose a third spinor with complex linearly independent components in the two copies of $\mathbb C^3$.  Its stabilizer in $\mathfrak{su}(3)$ is zero. These standard spinor restrictions also follow from the octonion model $\Delta_{16}=\mathbb O^2$. One can choose the successive spinors $(1,0)$, $(i,1)$, and $(j,\ell)$, with $i,j,\ell$ as in the three-vector octonion base. Orthogonalize the three spinors and use them, with distinct lengths, as the columns of $T:\mathbb R^3\to\Delta_{16}$.
\end{proof}

\begin{lemma}\label{lem828}
For $J=\Spin(9)\cdot \SO(3)\subset K$, write
$\mathfrak e_7/\mathfrak j=E\oplus F, \qquad E=\Delta_{16}\otimes\mathbb R^3, \quad F=(\mathbb R^9\otimes\mathbb R^3) \oplus\Delta_{16}\oplus\mathbb R^3.$
There is a vector of $F$ whose infinitesimal stabilizer is a root $\mathfrak{su}(2)\subset\mathfrak{so}(9)$, with
$\mathbb R^9|_{\mathfrak{su}(2)}=\mathbb R^5\oplus\mathbb H,$
the action on $\mathbb H$ being quaternionic scalar multiplication. Its external $\mathfrak{so}(3)$ stabilizer is zero.
\end{lemma}

\begin{proof}
In $\mathbb R^9\otimes\mathbb R^3$ take a full-rank $T:\mathbb R^3\to\mathbb R^9$ with distinct singular values. Its stabilizer is $\mathfrak{so}(6)$ fixing the three-dimensional image pointwise, and has zero component in the external factor. Under $\Spin(6)\simeq \SU(4)$,
$\Delta_{16}=(\mathbb C^4)_{\mathbb R} \oplus(\mathbb C^4)_{\mathbb R}.$
Choose independent unit vectors in the two complex summands. Their stabilizer is the $\SU(2)$ acting on the orthogonal complex two-plane.  The vector representation of $\Spin(6)$ is the real form of $\Lambda^2\mathbb C^4$. Its restriction to this $\SU(2)$ is $\mathbb R^2\oplus\mathbb H$. Adding the three directions fixed initially gives the asserted $\mathbb R^5\oplus\mathbb H$ decomposition. The last summand of $F$ does not change this stabilizer.
\end{proof}

\begin{proposition}
\label{prop829}
Every closed subgroup $H\subset J$ whose projection to $\Spin(9)$ preserves a proper nonzero subspace of $\mathbb R^9$ has only normal $E_7$-g.o. metrics.
\end{proposition}

\begin{proof}
Its Lie algebra lies in a conjugate of
$\mathfrak j_p=\mathfrak{so}(p)\oplus\mathfrak{so}(9-p) \oplus\mathfrak{so}(3),\qquad 1\le p\le4.$
The module $E$ in Lemma~\ref{lem828} is free for the full $\mathfrak j$ by Lemma~\ref{lem827}. We show that $F$ is free for $\mathfrak j_p$. Let $S$ be the $\SU(2)$ stabilizer constructed in the preceding lemma, and identify $\mathbb R^9=\mathbb R^5\oplus\mathbb H$ as there.  Choose independent $u_1,\ldots,u_p\in\mathbb R^5$ and vectors $w_i\in\mathbb H$ with $w_1\ne0$. The $p$-plane
$P=\operatorname{span}_{\mathbb R}\{u_i+w_i:1\le i\le p\}$
has injective projection to $\mathbb R^5$. If an element of $\mathfrak s$ preserves $P$, its action on $P$ must vanish, because it acts trivially on $\mathbb R^5$. It therefore kills $w_1$, forcing that element to be zero. Conjugate the vector of $F$ so that this $P$ becomes the prescribed $p$-plane of the block subgroup.  The stabilizer of the resulting vector intersects $\mathfrak j_p$ trivially. Thus both $E$ and $F$ are free for $\mathfrak j_p$ and remain free for $\mathfrak h$. The two-free-module criterion applies.
\end{proof}

\begin{proposition}
\label{prop830}
Every closed subgroup $H\subset\Spin(9)\cdot \SO(3) \subset E_7$ has only normal $E_7$-g.o. metrics.
\end{proposition}

\begin{proof}
Write $\mathfrak r\subset\mathfrak{so}(9)$ for its first projection. The reducible case is Proposition~\ref{prop829}. If the nine-dimensional real representation is irreducible and $\mathfrak r\ne\mathfrak{so}(9)$, then $\mathfrak r$ is semisimple. Its center would otherwise induce a complex structure on an odd-dimensional irreducible real space. Its complexification is also irreducible, since the real dimension is odd.  The list of irreducible modules of complex dimension nine and rank at most four gives just the principal $\mathfrak{su}(2)$ module $V_9$, or $V_3\otimes V_3$ for $\mathfrak{so}(3)\oplus\mathfrak{so}(3)$, besides the excluded full standard $\mathfrak{so}(9)$. This small-dimensional assertion can also be checked directly. Among simple algebras only $A_1$ and $B_4$ have an orthogonal nine-dimensional irreducible module. For a product the only nontrivial factorization is $9=3\cdot3$.

In the tensor case put $\mathbb R^9=V_3\otimes V_3$. Then
$\mathfrak{so}(9)/(\mathfrak{so}(3)\oplus\mathfrak{so}(3)) =(V_3\otimes V_5)\oplus(V_5\otimes V_3).$
Each summand has a free vector for the two factors. An injective map $\mathbb R^3\to\mathbb R^5$ with distinct singular values first kills one factor, and the spin-two action of a nonzero $\mathfrak{so}(3)$ element has only a one-dimensional kernel. Use the orthogonal modules
$E'=(V_3\otimes V_5)\oplus(\mathbb R^9\otimes\mathbb R^3), \qquad F'=(V_5\otimes V_3)\oplus(\Delta_{16}\otimes\mathbb R^3).$
In each, the first component kills the two internal factors. Then choose two independent coefficients in the external $\mathbb R^3$ to kill the remaining factor. Both are free for the full three-factor envelope.

For the principal nine-dimensional $A_1$,
$\mathfrak{so}(9)|_{A_1}=V_3\oplus V_7\oplus V_{11}\oplus V_{15}.$
Both $V_{15}$ and $V_{11}$ have generic free vectors for $A_1$. The same argument with
$E'=V_{15}\oplus(\mathbb R^9\otimes\mathbb R^3),\qquad F'=V_{11}\oplus(\Delta_{16}\otimes\mathbb R^3)$
gives hereditary normality for that two-factor envelope.

It remains to treat $\mathfrak r=\mathfrak{so}(9)$. Goursat's lemma gives $\mathfrak h=\mathfrak{so}(9)\oplus\mathfrak r_0$, with $\mathfrak r_0\subset\mathfrak{so}(3)$. The centralizer of this $B_4$ in $\mathfrak e_7$ is $\mathfrak{so}(4)$, and its centralizer envelope has quotient
$(\mathbb R^9\otimes\mathbb R^3\otimes1) \oplus(\Delta_{16}\otimes\mathbb R^4),$
of dimensions $27$ and $64$. For precision, this $B_4$ embedding is unique in $E_7$, with centralizer $A_1^1+A_1^2$. See \cite[Table~13, p.~256]{Mi2006}.
It may therefore be realized inside $\Spin(12)\cdot \SU(2)$ using $\mathbb R^{12}=\mathbb R^9\oplus\mathbb R^3$, where the displayed branching is immediate from the half-spin representation. The two-summand g.o. classification of Chen--Nikonorov gives normality on this envelope quotient, since its nonnormal list contains no ambient $E_7$ case. The second summand has a free vector. Use a map $\mathbb R^4\to\Delta_{16}$ with distinct singular values whose image contains the three spinors of Lemma~\ref{lem827}. The ideal-envelope descent lemma now proves normality whenever the isotropy contains this $B_4$ as an ideal. All possibilities for $\mathfrak r$ have been exhausted.
\end{proof}

\begin{theorem}
\label{thm831}
Every closed subgroup $H\subset F_4\cdot \SO(3)\subset E_7$ has only normal $E_7$-g.o. metrics.
\end{theorem}

\begin{proof}
If the first projection is the full $\mathfrak f_4$, the isotropy contains that simple factor by Goursat's lemma. Equivariance separates the $F_4$-trivial vertical complement from $W_{26}\otimes\mathbb R^3$. The latter is irreducible for $K$ and has a free vector, obtained from the three-vector base for $F_4$ on $W_{26}$. Scalar descent with a free vector gives normality.

Otherwise the first projection lies in one of the five maximal connected proper subgroups of $F_4$.
$\Spin(9),\quad \Sp(3)\cdot \Sp(1),\quad \SU(3)\cdot \SU(3),\quad G_2\times \SO(3),\quad A_{1,\mathrm{pr}}.$
The list and branchings are given in \cite[Section~4.3]{Sh2024}. The first case is Proposition~\ref{prop830}. The $G_2\times \SO(3)$ case is the hereditary three-factor envelope already proved in the tensor-chain section. Its quotient is $(7,3,3)\oplus(7,5,1)\oplus(1,5,3)$. For a principal $A_1$ in $F_4$, the two modules $V_{17}\otimes\mathbb R^3$ and $V_9\otimes\mathbb R^3$ are free for the two $A_1$ factors, by distinct-singular-value maps from $\mathbb R^3$ and the one-dimensional kernel of every nonzero operator on an odd-dimensional $A_1$ module.

For $R=\Sp(3)\cdot \Sp(1)$, use $W_{26}=W_{14}\oplus P_{12}$ and $\mathfrak f_4/\mathfrak r=P_{28}$ as in the $E_8$ argument. The representation $P_{12}=\mathbb H^3$ has the following three-vector base for $\Sp(3)\cdot \Sp(1)$.
$q_1=e_1,\qquad q_2=e_1i+e_2,\qquad q_3=e_1j+e_3.$
Indeed, from $Aq_s=q_sb$ and skew-adjointness one obtains
$[b,\langle q_1,q_2\rangle]=[b,i]=0,\qquad [b,\langle q_1,q_3\rangle]=[b,j]=0.$
Thus $b=0$. The three columns are a quaternionic basis and give $A=0$.  Therefore the two orthogonal modules
$E=P_{12}\otimes\mathbb R^3,\qquad F=(W_{14}\otimes\mathbb R^3)\oplus P_{28}$
are free.  For $E$ use the three-vector base just constructed. For $F$ use the three-vector base of $\Sp(3)$ on $W_{14}$, which first kills $\Sp(3)$ and the external $\SO(3)$. Any nonzero vector of $P_{28}$ then kills the remaining $\Sp(1)$.

Finally let $R=\SU(3)_a\cdot \SU(3)_b$, with $a$ the short-root factor.  Write $W_{26}=(\mathfrak{su}(3)_a,1)\oplus P_{18}$ and $\mathfrak f_4/\mathfrak r=P_{36}$ as before. The two modules
$E=P_{18}\otimes\mathbb R^3,\qquad F=(\mathfrak{su}(3)_a\otimes\mathbb R^3)\oplus P_{36}$
are free.  The three-matrix base on $P_{18}$ proves this for $E$. For $F$, the adjoint $\SU(3)$ module has a three-vector base, so its tensor factor kills $\SU(3)_a$ and the external $\SO(3)$. Two independent standard $\SU(3)_b$ coefficients in $P_{36}$ kill the remaining factor. In every case two free modules were constructed for the entire envelope, hence remain free for any subgroup, including coupled projections and central tori.  Disconnected isotropy introduces no further case at Lie algebra level.
\end{proof}

\paragraph{All subgroups of three further nonregular maximal subgroups}
\label{sec:small-exceptional-maximals}

In this section, we exclude all subgroups of the maximal $A_1A_1$ and $G_2A_1$ in $E_7$, and the maximal $A_2A_1$ in $E_8$. These conclusions concern arbitrary subgroups of the indicated maximal subgroup, rather than just the maximal isotropy itself.

\begin{theorem}\label{thm832}
Let $K\subset E_8$ be the maximal subgroup of type $A_2A_1$. For every closed $H\subset K$, every $E_8$-g.o. metric on $E_8/H$ is normal.
\end{theorem}

\begin{proof}
For this embedding the complex adjoint branching is
\begin{equation}\label{formula35}
 \mathfrak e_{8,\mathbb C}|_{A_2A_1}
 =V_{(1,1;0)}\oplus V_{(0,0;2)}\oplus V_{(1,1;6)}
 \oplus V_{(3,0;4)}\oplus V_{(0,3;4)}\oplus V_{(2,2;2)}.
\end{equation}
Here $(a,b;c)$ denotes the fundamental-weight coordinates. The quotient has orthogonal real irreducible summands of dimensions $56,100,81$. The branching is explicitly recorded in the official Sage reference manual, ``Branching rules'', under ``Miscellaneous''~\cite{Sa2026}. Its listed factor order is $A_1A_2$. Interchanging the factors gives \eqref{formula35}. Restriction of the $240$ roots and the eight zero Cartan weights gives the same six highest weights, each with multiplicity one. We use only the first two real summands, namely $U=\mathfrak{su}(3)\otimes_{\mathbb R}W_7$ and $V=(\operatorname{Sym}^3\mathbb C^3\otimes W_{5,\mathbb C})_{\mathbb R}$, where $W_j$ is the real irreducible $\mathfrak{su}(2)$-module of dimension $j$.

We show that both summands have a free vector for $\mathfrak k=\mathfrak{su}(3)\oplus\mathfrak{su}(2)$, whose dimension is eleven.  Write an infinitesimal element as $(B,D)$. If $D\ne0$, its eigenvalues on $W_{7,\mathbb C}$ are distinct. For each eigenvector of $\ad B$ at most one of these eigenvalues cancels its eigenvalue.  Thus the kernel on $U$ has dimension at most eight.  If $D=0$ and $B\ne0$, its dimension is at most $4\cdot7=28$, since a nonzero element of $\mathfrak{su}(3)$ has centralizer of dimension at most four.  The kernel codimension on $U$ is consequently at least $28>11$.

On the complex space underlying $V$, if $D\ne0$, the five distinct eigenvalues give a complex kernel of dimension at most ten.  Hence the real kernel has dimension at most twenty. If $D=0$ and $B\ne0$, the action of $B$ on $\operatorname{Sym}^3\mathbb C^3$ is nonzero and trace-free. It is diagonalizable, so it has at least two nonzero eigenvalues. Its complex rank is at least two.  After tensoring with $W_{5,\mathbb C}$ and passing to the underlying real space, the rank is at least twenty.  Thus the kernel codimension on $V$ is always at least $20>11$.

Lemma~\ref{lem23} gives the two free vectors. The hereditary form of Lemma~\ref{lem36} proves the assertion, including arbitrary central tori, diagonal subgroups and disconnected subgroups of $K$.
\end{proof}

\begin{theorem}\label{thm833}
Let $K\subset E_7$ be a maximal subgroup of type $A_1A_1$. For every closed $H\subset K$, every $E_7$-g.o. metric on $E_7/H$ is normal.
\end{theorem}
\begin{proof}
We work with Lie algebras.  A subalgebra of $\mathfrak{su}(2)\oplus\mathfrak{su}(2)$ is abelian, simple, the full algebra, or, up to interchanging factors, $\mathfrak{su}(2)\oplus\mathbb R$. Indeed, a nonsimple nonabelian proper subalgebra projects surjectively onto one factor, and the subdirect alternative with both projections surjective is either the full product or a simple graph.  The remaining residual subalgebra is a circle. Abelian and simple isotropy are covered by the known classifications. We give a direct argument for the full product. Its tangent space has odd dimension $127$, so it contains a real irreducible summand of odd dimension. Both simple factors act nontrivially on every irreducible summand. The normalizer of either factor is precisely $\mathfrak k$ by maximality, and a fixed tangent vector would centralize that factor. Thus the odd summand is $W_a\otimes W_b$ with odd $a,b\geq3$. A nonzero infinitesimal operator has kernel dimension at most $\max(a,b)$, since each nonzero $A_1$ operator has distinct weights. Its kernel codimension is at least $ab-\max(a,b)\geq6=\dim\mathfrak k$. Lemma~\ref{lem23} gives a free vector in that irreducible summand.

For completeness, this forces a $G$-g.o. cometric $C$ to be scalar. An equivariant spectral projection embeds the free irreducible module into some eigenspace $E_\lambda$ of $C$. For a different eigenspace $E_\mu$, the g.o. equation at $tx+y$, where $x\in E_\lambda$ and $y\in E_\mu$, gives
$(\lambda-\mu)t[x,y]=t[Z_t,x]+[Z_t,y].$
Consequently $[x,y]=u+v\in E_\lambda\oplus E_\mu$. If $x$ is free, the equation $[Z_t,x]=(\lambda-\mu)u$ determines $Z_t$ independently of $t$, while $[Z_t,y]=t(\lambda-\mu)v$ forces $v=0$. Density gives $[E_\lambda,E_\mu]\subset E_\lambda$. Invariance of $Q$ then shows that $\mathfrak k\oplus(E_\lambda^\perp\cap\mathfrak m)$ is a subalgebra. A second eigenvalue would make this a strict intermediate subalgebra, contradicting maximality.

Consider $\mathfrak h=\mathfrak s\oplus\mathbb R$ and write $\mathfrak k=\mathfrak s\oplus\mathfrak a$, with both factors simple. Maximality gives $\mathfrak n_{\mathfrak e_7}(\mathfrak s)=\mathfrak k$. This normalizer contains $\mathfrak k$ and cannot be all of $\mathfrak e_7$, since $\mathfrak s$ is not an ideal there. Consequently $\mathfrak k\ominus\mathfrak h$ is fixed by $\mathfrak s$, whereas $\mathfrak k^\perp$ has no fixed vector. Every invariant cometric preserves the two spaces. Descent makes its horizontal restriction normal. The isotropy representation of this semisimple $K$ has abelian principal stabilizer, by Benoist--Kobayashi's complete exceptional list.  The abelian-stabilizer descent criterion therefore makes the total metric normal. Corollary~\ref{cor33} finishes the proof for arbitrary closed $H$.
\end{proof}

\begin{theorem}\label{thm834}
Let $K\subset E_7$ be the maximal subgroup of type $G_2A_1$ whose tangent representation is $ (V_7\otimes W_5)\oplus(V_{27}\otimes W_3)$. Here $V_{27}=\operatorname{Sym}^2_0V_7$ and $W_j$ is the real $\SO(3)$-module of dimension $j$. Every closed subgroup $H\subset K$ has only normal $E_7$-g.o. metrics.
\end{theorem}
\begin{proof}
The quotient $E_7/K$ has precisely two real irreducible summands, of dimensions $35$ and $81$, so its $G$-g.o. metrics are normal by the two-isotropy-summand classification. Consider the projection of $\mathfrak h$ onto $\mathfrak g_2$. If it is contained in $\mathfrak{so}(4)$, the two free modules used in Proposition~\ref{prop84} still occur here, inside $(\mathfrak g_2/\mathfrak{so}(4))\oplus(V_7\otimes W_5)$. The proof uses only their representation structures and is unchanged by the additional $V_{27}\otimes W_3$ summand. Thus all subgroups of $\SO(4)\SO(3)$ are excluded.

The other possible projections are $G_2$, the standard $\SU(3)$, and principal $\SO(3)$.  Goursat's lemma leaves their products with $0$, $\SO(2)$ or $\SO(3)$, and a simple graph in the last case. The simple cases are already normal by the simple-isotropy classification.  For full $G_2$ projection with a circle, the vertical two-plane is $G_2$-trivial and the horizontal space has no fixed vector.  Therefore $C$ preserves the horizontal space, which becomes scalar by descent.  The rank-one test $X=v\otimes w\in V_7\otimes W_5$ and a nonzero vertical vector gives equality with the vertical coefficient, exactly as in Proposition~\ref{prop83}. Full $K$ is the already treated base.

For $\SU(3)\times \SO(2)$ use the two free modules in the circle case of Theorem~\ref{thm86}.  They lie entirely in $\mathfrak g_2\oplus\mathfrak{so}(3)\oplus(V_7\otimes W_5)$, so remain available in the present quotient. For $\SU(3)\times \SO(3)$ the vertical complement is $U_6$. It is trivial on the external factor, whereas both horizontal summands are nontrivial on that factor.  Hence there is no mixing, and descent makes the horizontal coefficient scalar. Testing $Y\in U_6$ and $X=e_0\otimes w$, where $e_0$ is the $\SU(3)$-fixed unit vector in $V_7$, gives equality of coefficients by the same orthogonal first-factor projection as before.

Finally, under principal $A_1\subset G_2$, $\mathfrak g_2/A_1=V_{11}$, $V_7=V_7$, and $V_{27}=V_{13}\oplus V_9\oplus V_5$. The latter equality follows by taking the trace-free symmetric square of the spin-three representation. None of the horizontal types equals $V_{11}$ or the trivial type. Consequently $C$ preserves the horizontal space for each product with a circle or with the full external $\SO(3)$. Descent makes that space scalar.  A generic free vector of $V_{11}$ and a rank-one tensor in $V_7\otimes W_5$ force the $V_{11}$ coefficient to agree with it. If the external complement is a two-plane, its coefficient agrees as well by the earlier circle test. These are all possibilities for the two projections.
\end{proof}

\section{Nonregular isotropy in exceptional groups}\label{sec:unified-exceptional-nonregular}
In this section, we complete the exceptional classification. We first prove the remaining centralizer cases and then use the complete embedding inventories in Appendices~\ref{app:high-rank-inventories} and~\ref{app:rank-one-inventories}. Every entry is accompanied by its normality argument, so the final reductions include arbitrary residual isotropy in a full centralizer.

\begin{theorem}\label{thm91}
Every $F_4$-g.o. metric on $F_4/H$, for an arbitrary closed subgroup $H$, is normal.
\end{theorem}
\begin{proof}
If $\mathfrak h$ is abelian or simple, Proposition~\ref{prop71} applies. Otherwise choose a simple ideal $\mathfrak s$ of $\mathfrak h$. Then $\mathfrak h=\mathfrak s\oplus\mathfrak r$ with $\mathfrak r\subset\mathfrak c_{\mathfrak f_4}(\mathfrak s)$. A regular $\mathfrak s$ is covered by Corollary~\ref{cor716}.

Suppose $\operatorname{rank}\mathfrak s\ge2$. Table~\ref{tab:f4-high-inventory} lists every embedded class and its centralizer. Zero centralizers reduce to simple isotropy. The $B_3$ case has a one-dimensional centralizer, so its only other possibility is maximal rank. The $G_2$ case is Theorem~\ref{thm86}. For $B_2$, a residual subalgebra of $A_1^2$ is regular or is a diagonal $A_1$. The latter is Lemma~\ref{lem811}. All other positive-centralizer entries are regular simple ideals.

For rank one, Table~\ref{tab:f4-a1-inventory} lists all fifteen classes. The two root classes have already been treated. Four classes have zero centralizer and reduce to simple isotropy. The remaining nine entries are covered by the envelope lemmas or by the pairs of free modules specified in the table and proved free in Lemma~\ref{lem88}. These modules lie in the complement of the full centralizer envelope, so every residual compact algebra is allowed. This proves the assertion for all isotropy Lie algebras. Corollary~\ref{cor33} gives the same conclusion for disconnected isotropy.
\end{proof}

\paragraph{\texorpdfstring{Nonregular isotropy in $E_6$. Centralizer envelopes}{Nonregular isotropy in E6. Centralizer envelopes}}
\label{sec:e6-nonregular}

All statements concern the prescribed transitive group $E_6$. We work at Lie algebra level with $Q=-B_{\mathfrak e_6}$. By Corollary~\ref{cor33}, passing from $H^0$ to $H$ leaves the $G$-g.o. cone unchanged. The following argument is useful precisely because normality for a maximal subgroup does not by itself imply normality for its subgroups.

\begin{proposition}\label{prop92}
Every $E_6$-g.o. metric on $E_6/H$ is normal if $\mathfrak h$ contains a simple ideal in any of the following embedded classes.
$G_2^1,\quad A_2^{2\prime\prime},\quad B_3^1, \quad A_1^2\ (\text{class }2A_1),\quad B_2^1.$
The superscripts and primes agree with Table~12 of Minchenko~\cite{Mi2006}. The assertion allows an arbitrary compact residual algebra in the centralizer of the specified ideal.
\end{proposition}

\begin{proof}
We check the hypotheses of Corollary~\ref{cor315} for three envelopes. Finite central quotients of the displayed groups have no effect on the argument.

\emph{The envelope $G_2\times \SU(3)$.} Here
$\mathfrak e_6=(\mathfrak g_2,1)\oplus(1,\mathfrak{su}(3)) \oplus(\mathbb R^7\otimes\mathbb R^8),$
where $\mathbb R^8$ is the adjoint $\SU(3)$-module. The horizontal module is irreducible, so its only invariant cometrics are scalar. It has a vector with zero infinitesimal stabilizer even under the larger group $\SO(7)\times \SO(8)$. View it as a real $7\times8$ matrix of full rank with seven distinct nonzero singular values. An infinitesimal stabilizer preserves its two Gram matrices. The seven nonzero simple eigenspaces and the one-dimensional kernel force both skew-symmetric infinitesimal factors to vanish. This treats both ideals of the envelope.

\emph{The envelope $\Spin(7)\times \SU(2)\times \U(1)$.} Use $\mathbb R^{10}=\mathbb R^7\oplus\mathbb R^3$ inside $\Spin(10)\times \U(1)\subset E_6$. Its horizontal representation is
$(\mathbb R^7\otimes\mathbb R^3) \ \oplus\ (\Delta_7\otimes\mathbb C^2)_{\mathbb R}, \qquad \dim_{\mathbb R}\Delta_7=8.$
These are two inequivalent irreducible real modules of dimensions $21$ and $32$. The classification of $G$-g.o. spaces with two isotropy summands, Chen--Nikonorov~\cite{ChNi2019}, Theorem~2, implies normality for this base. The only exceptional simple ambient entry permitting a nonscalar metric is $E_6/\Spin(10)$. On the second summand the action is a subgroup of $\SO(8)\times \U(2)\subset \SO(8)\times \SO(4)$ acting on real $8\times4$ matrices. Choose its four columns to be an orthogonal four-frame with trivial infinitesimal stabilizer in $\Spin(7)$, scaled to distinct lengths. Such a frame exists. In the octonion model take $1,i,j,\ell$. Successive stabilizers lie in $G_2,\SU(3),\SU(2)$, and the last $\SU(2)$ acts freely on unit vectors in the complementary quaternionic line. If $AX=XB$, preservation of $X^tX$ forces $B=0$. Then $A$ fixes this frame and hence vanishes. The two simple ideals have classes $B_3^1$ and $A_1^2$, respectively.

\emph{The envelope $\Spin(5)\times \Spin(5)\times \U(1)$.} Inside $\Spin(10)\times \U(1)$, split $\mathbb R^{10}=\mathbb R^5\oplus\mathbb R^5$. Its horizontal representation is
$(\mathbb R^5\otimes\mathbb R^5) \ \oplus\ (\mathbb C^4\otimes\mathbb C^4)_{\mathbb R},$
with a nonzero central circle weight on the second summand. Again these are two inequivalent real irreducible summands, of sizes $25$ and $32$, so the same two-summand classification makes the base normal. We verify that the second summand already has a locally free vector. Write its infinitesimal stabilizer equation as
$AX+XB^t+itX=0, \qquad A,B\in\mathfrak{sp}(2)\subset\mathfrak{su}(4),\quad t\in\mathbb R.$
For invertible $X$, the determinant, or the trace after multiplying by $X^{-1}$, gives $t=0$. The remaining equation implies that $B^t$ commutes with the positive Hermitian Gram matrix $X^*X$. Choose this matrix with distinct eigenvalues and choose its eigenbasis so that the symplectic form has all six upper-triangular entries nonzero. Such a basis exists. The real skew-orthogonal matrix
\begin{equation*}
 J'=\begin{pmatrix}
 0&a&b&c\\-a&0&-c&b\\-b&c&0&-a\\-c&-b&a&0
 \end{pmatrix},\qquad a^2+b^2+c^2=1,\quad abc\ne0,
\end{equation*}
is unitarily congruent to the standard symplectic form. A skew-Hermitian matrix commuting with the chosen Gram matrix is diagonal, say $i\operatorname{diag}(t_1,t_2,t_3,t_4)$, in this basis. The symplectic condition reads $(t_i+t_j)J'_{ij}=0$ for every $i<j$, and forces all $t_i=0$. Thus $B=0$ and then $A=0$. Each simple factor has class $B_2^1$.
\end{proof}

\begin{proposition}\label{prop93}
Every $E_6$-g.o. metric is normal if its isotropy has a simple ideal of class $A_2^{2\prime}$.
\end{proposition}
\begin{proof}
Inside the regular $\SU(3)^3$, this ideal is the diagonal of the first two factors with the same defining representations. Its centralizer is the third $\SU(3)$. For the envelope $L=\SU(3)_{\mathrm{diag}}\times \SU(3)$ the horizontal representation is
\begin{equation*}
 \mathfrak p=U\oplus V\oplus W,
 \qquad U=\mathfrak{su}(3)_{\mathrm{anti}},\quad
 V=(\operatorname{Sym}^2\mathbb C^3\otimes\mathbb C^3)_{\mathbb R},
 \quad W=(\Lambda^2\mathbb C^3\otimes\mathbb C^3)_{\mathbb R}.
\end{equation*}
Their dimensions are $8,36,18$, and they are mutually inequivalent real irreducibles. The anti-diagonal action on $\mathbb C^3\otimes\mathbb C^3$ is $A\otimes I-I\otimes A$. It interchanges the symmetric and alternating squares nontrivially. Consequently $[U,V]_W\ne0$ and $[U,W]_V\ne0$. Write the three coefficients of a $G$-g.o. cometric as $a_U,a_V,a_W$. For $u\in U$ and $v\in V$, project the g.o. equation at $u+v$ to $W$. Since the isotropy preserves $U$ and $V$, the compensating bracket has no $W$ component. Hence $(a_U-a_V)[u,v]_W=0$. Thus $a_U=a_V$. The same argument at $u+w$, projected to $V$, gives $a_U=a_W$. This argument works even when the coefficient on the third module coincides with one of the first two. Thus $G/L$ is normal.

The module $U\oplus W$ has a locally free vector for $L$. Identify $W$ with complex $3\times3$ matrices, with the usual left--right $\SU(3)$ action, up to dualizing one factor. A positive diagonal matrix with distinct entries has infinitesimal stabilizer the common diagonal Cartan. Choose the component in $U$ with every off-diagonal entry nonzero. Its centralizer in this Cartan is zero. Corollary~\ref{cor315} now applies.
\end{proof}

\begin{proposition}\label{prop94}
Every $E_6$-g.o. metric is normal if its isotropy has a simple ideal of class $B_2^2$.
\end{proposition}
\begin{proof}
Realize $\SO(5)\subset \SU(5)\subset \SU(6)$ by its real defining representation. Its centralizer in $E_6$ is $\SU(2)\times \U(1)$, where $\SU(2)$ is the other factor of $\SU(6)\times \SU(2)$. For $L=\SO(5)\times \SU(2)\times \U(1)$ one obtains three inequivalent irreducible real modules
$\mathfrak p=U\oplus V\oplus W, \quad U=\operatorname{Sym}_0(\mathbb R^5),\quad V=(\mathbb C^5)_{\mathbb R},\quad W=(\Lambda^2\mathbb C^5\otimes\mathbb C^2)_{\mathbb R},$
of dimensions $14,10,40$. The central circle has nonzero weights on $V$ and $W$ and acts trivially on $U$. A $G$-g.o. cometric is scalar on these three modules. Descent to the maximal-rank intermediate subgroup $K=\SU(5)\times \SU(2)\times \U(1)$ forces the coefficients on $V\oplus W=\mathfrak k^\perp$ to coincide, say with $b$. Write the coefficient on $U$ as $a$. If $a\ne b$, the horizontal g.o. equation on $y+x$, $y\in U$, $x\in V$, would imply
$(a-b)yx\in(\mathfrak{so}(5)+i\mathbb R I)x \quad\hbox{for every }y\in i\operatorname{Sym}_0(\mathbb R^5), \ x\in\mathbb C^5.$
Take $x=e_1$ and $y=i(E_{12}+E_{21})$. The left side is a nonzero imaginary multiple of $e_2$, whereas the right side lies in $\operatorname{span}_{\mathbb R}\{e_2,e_3,e_4,e_5,ie_1\}$. This contradiction proves $a=b$.

The summand $W$ identifies with real $10\times4$ matrices, with left $\SO(5)$ acting in its adjoint representation and right $\SU(2)\times \U(1)$ acting as a subgroup of $\SO(4)$. Choose four orthogonal vectors in $\mathfrak{so}(5)$ whose common centralizer is zero and scale them to distinct lengths. Such vectors exist. Two elements with zero common centralizer can be obtained from a regular Cartan element and a sum of nonzero simple root-plane components. Orthogonalizing their span does not change the common centralizer, and the frame can then be extended. The Gram-matrix argument forces the right infinitesimal factor to vanish, and the common-centralizer condition forces the left factor to vanish. Thus $W$ has a locally free vector and Corollary~\ref{cor315} applies.
\end{proof}

\begin{proposition}\label{prop95}
Every $E_6$-g.o. metric is normal if its isotropy has a simple ideal of class $A_2^{3\prime}$.
\end{proposition}
\begin{proof}
This is $\SU(3)$ acting by its real adjoint representation in $\SO(8)$, and its centralizer in $E_6$ is a two-dimensional torus. Thus $\mathfrak h=\mathfrak{su}(3)\oplus\mathfrak r$ with $\mathfrak r\subset\mathbb R^2$. Use the maximal-rank subgroup $K=\Spin(8)\times T^2$. The three triality modules of $\Spin(8)$ all restrict to the real adjoint $\SU(3)$-module $\mathbb R^8$. Therefore
$\mathfrak m=(\mathbb R^2\ominus\mathfrak r) \oplus (\operatorname{Sym}^3\mathbb C^3)_{\mathbb R} \oplus\bigoplus_{j=1}^3(\mathbb R^8\otimes\mathbb C_{\chi_j})_{\mathbb R}.$
The first two summands are respectively the trivial and the $20$-dimensional irreducible $\SU(3)$-modules. The last three have adjoint $\SU(3)$-type. Their torus weights $\chi_j$ span $(\mathbb R^2)^*$. In particular, every invariant $G$-g.o. cometric preserves $\mathfrak k^\perp$. Descent to $E_6/K$ makes it scalar, say $bI$, there. On the $20$-dimensional summand it is $aI$, and on $\mathbb R^2\ominus\mathfrak r$ it is an arbitrary positive self-adjoint operator $E$.

If $a\ne b$, use the g.o. equation on the sum of an arbitrary $y\in\mathfrak{so}(8)\ominus\mathfrak{su}(3)$ and a real vector $x$ in one of the three triality modules. The real part of the equation says $yx\in\mathfrak{su}(3)x$. The torus acts purely imaginarly and cannot affect this equality. Since $\mathfrak{so}(8)=\mathfrak{su}(3)\oplus (\mathfrak{so}(8)\ominus\mathfrak{su}(3))$, this would give $\mathfrak{so}(8)x=\mathfrak{su}(3)x$ for every real $x$. For a regular element $x$ of the adjoint module the dimensions are $7$ and $6$, a contradiction. Hence $a=b$.

Finally put $D=E-bI$ and take a vector whose components in all three charged triality modules are nonzero and real. For $y\in\mathbb R^2\ominus\mathfrak r$, the imaginary part of the g.o. equation gives $\chi_j(Dy-Z_{\mathfrak r})=0$ for $j=1,2,3$. The spanning property gives $Dy=Z_{\mathfrak r}$, and orthogonality forces $Dy=0$. Thus $C=bI$.
\end{proof}

\begin{proposition}\label{prop96}
Every $E_6$-g.o. metric is normal if its isotropy has a simple ideal of class $B_4^1$ or either class $B_2^3$.
\end{proposition}
\begin{proof}
In each case the centralizer is one-dimensional, so the isotropy is either the simple ideal alone or the full centralizer envelope. The simple case follows from the classification with simple isotropy. For $B_4^1=\mathfrak{so}(9)$ the full envelope has isotropy modules $\mathbb R^9\oplus(\Delta_9\otimes\mathbb C_\chi)_{\mathbb R}$, of real dimensions $9$ and $32$. For $B_2^3$, realize $\mathfrak{so}(5)\subset\mathfrak{so}(10)$ by its adjoint representation. Then $\mathfrak{so}(10)\ominus\mathfrak{so}(5)$ is the irreducible $35$-dimensional real module and the complex half-spin module of $\Spin(10)$ restricts to the irreducible complex $16$-dimensional $\Spin(5)$-module. Consequently the full envelope has two irreducible real isotropy modules of dimensions $35$ and $32$. The nonzero central circle weight on the second module makes it real irreducible. The two-summand classification excludes nonscalar $G$-g.o. metrics in both cases. The two outer-conjugate $B_2^3$ embeddings have the same branching and the same argument applies.
\end{proof}

\begin{proposition}\label{prop97}
Assume the previously proved exclusion for isotropy containing a regular root $A_1$ ideal in $E_6$. Every $E_6$-g.o. metric is normal if its isotropy has a simple ideal of class $A_2^5$, $A_3^2$, or $C_3^1$.
\end{proposition}
\begin{proof}
These three ideals lie in the $\SU(6)$ factor of $\SU(6)\times \SU(2)\subset E_6$, with defining representations $\operatorname{Sym}^2\mathbb C^3$, $\Lambda^2\mathbb C^4$, and $\mathbb C^6$ for $\Sp(3)$, respectively. Their centralizer is the external root $\SU(2)$. A proper connected residual subgroup of this centralizer is a circle or trivial. The trivial case is covered by the simple-isotropy classification, and the full $\SU(2)$ case by the stated root-ideal exclusion. It remains to treat the circle case. We use the two-free-submodule criterion. Two orthogonal invariant submodules each having zero principal infinitesimal stabilizer force a $G$-g.o. metric to be normal. This criterion requires no prior assumption that the cometric preserves the two submodules.

For $\SU(3)$ and $\SU(4)$ the exterior-power identities are
\begin{equation*}
 \Lambda^3(\operatorname{Sym}^2\mathbb C^3)
   =\operatorname{Sym}^3\mathbb C^3
      \oplus(\operatorname{Sym}^3\mathbb C^3)^*,
 \qquad
 \Lambda^3(\Lambda^2\mathbb C^4)
   =\operatorname{Sym}^2\mathbb C^4
      \oplus(\operatorname{Sym}^2\mathbb C^4)^*.
\end{equation*}
The $40$-dimensional horizontal representation splits under the external circle into two real $20$-dimensional modules. Each is the underlying real space of the corresponding displayed complex $10$-dimensional representation with nonzero scalar circle weight. Each is locally free. For ternary cubics, the vector $x_1^3+x_2^3+x_3^3$ has zero infinitesimal stabilizer under $\mathfrak{su}(3)+i\mathbb R I$. The six off-diagonal monomial coefficients force the off-diagonal matrix entries to vanish, and the three diagonal equations together with trace zero kill the remaining entries and the scalar. For complex symmetric $4\times4$ matrices, take a positive diagonal matrix with distinct entries. The equation $AD+DA^t+itD=0$, with $A\in\mathfrak{su}(4)$, first kills all off-diagonal entries of $A$, then its diagonal entries and $t$. Thus the required two free modules exist.

For $\Sp(3)$, use $\Lambda^3\mathbb C^6=\Lambda^3_0\mathbb C^6\oplus\mathbb C^6$. One free submodule is $(\Lambda^3_0\mathbb C^6\otimes\mathbb C_\chi)_{\mathbb R}$, of dimension $28$. It is the restriction of the isotropy representation of the compact symmetric space $F_4/(\Sp(3)\Sp(1))$. That representation has cohomogeneity $4$ and acting-group dimension $24$, hence zero principal infinitesimal stabilizer. Restriction to $\Sp(3)\U(1)$ remains locally free. A disjoint free submodule is
$\operatorname{Herm}_0(\mathbb H^3) \oplus(\mathbb C^6\otimes\mathbb C_\chi)_{\mathbb R} \oplus(\mathfrak{su}(2)\ominus\mathfrak u(1)),$
of dimensions $14+12+2$. Choose a diagonal quaternionic Hermitian matrix with three distinct eigenvalues, the quaternionic column $(1,1,1)^t$, and a nonzero vector in the last plane. The first component reduces the $\Sp(3)$ stabilizer to $\Sp(1)^3$. The last component kills the external circle. The middle component then kills all three $\Sp(1)$ infinitesimal factors. The two-free-submodule criterion proves the assertion.
\end{proof}

The complete rank-at-least-two list is Table~\ref{tab:e6-high-inventory}. Each row names one of the preceding proofs. In particular every nonregular simple ideal of rank at least two, with arbitrary residual isotropy in its centralizer, has only normal metrics.

\paragraph{The rank-one ideals}

\begin{lemma}
\label{lem98}
Let $\mathfrak s\simeq\mathfrak{su}(2)$ be embedded in $\mathfrak e_6$. Suppose its centralizer contains a simple ideal $\mathfrak q$ for which all $E_6$-g.o. metrics with an isotropy ideal $\mathfrak q$ have already been proved normal. Then the same assertion holds with $\mathfrak s$ in place of $\mathfrak q$.
\end{lemma}
\begin{proof}
For $\mathfrak l=\mathfrak s\oplus\mathfrak c(\mathfrak s)$, both $\mathfrak l$ and its derived algebra $\mathfrak d$ contain $\mathfrak q$ as an ideal. The hypothesis gives normality on $E_6/L$ and on $E_6/D$. Lemma~\ref{lem317} gives zero central projection of the principal $L$-stabilizer on $\mathfrak l^\perp$. The Benoist--Kobayashi~\cite{BeKo2021} classification of exceptional semisimple subalgebras with nonabelian generic stabilizer has, for ambient $E_6$, exactly the types $F_4$, $B_4$, and $D_5$. The algebra $\mathfrak d$ is a product containing the ideal $\mathfrak s$ and cannot be any of these simple algebras. Its generic stabilizer on $\mathfrak g/\mathfrak d$ is therefore abelian. Since $\mathfrak g/\mathfrak d=\mathfrak p\oplus\mathfrak z$ and $\mathfrak z$ is $D$-trivial, this is exactly the generic stabilizer needed on $\mathfrak p$. Adding the center preserves this property. The derived algebra of a compact $L$-stabilizer is semisimple and lies in the abelian stabilizer of $D$, hence is zero. The ideal-envelope descent theorem now applies, using zero central projection and the orthogonal-Cartan argument. This proves the claim.
\end{proof}

\begin{proposition}\label{prop99}
Every $E_6$-g.o. metric is normal if its isotropy has a rank-one ideal of Dynkin index $12,21,30$, or $60$.
\end{proposition}
\begin{proof}
The centralizers have dimensions $2,1,1,1$, respectively, and are abelian. Use the dimension convention $V_d=\operatorname{Sym}^{d-1}\mathbb C^2$. The branchings verified in Appendix~\ref{app:branching-checks} give the following pairs of orthogonal real submodules. Opposite complex weights are paired to take real forms.
\begin{equation*}
\begin{array}{c|c|c}
 \text{index}&W_1&W_2\\\hline
 12&(V_{5}\otimes\mathbb C_\alpha)_{\mathbb R}
       \oplus(V_{5}\otimes\mathbb C_\beta)_{\mathbb R}
    &(V_{3}\otimes\mathbb C_\alpha)_{\mathbb R}
       \oplus(V_{3}\otimes\mathbb C_\beta)_{\mathbb R}\\
 21&(V_{8}\otimes\mathbb C_\chi)_{\mathbb R}
    &(V_{6}\otimes\mathbb C_\chi)_{\mathbb R}\\
 30&(V_{8}\otimes\mathbb C_\chi)_{\mathbb R}
    &(V_{6}\otimes\mathbb C_\chi)_{\mathbb R}\\
 60&(V_{11}\otimes\mathbb C_\chi)_{\mathbb R}
    &(V_{5}\otimes\mathbb C_\chi)_{\mathbb R}
\end{array}
\end{equation*}
Here $\chi\ne0$, while $\alpha,\beta$ are linearly independent characters of the two-dimensional centralizer. Each $W_i$ is locally free by Lemma~\ref{lem26}. In the first row one component kills the $\SU(2)$ infinitesimal factor and its phase, and the other kills the remaining torus direction. This remains true after restricting to any residual subtorus. The two-free-submodule criterion proves normality. The index-$12$ characters in the calculator's centralizer-Cartan coordinates can be chosen as $(2,2)$ and $(-2,4)$. The other three rows have the displayed common nonzero character.
\end{proof}

\begin{proposition}\label{prop910}
Every $E_6$-g.o. metric is normal for every closed subgroup of the centralizer envelope of an index-$6$ rank-one subalgebra.
\end{proposition}
\begin{proof}
The branching calculation in Appendix~\ref{app:e6-a1-branching} gives this envelope the algebra $\mathfrak l=A_1^6\oplus A_1^6\oplus\mathbb R$. Under $L$, its orthogonal complement has branching
\begin{equation*}
 \mathfrak l^\perp
 =(\mathbb R^5\otimes\mathbb R^3)
 \oplus(\mathbb R^3\otimes\mathbb R^5)
 \oplus(\mathbb R^3\otimes\mathbb R^3)
 \oplus(V_{4}\otimes V_{2}\otimes\mathbb C_\chi)_{\mathbb R}
 \oplus(V_{2}\otimes V_{4}\otimes\mathbb C_\chi)_{\mathbb R},
 \qquad\chi\ne0.
\end{equation*}
The dimensions are $15+15+9+16+16=71$. We show that each of the last two modules is locally free for $L$. For the first, view a vector as a complex linear map $X:\mathbb C^2\longrightarrow V_{4}$ and choose its image to be
$P=\operatorname{span}_{\mathbb C}\{e_0,e_1+e_3\},$
where $e_0,e_1,e_2,e_3$ is the standard orthonormal weight basis of $V_{4}$, with weights $3,1,-1,-3$. If $A\in\mathfrak{su}(2)$ infinitesimally preserves $P$, the $e_1$ component of $Ae_0$ kills the off-diagonal coefficient, because $e_1\notin P$. Then the different weights of $e_1$ and $e_3$ force the Cartan coefficient to vanish. Thus $A=0$. In a stabilizer equation $\rho(A)X+X(B+itI)=0$, full rank of $X$ now gives $B+itI=0$, and trace zero gives $B=t=0$. Interchanging the two $\SU(2)$ factors proves the same assertion for the last module. These two orthogonal modules lie in $\mathfrak h^\perp$ for every $\mathfrak h\subset\mathfrak l$, and remain locally free after restriction. The two-free-submodule criterion completes the proof.
\end{proof}

Table~\ref{tab:e6-a1-inventory} lists all twenty rank-one classes and the argument for each. Its initial cases follow from the regular-ideal theorem and Table~\ref{tab:e6-high-inventory}. The index-two case is covered by Proposition~\ref{prop92}. Proposition~\ref{prop910} treats index six independently. The propagation at index nine uses the already treated index-three case. The remaining toral cases are Proposition~\ref{prop99}. The complete order is recorded with the table.

\begin{theorem}
\label{thm911}
Let $G$ be compact connected simple of type $E_6$ and let $H$ be any closed proper subgroup. Every $G$-g.o. metric on $G/H$ is normal, except when $\operatorname{Lie}(H)$ is the standard $\mathfrak{so}(10)$ subalgebra. In that case the complete family is the two-parameter family
$C=aI_{\mathbb R}\oplus bI_{\mathfrak p}, \qquad a,b>0, \qquad\mathfrak p=(\mathfrak{so}(10)\oplus\mathbb R)^\perp, \quad\dim\mathfrak p=32.$
The same family applies to disconnected $H$ by Corollary~\ref{cor33}.
\end{theorem}
\begin{proof}
If $\mathfrak h$ is abelian, the abelian-isotropy theorem gives normality. Otherwise choose a simple ideal of $[\mathfrak h,\mathfrak h]$. If it is regular, the regular-simple-ideal exclusion gives precisely the stated possibility. If it is nonregular of rank at least two, Table~\ref{tab:e6-high-inventory} and its row-by-row proofs give normality. If it has rank one, the complete $20$-row Table~\ref{tab:e6-a1-inventory} gives normality. This exhausts all compact isotropy Lie algebras, including arbitrary diagonal couplings in the residual centralizers. The exceptional family is the known simple-isotropy family, also obtained from the symmetric fibration through $\Spin(10)\U(1)$. Corollary~\ref{cor33} completes the proof for arbitrary closed isotropy without changing the parameter family.
\end{proof}

\paragraph{\texorpdfstring{All isotropy ideals of rank at least two in $E_7$}{All isotropy ideals of rank at least two in E7}}
\label{sec:e7-rank-ge-two}

We use the proved exclusion for a regular simple isotropy ideal, including a root $A_1$, and the hereditary normality results for the three maximal subgroups $G_2C_3$, $F_4A_1$, and $G_2A_1$ of $E_7$. Here hereditary means that every closed subgroup of the displayed maximal subgroup has only normal $E_7$-g.o. metrics. The argument below does not infer hereditary normality merely from normality of the maximal-subgroup quotient.

\begin{lemma}
\label{lem912}
Let $\mathfrak s\subset\mathfrak e_7$ be a nonregular simple subalgebra (of any rank). Suppose its centralizer contains a simple ideal $\mathfrak q$ such that every isotropy algebra having $\mathfrak q$ as an ideal has only normal $E_7$-g.o. metrics. Then the same assertion holds with $\mathfrak s$ in place of $\mathfrak q$.
\end{lemma}
\begin{proof}
Set $\mathfrak l=\mathfrak s\oplus\mathfrak c(\mathfrak s)$, $\mathfrak d=[\mathfrak l,\mathfrak l]$, $\mathfrak z=\mathfrak z(\mathfrak l)$, and $\mathfrak p=\mathfrak l^\perp$. Both $\mathfrak l$ and $\mathfrak d$ contain $\mathfrak q$ as an ideal, so the hypothesis gives normality for $E_7/L$ and $E_7/D$. The reverse-center lemma therefore gives $P_{\mathfrak z}(\mathfrak l_x)=0$ at a principal point of $\mathfrak p$.

Benoist--Kobayashi, Theorem~4.1 and Table~2, list precisely $D_6+A_1$, $D_6$, and $E_6$ as the semisimple subalgebras of $E_7$ with nonabelian generic stabilizer on the quotient module. The last two are simple and cannot equal the product $\mathfrak d$. In the first, both simple ideals are regular subalgebras of $E_7$, so $\mathfrak s$ cannot be an ideal of that embedded algebra. Thus $\mathfrak d$ has an abelian generic stabilizer on $\mathfrak g/\mathfrak d=\mathfrak p\oplus\mathfrak z$. The second summand is trivial, so the same holds on $\mathfrak p$. The principal $L$-stabilizer is also abelian. Its derived algebra is semisimple and lies in the abelian $D$-stabilizer, hence vanishes. Together with zero central projection, the ideal-envelope descent theorem proves the result for arbitrary residual isotropy in $\mathfrak c(\mathfrak s)$.
\end{proof}

\begin{proposition}
\label{prop913}
Every $E_7$-g.o. metric is normal if its isotropy has a simple ideal in one of the classes $A_2^9$, $G_2^3$, or $C_4^1$ of Minchenko's Table~13.
\end{proposition}
\begin{proof}
Each has centralizer $\mathbb R$ and lies in the standard $E_6\subset E_6T^1\subset E_7$. The isotropy is consequently either the simple ideal alone or its full circle extension. The simple case is covered by the simple-isotropy classification. For the extension, the isotropy representation is
$(\mathfrak e_6\ominus\mathfrak s) \oplus(V_{27}|_{\mathfrak s}\otimes\mathbb C_\chi)_{\mathbb R}, \qquad\chi\ne0,$
where $V_{27}$ is either complex minimal representation of $E_6$. The necessary restrictions are
\begin{equation*}
\begin{array}{c|c|c|c}
 \mathfrak s&\mathfrak e_{6,\mathbb C}|_{\mathfrak s}
       &V_{27}|_{\mathfrak s}&\text{real isotropy dimensions}\\\hline
 A_2^9&V_{(1,1)}\oplus V_{(4,1)}\oplus V_{(1,4)}
       &V_{(2,2)}&70+54\\
 G_2^3&V_{\omega_2}\oplus V_{\omega_1+\omega_2}
       &V_{2\omega_1}&64+54\\
 C_4^1&V_{2\omega_1}\oplus V_{\omega_4}
       &V_{\omega_2}&42+54.
\end{array}
\end{equation*}
For $G_2$, $\omega_1$ denotes the seven-dimensional fundamental weight. For $C_4$, $V_{\omega_2}=\Lambda^2_0\mathbb C^8$ and $V_{\omega_4}=\Lambda^4_0\mathbb C^8$. The first module in each row is real irreducible. In the $A_2$ row it is the real form of the pair of dual $35$-dimensional complex modules. The second is real irreducible because its complex $27$-dimensional factor is irreducible and carries a nonzero scalar circle character. The two real modules are inequivalent. Chen--Nikonorov's complete two-isotropy-summand classification has no nonnormal case with prescribed simple ambient group $E_7$. This proves normality.

The $A_2$ and $G_2$ restrictions of $V_{27}$ and the $A_2$ adjoint restriction are explicitly recorded in the official Sage branching-rules documentation under ``Miscellaneous'', including the outputs $V_{27}\mapsto A_2(2,2)$ and $V_{27}\mapsto G_2(2,0)$. The $G_2$ and $C_4$ adjoint restrictions are also directly listed in the explicit $E_6$ embedding table (entries $54$ and $109$). The symplectic row follows as well from the standard exterior-power realization of the symmetric pair $E_6/\Sp(4)$.
\end{proof}

\begin{proposition}
\label{prop914}
Every closed subgroup of the envelope $\SO(7)\U(1)\subset \SU(7)\U(1)\subset \SU(8)\subset E_7$ has only normal $E_7$-g.o. metrics. In particular, every isotropy having a simple ideal of class $B_3^2$ has only normal $G$-g.o. metrics.
\end{proposition}
\begin{proof}
Take the circle generator $\operatorname{diag}(iI_7,-7i)$ in $\mathfrak{su}(8)$. The standard $E_7/\SU(8)$ representation is the real form of $\Lambda^4\mathbb C^8$. Under $\SO(7)\U(1)$, the tangent module of the envelope quotient is
$\operatorname{Sym}_0(\mathbb R^7) \oplus(\mathbb C^7\otimes\mathbb C_8)_{\mathbb R} \oplus(\Lambda^3\mathbb C^7\otimes\mathbb C_{-4})_{\mathbb R}.$
Its dimensions are $27+14+70=111$. The sum of the first two modules is locally free. A real symmetric matrix with distinct eigenvalues kills the $\mathfrak{so}(7)$ stabilizer, and a nonzero vector in the charged module then kills the circle.

The last module is locally free as well. Let $\mathcal F$ be the seven triples of a Fano plane on $\{1,\ldots,7\}$, and take
\begin{equation*}
 \phi=\sum_{L\in\mathcal F} c_L e^L\in\Lambda^3\mathbb R^7,
 \qquad \{c_L:L\in\mathcal F\}=\{1,2,4,8,16,32,64\}.
\end{equation*}
Form the equivariant symmetric contraction matrix $M_{ij}(\phi)=\sum_{k,l}\phi_{ikl}\phi_{jkl}$. No pair of vertices belongs to two distinct Fano lines, so $M$ is diagonal. Its diagonal entries are twice the sums of the three $c_L^2$ incident at each vertex. These seven sums are distinct, because they are distinct sums of powers of $4$ with digits $0,1$. Every infinitesimal orthogonal stabilizer of $\phi$ commutes with this simple-spectrum symmetric matrix, and is therefore zero. At the real vector $\phi$, the circle derivative is purely imaginary whereas the $\mathfrak{so}(7)$ derivative is real. Hence the full $\SO(7)\U(1)$ stabilizer is zero. We have obtained two orthogonal invariant locally free modules. They remain locally free on restriction to every subgroup and lie in its tangent module. The two-free-submodule criterion proves the hereditary assertion. Minchenko's centralizer table identifies this $\SO(7)$ as class $B_3^2$ and gives its centralizer exactly as $\U(1)$.
\end{proof}

\begin{theorem}\label{thm915}
If a compact isotropy algebra in $\mathfrak e_7$ has a simple ideal of rank at least two, every $E_7$-g.o. metric is normal. Thus any remaining nonnormal $E_7$ possibility must have all simple isotropy ideals of type $A_1$.
\end{theorem}
\begin{proof}
Regular ideals are covered by Corollary~\ref{cor716}. Table~\ref{tab:e7-high-inventory} contains all twenty-seven nonregular classes, with their rows in Minchenko's exhaustive Table~7 and their full centralizers. Each row invokes a preceding envelope result, one of the four circle calculations, or Lemma~\ref{lem912}. The only nonregular factors used in propagation are $G_2^1$ and $A_1^{3\prime\prime}$. Their full envelopes have already been treated by Theorems~\ref{thm824} and~\ref{thm831}. Thus the argument includes every residual compact subalgebra.
\end{proof}

\paragraph{\texorpdfstring{All rank-one isotropy ideals in $E_7$}{All rank-one isotropy ideals in E7}}
\label{sec:e7-a1-exhaustion}

We next treat the rank-one ideals in $E_7$. We use the previously proved exclusion for regular simple ideals, normality for all isotropy ideals of rank at least two in $E_7$, and hereditary normality for the maximal subgroup $A_1^{24}A_1^{15}$. Superscripts denote Dynkin indices. A prime in a Bala--Carter name distinguishes the embedded class, not a new abstract Lie algebra. For representation formulas below, $V_d=\operatorname{Sym}^{d-1}\mathbb C^2$ has complex dimension $d$.

\begin{proposition}
\label{prop916}
Let $\mathfrak s\subset\mathfrak e_7$ be an $A_1$ in either Bala--Carter class $E_6(a_1)$ or $A_4+A_1$. Every $E_7$-g.o. metric on $E_7/H$ is normal whenever $\mathfrak h\subset\mathfrak s+\mathfrak c(\mathfrak s)$. In particular this holds whenever $\mathfrak s$ is an ideal of $\mathfrak h$, with arbitrary residual torus in its centralizer.
\end{proposition}
\begin{proof}
We exhibit two orthogonal invariant subspaces of the complement of the full envelope, each with a locally free vector for that envelope. The hereditary two-free criterion then applies to every subgroup of it, without assuming that its $G$-g.o. cometric preserves these two subspaces.

For $E_6(a_1)$ the Dynkin index is $84$ and the centralizer is one-dimensional. Use the regular $E_6T\subset E_7$, for which
$\mathfrak e_7=\mathfrak e_6\oplus\mathbb R \oplus (\mathbf {27}\otimes\mathbb C_1)_{\mathbb R}.$
The restriction of $\mathbf {27}$ to this $A_1$ is
\begin{equation}\label{formula36}
 \mathbf {27}\big|_{A_1^{84}}=V_{13}\oplus V_9\oplus V_5.
\end{equation}
Thus $(V_{13}\otimes\mathbb C_1)_{\mathbb R}$ and $(V_9\otimes\mathbb C_1)_{\mathbb R}$ are two required subspaces by Lemma~\ref{lem26}.

For $A_4+A_1$ the Dynkin index is $21$ and the centralizer is $\mathbb R^2$. Realize the $A_1$ diagonally in the principal $A_1\subset \SU(5)$ and the second factor of $\SU(6)\times \SU(2)\subset E_6$. Write $\mathbb C^6=\mathbb C^5_1\oplus\mathbb C_{-5}$. This defines the internal circle. The second circle is the center of $E_6T$ in $E_7$. The standard restrictions are
\begin{equation*}
 \mathfrak e_{6,\mathbb C}
 = (\mathfrak{sl}_6,1)\oplus(1,\mathfrak{sl}_2)
       \oplus(\Lambda^3\mathbb C^6\otimes\mathbb C^2),
 \qquad
 \mathbf {27}=(\Lambda^2\mathbb C^6,1)
       \oplus((\mathbb C^6)^*,\mathbb C^2).
\end{equation*}
Since $\Lambda^2V_5=V_7\oplus V_3$ and $V_5\otimes V_2=V_6\oplus V_4$, the latter gives
\begin{equation}\label{formula37}
 \mathbf {27}
 = V_{7,(2,1)}\oplus V_{3,(2,1)}\oplus V_{5,(-4,1)}
   \oplus V_{6,(-1,1)}\oplus V_{4,(-1,1)}\oplus V_{2,(5,1)}.
\end{equation}
Here the second displayed weight incorporates the external circle when the representation is viewed inside $\mathfrak e_7$. The term $\Lambda^3\mathbb C^5\otimes\mathbb C^2$ in $\mathfrak e_{6,\mathbb C}$ has internal weight $3$, and its restriction is $V_8\oplus V_6\oplus V_4\oplus V_2$. Consequently the complement of $\mathfrak s+\mathbb R^2$ contains orthogonal real submodules
$E=(V_{8,(3,0)})_{\mathbb R}\oplus(V_{7,(2,1)})_{\mathbb R}, \qquad F=(V_{6,(3,0)})_{\mathbb R}\oplus(V_{3,(2,1)})_{\mathbb R}.$
Choose the first component in each outside the eigenvector set of Lemma~\ref{lem26}, and the second component nonzero. The first component kills the $A_1$ and the internal circle. The second then kills the external circle. Each of $E,F$ is therefore locally free for the full envelope.
\end{proof}

The exact root calculation for both restrictions is given in Appendix~\ref{app:e7-a1-branching}.

Tables~\ref{tab:e7-a1-inventory} and~\ref{tab:e7-a1-inventory-2} list all forty-four rank-one classes and the complete propagation order. The two toral classes were treated in Proposition~\ref{prop916}. The mutually centralizing index-fifteen and index-twenty-four classes are covered simultaneously by Theorem~\ref{thm833}, independently of propagation.

\begin{theorem}\label{thm917}
If the isotropy algebra of a $E_7$-g.o. metric contains an ideal isomorphic to $\mathfrak{su}(2)$, the metric is normal. This allows an arbitrary compact residual algebra in the centralizer of that ideal and arbitrary finite components of the isotropy group.
\end{theorem}
\begin{proof}
Every embedded $\mathfrak{su}(2)$ occurs in exactly one of the $44$ rows of Tables~\ref{tab:e7-a1-inventory} and~\ref{tab:e7-a1-inventory-2}. Centralizer propagation handles every row marked $\mathrm P$ for arbitrary residual isotropy. The regular ideal theorem does the same for $\mathrm R$. For a row marked $\mathrm M$ or $\mathrm T$, the entire centralizer envelope is covered hereditarily. In a row marked $\mathrm S$ the zero centralizer forces the isotropy algebra to be precisely that $\mathfrak{su}(2)$, so the simple-isotropy classification applies. Corollary~\ref{cor33} treats disconnected isotropy with the same $G$-g.o. cone.
\end{proof}

\paragraph{\texorpdfstring{Three $\SO(14)$ circle envelopes in $E_8$}{Three \SO(14) circle envelopes in E8}}
\label{sec:e8-so14-circle}

The following criterion gives hereditary normality, rather than normality only for the envelope quotient. All representations and stabilizers in the proof are real unless explicitly complexified.

\begin{lemma}
\label{lem918}
Let $K\subset \SO(14)$ be a compact connected subgroup with $\dim K\le21$. Suppose every nonzero $A\in\mathfrak k$ has rank at least four on $\mathbb R^{14}$. Embed $L=K\cdot \SO(2)$ block-diagonally in $\SO(16)\subset E_8$, passing to finite covers or central quotients as necessary. Then every $E_8$-g.o. metric on $E_8/H$ is normal for every closed subgroup $H\subset L$.
\end{lemma}
\begin{proof}
Write $\mathfrak w=\mathfrak{so}(14)\ominus\mathfrak k$. Under $\mathfrak l=\mathfrak k\oplus\mathbb R$, the complement of $\mathfrak l$ in $\mathfrak e_8$ is
\begin{equation}\label{formula38}
 \mathfrak w\ \oplus\ (\mathbb C^{14}\otimes\mathbb C_2)_{\mathbb R}
 \ \oplus\ (\Delta_{14}^+\otimes\mathbb C_1)_{\mathbb R}.
\end{equation}
The first two terms come from $\mathfrak{so}(16)$. The last is its real half-spin complement in $\mathfrak e_8$. Here $\Delta_{14}^+$ is a \emph{complex} half-spin representation of complex dimension $64$. Its conjugate is $\Delta_{14}^-$. Thus the last summand of \eqref{formula38} has real dimension $128$, and the middle one has real dimension $28$. The weight normalization $2,1$ chooses the double cover of the $\SO(2)$ factor and has no effect on infinitesimal stabilizers.

First we show that $\mathfrak w$ has a locally free vector for $K$. For nonzero $A\in\mathfrak k$, let $r$ be the number of its nonzero rotation planes on $\mathbb R^{14}$, so $2\le r\le7$. Its zero eigenspace has dimension $14-2r$. If the multiplicities of its distinct nonzero absolute rotation frequencies are $m_1,\ldots,m_t$, then
\begin{equation*}
 \dim\mathfrak c_{\mathfrak{so}(14)}(A)
 =\binom{14-2r}{2}+\sum_jm_j^2
 \le \binom{14-2r}{2}+r^2\le49.
\end{equation*}
For $r=2,\ldots,7$, the upper bounds are $49$, $37$, $31$, $31$, $37$ and $49$, respectively. Hence
$\operatorname{rank}(\operatorname{ad}A|_{\mathfrak{so}(14)})\ge42.$
Writing $d=\dim\mathfrak k$ and using the invariant orthogonal splitting $\mathfrak{so}(14)=\mathfrak k\oplus\mathfrak w$ gives
$\operatorname{rank}(A|_{\mathfrak w})\ge42-d\ge d>d-1.$
The generic-tuple incidence criterion, with one vector, proves that $\mathfrak w$ has a vector with zero infinitesimal $K$-stabilizer. Consequently
$E=\mathfrak w\oplus(\mathbb C^{14}\otimes\mathbb C_2)_{\mathbb R}$
has a locally free vector for $L$. The first component kills $\mathfrak k$, and any nonzero second component kills the circle.

We next prove that $F=(\Delta_{14}^+\otimes\mathbb C_1)_{\mathbb R}$ is locally free for $L$. Conjugate $A\in\mathfrak{so}(14)$ to rotation coordinates $a_1,\ldots,a_7$. Its $64$ complex half-spin weights are
$\frac12(\epsilon_1a_1+\cdots+\epsilon_7a_7),\qquad \epsilon_j\in\{1,-1\},\quad \prod_j\epsilon_j=\varepsilon,$
for one fixed parity $\varepsilon$. If some $|a_i|\ne|a_j|$, flipping the $i$th and $j$th signs partitions these weights into $32$ pairs of distinct values. Thus every eigenvalue of $A$ has complex multiplicity at most $32$. Otherwise all $|a_i|$ are equal. If $A\ne0$, the multiplicities are the binomial coefficients $\binom7j$ in one parity class, so each is at most $35$. If $A=0$ and the circle component $t\ne0$, its scalar action has zero kernel. In all cases a nonzero element $(A,t)\in\mathfrak l$ therefore has complex kernel dimension at most $35$ on the charged half-spin module, and real rank at least
$2(64-35)=58>\dim\mathfrak l-1.$
The same incidence criterion makes $F$ locally free. The orthogonal $L$-invariant subspaces $E,F$ lie in $\mathfrak l^\perp$ and are both locally free. Their restrictions to every subalgebra of $\mathfrak l$ remain so. The hereditary two-free criterion now proves the assertion, without any assumption that a $G$-g.o. cometric preserves $E$ or $F$.
\end{proof}

\begin{theorem}\label{thm919}
Every $E_8$-g.o. metric is normal whenever its isotropy contains a simple ideal in one of the embedded classes
$B_2^7,\qquad G_2^4,\qquad C_3^2.$
In fact every closed subgroup of each of the full centralizer envelopes $B_2^7T$, $G_2^4T$, and $C_3^2T$ has only normal $E_8$-g.o. metrics.
\end{theorem}
\begin{proof}
Minchenko's exhaustive Table~8 places these classes in $D_7$. Table~14 gives their centralizer exactly as a one-dimensional torus. On the defining real $14$-dimensional module of $\SO(14)$, the three restrictions are, respectively,
$\operatorname{Sym}_0^2\mathbb R^5, \qquad \mathfrak g_2, \qquad (\Lambda^2_0\mathbb C^6)_{\mathrm{real}}.$
The last notation means the invariant real form of the complex $14$-dimensional $\Sp(3)$ representation, not its underlying $28$-dimensional real representation. These are the corresponding irreducible orthogonal embeddings into $\SO(14)$. Their group dimensions are $10,14,21$.

We verify the rank hypothesis explicitly. For $B_2$, the nonzero weights on the complexification of $\operatorname{Sym}_0^2\mathbb R^5$ are
$\pm2a,\quad\pm2b,\quad\pm a\pm b,\quad\pm a,\quad\pm b,$
with the zero weight having multiplicity two. If exactly one of $a,b$ vanishes, eight weights are nonzero. If both are nonzero, at least ten are nonzero. Thus the minimum nonzero operator rank is eight. For the adjoint representation of $G_2$, a nonzero Cartan element can vanish on at most one pair of roots. Including the two Cartan zero weights its kernel has dimension at most four. Thus its rank is at least ten. For $C_3$, the weights are
$\pm a\pm b,\quad\pm a\pm c,\quad\pm b\pm c,$
again with a zero weight of multiplicity two. If all $a,b,c$ are nonzero, each pair contributes at least two nonzero weights, so the rank is at least six. If just one or two coordinates are nonzero, the rank is at least eight. In particular all three representations satisfy the rank-four hypothesis of Lemma~\ref{lem918}, which proves the result.
\end{proof}

The embeddings and centralizers are given in \cite[Table~8, rows~29, 30, 32, and Table~14, pp.~257--258]{Mi2006}. The proof uses the restrictions in \eqref{formula38}, the incidence criterion and the hereditary two-free criterion. In particular, it does not require a classification of $G$-g.o. metrics on the three quotients by their full centralizer subgroups.

\paragraph{\texorpdfstring{Nonregular simple ideals of rank at least two in $E_8$}{Nonregular simple ideals of rank at least two in E8}}
\label{sec:e8-rank-two-reduction}

We next treat every isotropy containing a simple ideal of rank at least two in $E_8$. The rank-one ideals are treated after Theorem~\ref{thm927}. We keep the superscripts and primes of Minchenko's conjugacy-class notation.  A simple embedding is called \emph{closed} below when every $E_8$-g.o. metric whose isotropy contains that embedding as an ideal is normal.

\begin{lemma}
\label{lem920}
Let $\mathfrak s\subset\mathfrak e_8$ be a nonregular simple subalgebra.  If $\mathfrak c_{\mathfrak e_8}(\mathfrak s)$ has a simple ideal $\mathfrak q$ which is closed, then $\mathfrak s$ is closed.
\end{lemma}

\begin{proof}
Set $\mathfrak l=\mathfrak s\oplus\mathfrak c(\mathfrak s)$, $\mathfrak d=[\mathfrak l,\mathfrak l]$, $\mathfrak z=\mathfrak z(\mathfrak l)$ and $\mathfrak p=\mathfrak l^\perp$. Both $\mathfrak l$ and $\mathfrak d$ contain $\mathfrak q$ as an ideal, so normality is already known for $E_8/L$ and $E_8/D$. By Lemma~\ref{lem317}, the principal $L$-stabilizer on $\mathfrak p$ has zero projection to $\mathfrak z$.

By Benoist--Kobayashi, the only semisimple subalgebras of $E_8$ with nonabelian generic stabilizer on the adjoint quotient are the standard $E_7$ and $E_7+A_1$. The algebra $\mathfrak d$ contains the two distinct ideals $\mathfrak s$ and $\mathfrak q$, so the first is impossible. In the second, each simple ideal is regular in $E_8$, contrary to the assumption on $\mathfrak s$. Thus the generic $D$-stabilizer is abelian. Since $\mathfrak e_8/\mathfrak d=\mathfrak p\oplus\mathfrak z$ with a trivial $D$-summand, the same is true on $\mathfrak p$. Adding the central torus preserves abelianness. The derived algebra of a compact $L$-stabilizer is semisimple and lies in the abelian $D$-stabilizer, hence is zero. The ideal-envelope descent theorem applies.
\end{proof}

\begin{proposition}
\label{prop921}
For the standard embedding
$B_r+B_{7-r}\subset D_8\subset E_8,\qquad r=1,2,3,$
every $E_8$-g.o. metric whose isotropy contains either simple block as an ideal is normal. In particular, the classes $B_2^1,B_3^1,B_4^1,B_5^1,B_6^1$ are closed, as is the $A_1$ class of index two paired with $B_6^1$.
\end{proposition}

\begin{proof}
The vector module of $D_8$ splits as $\mathbb R^{2r+1}\oplus\mathbb R^{15-2r}$. Each block is the full centralizer of the other in $E_8$. The half-spin module has no fixed vectors under either block. The quotient of this envelope has two real irreducible summands,
$\mathbb R^{2r+1}\otimes\mathbb R^{15-2r}, \qquad \Delta_{16}^{+}\big|_{B_r+B_{7-r}}.$
The second is real irreducible of dimension $128$. Indeed, its complexification is the tensor product of the two odd-dimensional spin modules. Their real or quaternionic types match, so the product has a real form. The two summands have different dimensions. The two-isotropy-summand classification therefore gives normality on the base quotient. Its nonnormal list has no ambient $E_8$. The envelope is semisimple and is neither $E_7$ nor $E_7+A_1$, so its generic stabilizer is abelian by Benoist--Kobayashi. Ideal-envelope descent proves the assertion.
\end{proof}

\begin{proposition}
\label{prop922}
Every closed subgroup of the standard tensor subgroup $\Sp(2)\cdot \Sp(2)\subset \SO(16)\subset E_8$ has only normal $E_8$-g.o. metrics.  Consequently $B_2^{2\prime\prime}$ is closed.
\end{proposition}

\begin{proof}
The complexification of the real vector module is $\mathbb C^4\otimes\mathbb C^4$, where both factors are quaternionic.  Exterior squaring gives
$\mathfrak{so}(16)/(\mathfrak{sp}(2)+\mathfrak{sp}(2)) =(V_5\otimes\mathfrak{sp}(2)) \oplus(\mathfrak{sp}(2)\otimes V_5).$
These are orthogonal real modules of dimension $50$. Each has a free vector.  For example, choose an injective map $V_5\longrightarrow\mathfrak{sp}(2)$ with distinct singular values and image containing two elements with zero joint centralizer.  The stabilizer equation first kills the vector $\mathfrak{so}(5)$ action and then the adjoint action on the image.  Two such elements exist in every compact semisimple algebra. Take a regular Cartan element and a vector with nonzero components in all simple-root planes. The other summand is treated by interchanging the factors. The two-free-module criterion proves hereditary normality. The additional half-spin module is irrelevant.

For clarity, the tensor subgroup has two realizations in $E_8$ according to the choice of half-spin module of $D_8$. The restrictions, in complex notation, are
$(V_{14},1)\oplus(1,V_{14})\oplus(V_5,V_{10})\oplus(V_{10},V_5), \qquad (V_4,V_{16})\oplus(V_{16},V_4).$
Both have dimension $128$ and have real forms. They follow by restricting the half-spin weights to the vector weights $\pm e_i\pm f_j$ ($i,j=1,2$). For the second choice neither simple factor has half-spin fixed vectors, so its full $E_8$-centralizer is the other $\Sp(2)$. The embedding has Dynkin index two, and Minchenko's centralizer table identifies it uniquely as $B_2^{2\prime\prime}$. Thus every isotropy containing this class as an ideal lies in the tensor envelope just treated.
\end{proof}

\begin{proposition}
\label{prop923}
Every closed subgroup of either standard $\Sp(4)\cdot \Sp(1)\subset \SO(16)\subset E_8$ has only normal $E_8$-g.o. metrics.  In particular $C_4^{1\prime\prime}$ and the $A_1^{4\prime}$ in its centralizer are closed.
\end{proposition}

\begin{proof}
Let $U=\mathbb C^8$ be the standard symplectic module. The vector module of $\SO(16)$ has complexification $U\otimes\mathbb C^2$, and
$\mathfrak{so}(16)/(\mathfrak{sp}(4)+\mathfrak{sp}(1)) =W_{27}\otimes\mathbb R^3, \qquad W_{27}=\operatorname{Herm}_0(4,\mathbb H).$
This $81$-dimensional real module has a free vector. Choose a real diagonal matrix $D\in W_{27}$ with distinct entries, a Hermitian matrix $X\in W_{27}$ with all off-diagonal entries equal to $1$, and a Hermitian matrix $Y\in W_{27}$ whose $(1,2)$ and $(1,3)$ entries are $i$ and $j$. A skew quaternionic matrix commuting with $D$ is diagonal. Commuting with $X$ forces its entries to be a common imaginary quaternion, and commuting with $Y$ forces that quaternion to commute with $i$ and $j$, hence to vanish. Orthogonalize the three matrices and give them distinct lengths as columns of a tensor.  The distinct-singular-value equation kills the external $\mathfrak{so}(3)$ and then the internal $\mathfrak{sp}(4)$.

The two half-spin restrictions have real forms with respective complexifications
\begin{align*}
 &(\Lambda^4_0U\otimes1)
   \oplus(\Lambda^2_0U\otimes\operatorname{Sym}^2\mathbb C^2)
   \oplus(1\otimes\operatorname{Sym}^4\mathbb C^2),\\
 &(\Lambda^3_0U\otimes\mathbb C^2)
   \oplus(U\otimes\operatorname{Sym}^3\mathbb C^2).
\end{align*}
Their real dimensions are respectively $42+81+5$ and $96+32$. These are the usual symplectic half-spin decompositions, and can also be verified directly using the vector weights $\pm e_i\pm b$. In the first case there is a second copy of the free $81$-dimensional module, so the result follows immediately.

In the second case the $96$-dimensional real summand is free. The weights of $\Lambda^3_0U$, of complex dimension $48$, are
$\pm e_i\quad\text{with multiplicity two},\qquad \pm e_i\pm e_j\pm e_k\quad(i<j<k).$
An infinitesimal element is conjugate to $(a_1,a_2,a_3,a_4;b)$. If $b\ne0$, at most one of the two weights obtained by adding $\pm b$ to each of the $48$ internal weights can vanish. Thus its rank on the real $96$-module is at least $48$. If $b=0$, let $r$ be the number of nonzero $a_i$. For $r=1,2,3,4$, respectively, the displayed internal weights give rank at least $28,24,28,32$. For $r=2$, for instance, the two triples containing just one nonzero coordinate contribute $16$, and the single weights contribute $8$. For $r=3$, every triple contributes at least four nonzero weights, and the single weights contribute $12$. For $r=4$ these contributions are at least $16+16$. Tensoring with $\mathbb C^2$ doubles these bounds. Every nonzero infinitesimal operator therefore has real rank at least $48$, larger than $\dim(\mathfrak{sp}(4)+\mathfrak{sp}(1))=39$. The incidence-dimension criterion gives a free vector. The two-free-module criterion applies to the $81$- and $96$-dimensional summands.

In this second realization the half-spin module has no $\Sp(4)$-fixed vector.  Hence its full centralizer remains $\Sp(1)$, identifying the class as $C_4^{1\prime\prime}$, with centralizer $A_1^{4\prime}$, by Minchenko's table. For the other factor, the full centralizer is exactly $\Sp(4)$. The displayed quotient has no external $\Sp(1)$-fixed vectors. Thus the hereditary result closes both ideal classes.
\end{proof}

\begin{proposition}
\label{prop924}
Every closed subgroup of the full envelopes $B_3^{2\prime\prime}T$ and $D_4^{2\prime\prime\prime}T$ has only normal $E_8$-g.o. metrics.
\end{proposition}
\begin{proof}
Take $S=\Spin(7)$ acting on its real spin module $V=\mathbb R^8$, or $S=\SO(8)$ on its vector module. Embed $S\cdot \SO(2)\subset \U(8)\subset \SO(16)$ on $\mathbb R^{16}=V\otimes\mathbb R^2$. The orthogonal algebra quotient contains $\operatorname{Sym}_0^2V$ and $\mathfrak s\otimes\operatorname{Sym}_0^2\mathbb R^2$. Their circle weights are zero and two. Choose the realization in $E_8$ whose half-spin module is the odd exterior powers of $\mathbb C^8$. It contains a charged copy of $V\otimes\mathbb R^2$ of weight three, up to sign.  The two orthogonal modules
$E=\operatorname{Sym}_0^2V\oplus(V\otimes\mathbb R^2_3), \qquad F=\mathfrak s\otimes\mathbb R^2_2$
are free.  A symmetric matrix with simple spectrum kills all of $\mathfrak{so}(V)$, and a nonzero charged vector then kills the circle in $E$.  For $F$, take two adjoint vectors with zero common centralizer, orthogonalized and given distinct lengths. Their Gram matrix kills the circle and then they kill $\mathfrak s$. The two-free criterion applies. The odd half-spin module has no $S$-fixed vector. For $\SO(8)$ use $V$ and $\Lambda^3V$. For $\Spin(7)$ use $\Lambda^3\Delta_8=\Delta_8\oplus V_{48}$. Thus its full centralizer remains the displayed circle. The repeated realization has index two and identifies precisely the two stated classes in Table~14.
\end{proof}

\begin{proposition}
\label{prop925}
Every closed subgroup of the full envelope $A_2^4T_2\subset E_8$ has only normal $E_8$-g.o. metrics.
\end{proposition}
\begin{proof}
Use the real $\SO(16)$ vector decomposition
$\mathbb R^{16}=\mathfrak{su}(3)\oplus (\mathbb C^3)_{\mathbb R}\oplus\mathbb R^2$
with diagonal $\SU(3)$ on the first two terms.  The two circles act on the second and third terms.  Its index is $3+1=4$. The quotient $\mathfrak{so}(8)/\mathfrak{su}(3)$ is the underlying real $20$-module of $\operatorname{Sym}^3\mathbb C^3$. It has a free vector. In the infinitesimal stabilizer of $x_1^3+x_2^3+x_3^3$, the coefficients of $x_i^2x_j$ kill the off-diagonal entries, and those of $x_i^3$ kill the diagonal. The cross term $(\mathbb C^3)_{\mathbb R}\otimes\mathbb R^2$ splits into standard modules of circle weights $(1,1)$ and $(1,-1)$.  Its sum with the uncharged $20$-module is free. The cubic kills $\mathfrak{su}(3)$ and the independent charges kill both circles. A second orthogonal free module is
$(\mathfrak{su}(3)\otimes\mathbb R^2) \oplus(\mathfrak{su}(3)\otimes(\mathbb C^3)_{\mathbb R}).$
In its first term two adjoint vectors with zero common centralizer and distinct lengths kill $\mathfrak{su}(3)$ and the external circle.  A nonzero vector of the second kills the internal circle.  The two-free criterion applies. No additional centralizer comes from the half-spin module. Under $\Spin(8)\Spin(6)\Spin(2)$ its summands are products of a triality $8$-module, a half-spin $4$-module and a character. All three triality modules restrict to the adjoint $8$ of this $\SU(3)$, whereas the $4$ restricts to $3+1$ or its dual. None of their tensor products has a fixed vector. The full centralizer is therefore $T_2$, giving exactly $A_2^4$.
\end{proof}

\begin{proposition}
\label{prop926}
Every closed subgroup of the full envelope $B_2^{3\prime\prime}A_1^{10\prime\prime}\subset E_8$ has only normal $E_8$-g.o. metrics.  Both ideal classes are closed.
\end{proposition}
\begin{proof}
Inside the regular $\SU(5)\SU(5)$ embed $\Sp(2)$ diagonally, by its vector $5$ in the first factor and $4+1$ in the second. Its index is $2+1=3$.  The commuting circle generator $H$ has weights $1$ on $4$ and $-4$ on $1$. The adjoint branching
$(24,1)+(1,24)+(5,10)+(5^*,10^*)+(10,5^*)+(10^*,5)$
gives the following restrictions to $\Sp(2)$ and $H$.
\begin{equation*}
\begin{array}{c|l}
 H\text{-weight}&\Sp(2)\text{-modules}\\\hline
 0&2V_{10}+V_{14}+V_5+1\\
 \pm1&V_4+V_{16}+V_{20}\\
 \pm2&1+V_{10}+V_{14}+V_5\\
 \pm3&V_4+V_{16}\\
 \pm4&V_{10}\\
 \pm5&V_4.
\end{array}
\end{equation*}
Each $\pm$ entry occurs at both weights.  Here one uses $5\otimes5=1+10+14$, $5\otimes4=4+16$ and $10\otimes4=4+16+20$. The three invariant dimensions at weights $0,\pm2$ form the full centralizer $\mathfrak{su}(2)$. Regrouping the weights into $\SU(2)$ strings gives real forms of
\begin{align*}
 \mathfrak e_8={}&(10,1)+(1,3)+(4,6)+(16,4)\\
                &+(20,2)+(10,5)+(14,3)+(5,3).
\end{align*}
The dimensions sum to $248$. The sum of squared $H$-weights is $1200$. Dividing by $2h^\vee(E_8)=60$ and the squared length two of an $A_1$ coroot gives Dynkin index ten. The realization inside $A_4+A_3$ identifies the index-three factor as $B_2^{3\prime\prime}$ and its centralizer as $A_1^{10\prime\prime}$.

The real summand $(14,3)$ is free. Realize $V_{14}=\operatorname{Sym}_0^2\mathbb R^5$ and take a rank-three tensor of distinct singular values whose image contains a symmetric matrix with simple spectrum. It first kills the external $\mathfrak{so}(3)$ and then $\mathfrak{so}(5)$. The real summand $(10,5)$ is free by taking a rank-five tensor of distinct singular values whose image in $\mathfrak{so}(5)$ contains two vectors with zero common centralizer. The two-free criterion proves the assertion. Neither factor has fixed vectors in the displayed quotient, so the pair is the full ideal-centralizer envelope for either factor.  This computation also corrects the inconsistent group entry in Minchenko's table noted below.
\end{proof}

Table~\ref{tab:e8-high-inventory} lists all forty-four nonregular classes, obtained from Minchenko's exhaustive Table~8. The source cross-check and the two corrected centralizer entries are recorded immediately before the table. The orthogonal, symplectic and circle-envelope proofs above are independent of rank-one propagation. In particular they supply the $A_1^{4\prime}$ seed used by the $A_3^{2\prime\prime}$ row, as well as the later $A_1^{10\prime\prime}$ seed. The hereditary $A_2A_1$ theorem independently supplies the index-sixteen seed.

\begin{theorem}
\label{thm927}
Every $E_8$-g.o. metric is normal whenever the isotropy algebra contains a simple ideal of rank at least two.
\end{theorem}

\begin{proof}
Regular simple ideals were already excluded. The source cross-check and Table~\ref{tab:e8-high-inventory} cover every other rank-at-least-two class. The three $\SO(14)$ circle rows are proved independently in Theorem~\ref{thm919}. Any isotropy containing a specified ideal lies in its full ideal-centralizer envelope.  Thus the hereditary or propagation argument for an eliminated class applies regardless of the remaining isotropy, including central tori and diagonal identifications inside the centralizer.
\end{proof}

It remains to treat rank-one ideals in $E_8$.

\paragraph{\texorpdfstring{All rank-one isotropy ideals in $E_8$}{All rank-one isotropy ideals in E8}}
\label{sec:e8-a1-exhaustion}

We use the independently proved exclusion for every simple isotropy ideal of rank at least two in $E_8$, the regular-ideal theorem, and the hereditary two-free criterion. In particular, the $B_2^{3''}A_1^{10''}$ envelope is closed by its explicit $(\mathbf {14},V_3)$ and $(\mathbf {10},V_5)$ free modules. Proposition~\ref{prop926} treats both factors simultaneously. As in the preceding rank-one section, $V_d=\operatorname{Sym}^{d-1}\mathbb C^2$ denotes the complex $\SU(2)$-module of dimension $d$. A subscript is a circle weight.

We apply Lemma~\ref{lem920} from the rank-at-least-two section also to rank-one nonregular ideals. Its statement and proof impose no rank restriction. Its only nonabelian-generic-stabilizer exceptions are the regular $E_7A_1$ and the simple $E_7$, so neither creates an additional case in this propagation.

\begin{proposition}
\label{prop928}
Let $\mathfrak s\subset\mathfrak e_8$ be in one of the four Bala--Carter classes
$D_7(a_1)$, $D_7(a_2)$, $E_6(a_1)+A_1$ or $D_5+A_2$.
Every closed subgroup of $S\,C_{E_8}(S)^0$ has only normal $E_8$-g.o. metrics. Consequently every isotropy having $\mathfrak s$ as an ideal has only normal metrics.
\end{proposition}
\begin{proof}
In each case the centralizer is a single circle. We find two orthogonal invariant submodules in the complement of the full envelope, each locally free for $\mathfrak s\oplus\mathbb R$. The hereditary two-free criterion then proves the assertion for all of its subgroups. We use the elementary fact that $(V_d\otimes\mathbb C_q)_{\mathbb R}$ is locally free for $\mathfrak{su}(2)\oplus\mathbb R$ whenever $d\ge3$ and $q\ne0$. A generic vector is not an eigenvector of any nonzero element of $\mathfrak{su}(2)$, since the union of their complex eigenlines has real dimension at most four, smaller than $2d$.

For the first two cases use the regular $D_7T\subset E_8$. In a normalization of the circle character the complex adjoint restriction is
$\mathfrak e_{8,\mathbb C} =\mathfrak{so}_{14}\oplus\mathbb C \oplus(\Delta^+_{64})_1\oplus(\Delta^-_{64})_{-1} \oplus(\mathbb C^{14})_2\oplus(\mathbb C^{14})_{-2}.$
The $D_7$ nilpotent partitions of $D_7(a_1)$ and $D_7(a_2)$ are $(11,3)$ and $(9,5)$, respectively. Their vector restrictions are therefore
$\mathbb C^{14}=V_{11}\oplus V_3, \qquad \mathbb C^{14}=V_9\oplus V_5.$
The two charged real summands of the vector term provide the required locally free modules. The embedding indices are $(11(11^2-1)+3(3^2-1))/12=112$ and $(9(9^2-1)+5(5^2-1))/12=70$.

For $E_6(a_1)+A_1$ use the regular $E_6A_2\subset E_8$.
$\mathfrak e_{8,\mathbb C} = (\mathfrak e_6,1)\oplus(1,\mathfrak{sl}_3) \oplus(\mathbf {27},\mathbb C^3) \oplus(\mathbf {27}^*,(\mathbb C^3)^*).$
The $A_1$ is diagonal in the subregular $E_6(a_1)$ factor and a root $A_1\subset A_2$. The latter has $\mathbb C^3=(V_2)_1\oplus(V_1)_{-2}$ for its centralizing circle. The already verified restriction $\mathbf {27}|_{E_6(a_1)}=V_{13}\oplus V_9\oplus V_5$ then gives charged $V_{14}\oplus V_{12}$ from $V_{13}\otimes V_2$. Their underlying real modules are the required pair. The index is $84+1=85$.

For $D_5+A_2$ use the regular $D_5A_3\subset D_8\subset E_8$, with
$\mathfrak e_{8,\mathbb C} = (\mathbf {45},1)\oplus(1,\mathbf {15}) \oplus(\mathbf {10},\mathbf6) \oplus(\mathbf {16},\mathbf4) \oplus(\mathbf {16}^*,\mathbf4^*).$
The $A_1$ is diagonal in a principal $A_1\subset D_5$ and a principal $A_1\subset A_2\subset A_3$. Here $\mathbf4=(V_3)_1\oplus(V_1)_{-3}$ and $\mathbf {16}|_{A_1}=V_{11}\oplus V_5$. The latter half-spin restriction follows directly by evaluating the spin weights at the principal element, whose vector-positive coordinates are $8,6,4,2,0$. The half-spin character is $(x^{10}+x^8+\cdots+x^{-10})+(x^4+x^2+1+x^{-2}+x^{-4})$. Thus the term $(\mathbf {16},(V_1)_{-3})$ supplies charged $V_{11}$ and $V_5$, both locally free. The index is $60+4=64$.

All four circles displayed above centralize the indicated $A_1$, and the full centralizer has dimension one by the centralizer inventory. Hence no omitted circle acts trivially on the two chosen modules.
\end{proof}

Tables~\ref{tab:e8-a1-inventory}, \ref{tab:e8-a1-inventory-2} and~\ref{tab:e8-a1-inventory-3} list all sixty-nine rank-one classes and give the propagation order used in the proof. The eleven zero-centralizer rows reduce to simple isotropy. Proposition~\ref{prop928} treats the four purely toral rows. The other rows use previously proved centralizer ideals or the independent $B_2^{3\prime\prime}A_1^{10\prime\prime}$ envelope.

\begin{theorem}
\label{thm929}
Every $E_8$-g.o. metric on $E_8/H$ is normal if $\mathfrak h$ contains a simple ideal of rank one.
\end{theorem}
\begin{proof}
Choose such an ideal $\mathfrak s$. Then $\mathfrak h=\mathfrak s\oplus\mathfrak r$ with $\mathfrak r\subset\mathfrak c(\mathfrak s)$. The complete inventory in Tables~\ref{tab:e8-a1-inventory}--\ref{tab:e8-a1-inventory-3} and its propagation order cover every embedded class of $\mathfrak s$ and every residual $\mathfrak r$. The toral and hereditary cases also allow arbitrary subgroups of their full envelopes.
\end{proof}

\begin{corollary}
\label{cor930}
For every closed subgroup $H$ of a compact connected group with Lie algebra $\mathfrak e_8$, every $E_8$-g.o. metric on $E_8/H$ is normal. Equivalently, if $H\ne E_8$, its cometric relative to minus the Killing form is a positive scalar on $\mathfrak h^\perp$.
\end{corollary}
\begin{proof}
If $\mathfrak h$ is nonabelian, choose a simple ideal and apply either the rank-at-least-two theorem or the preceding theorem. Abelian $\mathfrak h$, including zero, is covered by the abelian-isotropy theorem. Corollary~\ref{cor33} extends the conclusion to every closed $H$.
\end{proof}

\begin{theorem}\label{thm931}
Let $G$ be a compact connected group with exceptional simple Lie algebra, and let $H\subsetneq G$ be closed. Every $G$-g.o. Riemannian metric on $G/H$ is normal, with precisely the following exception at Lie algebra level.
$(\mathfrak g,\mathfrak h)=(\mathfrak e_6,\mathfrak{so}(10)) \quad\text{in the standard embedding}.$
For this pair, with $Q=-B_{\mathfrak g}$, the complete family of $G$-g.o. cometrics is
$\mathfrak m=\mathbb R\oplus(\Delta^+_{10})_{\mathbb R},\qquad C=aI_{\mathbb R}\oplus bI_{(\Delta^+_{10})_{\mathbb R}}, \qquad a,b>0.$
The half-spin summand has real dimension $32$. The metric is normal if and only if $a=b$. These conclusions include disconnected $H$ and finite central quotients of the acting group.
\end{theorem}

\begin{proof}[Proof of Theorem~\ref{thm931}]
If $\mathfrak h$ is abelian, the abelian-isotropy theorem applies. Otherwise compact reductivity provides a nonzero simple ideal $\mathfrak s\triangleleft\mathfrak h$. If $\mathfrak s$ is regular in $\mathfrak g$, Corollary~\ref{cor716} gives normality or the stated $E_6/D_5$ exception.

For $G_2$ and $F_4$, arbitrary isotropy is already covered by Theorems~\ref{thm72} and~\ref{thm91}. Type $E_6$ is Theorem~\ref{thm911}. For $E_7$, if the chosen nonregular simple ideal has rank at least two, apply Theorem~\ref{thm915}. If it has rank one, apply Theorem~\ref{thm917}. For $E_8$, use Theorem~\ref{thm927} or Theorem~\ref{thm929}, respectively. Each of these results treats arbitrary residual compact isotropy inside the centralizer of the specified ideal. Thus no assumption of a product embedding or of an uncoupled residual subgroup is made.

Every positive pair in the exceptional $E_6/D_5$ family is g.o. by the simple-isotropy classification. Corollary~\ref{cor33} preserves the entire family for disconnected isotropy. Normality is exactly equality of the two coefficients.
\end{proof}

\section{The classification and the dimension of the cometric space}
\label{sec:unified-completion}

In this section, we combine the results for the simple Lie algebras and prove the parameter bound in the Introduction.

\begin{proof}[Proof of Theorem~\ref{thm11}]
For $\mathfrak g=\mathfrak{su}(N)$, Theorem~\ref{thm411} gives entries 1--3 of Table~\ref{tab:unitary-symplectic-list}, including the formula involving the projected central vector $w$. For $\mathfrak g=\mathfrak{sp}(n)$, Theorem~\ref{thm57} gives entries 4 and 5. The orthogonal classification in Theorem~\ref{thm615} gives entries 6--11 of Table~\ref{tab:orthogonal-exceptional-list}. Theorem~\ref{thm931} gives entry 12 and normality for all other exceptional pairs.

These statements include the necessity of the parameter conditions and the sufficiency for every positive parameter choice. Their proofs retain arbitrary equivariant cometrics until a stated argument proves the required invariant decomposition. Thus equivalent isotropy modules and coupled embeddings are included. The rank-one Lie algebras are covered by the unitary result, and the identifications $\mathfrak{so}(5)\cong\mathfrak{sp}(2)$ and $\mathfrak{so}(6)\cong\mathfrak{su}(4)$ give consistent parameter formulas.

By Corollary~\ref{cor33}, every $G$-g.o. cometric for $H^0$ descends to $H$. Passing to finite central covers or effective central quotients does not change the Lie-algebra equation. This proves the theorem for the stated global forms.
\end{proof}

To express the parameter bound intrinsically, let $\mathcal V(\mathfrak g,\mathfrak h)$ be the vector space of $Q$-self-adjoint, isotropy-equivariant operators $D$ on $\mathfrak m$ satisfying $[DX,X]\in[\mathfrak h,X]$ for every $X\in\mathfrak m$. It is a vector space because $[\mathfrak h,X]$ is linear for each fixed $X$. Its positive cone is precisely the cone of $G$-g.o. cometrics. Moreover, every $D\in\mathcal V(\mathfrak g,\mathfrak h)$ is a difference of positive operators in this cone. Indeed, $D+cI$ is positive for sufficiently large $c>0$. Hence $\mathcal V(\mathfrak g,\mathfrak h)$ is the linear span appearing in Corollary~\ref{cor12}.

\begin{proof}[Proof of Corollary~\ref{cor12}]
For every pair outside the tables, the positive cone consists of scalar operators, so its span has dimension one. In entry 1, the dimension is one when $w=0$ and two when $w\ne0$. Entry 2 has three independent scalar parameters. Entry 11 has parameter space $\operatorname{Sym}_2(\mathbb R)$, also of dimension three. Every remaining entry has two independent linear parameters. In particular, the central coefficient in \eqref{formula6} is a linear function of the other two coefficients. This proves the bound and its equality cases. The finite-isotropy assertion follows again from Corollary~\ref{cor33}.
\end{proof}

The parameter count concerns operators on the fixed tangent space. Quotienting by homotheties or by isometries is a different question. The formulas also show directly that every cometric in the classification has at most three eigenvalues.

\section{Weakly symmetric spaces and the choice of transitive group}
\label{sec:weak-symmetry-comparison}

In this section, we compare the classification with compact weakly symmetric spaces. We also give an example showing why the prescribed transitive group cannot be omitted from the g.o. condition.

A connected Riemannian manifold is weakly symmetric if, for every point $o$ and tangent vector $v\in T_oM$, an isometry fixes $o$ and sends $v$ to $-v$. This is equivalent to the condition that every two points can be interchanged by an isometry. For a presentation $M=G/H$, the term weakly symmetric pair is more restrictive. We allow a fixed isometry $\mu$ satisfying $\mu G\mu^{-1}=G$ and $\mu^2\in G$, and require the tangent-vector reversal to belong to the isotropy of the group generated by $G$ and $\mu$. This convention agrees with \cite[Definition~2.1 and Lemma~2.2]{Ng2000}. The connected group in this extension still has Lie algebra $\mathfrak g$.

Kr\"amer classified spherical pairs with compact connected simple first factor in \cite{Kr1979}. Nguyen proved that these are exactly the weakly symmetric pairs, allowing the fixed isometry $\mu$ above, and verified the tangent-vector reversal for every entry \cite[Theorem~1.1 and Table~1]{Ng2000}. Thus this list gives an independent comparison with the classification of $G$-g.o. metrics.

\begin{proposition}\label{prop111}
Every compact weakly symmetric pair with connected simple transitive group occurs in the preceding classification. For each pair $(G,H)$, all $G$-invariant Riemannian metrics are $G$-g.o. Up to local isomorphism and automorphisms of the acting Lie algebra, the pairs and the dimensions of their cometric cones are given in Table~\ref{tab:weak-symmetry-comparison}.
\end{proposition}

\begin{table}[ht]
\centering
\caption{Compact weakly symmetric pairs with simple transitive group. The dimension includes the overall scale.}
\label{tab:weak-symmetry-comparison}
\begin{tabular}{@{}p{0.84\linewidth}c@{}}
\hline
Pair & Dimension\\
\hline
Compact symmetric pairs with simple $G$ & $1$\\
$\SU(p+q)/(\SU(p)\SU(q))$, $p>q\ge1$ & $2$\\
$\SO(2n)/\SU(n)$, $n\ge3$ odd & $2$\\
$E_6/\Spin(10)$ & $2$\\
$\SU(2n+1)/\Sp(n)$, $n\ge2$ & $3$\\
$\SU(2n+1)/(\Sp(n)\mathrm \U(1))$, $n\ge2$ & $2$\\
$\SO(8)/\Spin(7)$, $\SO(7)/G_2$, $G_2/\SU(3)$ & $1$\\
$\SO(10)/(\Spin(7)\SO(2))$ & $2$\\
$\SO(9)/\Spin(7)$ & $2$\\
$\SO(8)/G_2$ & $3$\\
$\SO(2n+1)/\mathrm \U(n)$, $n\ge2$ & $2$\\
$\Sp(n)/(\Sp(n-1)\mathrm \U(1))$, $n\ge2$ & $2$\\
\hline
\end{tabular}
\end{table}

\begin{proof}
The nonsymmetric pairs are the five nonsymmetric families in \cite[Table~1]{Ng2000}. The first circle-bundle family belongs to the dominant standard unitary block in Theorem~\ref{thm411}. The other two circle-bundle families are $\SO(2n)/\SU(n)$ with odd $n$ and $E_6/\Spin(10)$. They occur in Theorems~\ref{thm615} and \ref{thm931}. The two symplectic-isotropy families in the unitary group occur in Theorem~\ref{thm411}. The Cayley-type families occur in Theorem~\ref{thm615}, and the final two families occur there and in Theorem~\ref{thm57}, respectively. The constant-curvature pairs and the symmetric pairs have irreducible isotropy and only the normal ray.

For each two-dimensional cone in the table, the isotropy representation has two inequivalent real irreducible summands. For $\SU(2n+1)/\Sp(n)$ it has three such summands. For $\SO(8)/G_2$ it is $\mathbb R^7\otimes\mathbb R^2$, and the cometric is $I_7\otimes P$ with arbitrary positive symmetric $P\in\operatorname{Sym}^2\mathbb R^2$. In every case the parameters in the preceding classification exhaust the full invariant cone. The three-dimensional cone in the last case is also described directly in \cite[Example~4]{Zi1996}.

The low-dimensional endpoints omitted from the table give no additional pairs. The case $\SU(3)/\Sp(1)=\SU(3)/\SU(2)$ is already in the first circle-bundle family and has two parameters. The cases $\SU(3)/(\Sp(1)\mathrm \U(1))$ and $\Sp(1)/\mathrm \U(1)$ are symmetric and have one parameter. Corollary~\ref{cor33} includes finite isotropy extensions and central coverings.
\end{proof}

The proposition concerns weakly symmetric pairs with prescribed connected group. It does not assert that every weakly symmetric metric invariant under a simple transitive subgroup $G$ is $G$-g.o. Yakimova's classification explicitly distinguishes a weakly symmetric reductive pair from the pair obtained by passing to its derived group \cite[Section~3]{Ya2004}. The following example proves this distinction directly.

\begin{example}\label{ex112}
There is a two-parameter family of weakly symmetric metrics on $\SU(6)/(\SU(3)\times\SU(3))$ for which only the normal ray is $\SU(6)$-g.o.
\end{example}
\begin{proof}
Put $G=\SU(6)$ and $H=\SU(3)\times\SU(3)$ in the standard block embedding. With respect to $Q=-B_{\mathfrak{su}(6)}$, write $\mathfrak m=\mathbb Rz\oplus\mathfrak p$, where
\begin{equation}\label{formula39}
 z=i\begin{pmatrix}I_3&0\\0&-I_3\end{pmatrix},\qquad
 P(X)=\begin{pmatrix}0&X\\-X^*&0\end{pmatrix},\qquad
 \mathfrak p=\{P(X):X\in M_3(\mathbb C)\}.
\end{equation}
For $a,b>0$, define the cometric $C=aI_{\mathbb Rz}\oplus bI_{\mathfrak p}$.

First we prove weak symmetry using a larger group. Set $G'=\mathrm \U(6)$ and $H'=\mathrm \U(3)\times\SU(3)$. Since $\det\colon H'\to\mathrm \U(1)$ is onto, $G$ acts transitively on $G'/H'$. Its isotropy is $G\cap H'=H$, so $G/H\cong G'/H'$. Under the tangent-space identification \eqref{formula39}, the isotropy element $(A,B)\in H'$ fixes $z$ and sends $P(X)$ to $P(AXB^*)$. Hence every cometric above is also $G'$-invariant. Entrywise complex conjugation induces an isometry $\mu$ fixing the origin, with $\mu_*(tz+P(X))=-tz+P(\overline X)$.

Fix a singular value decomposition $X=UDV^*$, where $U,V\in\mathrm \U(3)$ and $D$ is real diagonal. Choose a diagonal unitary matrix $R$ so that $\det(VRV^t)=1$. Set $A=-URU^t$ and $B=VRV^t$. Then $(A,B)\in H'$ and
\begin{equation}\label{formula40}
 A\overline X B^*=-URD R^*V^*=-UDV^*=-X.
\end{equation}
By \eqref{formula40}, the isometry induced by $(A,B)\mu$ fixes the origin and reverses $tz+P(X)$. It preserves each of the two summands in \eqref{formula39}. Thus every $a,b>0$ gives a weakly symmetric metric. This pair also appears in \cite[Table~1, row~1c]{Ya2004}.

Now suppose that this metric is $G$-g.o. for $G=\SU(6)$. Apply \eqref{formula11} to $Y=z+P(I_3)$. Any compensator has the form $Z=\operatorname{diag}(A_0,B_0)$ with $A_0,B_0\in\mathfrak{su}(3)$. Comparing the upper-right blocks, we obtain
\begin{equation}\label{formula41}
 A_0-B_0=2i(a-b)I_3.
\end{equation}
Taking the trace in \eqref{formula41} gives $a=b$. Conversely, $a=b$ is the normal metric and is $G$-g.o. Therefore changing the connected transitive group from $\mathrm \U(6)$ to $\SU(6)$ removes every nonnormal metric in this weakly symmetric family.
\end{proof}

\appendix

\section{Simple subalgebras of rank at least two}\label{app:high-rank-inventories}
In this section, we record the complete embedding lists used in the exceptional classification. We distinguish exhaustive simple-subalgebra lists from centralizer tables and retain the source corrections needed in type $E_8$.

\FloatBarrier
\paragraph{Simple subalgebras of rank at least two in $F_4$.}
Shan~\cite{Sh2024}, Sections~4.4--4.5, gives the complete list used in Table~\ref{tab:f4-high-inventory}. Each type other than $A_2$ has a unique embedding up to conjugacy. The four $A_2$ embeddings consist of the two mutually centralizing regular factors and two classes with zero connected centralizer. The two zero-centralizer classes need not be distinguished for the simple-isotropy argument.
\begin{table}[htbp]
\centering\small
\caption{Simple subalgebras of rank at least two in $F_4$}\label{tab:f4-high-inventory}
\begin{tabular}{lll}
Embedded class & Connected centralizer & Argument\\\hline
$B_4,D_4$ & $0$ & Simple isotropy\\
$A_3$ & $A_1$ & Regular simple ideal\\
$B_3$ & $\mathbb R$ & Maximal rank or simple isotropy\\
$C_3$ & $A_1$ & Regular simple ideal\\
$B_2$ & $A_1^2$ & Lemma~\ref{lem811} and regular isotropy\\
$G_2$ & $A_1$ & Theorem~\ref{thm86}\\
$A_2$ (long-root regular factor) & $A_2$ & Regular simple ideal\\
$A_2$ (short-root regular factor) & $A_2$ & Regular simple ideal\\
$A_2$ (two further classes) & $0$ & Simple isotropy
\end{tabular}
\end{table}

\FloatBarrier
\paragraph{Nonregular simple subalgebras of rank at least two in $E_6$.}
The exhaustive Table~6 of Minchenko~\cite{Mi2006}, cross-checked against the centralizers in Table~12 (printed p.~255), gives the following complete list of the $17$ nonregular linear-equivalence classes of simple subalgebras of rank at least two in $E_6$. The three classes marked by a star each split into two conjugacy classes exchanged by the outer automorphism. The proofs above cover both.
\begin{table}[htbp]
\centering\small
\caption{Nonregular simple subalgebras of rank at least two in $E_6$}\label{tab:e6-high-inventory}
$\begin{array}{c|c|l}
 \mathfrak s&\mathfrak c_{\mathfrak e_6}(\mathfrak s)&
 \text{normality argument for every }\mathfrak h\supset_{\rm ideal}\mathfrak s\\\hline
 A_2^{2\prime}&A_2^1&\text{Proposition~\ref{prop93}}\\
 A_2^{2\prime\prime}&G_2^1&\text{Proposition~\ref{prop92}}\\
 A_2^{3\prime}&\mathbb R^2&\text{Proposition~\ref{prop95}}\\
 A_2^{3\prime\prime}&0&\text{simple isotropy}\\
 A_2^5&A_1^1&\text{Proposition~\ref{prop97}}\\
 A_2^9\ (*) &0&\text{simple isotropy}\\
 B_2^1&B_2^1+\mathbb R&\text{Proposition~\ref{prop92}}\\
 B_2^2&A_1^1+\mathbb R&\text{Proposition~\ref{prop94}}\\
 B_2^3\ (*) &\mathbb R&\text{Proposition~\ref{prop96}}\\
 G_2^1&A_2^{2\prime\prime}&\text{Proposition~\ref{prop92}}\\
 G_2^3\ (*) &0&\text{simple isotropy}\\
 A_3^2&A_1^1&\text{Proposition~\ref{prop97}}\\
 B_3^1&A_1^2+\mathbb R&\text{Proposition~\ref{prop92}}\\
 C_3^1&A_1^1&\text{Proposition~\ref{prop97}}\\
 B_4^1&\mathbb R&\text{Proposition~\ref{prop96}}\\
 C_4^1&0&\text{simple isotropy}\\
 F_4^1&0&\text{simple isotropy}
\end{array}$
\end{table}
A zero centralizer forces $\mathfrak h=\mathfrak s$ because $\mathfrak s$ is an ideal. The regular-simple-ideal theorem supplies the root-ideal hypothesis in Proposition~\ref{prop97}. Thus every row applies to an arbitrary residual compact subalgebra in its centralizer.

Minchenko, Table~6, rows $1,2,15$, records the three outer-conjugate pairs. Each member is covered by the indicated proof.

\FloatBarrier
\paragraph{Nonregular simple subalgebras of rank at least two in $E_7$.}
The following is the complete list of nonregular rank-at-least-two simple ideals, with centralizers from Minchenko's Table~13. To check completeness independently of that centralizer table, the first column gives the corresponding row in Minchenko's exhaustive Table~7. Every row $1$--$27$ occurs exactly once. A label ``via $\mathfrak q$'' invokes Lemma~\ref{lem912}. ``whole envelope'' invokes one of the previously proved hereditary maximal-subgroup results. A zero centralizer forces the isotropy to be simple.
\begin{table}[htbp]
\centering\small
\caption{Nonregular simple subalgebras of rank at least two in $E_7$}\label{tab:e7-high-inventory}
$\begin{array}{r|c|c|l}
 \text{Table 7}&\mathfrak s&\mathfrak c(\mathfrak s)&\text{argument}\\\hline
 2&A_2^{2\prime}&A_2^1+A_1^{3\prime\prime}&\text{via regular }A_2^1\\
 3&A_2^{2\prime\prime}&G_2^1+\mathbb R&\text{via }G_2^1\\
 4&A_2^{3\prime}&3A_1^1&\text{via root }A_1^1\\
 5&A_2^{3\prime\prime}&A_1^{3\prime\prime}&\text{via }A_1^{3\prime\prime}\\
 6&A_2^{5\prime}&A_1^1+\mathbb R&\text{via root }A_1^1\\
 19&A_2^{5\prime\prime}&A_2^1&\text{via regular }A_2^1\\
 26&A_2^{6\prime}&0&\text{simple isotropy}\\
 27&A_2^{6\prime\prime}&0&\text{simple isotropy}\\
 13&A_2^9&\mathbb R&\text{Proposition~\ref{prop913}}\\
 1&A_2^{21}&0&\text{simple isotropy}\\
 7&B_2^1&B_3^1+A_1^1&\text{via root }A_1^1\\
 18&B_2^2&A_2^1+\mathbb R&\text{via regular }A_2^1\\
 14&B_2^3&A_1^1+\mathbb R&\text{via root }A_1^1\\
 8&G_2^1&C_3^{1\prime\prime}&G_2C_3\text{ whole envelope}\\
 20&G_2^2&A_1^7&G_2A_1\text{ whole envelope}\\
 15&G_2^3&\mathbb R&\text{Proposition~\ref{prop913}}\\
 16&A_3^{2\prime}&A_1^1+\mathbb R&\text{via root }A_1^1\\
 24&A_3^{2\prime\prime}&A_2^1&\text{via regular }A_2^1\\
 9&B_3^1&B_2^1+A_1^1&\text{via root }A_1^1\\
 25&B_3^2&\mathbb R&\text{Proposition~\ref{prop914}}\\
 10&C_3^{1\prime}&A_1^1+A_1^{3\prime\prime}&\text{via root }A_1^1\\
 21&C_3^{1\prime\prime}&G_2^1&G_2C_3\text{ whole envelope}\\
 11&B_4^1&A_1^1+A_1^2&\text{via root }A_1^1\\
 17&C_4^1&\mathbb R&\text{Proposition~\ref{prop913}}\\
 22&D_4^2&0&\text{simple isotropy}\\
 12&F_4^1&A_1^{3\prime\prime}&F_4A_1\text{ whole envelope}\\
 23&B_5^1&A_1^1&\text{via root }A_1^1
\end{array}$
\end{table}
The only nonregular factors used on the right are $G_2^1$ and $A_1^{3\prime\prime}$, whose full centralizer subgroups are covered by Theorems~\ref{thm824} and~\ref{thm831}. The four purely central one-dimensional envelopes are proved above. This exhausts every rank-at-least-two simple ideal and every compact residual algebra in its centralizer.

\FloatBarrier
\paragraph{Nonregular simple subalgebras of rank at least two in $E_8$.}
The complete list of nonregular simple subalgebras of rank at least two is Minchenko's Table~8, pp.~252--253. Its $43$ linear conjugacy classes give $44$ conjugacy classes because row~23 splits into two.  We cross-check their centralizers against Table~14, pp.~257--259, rather than treating that latter table as a complete list by itself. Two discrepancies must be recorded. First, Table~14 omits $B_6^1=\mathfrak{so}(13)$, although Table~8, row~34 explicitly includes it. The text on p.~239 also uses $D_6\subset B_6\subset E_8$. Its standard realization in $D_8$ gives centralizer $A_1^2$ and is covered by Proposition~\ref{prop921}. Second, the $B_2^{3\prime\prime}$ row of Table~14 gives Lie centralizer $A_1^{10\prime\prime}$ but group centralizer $\SL_4$.  These entries are incompatible. Proposition~\ref{prop926} independently computes the Lie centralizer as $A_1^{10\prime\prime}$. No exclusion is inferred from the erroneous group entry. All class labels in the following table refer to the Lie-algebra column and its cross-check with Table~8.

\FloatBarrier
\paragraph{Eliminated classes.}
The following table exhausts all $44$ classes. ``Propagation'' means Lemma~\ref{lem920}. The cited results apply to every isotropy containing the indicated factor as an ideal.
\begin{table}[htbp]
\centering\small
\caption{All $44$ nonregular rank-at-least-two classes in $E_8$}\label{tab:e8-high-inventory}
$\begin{array}{c|l}
 \text{simple-ideal classes}&\text{reason}\\\hline
 A_2^2,A_2^{3\prime},A_2^5,A_2^{6\prime}(1),
 A_2^{6\prime\prime},A_2^9,A_2^{21}
 &\text{propagation from a regular ideal}\\
 A_2^{3\prime\prime}&\text{propagation from }G_2^1\\
 A_2^{6\prime}(2)&\text{hereditary maximal }A_2A_1\\
 B_2^{2\prime},B_2^{3\prime},G_2^2,G_2^3,A_3^{2\prime},
 B_3^{2\prime},C_3^1,C_4^{1\prime},D_4^{2\prime}
 &\text{propagation from a regular ideal}\\
 G_2^1,F_4^1&\text{hereditary maximal }G_2F_4\\
 B_2^1,B_3^1,B_4^1,B_5^1,B_6^1
 &\text{odd orthogonal blocks}\\
 B_2^{2\prime\prime}&\Sp(2)\Sp(2)\text{ tensor envelope}\\
 C_4^{1\prime\prime}&\Sp(4)\Sp(1)\text{ tensor envelope}\\
 A_3^{2\prime\prime}&\text{propagation from }A_1^{4\prime}\\
 A_2^4&\text{Proposition~\ref{prop925}}\\
 B_2^{3\prime\prime}&\text{Proposition~\ref{prop926}}\\
 B_3^{2\prime\prime},D_4^{2\prime\prime\prime}
 &\text{Proposition~\ref{prop924}}\\
 B_2^7,G_2^4,C_3^2&\text{the three }\SO(14)\text{ circle envelopes}\\
 A_2^{10},B_2^4,B_2^{12},A_3^4,A_4^2,
 B_4^{2\prime},B_4^{2\prime\prime},D_4^{2\prime\prime},B_7^1
 &\text{zero centralizer, simple isotropy}
\end{array}$
\end{table}
Where Table~14 writes $A_4$ or $A_3$ without a superscript in a centralizer, these are the regular embeddings. The other classes of those types cannot have a centralizer containing the required $B_2$ factor.  Equivalently they are the complementary regular factors in $A_4+A_4$ and $D_5+A_3$.

\FloatBarrier
\paragraph{Application of centralizer propagation.}
The regular-simple-ideal theorem, the hereditary results for $G_2F_4$ and $A_2A_1$, and the orthogonal, symplectic and $\SO(14)$ circle calculations supply the hypotheses for centralizer propagation. The only row using a nonregular rank-one ideal is $A_3^{2\prime\prime}$, whose centralizer contains $A_1^{4\prime}$. Proposition~\ref{prop923} proves normality for this rank-one ideal directly. The other nonregular factor used in propagation is $G_2^1$, which is covered by Theorem~\ref{thm826}. Proposition~\ref{prop926} and Theorem~\ref{thm832} also prove the assertions for $A_1^{10\prime\prime}$ and $A_1^{16}$, respectively. These results are available before the rank-one classification.

\section{Rank-one subalgebras and the order of propagation}\label{app:rank-one-inventories}
In this section, we list every rank-one class used in the exceptional proof. The centralizers are those of the complete $\mathfrak{sl}_2$-triple. Each table identifies a preceding result, and the propagation order is specified before it is used.

\FloatBarrier
\paragraph{Rank-one subalgebras of $F_4$.}
Shan~\cite{Sh2024}, Sections~4.4--4.5 and~4.7, gives the fifteen classes in Table~\ref{tab:f4-a1-inventory}. The labels are Bala--Carter labels of the associated nilpotent orbit. The prose in Proposition~4.4.5(2) says seven $\SU(2)$ classes but lists eight partitions. Together with the seven $\PSU(2)$ classes, the list has fifteen entries. We use the listed classes rather than that numerical typo.
\begin{table}[htbp]
\centering\small
\caption{All rank-one classes in $F_4$}\label{tab:f4-a1-inventory}
\begin{tabular}{lll}
Class & Connected centralizer & Argument\\\hline
$A_1$ & $C_3$ & Regular simple ideal\\
$\widetilde A_1$ & $A_3$ & Regular simple ideal\\
$A_1+\widetilde A_1$ & $A_1^2$ & Lemma~\ref{lem89}(3)\\
$A_2$ & $A_2$ & Lemma~\ref{lem810}\\
$\widetilde A_2$ & $G_2$ & Theorem~\ref{thm86}\\
$A_2+\widetilde A_1$ & $A_1$ & Free modules $(V_4,V_2),(V_2,V_4)$\\
$B_2$ & $A_1^2$ & Lemma~\ref{lem89}(2)\\
$\widetilde A_2+A_1$ & $A_1$ & Free modules $(V_4,V_2),(V_2,V_4)$\\
$C_3(a_1)$ & $A_1$ & Free modules $(V_6,V_2),(V_4,V_2)$\\
$F_4(a_3)$ & $0$ & Simple isotropy\\
$B_3$ & $A_1$ & Theorem~\ref{thm86}\\
$C_3$ & $A_1$ & Free modules $(V_{10},V_2),(V_4,V_2)$\\
$F_4(a_2)$ & $0$ & Simple isotropy\\
$F_4(a_1)$ & $0$ & Simple isotropy\\
$F_4$ & $0$ & Simple isotropy
\end{tabular}
\end{table}
Here $(V_d,V_e)$ denotes the real form of $V_d\otimes V_e$. The required free vectors are proved in Lemma~\ref{lem88}. For $C_3(a_1)$, the branching also follows from $\mathbb C^6=V_4\oplus V_2$ and $\Lambda^3_0\mathbb C^6=V_6\oplus2V_4$. Each pair lies in the complement of the full centralizer envelope. Consequently every residual compact subalgebra, including a central torus, is covered by the same argument.

\FloatBarrier
\paragraph{Rank-one subalgebras of $E_6$.}
The complete rank-one table is as follows. The nilpotent labels and centralizer types are independently given by Arakawa--van Ekeren--Moreau~\cite{ArVaMo2024}, Table~13. The Dynkin indices and embedded centralizer factors agree with the explicit Milev calculator table. In the last column, ``via $\mathfrak q$'' means Lemma~\ref{lem98}.
\begin{table}[htbp]
\centering\small
\caption{All $20$ rank-one classes in $E_6$}\label{tab:e6-a1-inventory}
$\begin{array}{c|c|c|l}
 \text{index}&\text{Bala--Carter label}&\mathfrak c(\mathfrak s)&\text{argument}\\\hline
 1&A_1&A_5^1&\text{regular root ideal}\\
 2&2A_1&B_3^1+\mathbb R&\text{Proposition~\ref{prop92}}\\
 3&3A_1&A_2^{2\prime\prime}+A_1^1&\text{via }A_1^1\\
 4&A_2&2A_2^1&\text{via }A_2^1\\
 5&A_2+A_1&A_2^1+\mathbb R&\text{via }A_2^1\\
 6&A_2+2A_1&A_1^6+\mathbb R&\text{Proposition~\ref{prop910}}\\
 8&2A_2&G_2^1&\text{via }G_2^1\\
 9&2A_2+A_1&A_1^3&\text{via index }3\text{, already closed}\\
 10&A_3&B_2^1+\mathbb R&\text{via }B_2^1\\
 11&A_3+A_1&A_1^1+\mathbb R&\text{via }A_1^1\\
 12&D_4(a_1)&\mathbb R^2&\text{Proposition~\ref{prop99}}\\
 20&A_4&A_1^1+\mathbb R&\text{via }A_1^1\\
 21&A_4+A_1&\mathbb R&\text{Proposition~\ref{prop99}}\\
 28&D_4&A_2^{2\prime\prime}&\text{via }A_2^{2\prime\prime}\\
 30&D_5(a_1)&\mathbb R&\text{Proposition~\ref{prop99}}\\
 35&A_5&A_1^1&\text{via }A_1^1\\
 36&E_6(a_3)&0&\text{simple isotropy}\\
 60&D_5&\mathbb R&\text{Proposition~\ref{prop99}}\\
 84&E_6(a_1)&0&\text{simple isotropy}\\
 156&E_6&0&\text{simple isotropy}
\end{array}$
\end{table}
The regular-ideal theorem and Table~\ref{tab:e6-high-inventory} provide the initial cases. Proposition~\ref{prop92} treats index two, and Proposition~\ref{prop910} treats index six. The other rows follow by propagation in increasing index order. The only row requiring another nonregular rank-one class is index nine, whose centralizer contains the index-three class already covered through a regular root ideal.

\FloatBarrier
\paragraph{Rank-one subalgebras of $E_7$.}
Arakawa--van Ekeren--Moreau, Tables~14 and~15, give all $44$ nonzero $A_1$ classes in $E_7$ and their reductive centralizers. The table below records every row. The superscripts on simple centralizer factors give their indices in $E_7$. We omit the prime on index-three $A_1$ factors because the two classes satisfy the required normality assertion before either is used in propagation. For centralizer factors of rank at least two their precise embedded class remains the one specified by the source table. We use the result for \emph{all} those classes, and hence do not identify two embeddings merely from their abstract type and index.

The final column specifies the result used in each case. $\mathrm R$ is the regular-ideal theorem. $\mathrm P_{\ge2}$ is centralizer propagation through an ideal of rank at least two. $\mathrm P_j$ is propagation through a previously closed $A_1$ of index $j$. $\mathrm M$ is the already proved hereditary theorem for the maximal $A_1^{24}A_1^{15}$. $\mathrm T$ is Proposition~\ref{prop916}. The symbol $\mathrm S$ denotes the simple-isotropy theorem, applicable because the centralizer is zero. The propagation argument is precisely Lemma~\ref{lem912}. Its proof applies unchanged when the nonregular simple ideal has rank one.

\begin{table}[htbp]
\centering\small
\caption{Rank-one classes in $E_7$}\label{tab:e7-a1-inventory}
\begin{tabular}{llll}
Class & Index & Centralizer & Argument\\\hline
$A_1$ & $1$ & $D_6^1$ & $\mathrm R$\\
$2A_1$ & $2$ & $A_1^1+B_4^1$ & $\mathrm P_1$\\
$(3A_1)'$ & $3$ & $A_1^1+C_3^1$ & $\mathrm P_1$\\
$(3A_1)''$ & $3$ & $F_4^1$ & $\mathrm P_{\ge2}$\\
$4A_1$ & $4$ & $C_3^1$ & $\mathrm P_{\ge2}$\\
$A_2$ & $4$ & $A_5^1$ & $\mathrm P_{\ge2}$\\
$A_2+A_1$ & $5$ & $\mathbb R+A_3^1$ & $\mathrm P_{\ge2}$\\
$A_2+2A_1$ & $6$ & $A_1^1+A_1^2+A_1^6$ & $\mathrm P_1$\\
$A_2+3A_1$ & $7$ & $G_2^2$ & $\mathrm P_{\ge2}$\\
$2A_2$ & $8$ & $A_1^3+G_2^1$ & $\mathrm P_{\ge2}$\\
$2A_2+A_1$ & $9$ & $A_1^3+A_1^3$ & $\mathrm P_3$\\
$A_3$ & $10$ & $A_1^1+B_3^1$ & $\mathrm P_1$\\
$(A_3+A_1)'$ & $11$ & $A_1^1+A_1^1+A_1^2$ & $\mathrm P_1$\\
$(A_3+A_1)''$ & $11$ & $B_3^1$ & $\mathrm P_{\ge2}$\\
$A_3+2A_1$ & $12$ & $A_1^1+A_1^2$ & $\mathrm P_1$\\
$D_4(a_1)$ & $12$ & $A_1^1+A_1^1+A_1^1$ & $\mathrm P_1$\\
$D_4(a_1)+A_1$ & $13$ & $A_1^1+A_1^1$ & $\mathrm P_1$\\
$A_3+A_2$ & $14$ & $\mathbb R+A_1^1$ & $\mathrm P_1$\\
$A_3+A_2+A_1$ & $15$ & $A_1^{24}$ & $\mathrm M$\\
$A_4$ & $20$ & $\mathbb R+A_2^1$ & $\mathrm P_{\ge2}$\\
$A_4+A_1$ & $21$ & $\mathbb R^2$ & $\mathrm T$\\
$A_4+A_2$ & $24$ & $A_1^{15}$ & $\mathrm M$\\
\end{tabular}
\end{table}
\begin{table}[htbp]
\centering\small
\caption{Rank-one classes in $E_7$ (continued)}\label{tab:e7-a1-inventory-2}
\begin{tabular}{llll}
Class & Index & Centralizer & Argument\\\hline
$D_4$ & $28$ & $C_3^1$ & $\mathrm P_{\ge2}$\\
$D_4+A_1$ & $29$ & $B_2^1$ & $\mathrm P_{\ge2}$\\
$D_5(a_1)$ & $30$ & $\mathbb R+A_1^1$ & $\mathrm P_1$\\
$D_5(a_1)+A_1$ & $31$ & $A_1^8$ & $\mathrm P_8$\\
$(A_5)'$ & $35$ & $A_1^1+A_1^3$ & $\mathrm P_1$\\
$(A_5)''$ & $35$ & $G_2^1$ & $\mathrm P_{\ge2}$\\
$A_5+A_1$ & $36$ & $A_1^3$ & $\mathrm P_3$\\
$E_6(a_3)$ & $36$ & $A_1^3$ & $\mathrm P_3$\\
$D_6(a_2)$ & $38$ & $A_1^1$ & $\mathrm P_1$\\
$E_7(a_5)$ & $39$ & $0$ & $\mathrm S$\\
$A_6$ & $56$ & $A_1^7$ & $\mathrm P_7$\\
$D_5$ & $60$ & $A_1^1+A_1^2$ & $\mathrm P_1$\\
$D_5+A_1$ & $61$ & $A_1^2$ & $\mathrm P_2$\\
$D_6(a_1)$ & $62$ & $A_1^1$ & $\mathrm P_1$\\
$E_7(a_4)$ & $63$ & $0$ & $\mathrm S$\\
$E_6(a_1)$ & $84$ & $\mathbb R$ & $\mathrm T$\\
$D_6$ & $110$ & $A_1^1$ & $\mathrm P_1$\\
$E_7(a_3)$ & $111$ & $0$ & $\mathrm S$\\
$E_6$ & $156$ & $A_1^3$ & $\mathrm P_3$\\
$E_7(a_2)$ & $159$ & $0$ & $\mathrm S$\\
$E_7(a_1)$ & $231$ & $0$ & $\mathrm S$\\
$E_7$ & $399$ & $0$ & $\mathrm S$\\
\end{tabular}
\end{table}
Here $\mathrm P_2$ specifically means propagation through the already closed \emph{index-two} $A_1$ (the row $2A_1$).

The regular-ideal theorem first gives the index-one case. Centralizer propagation then treats index two, both index-three classes, and indices seven and eight through regular or rank-at-least-two centralizer ideals. These classes suffice for all subsequent propagation rows. The two rows with centralizers $A_1^{15}$ and $A_1^{24}$ are both covered by Theorem~\ref{thm833}. The maximal pair and its indices are given in \cite[Theorem~5.1, equation~(5.7)]{AdKaMoPaPe2018}.

The inventory is \cite[Tables~14--15, pp.~86--87]{ArVaMo2024}. The argument uses the regular-ideal theorem, Theorem~\ref{thm915}, Lemma~\ref{lem912}, Theorem~\ref{thm833}, and the two-free and simple-isotropy results. Each subgroup assertion is supplied by one of these results, rather than inferred from normality of a maximal quotient.

\FloatBarrier
\paragraph{Rank-one subalgebras of $E_8$.}
Arakawa--van Ekeren--Moreau, Tables~16 and~17, list all $69$ nonzero $A_1$ classes of $E_8$. There are $11$ zero-centralizer classes and exactly the four pure-toral classes just treated. The following table transcribes every row. Superscripts on simple centralizer factors are their indices in $E_8$, obtained from the coefficients of $k$ in the source. Embedded-class primes are omitted only when the argument applies to every class of the displayed type and index. The exceptional $B_2^{3''}A_1^{10''}$ seed is explicitly distinguished below.

The regular-ideal theorem treats index one. The class $2A_1$ has centralizer $B_6^1$, so Theorem~\ref{thm927} treats index two. The index-three class $3A_1$ has a root $A_1$ in its centralizer. The two index-four classes, $4A_1$ and $A_2$, are treated through $C_4$ and $E_6$, respectively. Similarly, $A_2+2A_1$ of index six is treated through $B_3$, and $A_2+3A_1$ of index seven through a root $A_1$. Among the index-ten classes, $A_3$ is treated through $B_5$, while $2A_2+2A_1$ is covered by Proposition~\ref{prop926}. The index-fifteen class $A_3+A_2+A_1$ has a root $A_1$ in its centralizer. These cases provide all rank-one factors used in the remaining propagation rows. The $B_6$ entry follows from \cite[Table~8, row~34]{Mi2006} and Proposition~\ref{prop921}, despite its omission from the centralizer Table~14.

In the last column, $\mathrm R$ denotes the regular-ideal theorem. $\mathrm P_{\ge2}$ denotes Lemma~\ref{lem920} applied to an already closed factor of rank at least two. $\mathrm P_j$ denotes that lemma applied to an already closed index-$j$ rank-one factor. $\mathrm B$ is the independent $B_2^{3''}A_1^{10''}$ hereditary theorem. $\mathrm T$ is Proposition~\ref{prop928}. The symbol $\mathrm S$ denotes the simple-isotropy theorem. The latter applies because a zero centralizer forces the isotropy to be exactly the given simple algebra.

\begin{table}[htbp]
\centering\small
\caption{Rank-one classes in $E_8$}\label{tab:e8-a1-inventory}
\begin{tabular}{lll}
Class & Centralizer & Argument\\\hline
$E_8$ & $\{0\}$ & $\mathrm S$\\
$E_8(a_1)$ & $\{0\}$ & $\mathrm S$\\
$E_8(a_2)$ & $\{0\}$ & $\mathrm S$\\
$E_8(a_3)$ & $\{0\}$ & $\mathrm S$\\
$E_8(a_4)$ & $\{0\}$ & $\mathrm S$\\
$E_7$ & $A_1^1$ & $\mathrm P_{1}$\\
$E_8(b_4)$ & $\{0\}$ & $\mathrm S$\\
$E_8(a_5)$ & $\{0\}$ & $\mathrm S$\\
$E_7(a_1)$ & $A_1^1$ & $\mathrm P_{1}$\\
$E_8(b_5)$ & $\{0\}$ & $\mathrm S$\\
$D_7$ & $A_1^2$ & $\mathrm P_{2}$\\
$E_8(a_6)$ & $\{0\}$ & $\mathrm S$\\
$E_7(a_2)$ & $A_1^1$ & $\mathrm P_{1}$\\
$E_6+A_1$ & $A_1^3$ & $\mathrm P_{3}$\\
$D_7(a_1)$ & $\mathbb R$ & $\mathrm T$\\
$E_8(b_6)$ & $\{0\}$ & $\mathrm S$\\
$E_7(a_3)$ & $A_1^1$ & $\mathrm P_{1}$\\
$E_6(a_1)+A_1$ & $\mathbb R$ & $\mathrm T$\\
$A_7$ & $A_1^4$ & $\mathrm P_{4}$\\
$D_7(a_2)$ & $\mathbb R$ & $\mathrm T$\\
$E_6$ & $G_2^1$ & $\mathrm P_{\ge2}$\\
$D_6$ & $B_2^1$ & $\mathrm P_{\ge2}$\\
$D_5+A_2$ & $\mathbb R$ & $\mathrm T$\\
\end{tabular}
\end{table}

\begin{table}[htbp]
\centering\small
\caption{Rank-one classes in $E_8$ (continued)}\label{tab:e8-a1-inventory-2}
\begin{tabular}{lll}
Class & Centralizer & Argument\\\hline
$E_6(a_1)$ & $A_2^1$ & $\mathrm P_{\ge2}$\\
$E_7(a_4)$ & $A_1^1$ & $\mathrm P_{1}$\\
$A_6+A_1$ & $A_1^7$ & $\mathrm P_{7}$\\
$D_6(a_1)$ & $A_1^1+A_1^1$ & $\mathrm P_{1}$\\
$A_6$ & $A_1^1+A_1^7$ & $\mathrm P_{1}$\\
$E_8(a_7)$ & $\{0\}$ & $\mathrm S$\\
$D_5+A_1$ & $A_1^1+A_1^2$ & $\mathrm P_{1}$\\
$E_7(a_5)$ & $A_1^1$ & $\mathrm P_{1}$\\
$E_6(a_3)+A_1$ & $A_1^3$ & $\mathrm P_{3}$\\
$D_6(a_2)$ & $A_1^1+A_1^1$ & $\mathrm P_{1}$\\
$D_5(a_1)+A_2$ & $A_1^6$ & $\mathrm P_{6}$\\
$A_5+A_1$ & $A_1^1+A_1^3$ & $\mathrm P_{1}$\\
$A_4+A_3$ & $A_1^{10}$ & $\mathrm P_{10}$\\
$D_5$ & $B_3^1$ & $\mathrm P_{\ge2}$\\
$E_6(a_3)$ & $G_2^1$ & $\mathrm P_{\ge2}$\\
$D_4+A_2$ & $A_2^2$ & $\mathrm P_{\ge2}$\\
$A_4+A_2+A_1$ & $A_1^{15}$ & $\mathrm P_{15}$\\
$D_5(a_1)+A_1$ & $A_1^1+A_1^8$ & $\mathrm P_{1}$\\
$A_5$ & $A_1^1+G_2^1$ & $\mathrm P_{1}$\\
$A_4+A_2$ & $A_1^{15}+A_1^1$ & $\mathrm P_{1}$\\
$A_4+2A_1$ & $\mathbb R+A_1^2$ & $\mathrm P_{2}$\\
$D_5(a_1)$ & $A_3^1$ & $\mathrm P_{\ge2}$\\
$2A_3$ & $B_2^2$ & $\mathrm P_{\ge2}$\\
\end{tabular}
\end{table}

\begin{table}[htbp]
\centering\small
\caption{Rank-one classes in $E_8$ (continued)}\label{tab:e8-a1-inventory-3}
\begin{tabular}{lll}
Class & Centralizer & Argument\\\hline
$A_4+A_1$ & $\mathbb R+A_2^1$ & $\mathrm P_{\ge2}$\\
$D_4(a_1)+A_2$ & $A_2^6$ & $\mathrm P_{\ge2}$\\
$D_4+A_1$ & $C_3^1$ & $\mathrm P_{\ge2}$\\
$A_3+A_2+A_1$ & $A_1^1+A_1^{24}$ & $\mathrm P_{1}$\\
$A_4$ & $A_4^1$ & $\mathrm P_{\ge2}$\\
$A_3+A_2$ & $\mathbb R+B_2^1$ & $\mathrm P_{\ge2}$\\
$D_4(a_1)+A_1$ & $A_1^1+A_1^1+A_1^1$ & $\mathrm P_{1}$\\
$A_3+2A_1$ & $A_1^2+B_2^1$ & $\mathrm P_{\ge2}$\\
$2A_2+2A_1$ & $B_2^3$ & $\mathrm B$\\
$D_4$ & $F_4^1$ & $\mathrm P_{\ge2}$\\
$D_4(a_1)$ & $D_4^1$ & $\mathrm P_{\ge2}$\\
$A_3+A_1$ & $A_1^1+B_3^1$ & $\mathrm P_{1}$\\
$2A_2+A_1$ & $A_1^3+G_2^1$ & $\mathrm P_{\ge2}$\\
$2A_2$ & $G_2^1+G_2^1$ & $\mathrm P_{\ge2}$\\
$A_2+3A_1$ & $A_1^1+G_2^2$ & $\mathrm P_{1}$\\
$A_3$ & $B_5^1$ & $\mathrm P_{\ge2}$\\
$A_2+2A_1$ & $A_1^6+B_3^1$ & $\mathrm P_{\ge2}$\\
$A_2+A_1$ & $A_5^1$ & $\mathrm P_{\ge2}$\\
$4A_1$ & $C_4^1$ & $\mathrm P_{\ge2}$\\
$A_2$ & $E_6^1$ & $\mathrm P_{\ge2}$\\
$3A_1$ & $A_1^1+F_4^1$ & $\mathrm P_{1}$\\
$2A_1$ & $B_6^1$ & $\mathrm P_{\ge2}$\\
$A_1$ & $E_7^1$ & $\mathrm R$\\
\end{tabular}
\end{table}

The full centralizer formulas and indices are given in \cite[Tables~16--17 and Appendix~A]{ArVaMo2024}. They are centralizers of the entire $\mathfrak{sl}_2$-triple and therefore give the compact reductive centralizers used here. They are not the larger nonreductive centralizers of a nilpotent element.

\section{Root and representation calculations}\label{app:branching-checks}
In this section, we verify the root deletions and representation restrictions used in the proof. Throughout, $V_d=\operatorname{Sym}^{d-1}\mathbb C^2$ has complex dimension $d$. Real forms and underlying real representations are distinguished explicitly.

\FloatBarrier
\paragraph{Maximal Levi subgroups.}
We first evaluate the central projection for every maximal Levi subgroup.  Number simple roots as in Bourbaki. The $E_r$ diagram has edges $1$--$3$--$4$--$5$--$\cdots$--$r$ and $2$--$4$. Let $L_i$ be the Levi obtained by deleting $\alpha_i$, let $m_i$ be its coefficient in the highest root $\theta$, and let $M_i$ be the root subgroup obtained by replacing $\alpha_i$ with $\alpha_0=-\theta$.  When $m_i>1$, $M_i$ is proper, semisimple, of full rank, and contains $L_i$. The following is the complete deletion calculation. The column ``method'' refers to an abelian-stabilizer overgroup (A), simple derived algebra and Corollary~\ref{cor75} (B), Proposition~\ref{prop710} (C), or the known nonzero $D_5$ entry (N).  A dash means that $m_i=1$, so replacement does not give a proper overgroup.
\begin{table}[htbp]
\centering\small
\caption{Maximal Levi subgroups and affine-node replacements}\label{tab:exceptional-levi-deletions}
$\begin{array}{c|c|c|c|c|c}
 \mathfrak g&i&m_i&[\mathfrak l_i,\mathfrak l_i]&\mathfrak m_i&\text{method}\\\hline
 G_2&1&3&A_1&A_2&B\\
    &2&2&A_1&A_1+A_1&A\\\hline
 F_4&1&2&C_3&C_3+A_1&A\\
    &2&3&A_2+A_1&A_2+A_2&A\\
    &3&4&A_2+A_1&A_3+A_1&A\\
    &4&2&B_3&B_4&B\\\hline
 E_6&1&1&D_5&-&N\\
    &2&2&A_5&A_5+A_1&A\\
    &3&2&A_4+A_1&A_5+A_1&A\\
    &4&3&A_2+A_2+A_1&A_2+A_2+A_2&A\\
    &5&2&A_4+A_1&A_5+A_1&A\\
    &6&1&D_5&-&N\\\hline
 E_7&1&2&D_6&D_6+A_1&B\\
    &2&2&A_6&A_7&A\\
    &3&3&A_5+A_1&A_5+A_2&A\\
    &4&4&A_3+A_2+A_1&A_3+A_3+A_1&A\\
    &5&3&A_4+A_2&A_5+A_2&A\\
    &6&2&D_5+A_1&D_6+A_1&C\\
    &7&1&E_6&-&B\\\hline
 E_8&1&2&D_7&D_8&A\\
    &2&3&A_7&A_8&A\\
    &3&4&A_6+A_1&A_7+A_1&A\\
    &4&6&A_4+A_2+A_1&A_5+A_2+A_1&A\\
    &5&5&A_4+A_3&A_4+A_4&A\\
    &6&4&D_5+A_2&D_5+A_3&A\\
    &7&3&E_6+A_1&E_6+A_2&A\\
    &8&2&E_7&E_7+A_1&B
\end{array}$
\end{table}
For completeness, the added affine node meets node $2,1,8$ in $E_6,E_7,E_8$, respectively.  For $F_4$, it meets node $1$. Nodes $1,2$ are long and $3,4$ short.  For $G_2$, node $1$ is short and the affine node meets node $2$. These edges and the displayed marks determine every entry.

\FloatBarrier
\paragraph{Fixed vectors for $D_4$ in the fundamental $E_7$ representations.}
\begin{proposition}\label{propC1}
Every fundamental representation of the simply connected compact group $E_7$ contains a nonzero vector fixed by the standard regular $\Spin(8)$ subgroup.
\end{proposition}
\begin{proof}
Let $S=\Spin(8)$ denote this subgroup. We denote the indicated $E_7$, $D_6$ and $D_4$ modules by $U$ with dimension subscripts, and retain $V_d$ for the $A_1$ factor. We check the seven fundamental representations, whose dimensions are $56,133,912,1539,8645,27664,365750$. Under the standard $D_4A_1^3$ subgroup,
\begin{align*}
 U_{56}|_S&=2U_{8v}\oplus2U_{8s}\oplus2U_{8c}\oplus8\mathbb C,\\
 U_{133}|_S&=U_{28}\oplus4U_{8v}\oplus4U_{8s}
                       \oplus4U_{8c}\oplus9\mathbb C.
\end{align*}
In particular their fixed-space dimensions are $8$ and $9$. The three triality modules are inequivalent orthogonal complex representations.  Consequently
$\dim(\Lambda^2U_{56})^S=\binom82+3=31.$
The extra three dimensions are the invariant alternating forms on the multiplicity spaces of the three doubled orthogonal modules. There are no invariant cross terms between inequivalent triality types.

Use the decompositions
\begin{align*}
 \Lambda^2U_{56}&=\mathbb C\oplus U_{1539},\\
 \Lambda^3U_{56}&=U_{56}\oplus U_{27664},\\
 \Lambda^4U_{56}&=\mathbb C\oplus U_{1539}\oplus U_{365750},\\
 \Lambda^2U_{133}&=U_{133}\oplus U_{8645}.
\end{align*}
They can be checked without a large character calculation. For $k=2,3,4$, exterior products of the first $k$ weight vectors along the unbranched chain starting at the minuscule end node give highest weights $\omega_6,\omega_5,\omega_4$, respectively, in Bourbaki numbering.  Their Weyl dimensions are $1539,27664,365750$.  The invariant symplectic form embeds $\Lambda^{k-2}U_{56}$ into $\Lambda^kU_{56}$, and
$\binom{56}{k}-\binom{56}{k-2} =1539,\ 27664,\ 365750 \quad(k=2,3,4),$
so these exhaust the primitive parts. For the last identity, the bracket onto $\mathfrak e_7$ is surjective. The wedge of a highest-root vector with the adjacent root vector has highest weight $\omega_3$ and Weyl dimension $8645=\binom{133}{2}-133$. It follows that
\begin{equation*}
 \begin{array}{c|rrrrrr}
 V&56&133&1539&27664&365750&8645\\ \hline
 \dim V^S\text{ is at least}&8&9&30&48&39&27.
 \end{array}
\end{equation*}
For example the fourth exterior power has at least $\binom84=70$ fixed vectors, and subtracting $30+1$ leaves at least $39$.

Finally, the standard $D_6A_1$ branching contains
$U_{912}|_{D_6A_1} =(U_{12}\otimes V_2)\oplus(U_{220}\otimes V_2) \oplus(U_{32}\otimes V_3)\oplus(U_{352}\otimes\mathbb C).$
The half-spin and vector-spinor chirality labels do not affect the argument.  Since $U_{12}|_{D_4}=U_{8v}\oplus4\mathbb C$, the first summand alone contributes eight fixed vectors. This branching is given in \cite[Table~52, p.~112]{Sl1981}. All seven fundamental representations therefore have nonzero fixed spaces.
\end{proof}
\FloatBarrier
\paragraph{Five rank-one branchings in $E_6$.}\label{app:e6-a1-branching}
We verify the five branchings used in Propositions~\ref{prop99} and~\ref{prop910} by classical restrictions and by roots. We write $R_d$ for the real form of $V_d$ when $d$ is odd. For index $6$, embed $\SO(3)\times \SO(3)$ in $\SO(10)$ through $\mathbb R^{10}=(R_3\otimes R_3)\oplus\mathbb R$. Each factor has ambient index $6$, and the complex half-spin module restricts as $(V_4\otimes V_2)\oplus(V_2\otimes V_4)$. For index $12$, embed the diagonal $\SU(2)$ in $\Spin(5)\times\Spin(3)\subset\Spin(8)$ with vector representation $R_5\oplus R_3$. Both half-spin representations restrict as $V_4\otimes V_2=V_5\oplus V_3$. The three charged triality modules in $E_6/(\Spin(8)T^2)$ therefore each contain the required $V_{5}\oplus V_{3}$. For index $21$, take the diagonal of the principal $\SU(2)\subset \SU(5)$ (index $20$) and the external root $\SU(2)$ in $\SU(6)\SU(2)$. Then $\Lambda^2V_5\otimes V_2=(V_7\oplus V_3)\otimes V_2 =V_8\oplus V_6\oplus V_4\oplus V_2$, with nonzero central charge. For index $30$, use $\mathbb R^{10}=R_7\oplus R_3$. The spin module of $\Spin(7)$ restricts as $V_7\oplus V_1$, so the half-spin module of $\Spin(10)$ becomes $(V_7\oplus V_1)\otimes V_2=V_8\oplus V_6\oplus V_2$. For index $60$, the principal $\SU(2)$ in $\Spin(9)$, with vector $R_9\oplus\mathbb R$ in $\SO(10)$, restricts the half-spin module as $V_{11}\oplus V_5$. All these spin restrictions follow immediately by evaluating the standard half-spin weights. The indices are checked by $\operatorname{index}(R_d\subset \SO(N))=d(d^2-1)/12$ and additivity. For the index-$6$ tensor embedding the multiplicity of $R_3$ is $3$. These restrictions give the modules used in the two-free arguments.

The branchings in Propositions~\ref{prop99} and~\ref{prop910} can also be computed directly from roots. Use the $E_6$ Cartan matrix with edges $13,34,24,45,56$, and generate its $72$ roots by simple reflections. The following are simple-coroot coordinates of the $A_1$ element $H$ and a basis of its centralizer Cartan used in the calculation.
\begin{table}[htbp]
\centering\small
\caption{Coroot and centralizer coordinates for five $E_6$ rank-one classes}\label{tab:e6-a1-root-check}
$\begin{array}{c|c|l}
\text{index}&H&\text{centralizer Cartan basis}\\\hline
6&(3,4,6,8,6,3)&(1,-1,1,-1,-1,0),\ (1,-2,0,-2,0,1)\\
12&(4,6,8,12,8,4)&(0,0,1,0,-1,0),\ (1,0,0,0,0,-1)\\
21&(6,8,11,15,11,6)&(2,0,1,3,-1,-2)\\
30&(7,10,13,18,13,7)&(1,0,-1,0,1,-1)\\
60&(10,14,18,26,18,10)&(1,0,2,0,-2,-1)
\end{array}$
\end{table}
For a root $\gamma$, evaluate $\gamma(H)$ and its values on the centralizer basis, add six zero Cartan weights, and subtract $A_1$ characters $t^m+t^{m-2}+\cdots+t^{-m}$ in each central weight space. This yields exactly the charged modules displayed above. In direct evaluation coordinates the two index-$12$ characters are $(1,1),(-1,2)$, and the common nonzero character in the index-$21,30,60$ rows is $3$. The calculator uses twice these centralizer weight coordinates. For index $6$, the second $A_1$ coroot is $(1,-2,0,-2,0,1)$ and a central generator is $(1,0,2,0,-2,-1)$. Subtraction of the two $A_1$ characters gives
\begin{equation*}
 \begin{aligned}
 \mathfrak e_{6,\mathbb C}
 ={}&(V_{3},V_{1})_0\oplus(V_{1},V_{3})_0\oplus(V_{1},V_{1})_0\\
 &\oplus(V_{5},V_{3})_0\oplus(V_{3},V_{5})_0\oplus(V_{3},V_{3})_0\\
 &\oplus(V_{4},V_{2})_{3}\oplus(V_{4},V_{2})_{-3}
 \oplus(V_{2},V_{4})_{3}\oplus(V_{2},V_{4})_{-3}.
 \end{aligned}
\end{equation*}
This calculation determines the weights and multiplicities, and verifies the dimension sums and the normalization $(H,H)=2\,\mathrm{index}$ by integer and rational arithmetic.

\FloatBarrier
\paragraph{Two rank-one branchings in $E_7$.}\label{app:e7-a1-branching}
Let $A$ be the $E_7$ Cartan matrix with edges $13,34,42,45,56,67$. All roots have squared length $2$. The first six simple roots form $E_6$. The element
$h=(12,16,22,30,22,12,0)$
is written in the simple-coroot basis. Its labels on $E_6$ are $(2,2,2,0,2,2)$, the diagram of $E_6(a_1)$, and $h^tAh/2=84$. The $27$ roots whose coefficient at the last simple root is $1$ have $h$-character
$(x^{12}+x^{10}+\cdots+x^{-12}) +(x^8+x^6+\cdots+x^{-8}) +(x^4+x^2+\cdots+x^{-4}),$
which proves \eqref{formula36}. This can also be verified by generating the $126$ roots by simple reflections. Only integer arithmetic is needed. For $h=(6,8,11,15,11,6,0)$ and the internal central direction $z=(2,0,1,3,-1,-2,0)$ the same computation gives \eqref{formula37}, with internal weights $2,-4,-1,5$. This proves the charged restrictions used in Proposition~\ref{prop916}.

\FloatBarrier
\medskip
\noindent\textbf{AI Assistance.} AI was used under the authors' direction to assist with mathematical exploration, drafting, literature searches, and internal review of this work. The authors are responsible for independently checking the mathematical arguments and references in the final manuscript, make all final decisions concerning its contents, and take full responsibility for the paper.


\begin{thebibliography}{99}

\bibitem{AdKaMoPaPe2018}
D. Adamovi\'c, V. G. Kac, P. M\"oseneder Frajria,
P. Papi and O. Per\v se,
\emph{On the classification of non-equal rank affine conformal
embeddings and applications}, Selecta Math. \textbf{24} (2018),
2455--2498, equation (5.7),
\url{https://arxiv.org/abs/1702.06089}.

\bibitem{AlAr2007}
D. V. Alekseevsky and A. Arvanitoyeorgos, \emph{Riemannian flag manifolds with homogeneous geodesics}, Trans. Amer. Math. Soc. \textbf{359} (2007), no. 8, 3769--3789, \url{https://doi.org/10.1090/S0002-9947-07-04277-8}.

\bibitem{AlNi2009}
D. Alekseevsky and Yu. G. Nikonorov,
\emph{Compact Riemannian manifolds with homogeneous geodesics},
SIGMA \textbf{5} (2009), 093,
\url{https://arxiv.org/abs/0904.3592}.

\bibitem{ArVaMo2024}
T. Arakawa, J. van Ekeren and A. Moreau,
\emph{Singularities of nilpotent Slodowy slices and collapsing levels
of W-algebras}, Forum Math. Sigma \textbf{12} (2024), e95,
\url{https://doi.org/10.1017/fms.2024.81}.

\bibitem{Ar2017}
A. Arvanitoyeorgos, \emph{Homogeneous manifolds whose geodesics are orbits. Recent results and some open problems}, Irish Math. Soc. Bull. \textbf{79} (2017), 5--29, \url{https://www.irishmathsoc.org/bull79/Arvanitoyeorgos.pdf}.

\bibitem{ArQiWaZh2019}
A. Arvanitoyeorgos, H. Qin, Y. Wang and G. Zhao, \emph{Riemannian generalized C-spaces with homogeneous geodesics}, Filomat \textbf{33} (2019), no. 4, 1117--1124, \url{https://doi.org/10.2298/FIL1904117A}.

\bibitem{ArSoSt2021a}
A. Arvanitoyeorgos, N. P. Souris and M. Statha, \emph{Geodesic orbit metrics in a class of homogeneous bundles over real and complex Stiefel manifolds}, Geom. Dedicata \textbf{215} (2021), 31--50, \url{https://doi.org/10.1007/s10711-021-00639-6}.

\bibitem{ArSoSt2021b}
A. Arvanitoyeorgos, N. P. Souris and M. Statha, \emph{Geodesic orbit metrics in a class of homogeneous bundles over quaternionic Stiefel manifolds}, J. Geom. Phys. \textbf{165} (2021), 104223, \url{https://doi.org/10.1016/j.geomphys.2021.104223}.

\bibitem{ArWa2017}
A. Arvanitoyeorgos and Y. Wang, \emph{Homogeneous geodesics in generalized Wallach spaces}, Bull. Belg. Math. Soc. Simon Stevin \textbf{24} (2017), no. 2, 257--270, \url{https://doi.org/10.36045/bbms/1503453709}.

\bibitem{ArWaZh2017}
A. Arvanitoyeorgos, Y. Wang and G. Zhao, \emph{Riemannian g.o. metrics in certain M-spaces}, Differential Geom. Appl. \textbf{54, Part A} (2017), 59--70, \url{https://doi.org/10.1016/j.difgeo.2016.12.005}.

\bibitem{ArWaZh2018}
A. Arvanitoyeorgos, Y. Wang and G. Zhao, \emph{Riemannian M-spaces with homogeneous geodesics}, Ann. Global Anal. Geom. \textbf{54} (2018), 315--328, \url{https://doi.org/10.1007/s10455-018-9603-7}.

\bibitem{BeElMa2009}
P. Benito, A. Elduque and F. Mart\'in-Herce, \emph{Irreducible Lie--Yamaguti algebras}, J. Pure Appl. Algebra \textbf{213} (2009), 795--808, \url{https://doi.org/10.1016/j.jpaa.2008.09.003}.

\bibitem{BeKo2021}
Y. Benoist and T. Kobayashi,
\emph{Tempered homogeneous spaces III}, J. Lie Theory \textbf{31} (2021), 833--869,
\url{https://arxiv.org/abs/2009.10389}.

\bibitem{BeNi2020}
V. N. Berestovskii and Yu. G. Nikonorov, \emph{Riemannian Manifolds and Homogeneous Geodesics}, Springer Monographs in Mathematics, Springer, Cham, 2020, \url{https://doi.org/10.1007/978-3-030-56658-6}.

\bibitem{BeKoVa1997}
J. Berndt, O. Kowalski and L. Vanhecke, \emph{Geodesics in Weakly Symmetric Spaces}, Ann. Global Anal. Geom. \textbf{15} (1997), no. 2, 153--156, \url{https://doi.org/10.1023/A:1006565909527}.

\bibitem{ChChCh2025}
C. Chen, H. Chen and Z. Chen, \emph{Geodesic orbit metrics on homogeneous spaces arising from generalized flag manifolds with two irreducible summands}, J. Geom. Phys. \textbf{217} (2025), 105641, \url{https://doi.org/10.1016/j.geomphys.2025.105641}.

\bibitem{ChChDe2018}
H. Chen, Z. Chen and S. Deng, \emph{Compact simple Lie groups admitting left-invariant Einstein metrics that are not geodesic orbit}, C. R. Math. \textbf{356} (2018), no. 1, 81--84, \url{https://doi.org/10.1016/j.crma.2017.11.018}.

\bibitem{ChChWo2018}
H. Chen, Z. Chen and J. A. Wolf, \emph{Geodesic orbit metrics on compact simple Lie groups arising from flag manifolds}, C. R. Math. \textbf{356} (2018), no. 8, 846--851, \url{https://doi.org/10.1016/j.crma.2018.06.004}.

\bibitem{ChChYaZh2023}
H. Chen, Z. Chen, Z. Yan and F. Zhu, \emph{Invariant geodesic orbit metrics on certain compact homogeneous spaces}, Manuscripta Math. \textbf{172} (2023), 651--668, \url{https://doi.org/10.1007/s00229-022-01416-9}.

\bibitem{ChChZh2021}
H. Chen, Z. Chen and F. Zhu, \emph{Geodesic orbit metrics on homogeneous spaces constructed by strongly isotropy irreducible spaces}, Sci. China Math. \textbf{64} (2021), no. 10, 2313--2326, \url{https://doi.org/10.1007/s11425-020-1818-y}.

\bibitem{ChChZh2026}
H. Chen, Z. Chen and F. Zhu,
\emph{Invariant geodesic orbit metrics on homogeneous spaces
with intermediate subgroups}, J. Lie Theory \textbf{36} (2026), 9--21,
\url{https://doi.org/10.5802/jolt.1415}.

\bibitem{ChNi2019}
Z. Chen and Yu. G. Nikonorov,
\emph{Geodesic orbit Riemannian spaces with two isotropy summands. I},
Geom. Dedicata \textbf{203} (2019), 163--178,
\url{https://arxiv.org/abs/1704.01913}.

\bibitem{ChNiNi2023}
Z. Chen, Y. Nikolayevsky and Yu. G. Nikonorov,
\emph{Compact geodesic orbit spaces with a simple isotropy group},
Ann. Global Anal. Geom. \textbf{63} (2023), 7,
\url{https://arxiv.org/abs/2009.06329}.

\bibitem{DaZi1979}
J. E. D'Atri and W. Ziller, \emph{Naturally reductive metrics and Einstein metrics on compact Lie groups}, Mem. Amer. Math. Soc. \textbf{18} (1979), no. 215, iii+72 pp., \url{https://doi.org/10.1090/memo/0215}.

\bibitem{DuKoNi2004}
Z. Du\v{s}ek, O. Kowalski and S. \v{Z}. Nik\v{c}evi\'{c}, \emph{New examples of Riemannian g.o. manifolds in dimension 7}, Differential Geom. Appl. \textbf{21} (2004), no. 1, 65--78, \url{https://doi.org/10.1016/j.difgeo.2004.03.006}.

\bibitem{Go1996}
C. S. Gordon, \emph{Homogeneous Riemannian Manifolds Whose Geodesics Are Orbits}, in \emph{Topics in Geometry: In Memory of Joseph D'Atri} (S. Gindikin, ed.), Progr. Nonlinear Differential Equations Appl., vol. 20, Birkh\"auser, Boston, 1996, 155--174, \url{https://doi.org/10.1007/978-1-4612-2432-7_4}.

\bibitem{GoNi2018}
C. S. Gordon and Yu. G. Nikonorov, \emph{Geodesic orbit Riemannian structures on $\mathbb R^n$}, J. Geom. Phys. \textbf{134} (2018), 235--243, \url{https://doi.org/10.1016/j.geomphys.2018.08.018}.

\bibitem{GrGrNe2020}
B. Grajales, L. Grama and C. J. C. Negreiros, \emph{Geodesic orbit spaces in real flag manifolds}, Comm. Anal. Geom. \textbf{28} (2020), no. 8, 1933--2003, \url{https://doi.org/10.4310/CAG.2020.v28.n8.a7}.

\bibitem{GrRoSa2025}
B. Grajales, G. Rond\'on and J. Saavedra, \emph{Riemannian geometry of $G_2$-type real flag manifolds}, Eur. J. Math. \textbf{11} (2025), 59, \url{https://doi.org/10.1007/s40879-025-00848-3}.

\bibitem{Ka1983}
A. Kaplan, \emph{On the geometry of groups of Heisenberg type}, Bull. London Math. Soc. \textbf{15} (1983), 35--42, \url{https://doi.org/10.1112/blms/15.1.35}.

\bibitem{Ko2009}
A. Kollross,
\emph{Low cohomogeneity and polar actions on exceptional compact Lie groups},
Transform. Groups \textbf{14} (2009), 387--415,
\url{https://doi.org/10.1007/s00031-009-9052-2}.

\bibitem{Ko1956}
B. Kostant, \emph{On differential geometry and homogeneous spaces. II}, Proc. Natl. Acad. Sci. USA \textbf{42} (1956), no. 6, 354--357, \url{https://doi.org/10.1073/pnas.42.6.354}.

\bibitem{Ko1973}
B. Kostant,
\emph{On convexity, the Weyl group and the Iwasawa decomposition},
Ann. Sci. \'Ecole Norm. Sup. \textbf{6} (1973), 413--455,
\url{https://doi.org/10.24033/asens.1254}.

\bibitem{KoVa1991}
O. Kowalski and L. Vanhecke, \emph{Riemannian manifolds with homogeneous geodesics}, Boll. Un. Mat. Ital. B (7) \textbf{5} (1991), 189--246.

\bibitem{Kr1979}
M. Kr\"amer, \emph{Sph\"arische Untergruppen in kompakten zusammenh\"angenden Liegruppen}, Compositio Math. \textbf{38} (1979), no. 2, 129--153.

\bibitem{La1999}
J. Lauret, \emph{Modified H-type groups and symmetric-like Riemannian spaces}, Differential Geom. Appl. \textbf{10} (1999), no. 2, 121--143, \url{https://doi.org/10.1016/S0926-2245(99)00002-9}.

\bibitem{LeMeToWo1998}
E. Lerman, E. Meinrenken, S. Tolman and C. Woodward,
\emph{Non-abelian convexity by symplectic cuts},
Topology \textbf{37} (1998), 245--259,
\url{https://arxiv.org/abs/dg-ga/9603015}.

\bibitem{LuTaXu2025}
W. Luo, J. Tan and N. Xu, \emph{New non-geodesic orbit Einstein metrics on $\Sp(n)$}, Manuscripta Math. \textbf{176} (2025), 13, \url{https://doi.org/10.1007/s00229-025-01614-1}.

\bibitem{Mi2006}
A. N. Minchenko,
\emph{The semisimple subalgebras of exceptional Lie algebras},
Trans. Moscow Math. Soc. (2006), 225--259,
\url{https://pi.math.cornell.edu/~andreim/semsimC.pdf}.

\bibitem{Ng2000}
H. D. Nguyen, \emph{Compact weakly symmetric spaces and spherical pairs}, Proc. Amer. Math. Soc. \textbf{128} (2000), no. 11, 3425--3433, \url{https://doi.org/10.1090/S0002-9939-00-05425-3}.

\bibitem{Ni2020}
Y. Nikolayevsky, \emph{Geodesic orbit and naturally reductive nilmanifolds associated with graphs}, Math. Nachr. \textbf{293} (2020), no. 4, 754--760, \url{https://doi.org/10.1002/mana.201800456}.

\bibitem{NiNi2019}
Y. Nikolayevsky and Yu. G. Nikonorov, \emph{On invariant Riemannian metrics on Ledger--Obata spaces}, Manuscripta Math. \textbf{158} (2019), no. 3--4, 353--370, \url{https://doi.org/10.1007/s00229-018-1029-9}.

\bibitem{NiZi2026a}
Y. Nikolayevsky and W. Ziller, \emph{Non-singular geodesic orbit nilmanifolds}, Ann. Mat. Pura Appl. \textbf{205} (2026), 1939--1973, \url{https://doi.org/10.1007/s10231-026-01661-9}.

\bibitem{Ni2013a}
Yu. G. Nikonorov, \emph{Geodesic orbit manifolds and Killing fields of constant length}, Hiroshima Math. J. \textbf{43} (2013), no. 1, 129--137, \url{https://doi.org/10.32917/hmj/1368217953}.

\bibitem{Ni2013b}
Yu. G. Nikonorov, \emph{Geodesic orbit Riemannian metrics on spheres}, Vladikavkaz. Mat. Zh. \textbf{15} (2013), no. 3, 67--76, \url{https://www.mathnet.ru/eng/vmj473}.

\bibitem{Ni2017}
Yu. G. Nikonorov, \emph{On the structure of geodesic orbit Riemannian spaces},
Ann. Global Anal. Geom. \textbf{52} (2017), 289--311,
\url{https://arxiv.org/abs/1611.01050}.

\bibitem{Ni2019}
Yu. G. Nikonorov, \emph{On left-invariant Einstein Riemannian metrics that are not geodesic orbit}, Transform. Groups \textbf{24} (2019), no. 2, 511--530, \url{https://doi.org/10.1007/s00031-018-9476-7}.

\bibitem{Ni2024a}
Yu. G. Nikonorov, \emph{On geodesic orbit nilmanifolds}, J. Geom. Phys. \textbf{203} (2024), 105257, \url{https://doi.org/10.1016/j.geomphys.2024.105257}.

\bibitem{Ri1984}
C. Riehm, \emph{Explicit Spin representations and Lie algebras of Heisenberg type}, J. London Math. Soc. (2) \textbf{29} (1984), 49--62, \url{https://doi.org/10.1112/jlms/s2-29.1.49}.

\bibitem{Sa2026}
SageMath developers, \emph{Branching rules},
SageMath Reference Manual, ``Miscellaneous'',
\url{https://doc.sagemath.org/html/en/reference/combinat/sage/combinat/root_system/branching_rules.html}.

\bibitem{Sh2024}
Y. Shan, \emph{Level one automorphic representations of an anisotropic exceptional group over $\mathbb Q$ of type $F_4$}, preprint, arXiv:2407.05859 (2024), to appear in M\'em. Soc. Math. Fr., \url{https://arxiv.org/abs/2407.05859}.

\bibitem{Sl1981}
R. Slansky,
\emph{Group theory for unified model building}, Phys. Rep.
\textbf{79} (1981), 1--128, Table 52,
\url{https://doi.org/10.1016/0370-1573(81)90092-2}.

\bibitem{So2018}
N. P. Souris, \emph{Geodesic orbit metrics in compact homogeneous manifolds with equivalent isotropy submodules}, Transform. Groups \textbf{23} (2018), no. 4, 1149--1165, \url{https://doi.org/10.1007/s00031-017-9464-3}.

\bibitem{So2021}
N. P. Souris, \emph{On a class of geodesic orbit spaces with abelian isotropy subgroup}, Manuscripta Math. \textbf{166} (2021), 101--129, \url{https://doi.org/10.1007/s00229-020-01236-9}.

\bibitem{So2022}
N. P. Souris, \emph{Geodesic orbit spaces of compact Lie groups of rank two}, Geom. Dedicata \textbf{216} (2022), 1, \url{https://doi.org/10.1007/s10711-021-00668-1}.

\bibitem{So2025a}
N. P. Souris, \emph{Einstein Lie groups, geodesic orbit manifolds and regular Lie subgroups}, Commun. Contemp. Math. \textbf{27} (2025), 2350068, \url{https://doi.org/10.1142/S0219199723500682}.

\bibitem{So2025b}
N. P. Souris, \emph{Geodesic orbit metrics on the compact Lie group $G_2$}, Proc. Amer. Math. Soc. (2025), published online 23 October 2025, \url{https://doi.org/10.1090/proc/17436}.

\bibitem{So2026}
N. P. Souris,
\emph{Compact geodesic orbit manifolds with regular stabilizer},
Ann. Mat. Pura Appl. \textbf{205} (2026), 689--702,
\url{https://doi.org/10.1007/s10231-025-01618-4}.

\bibitem{Ta1999}
H. Tamaru, \emph{Riemannian G.O. spaces fibered over irreducible
symmetric spaces}, Osaka J. Math. \textbf{36} (1999), 835--851,
\url{https://doi.org/10.18910/10200}.

\bibitem{To2010}
G. Z. T\'{o}th, \emph{On Lagrangian and Hamiltonian systems with homogeneous trajectories}, J. Phys. A Math. Theor. \textbf{43} (2010), no. 38, 385206, \url{https://doi.org/10.1088/1751-8113/43/38/385206}.

\bibitem{XuTa2019}
N. Xu and J. Tan, \emph{On left-invariant Einstein metrics that are not geodesic orbit}, C. R. Math. \textbf{357} (2019), no. 7, 624--628, \url{https://doi.org/10.1016/j.crma.2019.07.003}.

\bibitem{Ya2004}
O. S. Yakimova, \emph{Weakly symmetric Riemannian manifolds with reductive isometry group}, Sb. Math. \textbf{195} (2004), no. 4, 599--614, \url{https://doi.org/10.1070/SM2004v195n04ABEH000817}.

\bibitem{Zi1996}
W. Ziller, \emph{Weakly symmetric spaces}, in \emph{Topics in Geometry: In Memory of Joseph D'Atri} (S. Gindikin, ed.), Progress in Nonlinear Differential Equations and Their Applications, vol. 20, Birkh\"auser, Boston, 1996, 355--368.

\end{thebibliography}
\end{document}